%% file: main.tex
\documentclass[alpha-refs]{wiley-article}
\usepackage{color}
\usepackage{ulem}

\DeclareMathOperator{\UOT}{UOT_{\varepsilon}}
\DeclareMathOperator{\Sink}{S_{\varepsilon}}

\usepackage{siunitx}
\usepackage{lineno}
\usepackage{setspace}
\usepackage{tabularx}
\usepackage[page]{appendix}
\usepackage{listings,xcolor}
\definecolor{oxfordblue}{rgb}{0.0, 0.1, 0.5}
\usepackage[colorlinks=true, allcolors=oxfordblue]{hyperref}
\usepackage{float}
\usepackage{csvsimple}
\usepackage{subcaption}
\usepackage{graphicx} 
\usepackage{booktabs}
\usepackage{algorithm}
\usepackage{algorithmic}

\usepackage{multicol}
\usepackage[english]{babel}
\usepackage[T1]{fontenc}

\makeatletter
\csvset{
  autotabularcenter/.style={
    file=#1,
    after head=\csv@pretable\begin{tabular}{|*{\csv@columncount}{c|}}\csv@tablehead,
    table head=\hline\csvlinetotablerow\\\hline,
    late after line=\\,
    table foot=\\\hline,
    late after last line=\csv@tablefoot\end{tabular}\csv@posttable,
    command=\csvlinetotablerow},
  autobooktabularcenter/.style={
    file=#1,
    after head=\csv@pretable\begin{tabular}{*{\csv@columncount}{c}}\csv@tablehead,
    table head=\toprule\csvlinetotablerow\\\midrule,
    late after line=\\,
    table foot=\\\bottomrule,
    late after last line=\csv@tablefoot\end{tabular}\csv@posttable,
    command=\csvlinetotablerow},
}
\makeatother

\newcommand{\csvautobooktabularcenter}[2][]{\csvloop{autobooktabularcenter={#2},#1}}

\newcommand{\norm}[1]{\lVert #1 \rVert}

\newcommand{\defeq}{\vcentcolon=}

\title{Examining Entropic Unbalanced Optimal Transport and Sinkhorn Divergences for Spatial Forecast Verification}

\author[1,2]{Jacob J. M. Francis}
\author[2]{Colin J. Cotter}
\author[3]{Marion P. Mittermaier}

\affil[1]{NERC Doctoral Training Programme in Science and Solutions for a Changing Planet, Imperial College London, London, SW7 2AZ, UK}
\affil[2]{Department of Mathematics, Imperial College London, London, SW7 2AZ, UK}
\affil[3]{Met Office, Exeter, Devon, EX1 3PB, UK}

\corremail{jacob.francis18@imperial.ac.uk}

\fundinginfo{This work was supported by the Natural Environment Research Council [grant number NE/S007415/1] and the CASE Studentship "Transport methods for the verification of Numerical Weather Prediction (NWP) forecasts" at the Met Office}

\begin{document}

\maketitle

\begin{abstract}
An optimal transport (OT) problem seeks to find the cheapest mapping between two distributions with equal total density, given the cost of transporting density from one place to another. 
Unbalanced OT allows for different total density in each distribution. 
This is the typical setting for precipitation forecast and observation data, when considering the densities as accumulated rainfall, or intensity.
True OT problems are computationally expensive, however through entropic regularisation it is possible to obtain an approximation maintaining many of the underlying attributes of the true problem. 
In this work, entropic unbalanced OT and its associated Sinkhorn divergence are examined as a spatial forecast verification method for precipitation data.
The latter being a novel introduction to the forecast verification literature. 
It offers many attractive features, such as morphing one field into another, defining a distance between fields and providing feature based optimal assignment.
This method joins the growing research by the Spatial Forecast Verification Methods Inter-Comparison Project (ICP) which aims to unite spatial verification approaches. 
After testing this methodology’s behaviour on numerous ICP test sets, it is found that the Sinkhorn divergence is robust against the common double penalty problem (a form of phase error), on average aligns with expert assessments of model performance, and allows for a variety of novel pictorial illustrations of error. It provides informative summary scores, and has few limitations to its application.
Combined, these findings place unbalanced entropy regularised optimal transport and the Sinkhorn divergence as an informative method which follows geometric intuition. 
\end{abstract}

%% For editing and Reviewing; Line numbers and double spacing
% \linenumbers
%% \doublespacing  

% \tableofcontents

%% Main Tex
\input{content/1_introduction_discrete}
\input{content/2_ot_discrete}
\input{content/3_data_description}
\input{content/4_results}
\input{content/5_conclusion}

\bibliography{bib/references}

%% appendix
\input{appendix/appendix.tex}

\end{document}

%% file: content/1_introduction_discrete.tex
\newpage
\section{Introduction}

Without forecast verification, it would be impossible to accurately assess a precipitation, or more generally a Numerical Weather Prediction (NWP), model's overall \textit{goodness}. 
As observed by \citet{murphy_1992}, there are three key factors that contribute to this assessment: consistency, quality, and value. 
Since Murphy's essay, there have been significant developments in the accuracy of weather models over longer time periods and in their increased spatial resolution. 
However, as forecasters move to higher resolution models, traditional metrics of quality have struggled to capture the improvements in model accuracy which are believed to be present \citep{mass_et_al_2002}.
This development has been coupled with more readily available spatial and temporal observational data \citep{hou_et_al_2014}. 
Though, this expansion in data has complicated the determination of quality. While the value of a precipitation forecast is subjective to the user, quality should be objectively quantified, either numerically, graphically, or pictorially, to show the correspondence between the forecast and the observation.
 
It has been observed that traditional metrics, based on point-wise evaluations, fail to accurately score increasingly complex data \citep{brown_et_al_2012}.
One significant complication is the \textit{double-penalty problem}; in this scenario, a forecast predicts the correct spatial extent (shape) and intensity but in the wrong location. 
Point-wise, or \textit{local}-metrics, penalise this forecast twice — once for failing to detect the correct event and once for predicting an event where there was none.
Despite this, local-metrics have undergone many years of research and are a foundation of operational forecasting, meaning that users are often more familiar with reading and interpreting their outputs.
What is missing from these local metrics is the consideration of spatial relationships. 
The emphasis in modern precipitation verification now lies in incorporating this spatial information for greater diagnostic capacity. To maximise impact, this needs to be combined with education around these novel metrics. This will enable users to utilise them in a meaningful and streamlined way.
 
To address the double penalty problem and incorporate non-local information, such as spatial patterns, structure, and covariance, verification practitioners have researched spatial forecast verification methods under the ICP (Spatial Forecast Verification Inter-comparison Project, \citet{gillelandetal_2009a, gilleland_et_al_2010}). 
This collective effort led to the rigorous study, creation, and adaptation of many spatial methods, often leveraging principles from computer vision and image processing in a new guise. 
From their work and that of \citet{dorninger_et_al_2018}, five overarching categories were established, broadly encompassing most identified methods (which are not limited to precipitation verification). 
These are\footnote{Table 1 of \citet{gillelandetal_2009a} presents a more complete picture, though an update-to-date and full reference to all ICP publications is available at \citet{eric_reference_list}}: (i) \textit{distance metrics}, such as Centroid Distance (CDST), Baddeley's Distance (BDEL), and Mean Error Distance (MED)  \citep{gilleland_et_al_2019}; (ii) \textit{neighbourhood methods}, including Fraction Skill Score (FSS) and others as reviewed by \citet{ebert_2008}; (iii) \textit{scale separation schemes}, such as intensity-scale techniques using Fourier or wavelet decomposition \citep{casati_ross_stephenson_2004}; (iv) \textit{feature-based methods} for example, the Method for Object-based Diagnostic Evaluation (MODE, \citet{davis_et_al_2006, davis_et_al_2009}) and Structure, Amplitude, and Location (SAL) \citep{wernli_2008}; and (v) \textit{field deformation techniques} like Optical Flow (OF) \citep{keil_craig_2007,keil_craig_2009, marzban_sandgathe_2010} or the image warping statistic (IWS) \citep{gilleland_2010}.
Through the ICPs combined effort there is now  centralised resources for metrics, their references \citep{eric_reference_list}, their implementation \citep{cran_spatialvx}, as well as datasets for testing any novel methodologies \citep{ral_icp_webpage}.
 
This present work shall introduce entropy regularised unbalanced optimal transport (UOT, or \(\UOT\) denoting the score) and its unbalanced Sinkhorn divergence \((\Sink)\) counterpart as alternative spatial forecast verification methods.
Specifically, whilst optimal transport (OT) and UOT have been studied by \citet{benamou_2003, skok_2023, nishizawa_2024}, this is the first presentation of the Sinkhorn divergence for precipitation verification (to our knowledge). 
These approaches draw on recent advances in computer vision and machine learning, and are rooted in rigorous mathematical concepts, including the related Wasserstein distance, balanced and partial optimal transport and the Gaussian Hellinger-Kantorovich distance \citep{peyre_cuturi_2018, marzban_sandgathe_2010, sejourne_et_al_2021, liero_mielke_savare_2016}.

The UOT methodology incorporates a unit-preserving distance, an approximate transport plan with diagnostic vectors, and a tunable parameter that defines the scale of relation between features. 
Therefore it does not directly fit into one of the five above categories of \citet{dorninger_et_al_2018}.
Instead, the UOT methodology can bridge multiple categories.
% In certain scenarios 
There are parallels between OT and OF methods, where both morph fields, yet through different methodologies.
% Equally, a
As shall be seen through the ICP idealised cases, UOT inherently matches features and those that are separated too far will not be matched: this then mimics a feature based approach.
% Similarly, t
The raw cost is unit preserving and hence a suitable distance metric.
Additionally, there is the potential, through numerical efficiency considerations, to capture information at sequentially finer scales, imitating a neighbourhood method.
 
By no means is this the first application of OT or the Wasserstein distance in meteorology. 
These concepts have been applied in data assimilation \citep{liu_frank_2022, bocquet_et_all_2024, feyuex_vidard_nodet_2018} , wind speed model aggregation \citep{papayiannis_et_al_2018}, hydrological misfit functions \citep{magyar_sambridge_2023}, and as a pollutant metric \citep{farchi_et_al_2016, vanderbecken_et_al_2023, hyan_et_al_2022}.
These studies range from testing 1D cases where the precise Wasserstein problem is simple to solve, to more complex 2D fields utilising both static and dynamic OT.
More pressingly, \citet{benamou_2003, skok_2023, skok_lledo_2024, nishizawa_2024} used OT theory as a metric for precipitation data in some capacity. 
\citet{benamou_2003} investigated an L\(_2\) penalised UOT problem without entropy.
While \citet{nishizawa_2024} studied the same entropy regularised UOT object as investigated in the present work, and \citet{skok_2023} introduced a related random-nearest neighbour OT type verification methodology.  
However, to the best of our knowledge, none have considered the bias effect of the entropic regularisation and debiasing of \(\Sink\).
% , or focussed on UOT's use as a score in its own right.

Specifically, through this investigation the following questions shall be answered:
\begin{enumerate}[label=Q\arabic*]
    \item\label{q1} How does UOT behave in the context of the double penalty problem? Is there a diagnostic advantage in using the Sinkhorn Divergence? (Sections \ref{section:the_perfect_forecast},  \ref{section:phase_error},  \ref{section:perturbed_fake_cases} )

   \item\label{q2}  How does UOT perform geometrically? Can it handle object rotations, detect zero-points and boundaries, or determine whether one field is a subset of another? (Sections \ref{section:boundary_position_effect},  \ref{section:rotation_scale})
   
    \item\label{q3} Can the score alone be diagnostic of translation error, spatial bias and intensity imbalances? Are additional tools such as visualisations of the transport vectors and the marginals necessary for diagnosis of error, or are these a useful aid in their own right? (Sections  \ref{section:unbalance_circles},  \ref{section:unbalanced_ellipses}, \ref{section:perturbed_fake_cases}, \ref{section:spring} )
 
    \item\label{q4} Given one tuneable parameter, the reach (or \(\rho\)), how does changing it affect the methods behaviour and associated tools, such as the optimised marginals of the plan? (Sections \ref{section:phase_error},   \ref{section:unbalance_circles}, \ref{section:perturbed_fake_cases} )
    %  old text: Does the tunable reach parameter define a strict skill threshold (or skill scale), and how does it affect the marginals of the returned plan?

    \item\label{q5} How does changing the type of marginal penalisation — total variation or Kullback-Leibler divergence — affect the behaviour of the score? Does this choice alter our understanding of \(\rho\) and the scores’ sensitivity to intensity imbalance? (Sections  \ref{section:phase_error},  \ref{section:unbalance_circles},  \ref{section:unbalanced_ellipses}, \ref{section:scatterholenoisy} , \ref{section:perturbed_fake_cases} )
    % olf text: Is there a difference in the diagnostic behaviour for different marginal penalisation of the UOT problem? Does this affect the influence of the tunable parameter and the sensitivity to any imbalance in intensity?

    \item\label{q6} Is there a way a forecaster can preferentially forecast to improve their UOT score, i.e. hedging and increasing overlap? (Sections  \ref{section:unbalance_circles},  \ref{section:unbalanced_ellipses})

    \item\label{q7} How does UOT perform in extreme, edge or null cases? Is it sensitive to small amounts of noise in the data? Is it always defined or are there situations where it is no longer defined? (Sections \ref{section:scatterholenoisy},  \ref{section:edge})
    
    \item\label{q8} Given real intensity data, are any properties from the simplified setting retained or lost? (Sections  \ref{section:perturbed_fake_cases} )

    \item\label{q9} Does the spatial metric achieve better alignment with human subjective evaluation, compared to traditional local metrics? (Sections \ref{section:unbalanced_ellipses},  \ref{section:spring})
\end{enumerate}
 
The rest of this article is organised as follows. In Section \ref{section:optimal_transport} a brief introduction to the background of OT is given (those familiar with the concept may wish to skip this).
Section \ref{section:methodology} presents the methodology of \(\UOT\) and \(\Sink\), with a brief note on implementation in Section \ref{section:implmentation}.
In Section \ref{section:data} the data is introduced, including the ICP's binary idealised cases in Section \ref{section:binary_data}, as well as the real textured intensity cases in Section  \ref{section:intensity_data} and the MesoVICT core case in \ref{section:vera_core_case_data}. 
Then, Section \ref{section:results} confers the first results for the simple set of idealised cases. These comprise: the perfect forecast (Section
 \ref{section:the_perfect_forecast}), boundary and relative position effects (Section
 \ref{section:boundary_position_effect}), sensitivity to phase error (translational error) (Section
 \ref{section:phase_error}), rotation and scale in the balanced setting (Section
 \ref{section:rotation_scale}), unbalanced circular (Section
 \ref{section:unbalance_circles}) and elliptical cases (Section
 \ref{section:unbalanced_ellipses}). Finishing the geometric cases with noisy, scattered and hole cases (Section
 \ref{section:scatterholenoisy}), and edge and extreme ones (Section
 \ref{section:edge}).
Real textured intensity cases are then tested in Section \ref{section:real_intensities}, with the Perturbed, Spring 2005, and MesoVICT cases (Sections  \ref{section:perturbed_fake_cases},  \ref{section:spring},  \ref{section:vera_mesovict}). Finally, Section \ref{section:conculsion} presents a summary and outlook.

\subsection{Optimal Transport}\label{section:optimal_transport}

OT was originally devised by \citet{monge1781} as a logistics problem, with the goal of optimising the construction of wartime trenches. 
From its inception, to economic interpretation and development by \citet{kantoro_1958}, and to its more recent boom in research and application in machine learning due to \citet{cuturi_2013}'s regularised version, OT has found use in a vast array of fields, from computer vision, machine learning and signal processing to fluid dynamics, economics, biology, and, more recently, spatial forecast verification \citep{Santambrogio_2015, benamou_brenier_1998, kolouri_et_al_2017, gramfort_peyre_cuturi_2015, benamou_2003, nishizawa_2024, skok_2022}.
 
By definition, an OT problem aims to find a minimising map between two density distributions, provided with a cost of transporting density from one place to another. 
In the meteorological context this could be 
any quantity described as an intensity, e.g. accumulated rainfall.
Mathematically and computationally, the space of possible maps to minimise over is unwieldy, and progress is made by relaxing the condition that
the mapping is one-to-one.
Instead of transporting mass from just one point to another, now mass from any point is allowed to be split and distributed to many points.

This relaxation is known as Kantorovich's (balanced) OT problem and is defined as follows (for full details starting from Monge's formulation, see \citet{villani_2009}).
Given two fields, which are represented by discrete data \(\{\alpha_i\}_{i=1}^{N}\),  \(\{\beta_j\}_{j=1}^{M}\) on grids \(x_i \in \mathcal{X} \subset \mathbb{R}^2  \) and \(y_j \in \mathcal{Y} \subset \mathbb{R}^2   \) which do not necessarily have to be the same grid. 
The density is constructed as \(\alpha = \sum_{i} \alpha_i \delta_{x_i}\), \(\beta = \sum_{j} \beta_j \delta_{y_i}\), with the condition (for balanced OT) that \(\sum_{i}\alpha_i = \sum_j \beta_j = 1\).
Notationally, the Dirac shorthand, \(\delta_{x_i} = \delta(x-x_i)\) allows us to define a grid through summations, and \(\alpha_i \delta_{x_i}\) says there is weight \(\alpha_i\) at grid point \(x_i\)\footnote{Formally, \(\delta(x) = \infty\) if \(x=0\) and zero otherwise, such that \(\int_{-\infty}^{\infty} \delta(x) dx = 1\)}.
The Kantorovich balanced OT problem searches for an optimal joint density, \(\pi = \sum_{i,j}\pi_{i,j}\delta_{x_i}\delta_{y_j} \), such that it has left marginal equal to \(\alpha\) and right marginal equal to \(\beta\), i.e. \(\pi_{1,j} \defeq \sum_i \pi_{i,j} = \beta_j, \  \pi_{0,i} \defeq\sum_j \pi_{i,j} = \alpha_i \). The transport plan $\pi_{i,j}$ tells us how much mass is transported from $x_i$ to $y_j$,
allowing for the mass in $x_i$ to be split and transported to multiple points. The transport plan is optimal if it minimises the associated cost or energy,
\begin{align}
    \sum_{i,j} c_{i,j}\pi_{i,j} .\label{eq:eq1}
\end{align}
Here \(c_{i,j} = c(x_i, y_j)\) for some function \(c\) which  defines a cost.
This cost function should characterise the cost of moving mass from point \(x_j\in\mathcal{X}\) to \(y_j\in\mathcal{Y}\).
In the setting of precipitation data, the `mass' is considered as accumulated rain volume, or intensity, for a given mesh-point or region, i.e. the weights \(\alpha_i, \beta_j\).
The p-Wasserstein distance is defined by using an \(L_p\)-norm cost function (defined via \(\norm{X}_{p} = (\sum_{i} (x_i)^p)^{\frac{1}{p}}\) for some vector \(X\)) \citep{peyre_cuturi_2018}.
For precipitation data, let us define the cost of moving a certain volume of rain via the squared Euclidean norm or L\(_2\)-norm with a factor of \(1/2\) so that transport vectors align with barycentres (Section \ref{section:barycentric_maps}).
That is, the cost function throughout shall be defined as \(c(x,y) = \frac{1}{2}\norm{x-y}_{2}^2\).
 
In the 1D setting, solving \eqref{eq:eq1} becomes a simple sorting problem. However, in 2D solving the true Wasserstein distance is non-trivial. 
For instance, given discrete densities the cost, \(c_{i,j}\), defines a matrix, and the associated optimisation problem is a linear program.
Subsequentially, it suffers from the \textit{curse of dimensionality}.
Solving this problem becomes inhibitive for a hundred or so points. 
Yet for NWP data, grids typically consist of many of thousands of points, even into the millions.
 
To overcome this dimensionality problem, \citet{cuturi_2013} proposed a further relaxation, known as \textit{entropic regularisation}.
This corresponds to the inclusion of entropy, or a Kullback-Leibler (KL) penalty term, (\ref{eq:kl}), between some reference and the plan.
It has received much interest and has been shown to be capable of scaling to millions of points on modern GPU hardware \citep{feydy_sejourne_vialard_2018, geomloss_webpage}.
Beneficially, whilst any relaxation of a problem will introduce error, the error (or \textit{bias}) in entropic regularisation has a known \textit{debiased} counterpart \citep{feydy_sejourne_vialard_2018}. 
This error-corrected formulation of OT is known as the Sinkhorn divergence (\(\Sink\)). It requires the computation of two extra OT problems which are cheap to calculate and capable of being solved in parallel to the main computation.
 
Unfortunately, to the detriment of application to real forecast-observation pairs, the definition given in \eqref{eq:eq1} assumes the two input densities have equal total mass (hence the term "balanced OT"), i.e. it is assumed that the two fields are unbiased (\(\sum_i \alpha_i = \sum_j \beta_j\)).
This is generally an unrealistic situation and so prevents the use of real rain volumes and intensities (or spatial coverage).
One solution is to scale the forecast and observation fields to have the same total mass (or total accumulated rainfall).
This mimics the idea of applying a percentile threshold which effectively removes the bias, e.g.
\citet{mittermaier_2013}.
However, this will lead to a loss of information: a model could be significantly over-forecasting intensities, which is undetected if the volumes are forced to be equal.
Instead, a further relaxation to the optimisation is introduced. This generalised problem corresponds to the static formulation of unbalanced optimal transport (UOT, \citet{liero_mielke_savare_2016}) (as opposed to the dynamic (unbalanced) formulation of \citet{benamou_brenier_2000}). 
In this setting, the input fields are allowed to have different total mass and the marginals of the optimised plan should be \textit{close} to the input fields. 
Moreover, the definition of \textit{close} is open to choice provided it follows a suitable divergence  (see \citet{sejourne_et_al_2021}).
UOT then compares two fields by a combination of transport and 
pointwise modification of the marginals.

In this work two penalisations are investigated; KL divergence and total variation (TV).
These are defined given discrete densities \(\lambda = \sum_k \lambda_k \delta_{z_k}, \gamma = \sum_k \gamma_k \delta_{z_k}\), for some grid points \(z_k\), as
\begin{align}
    KL(\lambda | \gamma ) &\defeq \sum_{k}\lambda_k (\log(\frac{\lambda_k}{\gamma_k}) - 1) + \sum_k \gamma_k , \label{eq:kl} \\
    TV(\lambda | \gamma ) &\defeq \sum_{k} |\frac{\lambda_k}{\gamma_k} - 1| \gamma_k .\label{eq:tv}
\end{align}
Given that wherever \(\gamma_k\) is zero then \(\lambda_k\) is also zero, and recall that \(\lim_{x\to 0} x\log(x)= 0\).
These divergences are defined as \(+\infty\) for negative weights, enforcing positivity by construction.
 
By directly viewing marginals of \(\pi\) or calculating the error between the input and marginal, various diagnostic capabilities can be identified.
These diagnostics are also linked to a reach parameter, resembling a skill-scale parameter, which we discuss later.
Through the joint density it is possible to calculate an approximate map from the observation to the forecast and vice versa. In a similar fashion to OF, these are not one-to-one mappings and each provides different information.
Given these transport vectors one can calculate average statistics, such as average direction and magnitude of transport, or to consider the distribution as a whole. 
This form of information is considered in the precipitation attribution distance (PAD) methodology \citep{skok_2023}.
As well as \citet{ebert_mcbridge_2000,gilleland_2013, gilleland_2010},  who use different methods to also obtain vector fields describing the alignment between the forecast and observation.
It is clear that united under the same framework of OT are an array of tools.
Rather than just point-wise considerations being evaluated and averaged, the underlying optimisation is taking into account spatial and intensity information simultaneously.

Throughout this work, no specific comparator metric is selected, though it is worth noting some similarities that OT holds with other metrics.
% A full literature review and comparison is out-of-scope for this current work.
The IWS provides a deformation technique that takes into account affine, rotational and non-linear transformations.
It comprises an intensity error term after morphing the observation.
In contrast, UOT does not allow for rotation, however has the benefit of inherently accounting for intensity error and dispersion of mass.
The IWS methodology also assumes the choice of control points which define the deformation — these are sparsely sampled on a  regular grid without feature selection. 
Our UOT methodology, also does not require feature selection (unlike MODE or contiguous rain areas, CRA \citep{ebert_mcbridge_2000}), however nor does it require control points to be selected.

OF techniques, also typically find an optimal deformation before intensity errors are studied.
For a comparison of computational expensive we consider the LK-OF method \citep{lucas_kanade_1981} which relies on the underlying assumptions that the fields are continuous as local spatial derivatives are taken, and that changes in derivatives are small between the two fields. These assumptions are not always true and are not required for OT. 
\citet{marzban_sandgathe_2010} extend the LK-OF method by including higher-order terms in the Taylor expansion, creating a non-linear optimisation problem that they tackle through the BFGS algorithm, which has a per-iteration complexity of \(O(n^2)\) (n being the number of grid points).
Alternatively, \citet{keil_craig_2007} consider displacement vectors found through pyramidal matching. 
The pyramid matching has complexity \(O(nL)\) \citep{grauman_darrell_2005}  (where L in the number of layers in the pyramid), then \citet{keil_craig_2007}'s method searches each layer for an optimal local pairing.
This process is similar to solving the Earth Mover's Distance \citep{rubner_et_al_2000}, a related optimal transport problem, though without entropic regularisation.  
Both OF methods have one tunable parameter which defines the window in which to look for perturbations, or the number of layers in the pyramid.
Entropy regularised UOT has a complexity of \(O(n^2 \log n)\) (provided the entropic parameter is suitably chosen, see Supplementary Materials \ref{appendix:implementation}), and has one tunable parameter relating to locality \citep{berman_2020, merigot_thibert_2020}.
However, unlike the LK-OF method, which handles non-matches through local squared errors, UOT allows for the destruction of mass when no match is possible or when matches are too distant.

It is clear then that UOT offers a comprehensive and integrated solution to spatial forecast verification, being multifaceted in its dealing with spatial and intensity error, computationally feasible and on a similar order to current methods, and parsimonious with some necessary choice.

%% file: content/2_ot_discrete.tex
\section{ Unbalanced Optimal Transport and Sinkhorn Divergence}\label{section:methodology}

To consolidate notation, consider a verifier is given a forecast density \(\mu_F = \sum_{j}\mu_{F,j}\delta_{X_{F,j}}\) and an observation density \(\mu_O = \sum_{i}\mu_{O,i}\delta_{X_{O,i}}\) both representing accumulated precipitation, over some region given by grids \(X_{F,j}\in\Omega_O\subset\mathbb{R}^2,\ X_{O,i}\in\Omega_F\subset\mathbb{R}^2\).
In general, these do not have to be the same grid.
The goal is to assess the quality of agreement between these fields, using the machinery of UOT and in two variants. 
The first variant uses KL relaxation of the marginal mass, the second variant TV (see \eqref{eq:tv} and \eqref{eq:kl}).
These were picked for their distinct characteristics: TV favours edges, filling the weighting in the nearest support, creating sharp cut-offs. In contrast, KL favours filling all the support smoothly with some mass everywhere, at the expense of over or under distributing mass.
For full details of UOT and the associated Sinkhorn divergence, see \citet{sejourne_et_al_2021} and \citet{chizat_peyre_et_al_2016}. Here a high level sketch will be given, providing only the necessary points for our methodology.
Let us define the UOT problem as
\begin{align}
    \frac{\UOT(\mu_O, \mu_F | \rho)}{\rho} =  \min_{\pi} \quad  \underbrace{\sum_{i,j} \left(\frac{\norm{X_{O,i} -X_{F,j}}}{\sqrt{2\rho}}\right)^2 \pi_{i,j} + \frac{\varepsilon}{\rho} KL(\pi | \mu_{O} \otimes \mu_{F})}_{\text{regularised transport cost}} + \underbrace{D(\pi_0 | \mu_{O}) + D(\pi_1 | \mu_{F})}_{\text{marginal mass imbalance cost}}, \label{eq:uot_general_cost}
\end{align}
where  entropic regularisation parameter \(\varepsilon>0\), tunable parameter \(\rho>0\), and where \(D\) is either the KL divergence or TV and the left and right marginals are denoted \(\pi_0\)  and \(\pi_1\), respectively.
The notation \(\mu_{O} \otimes \mu_{F}\) defines the outer product, i.e. \((\mu_{O}\otimes \mu_{F})_{ij} = \mu_{O,i}
\mu_{F,j}\). The entropic regularisation enforces positivity of the plan, and provides access to Sinkhorn's iterative algorithm through a dual associated problem. 
While the marginal penalty terms allow for the use of non-balanced, and hence realistic datasets.
See Supplementary Materials \ref{appendix:supplementary_material}.

Notice that changing \(\rho\) changes the balance of transport cost and marginal error penalty. 
That is, increasing \(\rho\) penalises the marginal mass imbalance cost, so that the inputs and plan's marginals must match as close as possible.
Whilst a smaller parameter permits more error in the marginals, thus greater differences between the input and marginal.
Crucially, this is a trade-off between marginal mass imbalance and transport cost.
% If the data were balanced then it possible that as \( \rho \rightarrow \infty, \UOT \rightarrow OT_{\varepsilon}\), whilst as \( \rho \rightarrow 0, \UOT \rightarrow m(\mu_0) + m(\mu_F)\). Zero limit is wron,g whilst is would be true without regulairsation, need to think a bit more when regularisation is included
Writing the unbalanced problem in the manner of \eqref{eq:uot_general_cost} highlights the distance dependence on the tunable parameter \(\rho\).
If two points \(X_{O,i}\) and \(X_{F,j}\) are too far apart, i.e. the distance between them is greater than \(\mathcal{O}(\sqrt{2\rho})\), then it is cheaper to pay to destroy mass than to move it.
If this is the case \(\pi_{i,j}\) is set close to, or equal to zero. 
On the other hand, if the points are separated by less than this, it is feasible to transport the point and not pay to destroy the associated mass.
\citet{sejourne_et_al_2021} refer to this relationship in \(\rho\) as the \textit{reach}. 
Hence, let us define the reach as \(\sqrt{2\rho}\) with units of distance.
These two views of \(\rho\), as a penalty parameter or reach distance, are central and will be drawn upon in Section \ref{section:results}.
% \textcolor{red}{rho infiity and zero limits for balanced data.}

By way of an overview, the iterative Sinkhorn optimisation is undertaking the following process; the error between the observation marginal (without loss of generality) and true observation is found, defined by KL or TV. 
This error is then mapped through the proposed plan and used to optimise the forecast marginal, and implicitly the plan.
Now the forecast marginal and input error is calculated, mapped, and used to update the observation marginal and implicitly the plan.
This iterative process is repeated, minimising the transport plan and adapting to the density inputs.
By controlling the contributions of each term in \eqref{eq:uot_general_cost}, through \(\rho\), it is possible to prioritise making the marginal and input densities more similar or prioritise cheaper transport.
% Again, this highlights the alternative optimisation viewpoint of \(\rho\) as a marginal penalty parameter.

\subsection{Sinkhorn Divergence}\label{section:sinkhorndivergence}
Through entropic regularisation, transport plans become increasingly diffusive around the true non-entropic plan. 
This is known as \textit{blurring} of the plan, which is a form of bias.
This bias manifests as the transport plan (and transport vectors) being distorted towards centres of mass \citep{feydy_sejourne_vialard_2018}. It also results in \(\UOT(\mu, \mu) \neq 0\) for arbitrary \(\mu\).
A significant contribution of this work shall be the study of the
unbalanced Sinkhorn divergence, a debiased counterpart to the entropy regularised biased problem, as a spatial forecast metric.

In \eqref{eq:sinkhorn_definition} the post-debiasing technique is introduced, where the \textit{symmetric} debiasing terms remove the associated entropic error up to the next order in \(\varepsilon\). 
\citet{sejourne_et_al_2021} define the unbalanced Sinkhorn divergence as
\begin{align}
    \Sink(\mu_O, \mu_F | \rho)= \UOT(\mu_O, \mu_F | \rho) \underbrace{- \frac{1}{2} \UOT(\mu_O, \mu_O | \rho) - \frac{1}{2} \UOT(\mu_F, \mu_F | \rho)}_{\text{symmetric debiasing terms}} + \underbrace{\frac{\varepsilon}{2}(m(\mu_O) + m(\mu_F))^2}_{\text{mass imbalance term}}, \label{eq:sinkhorn_definition}
\end{align}
where \(m(\cdot)\) measures the total mass of the density.
This corrects the transport vectors, and gives provably that \(\Sink(\mu, \mu) = 0\). In fact, \(\Sink\) is a pseudo-metric, being positive and definite\footnote{\(d(x, y)=0  \iff x=y \) for some function \(d\).}, though lacks the triangle inequality which is require for a true metric \citep{sejourne_et_al_2021}.
Further, since these terms are post-debiasing steps, each \(\UOT\) problem is solved separately, and for efficiency they may be solved in parallel. 

\subsubsection{Transport vectors}\label{section:barycentric_maps}
Alongside the scores (interchangeably the costs) returned from \(\Sink\) and \(\UOT\), there is the possibility to calculate transport vectors, which are found through the \textit{barycentric projection}.
These barycentric projections are approximations to transport maps.
They are formed through the joint density, and allow for the creation of transport vectors between the two input densities. 
That is, a vector from an initial point to the average location it is sending mass.
The definition of the barycentric projection is
\begin{align}
    X_{O,i} \in \Omega_{O} \rightarrow \frac{\sum_j \pi_{i,j}X_{F,i}}{\pi_{0,i}} \in \Omega_{F}, \\
    X_{F,j} \in \Omega_{F} \rightarrow \frac{\sum_i \pi_{i,j}X_{O,i}}{\pi_{1,j}} \in \Omega_{O}.
\end{align}
Notice, if the total mass sums to 1, so that $\pi_{i,j}$ defines a discrete joint probability distribution, they correspond to conditional expectations.

As stated, this quantity is not debiased, nor can \(\pi\) be debiased through its primal interpretation. 
However, \citet{feydy_sejourne_vialard_2018} showed that the gradient of \(\Sink\) in space is precisely the desired debiased quantity. Given the smoothness and differentiability of the half squared Euclidean cost, the gradient corresponds to the debiased barycentric projection, and is given as, 
\begin{align}
    X_{O,i} \in \Omega_{O} \rightarrow  \frac{\sum_j \pi_{i,j} X_{F,j}}{\pi_{0,i}}  - \left(X_{O,i} -  \frac{\sum_k \pi^{sym,O}_{i,k} X_{O,k}}{\pi^{sym,O}_{0,i}}\right), \\ 
    X_{F,j} \in \Omega_{F} \rightarrow  \frac{\sum_i \pi_{i,j} X_{O,i}}{\pi_{1,j}} - \left(X_{F,i} - \frac{\sum_k \pi^{sym,F}_{k,j} X_{F,k}}{\pi^{sym,F}_{1,j}}\right), 
\end{align}
where \(\pi^{sym, O/F}\) denotes the symmetric plans for solving the problems \(\UOT(\mu_O, \mu_O)\) and \(\UOT(\mu_F, \mu_F)\),  respectively.
We use the (debiased) barycentric projection to define transport vectors for each point. 
Figure \ref{fig:c1c1_biased_debiased} illustrates the significance of this debiasing to the transport vectors and is discussed in more detail in Section \ref{section:the_perfect_forecast}. 

Given the underlying transport vectors, the mean average transport magnitude and direction (ATM and ATD) is formed via averaging all vectors and calculating its magnitude and direction. 
Note that points with zero mass have undefined transport vectors and are left out of the averaging, and in general are ignored.
Moreover, like the OF based DAS verification methodology of \citet{keil_craig_2007, keil_craig_2009}, it is possible to calculate the approximate map in both a forward (observation to forecast) and inverse direction.
Both of these will be reported, with differences used as a diagnostic tool.

%%%%%%%%%%%%%%%%%%%%
\subsection{Implementation}\label{section:implmentation}
For more details on the Sinkhorn algorithm and implementation, see Supplementary Materials \ref{appendix:supplementary_material} or \citep{uot_own_implementation} for the Python implementation.
For our purposes, it is sufficient to know that \(\varepsilon \sim 1/\sqrt{N}\) and is fixed through the grid resolution to maintain convergence \citep{merigot_thibert_2020}. 
While efficiency was not the primary focus of this work, the Sinkhorn algorithm implementation was accelerated using two straightforward techniques: tensorisation (Remark 4.17 in \citet{peyre_cuturi_2018}) and \(\varepsilon\)-annealing  \citep{schmitzer_2019, chizat_2024}, which are detailed in Section \ref{appendix:implementation}.
With these optimisations the Sinkhorn algorithm for the ICP binary cases (40,000 points, Section \ref{section:binary_data}) required an average of 2.21 seconds, while the perturbed and Spring 2005 cases (301,101 points, Section \ref{section:intensity_data}) took on average 36.6 seconds.
These were ran on one NVIDIA GeForce RTX 3090 GPU, with a convergence tolerance of \(10^{-12}\). Note that the runtime depended on the type of penalisation used (e.g. TV regularisation was seen to be typically slower) and the chosen reach parameter.

Additionally, because physical data often exists on large scales and the underlying algorithm operates in exponential and logarithmic domains (which are generally more stable but can still result in numerical overflows and underflows on such large scales), a scaling is introduced to remove the dimensions of grid length and mass (or total intensity). 
Specifically, given a domain or region of interest, the typical length scale is defined as \(L\), and the total rain volume (or mass) is defined as \(M\).  
These dimensions can be reintroduced by simple multiplication of \(M \cdot L^2\), or the cost can be left in a dimensionless form. 

The length scale, \(L\), is typically straightforward to determine from the grid. 
However, the mass \(M\) may not be; it can be estimated from climatological data for the region or using a climatological average.
In this work, the suboptimal in-sample average is utilised, although an out-of-sample value would be preferable for future studies.
Throughout, for interpretability reasons the mass scaling is not reintroduced. 
To recover the actual cost, one would need to multiply the metric costs by \(M\).

This now provides us with a range of different scores to consider: the UOT energy in two flavours (\(\UOT^{KL}\) and \(\UOT^{TV}\)) and novelly each of their debiased Sinkhorn divergence counterparts (\(\Sink^{KL}\) and \(\Sink^{TV}\)), the debiased ATM and ATD in the forward and inverse direction (\(KL/KL^{-1}\) and \(TV/TV^{-1}\)), as well as the decomposition of the cost and various visualisations (including transport vectors, viewing the returned marginals, 2D histograms of the underlying transport vectors and decomposition of the costs into transport and marginal mass imbalance, e.g. Figures  \ref{fig:c1c1_biased_debiased},   \ref{fig:c1c3_rho_dependece_example},  \ref{fig:2dhist},  \ref{fig:perturbedspreadofcases}, respectively).
These are intended to illustrate the possibilities behind UOT, and demonstrate the capacity of this methodology.

%% file: content/3_data_description.tex
\section{Data Set Description}\label{section:data}

In this investigation, we examine UOT's behaviour through the binary geometric test cases of \citet{gilleland_et_al_2019}, the perturbed and Spring 2005 sets of \citet{ahijevychetal_2009, Kain_et_al_2008} and the core case (case 1) from the MesoVICT real test sets \citep{dorninger_et_al_2013}.
The two-step-intensity cases of \citet{ahijevychetal_2009} are covered in \citet{nishizawa_2024}.
Rather than investigating the test cases in chronological order of the studies, first the simple binary cases are examined, before looking at the inclusion of intensity and concluding with real textured precipitation forecasts.
In this way results are built up from the most simple cases.

\subsection{ICP Binary Geometric Cases}\label{section:binary_data}

The idealised geometric binary test cases were created by \citet{gilleland_et_al_2019}, under the MesoVICT project (Mesoscale Verification Intercomparison over Complex Terrain).
In the preliminary study, distance metrics were tested against various combinations of the 50 plus cases which make up the test set.
The chosen shapes are intended to replicate typical fields found over more complex terrain and at higher resolution, under a thresholding strategy.
This is where an event is defined as being above some threshold of precipitation or in a certain quantile of that forecast, then assigned one or zero if the event occurred, or not.
In this fashion, only the shape within the domain changes, not the intensities.
The benefit of studying these binary cases first is to distinguish amplitude (or intensity) errors from transport and shape or extent (e.g. spatial) errors.
 
These binary cases lie on a 200 by 200 regular grid, with no underlying physical units.
Since UOT operates in density space, it has the capacity to squeeze and stretch shapes and thus over- or under-fill grid points from their input intensities. 
With this in mind, a thresholded interpretation is not wholly appropriate. 
Instead, these are treated as examples of areal extents changing with, simple, uniform intensity within the shape bounds.
The full set of possible cases is illustrated in Figures 2, 3, 5, 7, 9, in \citet{gilleland_et_al_2019}, along with their interpretation.
Only a subset is considered here — see Table \ref{table:binary_cases_here}.

\begin{table}
\centering
\footnotesize

\begin{minipage}{0.48\linewidth}
    \centering
    \begin{tabular}{c|p{5cm}}
    \hline
    \textbf{Case} & \textbf{Description} \\ \hline
    C1 &  Radius 20 circle, centred at (100, 100) (base case) \\ 
    C2 & C1 centred at (140, 100) \\ 
    C3 & C1 centred at (180, 100) \\ 
    C4 & C1 centred at (140, 140)  \\ 
    C6 &  Two radius 20 circles centred at (100, 140) and (100, 60)\\ 
    C7 & C6 with lower circle shifted (40, \(\rightarrow\)) \\ 
    C8 &  C6 with lower circle shifted (80, \(\rightarrow\) )\\ 
    C9 &  Radius 60 circle, centred at (100, 100) \\ 
    C11 &  Union of C1, C3 and C4 \\ 
    % C12 & Two radius 20 circles centred at (120, 160) and (80, 40) \\
    C13 & One radius 8 circle centred at (75, 25) and one radius \(8/\sqrt{2}\) centred at (88, 180) \\ 
    C14 & C13 with upper circle shifted (25, \(\rightarrow\)) and lower circle shifted (50, \(\rightarrow\)) (mimicking mirroring in x=100)\\
    E1 & Vertical large ellipse \\ 
    E2 & E1 at \(45^{\circ}\) \\ 
    E3 & E1 at \(90^{\circ}\) \\ 
    E4 & E1 at \(135^{\circ}\) \\ 
   
    \end{tabular}
\end{minipage}
\hfill
\begin{minipage}{0.48\linewidth}
    \centering
    \begin{tabular}{c|p{5cm}}
    \hline
    \textbf{Case} & \textbf{Description} \\ \hline
     E6 & Small ellipse at \(45^{\circ}\) angle \\ 
    E7 & Small horizontal ellipse \\ 
    E9 & E1 shifted (25, \(\rightarrow\))\\ 
    E10 & E2 shifted (15, \(\rightarrow\)) and (20, \(\downarrow\)) \\ 
    E11 & E3 shifted (25, \(\downarrow\))\\ 
    % E12 & E4 shifted (15, \(\rightarrow\)) and (20, \(\uparrow\)) \\ 
    E14 & E6 shifted (8, \(\downarrow\)) and (10, \(\uparrow\))\\ 
    % E19 & Three vertical ovals centred at (100, 40), (100, 55),
% and (125, 75), scaled by (40, 5),
% (35, 5), and (25, 5), respectively \\ 
    % E20 & E19 smoothed using a disk kernel with radius 12 \\ 
    H1 & Compliment of C1 \\ 
    H2 & Compliment of C2 \\
    P2 & Full field, all ones. \\ 
%     P3 & One point at (1, 1) \\ 
%     P4 & One point at (200, 200) \\ 
%     P5 & One point at (100, 100) \\ 
%     P6 & Four points in the four corners\\ 
%     P7 & Four points near the middle of each boundary
% at (1, 100), (100, 1), (200, 100) and (100, 200) \\ 
    N1 & C1, but with 0.1\% frequency randomly generated
noise \\ 
    N2 & C4, but with 0.1\% frequency randomly generated
noise \\ 
    N3 & Union of C4 and P5 \\
    N4 & Union of C4 and P3 \\
    S1 &  Scattered events (frequency 5\%) inside a
circular envelope with radius 35 centred at (50, 100) \\ 
    S2 & As in S1, but a different realisation \\ 
    S3 & As in S1 and S2, but with the envelope shifted (100, \(\rightarrow\)) and a different realisation \\
    \end{tabular}
\end{minipage}
\caption{Table of idealised cases discussed in this report. For full description of all possible cases, see \citet{gilleland_et_al_2019}. The domain in a 200 by 200 grid, indexing from (1,1) to (200,200), using coordinates (horizontal, vertical). For the ellipse cases, the ratio of the major vs minor axis is 5:1 with dimensions of the ellipses either 100 by 20 (large ellipse) or 25 by 5 (small ellipse). All ellipses are centred at (100, 100) unless translated. 
Notation; C circles, E ellipses, P points, S scattered, N noisy and H hole cases.}\label{table:binary_cases_here}
\end{table}

%%%%%%%%%%%%%%%%%%%%%%%%%%%
\subsection{ICP Intensity Cases}\label{section:intensity_data}
The second dataset probes UOT's behaviour with more realistic textures and intensities.
Inherently, UOT operates with densities, therefore the framework remains unchanged. 
However, the data now carries true units of intensity.
These cases originate from \citet{ahijevychetal_2009}, which was the initial comparative study for spatial metrics under the ICP; the MesoVICT sets emerged a decade later.
This integral dataset consists of three subsets,  (i) two step idealised geometric intensity cases, (ii) a single real forecast over the US where the forecast-observation pairs are synthetic perturbations of the observation, and (iii) real forecast-observation pairs for 3 models in Spring 2005.
However, as mentioned, analysis of the geometric intensity cases, (i), is not repeated here.
(ii) and (iii) shall be referred to as the perturbed (notationally fake00X) cases, and the Spring 2005 cases, respectively.
 
All three datasets use the same 4km grid with size 601 by 501 (501 in Latitude, 601 in Longitude), where the real data corresponds to forecasts over the Rockies during Spring 2005.
For our exploratory investigation, the fields are placed on a regular uniform grid of the same size, i.e. a 605 by 501 box, rather than the true latitude-longitude grid.
This allows a more principled examination of the methodology and a true distinction of grid effects. 
Since OT allows for two different meshes to be compared, grid effects and imperviousness to the grids themselves will be explored in future work.
 
For a full description of the perturbed cases, consult Table 2 in \citet{marzban_et_al_2009} and there-in.
In brief, the observation is one realisation (fake000), then forecasts fake001-005 start with a 5 point shift south and 3 point shift east, where sequential cases double in both directions.
fake006 and fake007 are at the same location as fake003, but with 1.5 times (multiplicative) the intensity of fake003, and minus 0.05 inch (1.27 mm) (additive) the intensity, respectively. 
Importantly, all these cases are not balanced, because points shifted out of the domain are removed and new points set to zero.

The Spring 2005 cases consist of predictions from three models plus observations (stage II precipitation analysis) \citep{lin_mitchell_2005} of 60 min precipitation accumulations. 
The models are wrf4ncar, wrf4ncep, and wrf2caps, which are run at 4km, 4km and 2km respectively. 
They are then interpolated to the courser resolution 4km grid.
Experts were asked to subjectively score the models for each valid time. 
Here, as in \citet{keil_craig_2009}, 9 time stamps are presented and the authors' Table 3 is used for the average expert scores.
An important finding through their work was the lack of alignment between subjective scoring and traditional metrics.
In contrast, UOT methodology restores, on average, this association with subjective scoring, as shown in Section \ref{section:real_intensities}.
With the caveat that this is a small sample size and experts scores on average showed small difference, and thus were not strongly discerning of the best model.

%%%%%%%%%%%%%%%%%%%%%%%%%%%
\subsection{MesoVICT Core Case}\label{section:vera_core_case_data}

The last case compares two models against the Vienna Enhanced Resolution analysis (VERA) \citep{steinacker_et_al_2000} over the Alpine region during the period 20-22 June 2007. 
The two NWP models are the Swiss model COSMO-2 (at 2.2KM resolution, denoted CO2) and the Canadian high-resolution model GEM-LAM (at 2.5km resolution, denoted CMH).
Both models are interpolated onto the VERA grid at 8km.
Data is given at accumulation periods 1, 3, 6, and 12 hours, denoted: AC01, AC03, AC06, and AC12.
See \citet{dorninger_et_al_2013} for technical details, and Supplementary Materials \ref{appendix:real_texture} for more on data preparation.

This thorough testing aims to provide insight into how UOT performs with intensity and real textured rain, assessing its ability to diagnose pure transport alongside unequal weighting of regions of rain.
The evaluation also examines how much can be deciphered from the scores, with or without further visualisation, and whether model rankings align with the subjective scoring from the original ICP work.
 
Note all data used is available through \citet{ral_icp_webpage}.

%% file: content/4_results.tex
\section{Results and Analysis}\label{section:results}

Through-out the proceeding section, recall that the units of intensity are omitted to aid interpretability of the scores.
Hence, one may imagine scaling the scores by \(M\) if full units are desired.
For reference, for the binary cases \(L=200\), the textured cases, L=601, and the core case L=1000. Any cases which are not presented in the following section but are presented in \citet{gilleland_et_al_2019} can be found in Section \ref{appendix:extra_figures} and Section \ref{appendix:all_cases}.

Notationally, comparisons between two fields, e.g. C1 and C2, are written as C1C2 (or fake000fake001 for the intensity cases). This means that C1 is treated as the observation field and C2 the forecast.
For each case the UOT score or costs, in KL and TV flavours, are calculated. Both of these are also considered in their biased (p=\(\UOT\)) and debiased (se=\(\Sink\)) forms.
We write \(\UOT^{KL/TV}\) and \(\Sink^{KL/TV}\), whilst UOT is used when referring to the general methodology.
Additionally, the mean average transport magnitude (ATM) and direction (ATD) are reported in both flavours and in the forward (KL/TV) or inverse (KL\(^{-1}\)/TV\(^{-1}\)) directions, i.e. observation to forecast (C1 \(\rightarrow\) C2) or vice versa. 
The UOT\(_{\varepsilon}\) cost (and \(\Sink\)) reported throughout is equation \eqref{eq:uot_general_cost} multiplied by \(\rho\).

\subsection{ICP Binary Geometric cases}\label{section:bindary_results}

See Table \ref{table:binary_cases_here} for the cases used in this section, or \citet{gilleland_et_al_2019} for their full description and interpretation.
For these binary cases \(\varepsilon =0.005L^2\), following  \citet{merigot_thibert_2020}. This is set by the grid structure, not the user.
Recall that the units of the objective costs are in (half) squared Euclidean space, however the ATM is in units of distance (grid points). For example, C1C2 (Figure \ref{fig:rotation_and_orientation}) has an expected transport cost of \(800\), which corresponds to a 40-point shift, as \(\frac{1}{2} \cdot 40^2 = 800\). 
This cost is then scaled by the relative total rainfall/areal extent, or mass. The average mass is \(M=1873.5\), and C1 has a mass of \(1345\), so in Figure \ref{fig:rotation_and_orientation}, the cost is \(574 \sim 800 \cdot 1345 / 1873.5\).

Rather than jumping straight to unbalanced translated cases, intuition is first built up for the reader. First dealing with the prefect forecast, accounting for boundary effects and position, checking resistance to the double penalty problem, rigorously probing the behaviours of the reach parameter, and investigating rotation and scale, all in the balanced setting first. Only then are the unbalanced cases investigated.

\subsubsection{The perfect forecast}\label{section:the_perfect_forecast}

\begin{figure}[h!]
    \centering
    \includegraphics[width=0.8\linewidth, ]{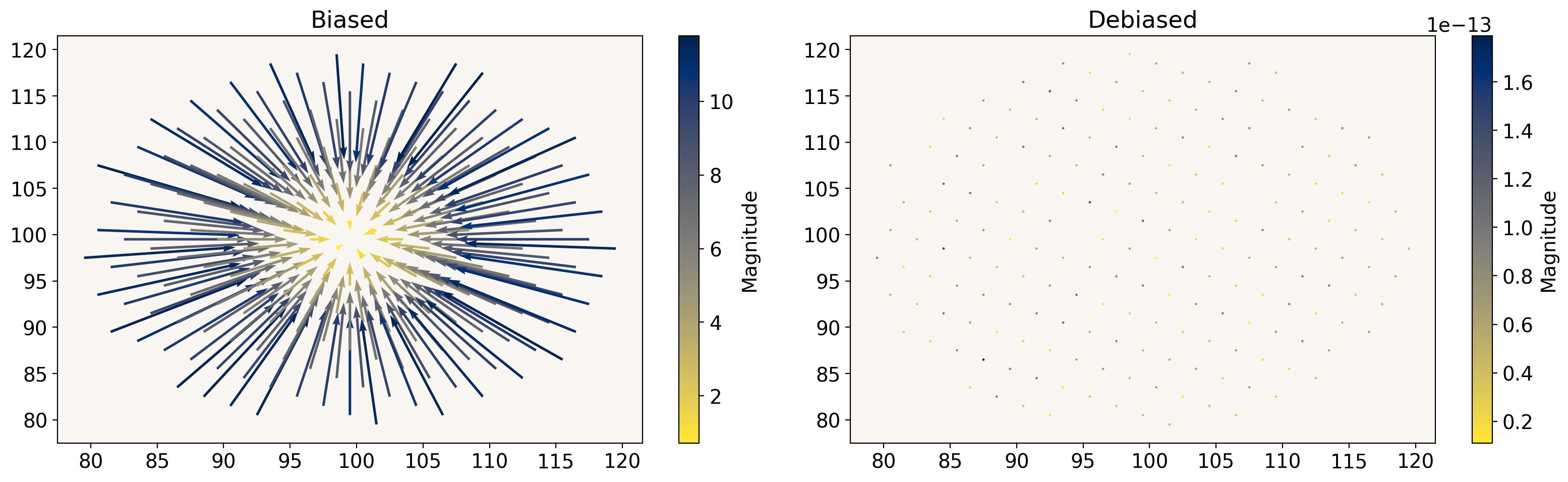}
    \caption{Biased vs debiased transport vector illustration with KL marginal penalty for the C1C1 case which is expected to have no transport (cropped to region of interest). The TV flavour similarly suffers from bias. Left: Biased UOT transport vectors which cause a contraction of the support, Right: Debiased UOT transport vectors. A regular sample of the vectors are shown to prevent overcrowding. Notice the scale of the magnitude bar for the debiased cases. \(\varepsilon = 0.005L^2,\ \rho = L^2\).}
    \label{fig:c1c1_biased_debiased}
\end{figure}

By considering a perfect forecast-observation pair, C1C1, this section highlights the importance of debiasing and responds to \ref{q1}. In Figures \ref{fig:c1c1_biased_debiased} and  \ref{fig:translation}, \(\Sink^{KL/TV}\) returns zero, whilst \(\UOT\) is dominated by entropic error as expected by \citet{sejourne_et_al_2021}. 
% This highlights \(\Sink\) as a pseudo-metric.
Figure \ref{fig:c1c1_biased_debiased} illustrates how this bias distorts the transport vectors, and how through debiasing this error is corrected. 
Now, the ATM remains zero, though the ATD may exhibit small non-zero values from machine precision.

\subsubsection{Boundary and relative position effects}\label{section:boundary_position_effect}

\begin{figure}[h!]
    \centering
    \includegraphics[width=1.0\linewidth, trim=70 70 70 70, clip ]{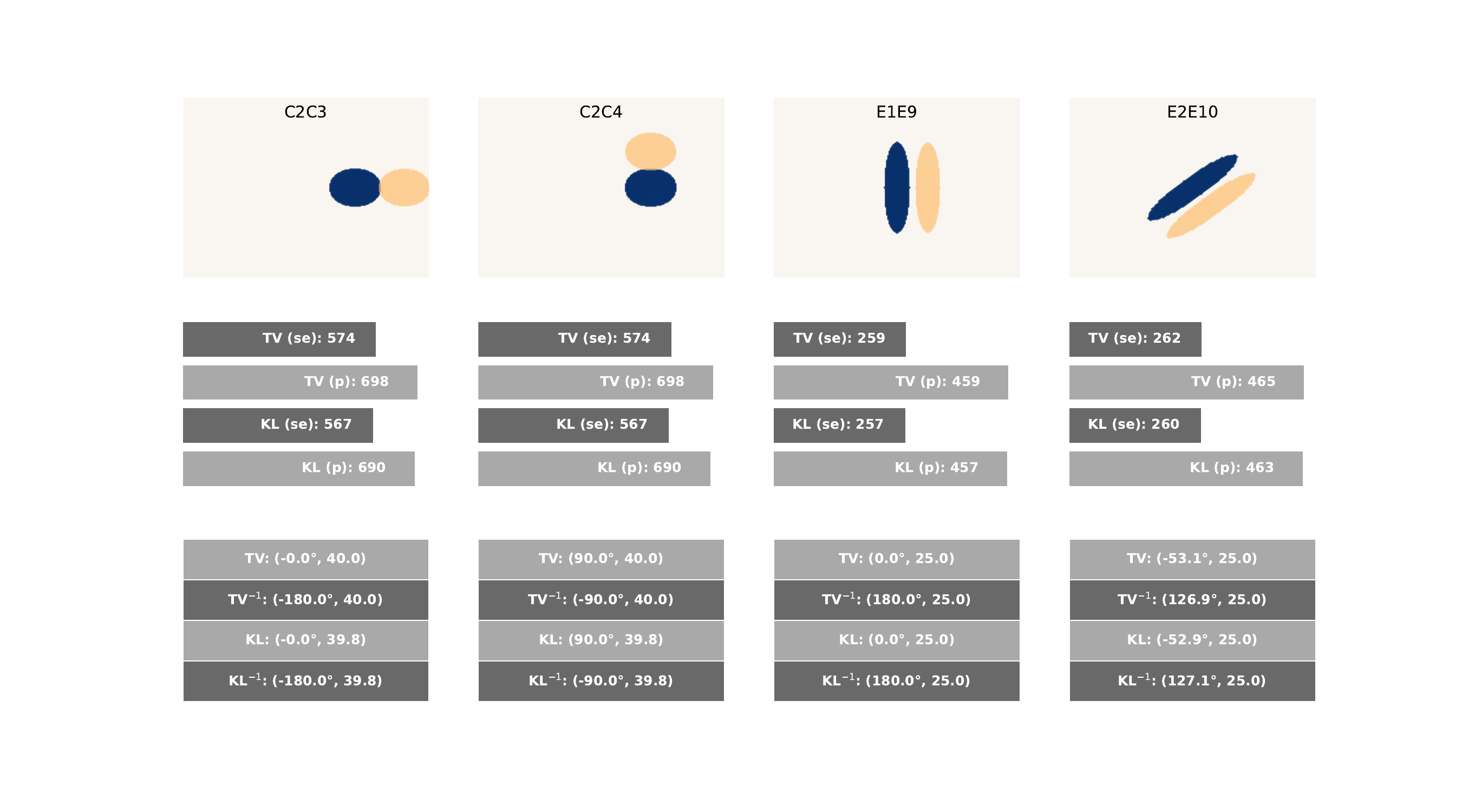}
    \caption{Cases illustrating robustness to orientation and rotation of events on a regular grid. Notice that for E2E10 the -53.1 degree ATD (and its inverse) appears from these events being 15 grid points east and 20 south which preserves the correct separation after rotation of the observation. The top four horizontal bars display; \(\Sink^{TV}, \UOT^{TV}, \Sink^{KL}, \UOT^{KL}\). The lower table presents the mean (ATD, ATM) in both flavours, and with the forward and inverse debiased vectors. The blue (darker) colour indicates observations, while the pale orange (lighter) represents forecasts. \(\varepsilon = 0.005L^2,\ \rho = L^2\).}
    \label{fig:rotation_and_orientation}
\end{figure}

The next question investigates whether UOT is sensitive to proximity to the boundary, or features' relative position and orientation (addressing \ref{q2}).
Good behaviour would be insensitive to these factors.
Through C1C2, C2C3 and C2C4 (Figures \ref{fig:rotation_and_orientation} and \ref{fig:translation}) all three cases have equal scores and the ATD shows the expected transport direction.
Simultaneously, cases C2C3 vs C2C4 and E1E9 vs E2E10 (Figure \ref{fig:rotation_and_orientation}) show that pairs at \(90^{\circ}\) rotation are equal, with the expected ATD.
However,  E1E9 (1553 points) and E2E10 (1573 points) differ slightly in \(\Sink\), although their ATM and ATD is correct.
This sensitivity is a discretisation effect on the regular grid. 
Despite this, UOT demonstrates robustness to the boundary, which is not true for all distance metrics (see Baddeley’s metric in \citet{gilleland_et_al_2019}), and some sensitivity to rotation on the regular grid, but this is not detrimental for equivalent shifts.

\subsubsection{Sensitivity to phase error for balanced cases}\label{section:phase_error}
\begin{figure}[h!]
    \centering
    \includegraphics[width=1.0\linewidth, trim=70 70 70 70, clip]{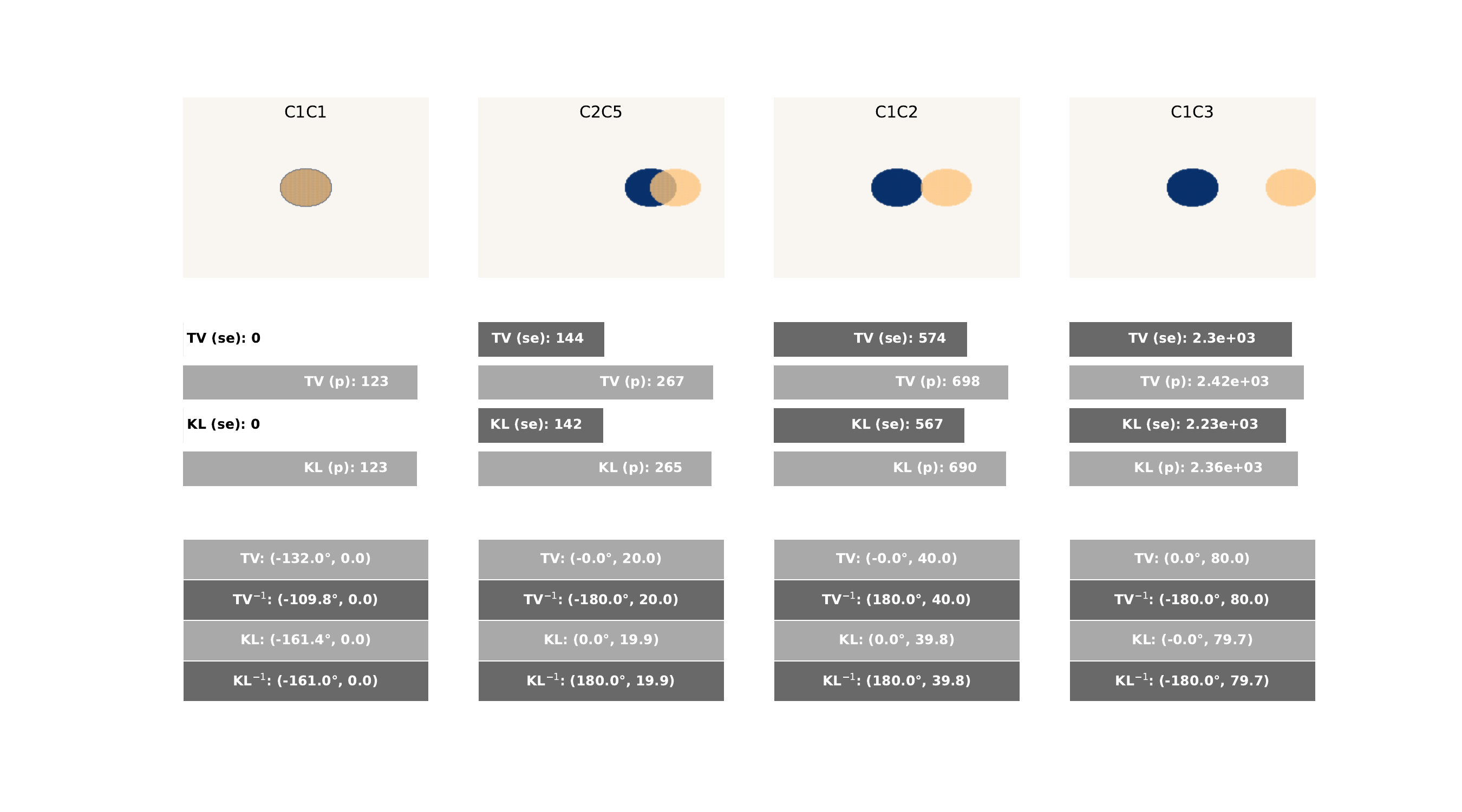}
    \caption{Cases which illustrate the similarity to \(\frac{1}{2}W^2_2\), and thus robust to the double penalty problem. Notice the quadratic translation relationship in the cost function (doubling the shift quadruples the cost). The top four horizontal bars display; \(\Sink^{TV}, \UOT^{TV}, \Sink^{KL}, \UOT^{KL}\). The lower table presents the mean (ATD, ATM) in both flavours, and with the forward and inverse vectors. The blue (darker) colour indicates observations, while the pale orange (lighter) represents forecasts. \(\varepsilon = 0.005L^2,\ \rho = L^2\). }
    \label{fig:translation}
\end{figure}

Through the  cases C1C1, C1C2, C1C5, and C1C3 (Figure \ref{fig:translation}), transport and sensitivity to the double penalty problem is scrutinised by comparison to the known Wasserstein-2 (\(\frac{1}{2}W^2_2\)) distance which corresponds to the balanced, non-entropic version of UOT\(_{\varepsilon}\).
This answers \ref{q1}, as \(\frac{1}{2}W^2_2\) is known to diagnose the shift error and optimal cost corresponding to this shift \citep{benamou_2003}.
It is expected, and shown, that the methodology is capable of diagnosing this translation error.
In Figure \ref{fig:translation}, \(\Sink\) correctly diagnoses only the transport error, maintaining a quadratic cost scaling with the shifted distance.
This aligns more with \(\frac{1}{2}W^2_2\) in \(\Sink^{TV}\), than in \(\Sink^{KL}\), although it will be shown to be \(\rho\) and mass-balance dependent. 
The ATM is correctly diagnosing the shift. 
Hence \(\Sink\) and the ATM is not susceptible to the double penalty problem in the balanced mass setting. 
 
\begin{figure}[h!]
    \centering

    \begin{subfigure}[t]{0.48\textwidth}
        \centering
 \includegraphics[width=0.8\textwidth, trim=40 30 50 50,clip]{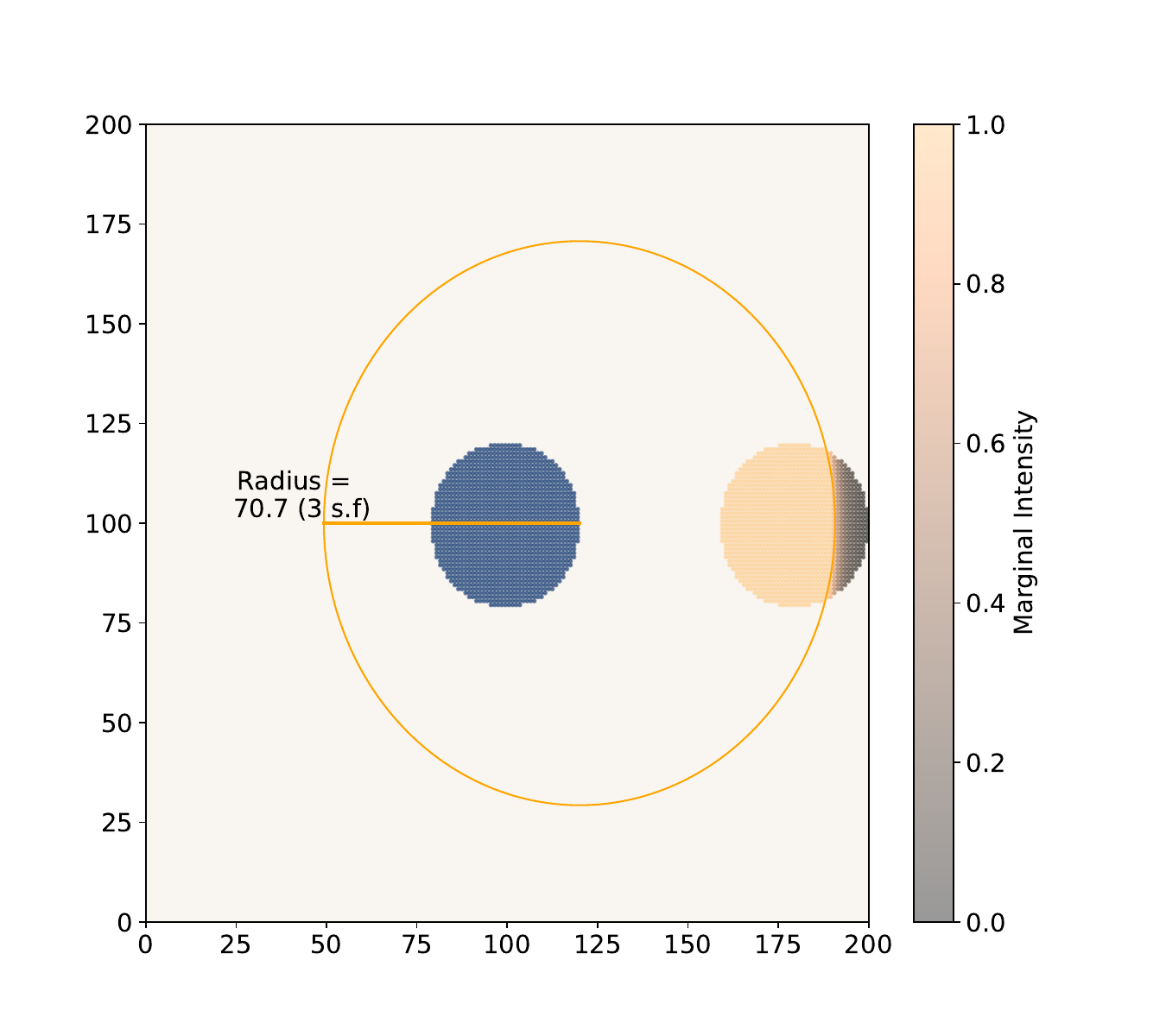}
        % \label{fig:tv_rho_dependece_example}
    \end{subfigure}
    \quad 
    \begin{subfigure}[t]{0.48\textwidth}
        \centering
        \includegraphics[width=0.8\textwidth, trim=40 30 50 50,clip]{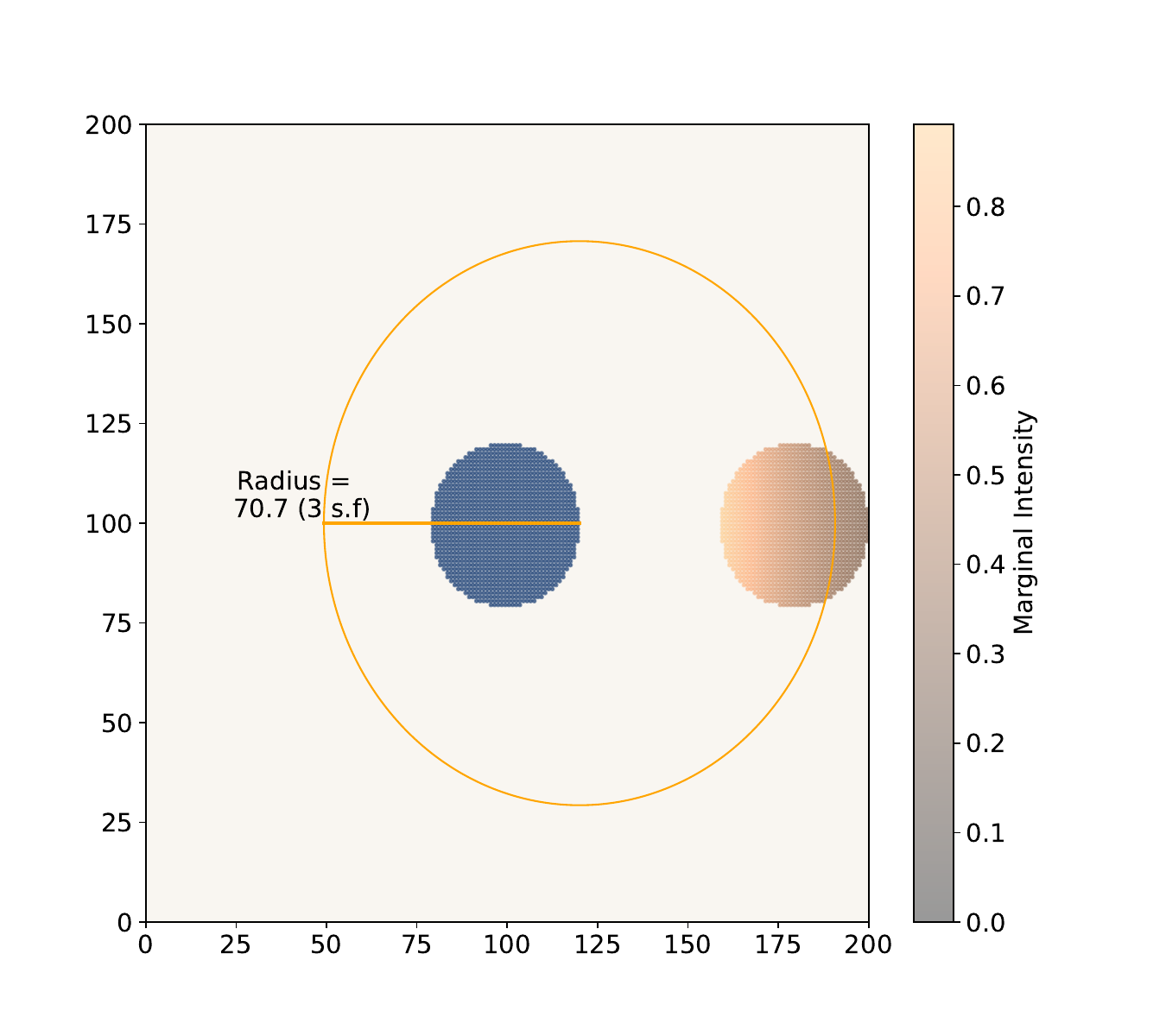}
        % \label{fig:kl_rho_dependece_example}
    \end{subfigure}
    \caption{ Illustration of the dependence on the TV and KL penalties on \(\rho\) for case C1C3. In both figures, the central observation C1 is shown in blue (darker), while the forecast is shifted 80 points to the right (C3). Potted in pale orange to grey gradient is the returned marginal, \(\pi_1\). This shows where mass has been paid. Left: the TV penalty, mass too far from the target is destroyed with a strict dependence on the reach parameter, resulting in a 14.4\% reduction in total intensity. Right: the KL penalty displays a smoother dependence on the reach, although with a 45.9\% reduction in total intensity. The radius represents the reach for \(\rho = 2^{-4} L^2\), beyond which it is cheaper to destroy mass. Both figures were generated using \(\rho= 2^{-4} L^2\) and \(\varepsilon=0.005 L^2\).
}\label{fig:c1c3_rho_dependece_example}

\end{figure}

\begin{figure}[h!]
    \centering
    \includegraphics[width=0.9\linewidth]{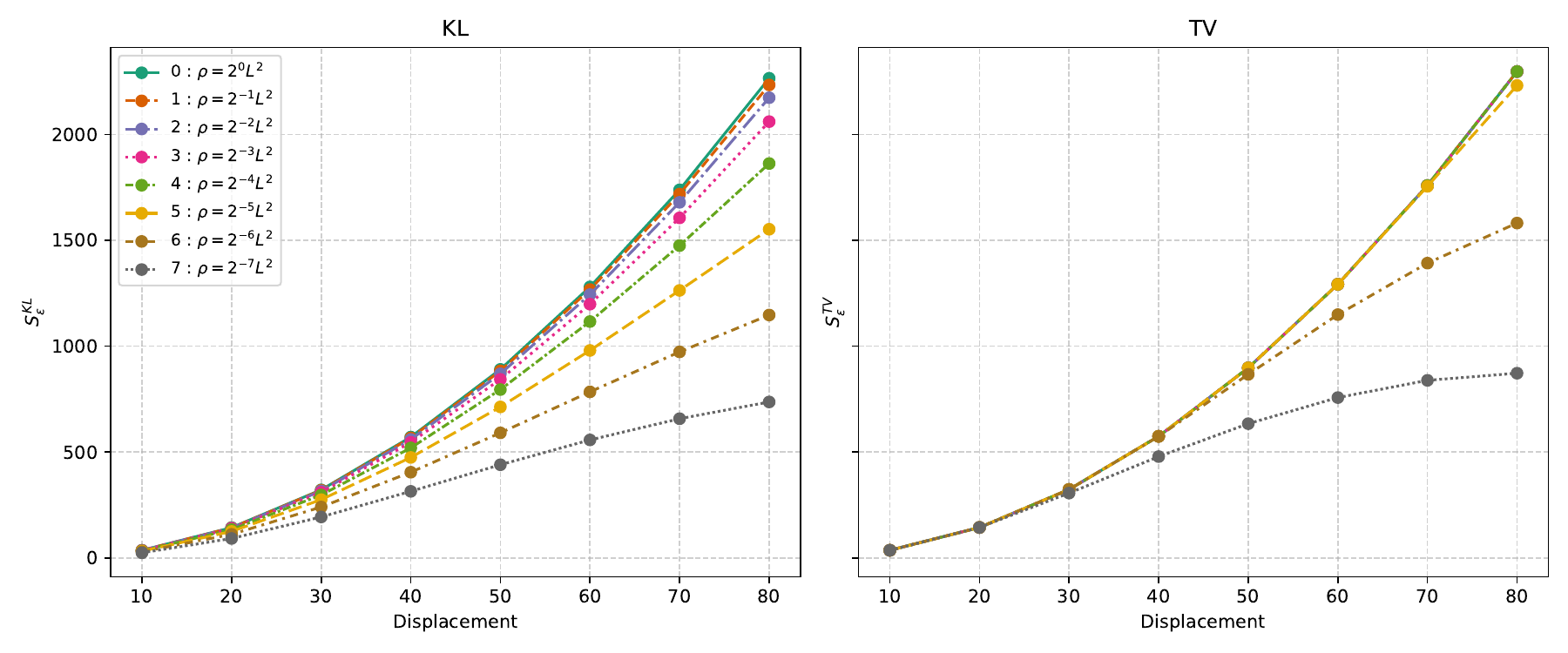}
    \caption{ Reach sensitivity study for case C1 vs C1 displaced east. Notice that a quadratic relation is maintained in the displacement for sufficiently large \(\rho\). However, this changes as \(\rho\) decreases, and it becomes cheaper to destroy mass. \(\Sink^{TV}\) is more robust with a sharp drop-off in the reach dependence. While \(\Sink^{KL}\) is a smoother transition.  Left: \(\Sink^{KL}\), Right: \(\Sink^{TV}\). \(\varepsilon = 0.005 L^2 \) .}\label{fig:transport_vs_rho}
\end{figure}

To test how the reach parameter effects our results, answering \ref{q4} and \ref{q5},  consider cases at 10-grid-point increment shifts\footnote{All shifts are described from the centre of mass, and east from C1. These new cases supplement C1C1, C2C5/C3C5, C1C2/C2C3, C1C5, C1C3 which are at 20 increments.} and across 8 values of \(\rho =  2^{-i}L^2 : i \in \{-6, \ldots, 1\}\).
Figure \ref{fig:transport_vs_rho} illustrates the returned \(\Sink^{TV/KL}\).
For reach \(\sim L\) (\(\, \rho\sim L^2\)) the expected trend in the shift is maintained.
However, as \(\rho\) decreases, allowing more relaxation in the marginals, it becomes cheaper to destroy mass and change the marginals than to pay the transport cost. Figure \ref{fig:8case_2} similarly shows for both flavours at low reach it becomes cheaper to destroy mass rather than transport, correcting  only local (small) displacement.
In fact, in this simple setting, one can predict at what value of \(\rho\) the cost will change.
Consider the 80 point separation case (C1C3).
The pure transport cost would be \(\frac{1}{2}(80)^2 \cdot 1345/1873.5 \sim 2300 \) and the value at which the \(\Sink^{TV}\) cost changes from this is between \(\rho =2^{-3}L^2\) and \(\rho=2^{-4}L^2\).
That is, when the reach goes from 100 to \(\sim70\) which is below the 80 point separation. Figure \ref{fig:c1c3_rho_dependece_example} visualises this transition through the returned marginals.
It demonstrates that \(\Sink^{TV}\) is more robust in the reach parameter (destroying mass with a sharp edge at the limit of the reach), whilst KL is a smoother penalty term.
In summary, for the balanced scenarios, \(\Sink\) reliably diagnoses transport error within the reach radius. 
Distant features can be removed at a  reduced cost depending on how \(\rho\) is set.
Marginal visualisations can clarify which regions remain associated.

Investigation of the balanced setting with multiple features present, demonstrates cases with multiple associated features, e.g. one-hit-one-miss scenarios. The results follow a similar trend to the above discussion.
However, each feature has its own reach radius; a notion we hypothesis to follow human evaluation of a forecast. See Supplementary Materials \ref{appendix:multiple_features}.

\subsubsection{Rotation and scale for balanced cases}\label{section:rotation_scale}

\begin{figure}[h]
    \centering
    \includegraphics[width=0.75\linewidth]{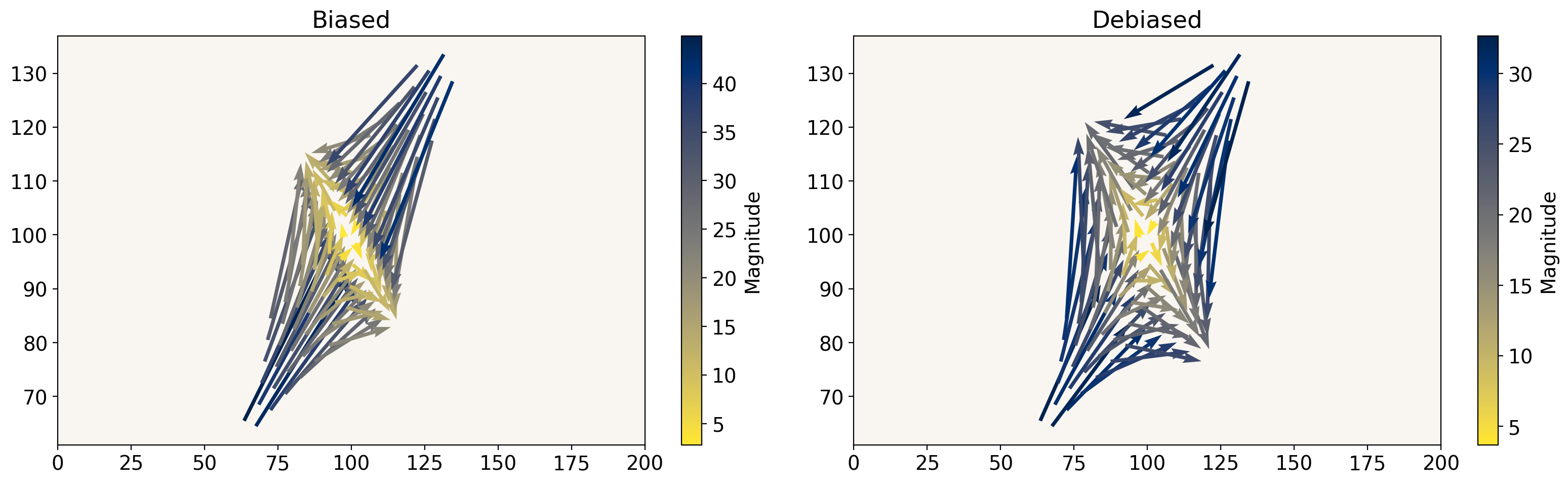}
    \caption{Biased vs debiased transport vector illustration with KL marginal penalty for the E2E4 case. Here the squeezing and stretching behaviour for rotated objects is demonstrated. Rather than rotation it shows aspect-ratio correction. Left: Biased UOT transport vectors, Right: Debiased UOT transport vectors. A regular sample of the vectors are shown to prevent overcrowding. \(\varepsilon = 0.005L^2,\ \rho = L^2\).}
    \label{fig:e2e4_transportvectors_kl}
\end{figure}

A good attribute for a spatial verification methodology would be to diagnose rotational error, e.g. in MODE or IWS.
However, balanced cases E1E4 vs E2E4 (Figure \ref{fig:rotation_ellipse}) highlight a limitation of UOT, since it is not built to rotate objects (responding to \ref{q2}).
By construction, the transport plan will stretch and squeeze the shapes to match (Figure \ref{fig:e2e4_transportvectors_kl} and Figure \ref{fig:e2e4_transportvectors_tv}), mimicking the behaviour of quasi-rotational (or aspect ratio) error, similar to geom004 \citep{ahijevychetal_2009}.
Though we notice, \(\Sink^{KL}\) scores the 45 degrees rotation better than the 90 degree (unlike \(\Sink^{TV}\)), thus potentially following a more subjective evaluation.
The ATM and ATD demonstrate a downfall in taking the average due to the symmetry of the stretching and squeezing which leads to cancellation. 
Figure \ref{fig:e2e4_transportvectors_kl} shows this effect. Notably, unlike CDST, \(\Sink^{TV/KL}\) is able to discern the two different rotations.
The relative scale of events, investigated through E6E14 and E2E10 (Figure  \ref{fig:scaled_Cases}), follows the expected size of events, larger events scoring higher if there is translational error present. More details are given in Section \ref{appendix:scale_rotation}.

\subsubsection{Unbalanced circle scenarios}\label{section:unbalance_circles}

Changing focus to the unbalanced setting, where cases have different total rain coverage, or mass. Cases C1C6 and C2C11 against C1C2  (Figure \ref{fig:unbalanced_extent} and  \ref{fig:translation}) demonstrate over- and under- forecasting, yet without overlap and deals with \ref{q6}. 
A method should differentiate, and correctly rank these, unlike say the Hausdorff distance and MED.
UOT is shown to rank them in the expected order, with C2C11 (worse), C1C6, C2C1 (best). 
\(\Sink^{TV}\) returned more severe penalties.
C1C6 returns the expected ATM with equal and opposite transport, while C2C11’s opposing east-west shifts cancel, leaving north-south transport at one-third of C1C2 due to spatial bias (with three times the mass).
This means that a forecaster is  encouraged to match total rain coverage.
% In contrast, the Hausdorff distance and MED failed to distinguish these cases.
 
To answer \ref{q4} in the unbalanced setting, the reach parameter's influence is examined through cases C1C6, C1C7, C1C8 in Figure \ref{fig:unbalanced_displcement_rho} (also Figures \ref{fig:unbalanced_extent} and  \ref{fig:unbalanced_reach}).
The score should, and does, degrade as the lower event in the observations moves further away.
Again,  \(\Sink^{TV}\) is more harsh, being an order of magnitude larger than \(\Sink^{KL}\).
Noticeably, in comparison to Figure \ref{fig:transport_vs_rho}, the dependence is rotated in \(\rho\) and displacement, i.e. larger \(\rho\) enforces more balanced mass which is expensive in the unbalanced setting and this dominates the cost.
However there is still an increase in the cost as the southern event becomes more displaced.
Fortunately, the ATM picks up on the increased transport required. 
In fact, closer examination shows that the KL flavour favours the closer event, which appears as a skewed ATD (Figure \ref{fig:unbalanced_reach}).
 
\begin{figure}[h!]
    \centering

    \begin{subfigure}[t]{0.48\textwidth}
        \centering
        \includegraphics[width=0.8\textwidth, trim= 40 30 50 50, clip]{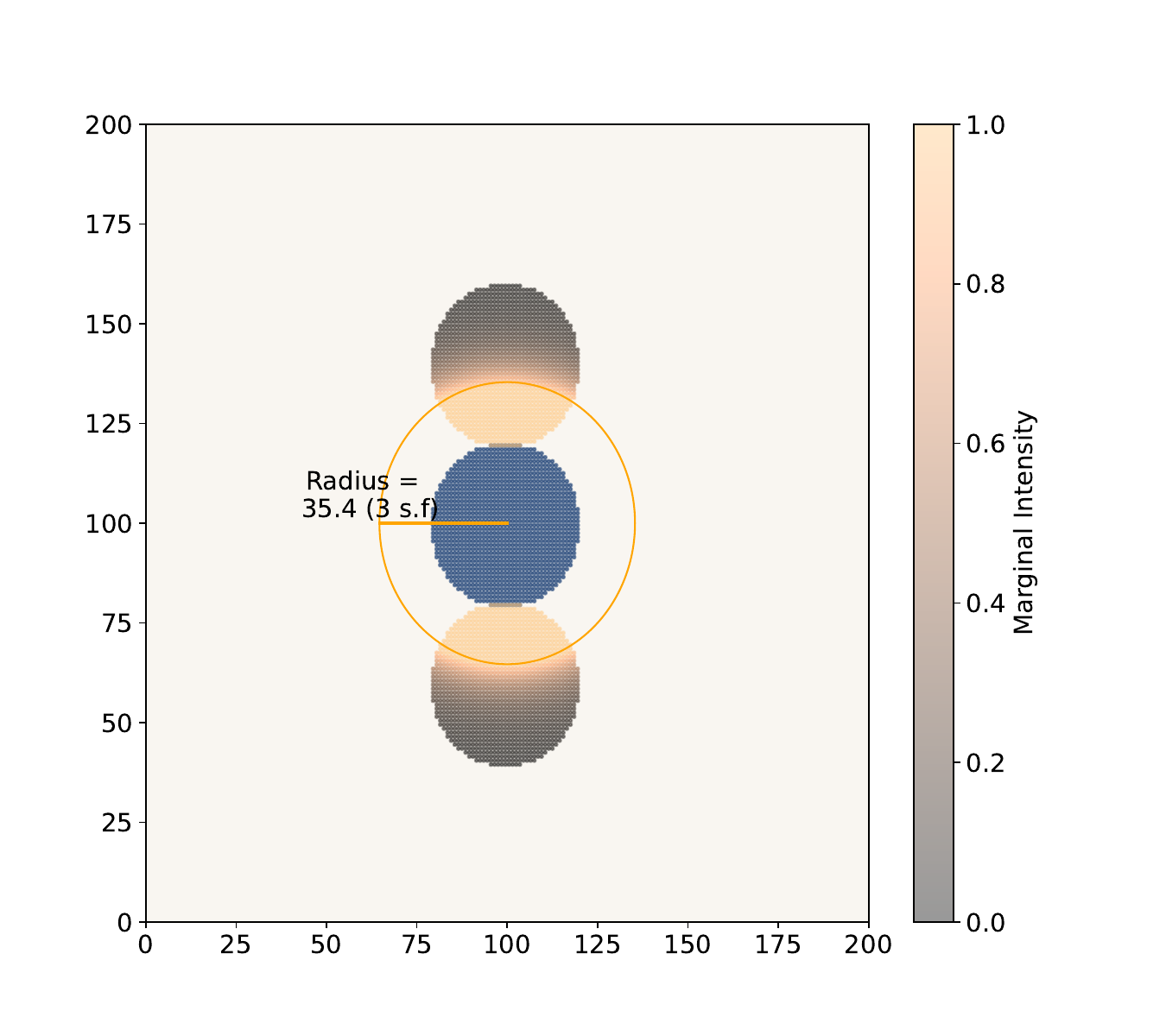}
    \end{subfigure}
    \quad 
    \begin{subfigure}[t]{0.48\textwidth}
        \centering
        \includegraphics[width=0.8\textwidth, trim= 40 30 50 50, clip]{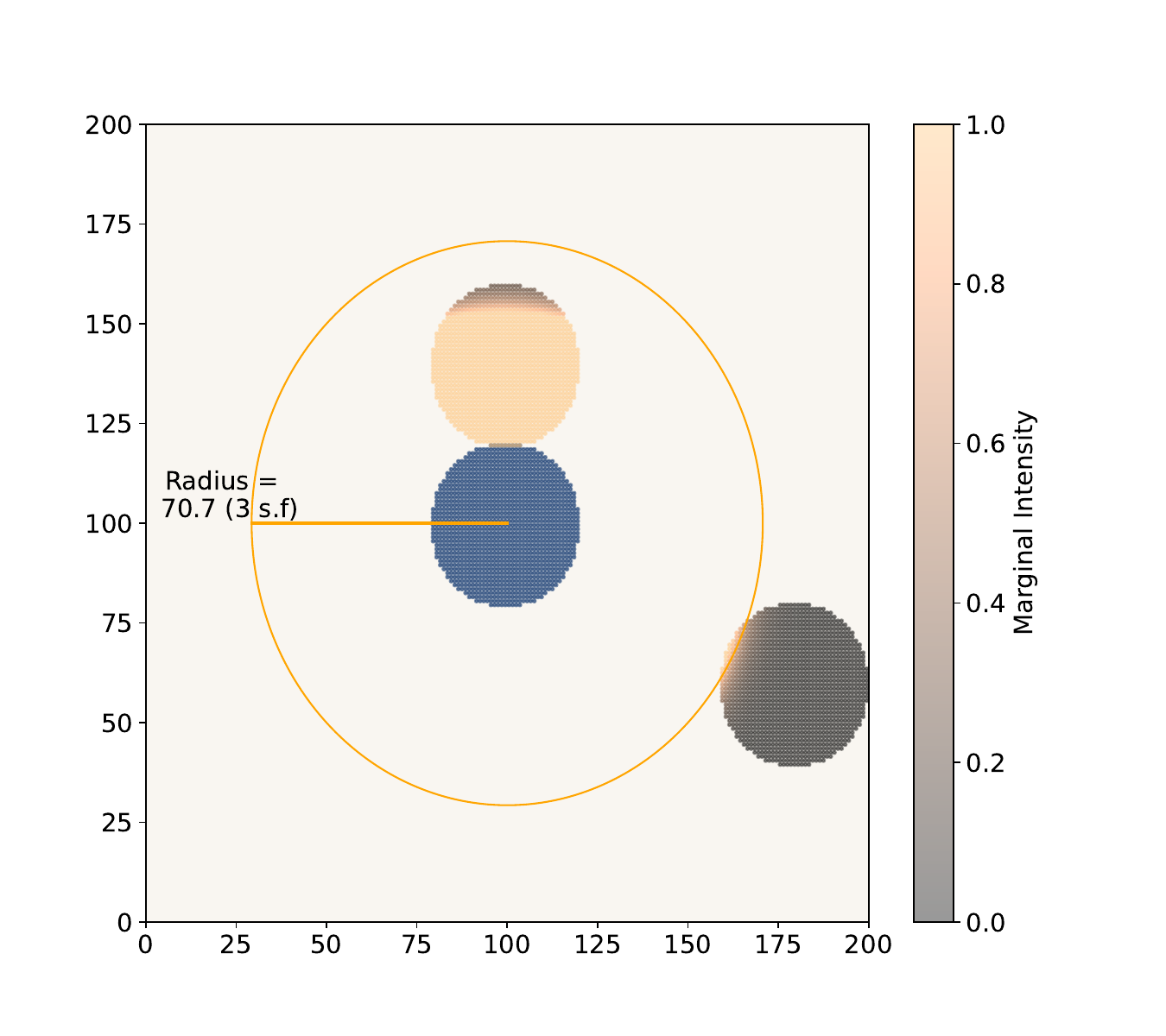}
    \end{subfigure}
    \caption{Illustration of unbalanced displacement and \(\rho\) sensitivity with TV penalty. In both figures, the central observation C1 is shown in blue (darker), while the forecast is C6 (left) or C8 (right).
    Marginals \(\pi_1\) are plotted in pale orange to grey shades, not the forecast, highlighting regions where mass is destroyed. Left: there is again a strict relationship to the reach parameter although the mass is shared within half the reach. Hence the radius shown is equal to half the reach. Right: the closer event is assigned more intensity whilst the distant feature in this unbalanced setting is given very little. The radius is equal to the reach.  The figures were generated for \(\rho= 2^{-4} L^2, \varepsilon=0.005L^2\).}
    \label{fig:comparison_tv_marginal}
\end{figure}

The marginals (Figure \ref{fig:comparison_tv_marginal} in TV and  \ref{fig:comparison_kl_marginal} in KL) show that for the symmetric case of C1C6 mass is distributed evenly within half the reach radius (answering \ref{q3} and \ref{q5}).
When we break this symmetry, closer features receive more mass, indicating that while the reach’s linear interpretation breaks down in the unbalanced setting, its preference for proximity remains. 

\subsubsection{Unbalanced ellipses and hedging}\label{section:unbalanced_ellipses}

\begin{figure}[h!]
    \centering
    \includegraphics[width=0.75\linewidth, trim=70 70 70 70, clip]{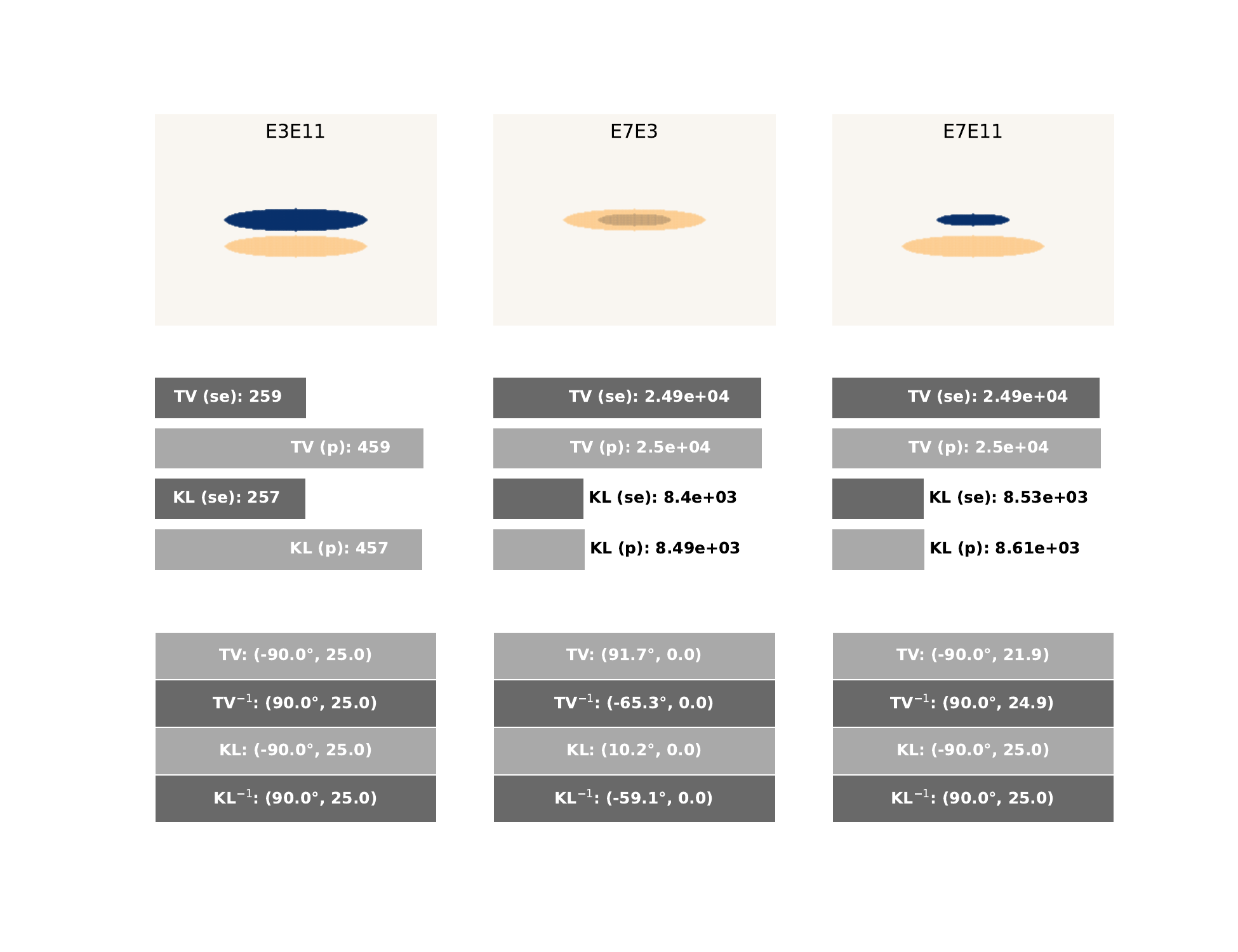}
    \caption{Comparison of subjective cases, exploring both balanced and unbalanced settings. Note that between  E7E3 and E7E11 the TV scores' do increase, but below the shown significant figures. The top four horizontal bars display: \(\Sink^{TV}, \UOT^{TV}, \Sink^{KL}, \UOT^{KL}\). The lower table presents the mean (ATD, ATM) in both flavours, and with the forward and inverse vectors. The blue (darker) colour indicates observations, while the pale orange (lighter) represents forecasts. \(\varepsilon = 0.005L^2,\ \rho = L^2\).}
    \label{fig:hedging_forecast}
\end{figure}

Through the comparison of cases E3E11, E7E3, E7E11 (Figure \ref{fig:hedging_forecast}) 
common scenarios with a model over- or under- forecasting  areal extent, with overlap, can be tested.
They also address \ref{q3}, \ref{q5}, \ref{q6}, and \ref{q9}, probing how transport and spatial bias interplay, along with subjective assessment as different users may favour cases differently.
UOT scores them across flavours in ascending order: E3E11, E7E3, E7E11\footnote{For the TV flavours this discrepancy is below the shown significant figures.}.
E3E11 follows the expected \(\frac{1}{2}W_2^2\) cost, whereas the penalty paid for over-/under forecasting, a form of mass imbalance, dominates the other two  cases. 
Compellingly there is only a small difference in scores between E7E3, E7E11.
Through the cost decomposition in Figure  \ref{fig:ellipse_cost_distintergartion_across_rho} and Section  \ref{appendix:decomposition}, the extra transport in E7E11 is diagnosable.
% Here the transport cost and marginal mass imbalance terms, \ref{eq:uot_general_cost}, can be distinguished, i.e. .
We revisit these diagrams later for the perturbed cases.
This provides utility as forecasters are not encouraged to increase overlap to hedge scores.
Instead, matching shape and extent reduces spatial bias and is a preferred strategy for minimising the UOT scores. 
This is in opposition to the distance metrics considered in \citet{gilleland_et_al_2019} which could all be hedged.

This analysis comes with the caveat that the forecast value is a subjective assessment, since some users may prefer E7E3 or E7E11.
If this were the case, and absolute transport were more important to a user than intensity and spatial bias, then this user could use the ATM . 
Now, E7E3 is favoured since less transport is required, (note the average has cancelled with equal and opposite transport vectors\footnote{Some readers may be interested to know that if E3 is given the same total intensity as E7, so that it becomes a balanced case again. 
Then E7E3 obtains a lower score due to less transport than E3E11.}).

Finally, a limitation of \(\Sink\) is that it cannot discern subsets (\ref{q2}), a useful skill for preventing hedging and understanding feature location. See \ref{appendix:subset_section}.
This is in contrast to MED which is able to diagnose if an event was a subset of another.
Until now only simple circles and ellipses have been considered, however a large benefit of the idealised cases is the opportunity to showcase a methodology's capacity (or deficiency) to cope with non-trivial fields.

\subsubsection{Noisy, scattered and hole cases}\label{section:scatterholenoisy}

Through cases S1S2, S1S3 (Figure \ref{fig:scattered_hole}), we address \ref{q5} and \ref{q7}, and expose how a methodology copes with random scattered events contained within an envelope, imitating the occurrence of a collection of small showers.
Note that the events do not have balanced mass.
A good score should still diagnose the spatial shift in S1S3, and return small error for S1S2 as they share an envelope.
For S1S2, \(\Sink^{KL}\) correctly returns a low score.
Conversely, \(\Sink^{TV}\) penalises these events harshly, which from the above analysis is due to mass imbalance.
The ATM, in both flavours, indicates 1.3 on average (2.1 median average), aligning with MED and the FSS, although these were able to diagnose closer to the average nearest neighbour transport distance between the two realisations, which is 2.35 \citep{gilleland_et_al_2019}.
This trend continues in S1S3, where it scores worse than S1S2 due to a shift between the scattered events, and the ATM diagnoses the shift.
In  summary, \(\Sink^{KL}\) is the more robust flavour to total rain fall imbalance, and with greater capacity to cope with random scatter events.

Noisy cases C1N3, C1N4, N1N2 versus C1C4 (Figure \ref{fig:noise}) reflects that observations may not be clean, and spurious points may remain even after cleaning.
The score is expected to be insensitive to small noise, and to degrade as the noise increases.
Without noise, C1C4 returns the expected cost.
After introducing noise, both flavours remain close to the known cost, ATM, and ATD, with \(\Sink^{KL}\) being more robust as expected. 
Yet both are still capable of diagnosing transport in the presence of a small amount of noise.
Notice for N1N2 there is still degradation of the score, below the shown signifiant figures, however this small amount of noise has not dominated the cost. 
The full ordering, from best to worst, for \(\Sink\) is C1C4, C1N3, C1N4, N1N2, hence follows intuition.

Finally, case H1H2 against C1C2 demonstrates that the methodologies does reward the increased spatial match. However the spatial displacement is equal, and this is seen in the ATM. For more details for these cases see \ref{appendix:hole_edge}.

\subsubsection{Edge and extreme cases}\label{section:edge}

Before moving onto real intensity forecasts and to answer \ref{q7}, attention is required for edge or extreme cases, including: very little to no precipitation or full precipitation over the whole domain.
How a score should respond to these cases is not always clear, though P1P1, P2P2 should and do return zero (Figure \ref{fig:new_new_paper_pcase_1}). However, P1 or P2 against another forecast is less clear.
Purely from the optimisation interpretation, the null case has only one feasible plan that has a null marginal — the zero plan.
In this scenario, all mass is destroyed in the non-zero field, hence the \(\Sink^{TV}/\UOT^{TV}\) can be defined as the total mass of the non-zero field (which can be zero). Whilst \(\Sink^{KL}/\UOT^{KL}\) is zero .
In neither flavour do the ATD and ATM make sense with null fields.
Additional, simple extremal cases are explored in \ref{appendix:hole_edge}.

To aid comparison with distance metrics, recalling \(\Sink\) is a pseudo-metric, we provide our evaluation of the criteria given in \citet{gilleland_2021} Table 3. 
1) \textit{Good pathological handling}, given there is some available mass then scores are defined. Null and full cases are appropriated scored. Extreme imbalance however may cause issues with numerical stability. 
2) \textit{No position effects}, this was explored in section \ref{section:boundary_position_effect} and UOT is insensitive to the relative position to boundary and rotation on the regular grid. 
3) \textit{Sensitivity to frequency bias}, is also answered in the affirmative, since  C1C2, C2C6 and C2C11 score in the best to worse order. 
4) \textit{Useful for rare events}, considering C13C14, UOT scores this intuitively, since given sufficiently large reach the ATM scored the average transport, while the \(\Sink\) is relatively small. Moreover, illustrations of the marginals show dissociation for a sufficiently small reach parameter. Although, we provide the same caveat as \citet{gilleland_2021}; that results are dependent on \(\rho\). 
5) \textit{Reward partial perfect match}, cases with one-hit-one-miss score intuitively again. C6C7 vs C1C2, Figure \ref{fig:rho_balanced_reach}. Only the miss is transported. The mean ATM does suffer here due to zero average transport for the hit, yet carrying half the mass. However, UOT is able to provide zero scores for a complete hit.
6) \textit{Correct penalisation despite partial match}, asked if we can differentiate between C1C9 and C1C10 — which we do (Figures  \ref{fig:subsets} and \ref{fig:new_new_paper_circles_0}).

To conclude the binary geometric shapes, UOT has shown insensitivity to the relative boundaries position. 
It can diagnose pure transport depending on the reach parameter, and hence it is robust to the double penalty problem.
It is crucial that the debiased Sinkhorn divergence and transport plans are used, though decomposition of the cost may offer deeper diagnostics. 
See Figure \ref{fig:cost_disintergration}.
Markedly, UOT prevents hedging and does not prioritise overlapping events.
Instead, matching spatial extent is scored better by the metric. 
UOT is not able to pick up on purely rotational errors and fails to provide useful information under certain extreme cases, such as the null cases.
Moreover, the cost itself is not diagnostic of subsets.
In general, \(\Sink^{KL}\) appears to offer a more robust penalisation in the unbalanced setting.
However, if (like \citet{skok_2023}) both fields are normalised to the same mass, then \(\Sink^{TV}\) offers attractive features, especially when considering the link to the reach and the returned marginals of the plan.
After this thorough exploration in a binary setting, let us turn to the real intensity cases.

\subsection{Real Intensity Cases}\label{section:real_intensities}
Given the observations above, the natural proceeding question is if these behaviours and properties extend into real intensity textures and shapes.

\subsubsection{Fake (Perturbed) cases}\label{section:perturbed_fake_cases}
Through these perturbed cases questions \ref{q1}, \ref{q3}, \ref{q4}, \ref{q5} are examined, re-confirming previous trends found in the simple cases.
\ref{q8} is then answered showing various visualisations are often required to separate intensity and transport errors.
The average total accumulation, M, is chosen as the total mass of fake000 (356100.0) and \(\varepsilon=0.001 L^2\).
As in the binary cases, the \(\Sink^{KL/TV}\) score is expected to quadruple between cases fake001-fake005 — irrespective of the intensity now being taken into consideration — and the ATM and ATD should follow the known displacement (with the caveat that these cases are not perfectly balanced). 
Briefly, we note that fake000fake000 achieved zero \(\Sink\) and ATM as expected.

\begin{figure}[h!]
    \centering
    \includegraphics[width=0.9\linewidth]{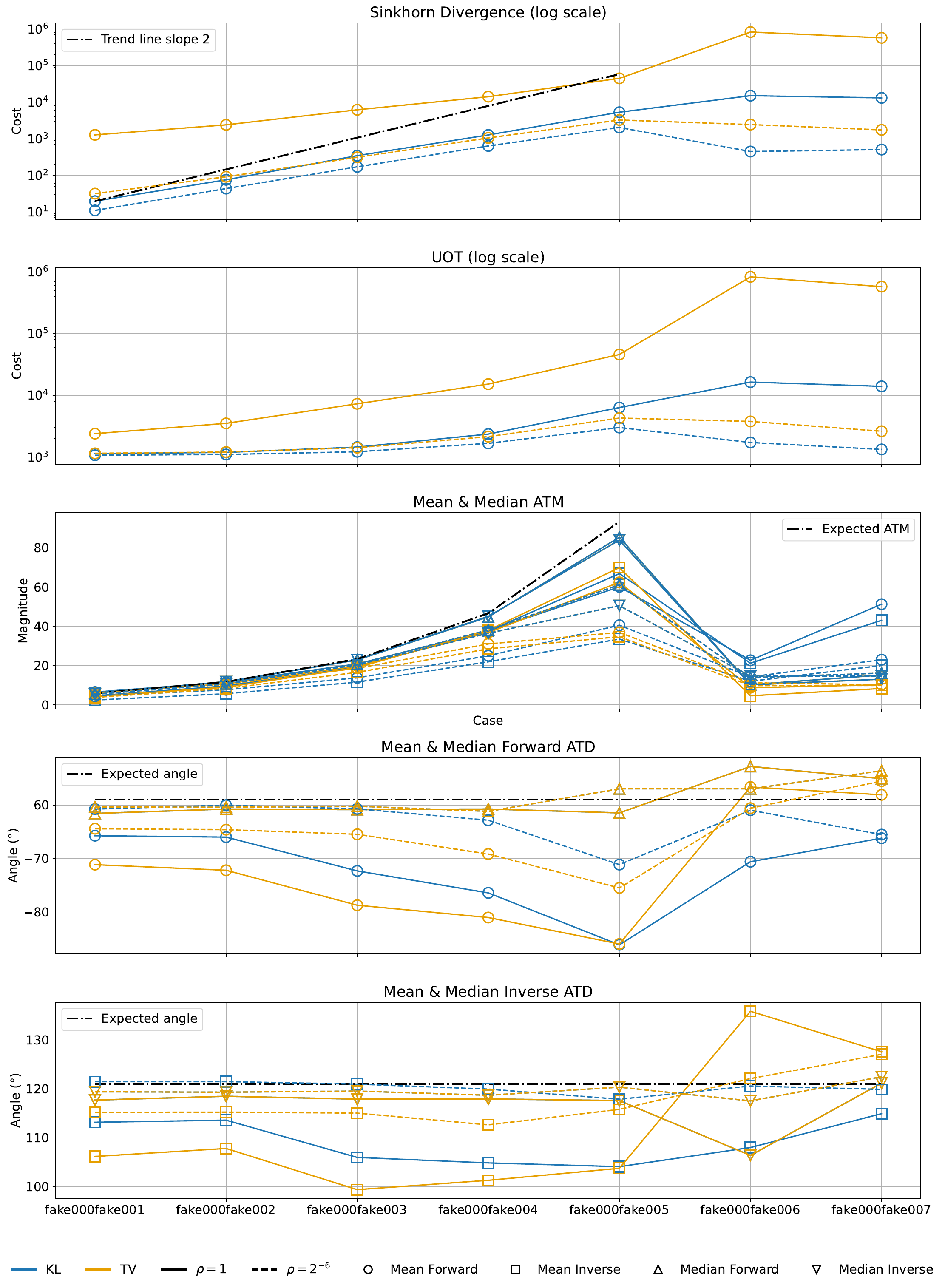}
    \caption{Perturbed cases with real textured intensities. All cases are a changed version of fake000. 
    From top of bottom is displayed, the \(\Sink^{TV/KL}\), \(\UOT^{TV/KL}\), ATM in both flavours, directions and by taking the mean and median, the forward ATD (mean and median, and both flavours) and inverse ATD (mean and median, and both flavours). 
    The key at the bottom of the figure differentiates the values of \(\rho\) tested, flavours, and direction. \(\Sink\) and \(\UOT\) do not have a direction, so are both plotted with circles.\(\varepsilon = 0.001L^2\). the expected ATM and ATD are plotted, whilst a trend line corresponding to quadratic growth is plotted for  \(\Sink^{TV/KL}\).
    }
    \label{fig:all_perturbed_cases}
\end{figure}

Firstly, consider fake000 against fake001-fake005, Figure \ref{fig:all_perturbed_cases} (with a detail breakdown in Figures \ref{fig:perturbed_cases},  \ref{fig:perturbed_cases_smallrho})  where the expected \(\Sink\) score for fake000fake001 is approximately \( 16.9\sim 0.5\cdot(3^2 + 5^2) \cdot 354859 / 356100\) (the exact cost is not known as they are unbalanced).
\(\Sink^{KL}\) is appearing more robust, however neither score maintains the quadratic relation well, due to changing mass imbalance.
As previous seen, \(\Sink^{TV}\) is overly penalising the forecasts, because the marginal penalty terms dominate the costing (Figure \ref{fig:perturbedspreadofcases}). 
This sensitivity is reduced for a smaller reach parameter as transport is kept local.

For both flavours the ATM starts close to the expected values, before diverging at larger translations. 
However through the median average the ATM stays closer to the expected value (breakdown in Figures \ref{fig:fake_cases_median_rho1} and  \ref{fig:fake_cases_median_rhominus6}).
Moreover, for smaller reach and displacement the ATM is closer again to the known perturbation, though at the larger shifts this diagnostic fails.
The reach is too small.
The ATD similarly improves for smaller reach and by taking the median average.
This reconfirmed \(\Sink^{KL/TV}\)s' expected behaviour, yet with non-uniform intensities.

\begin{figure}[h!]
    \centering
    \includegraphics[width=0.98\linewidth, trim= 0 0 0 70,clip]{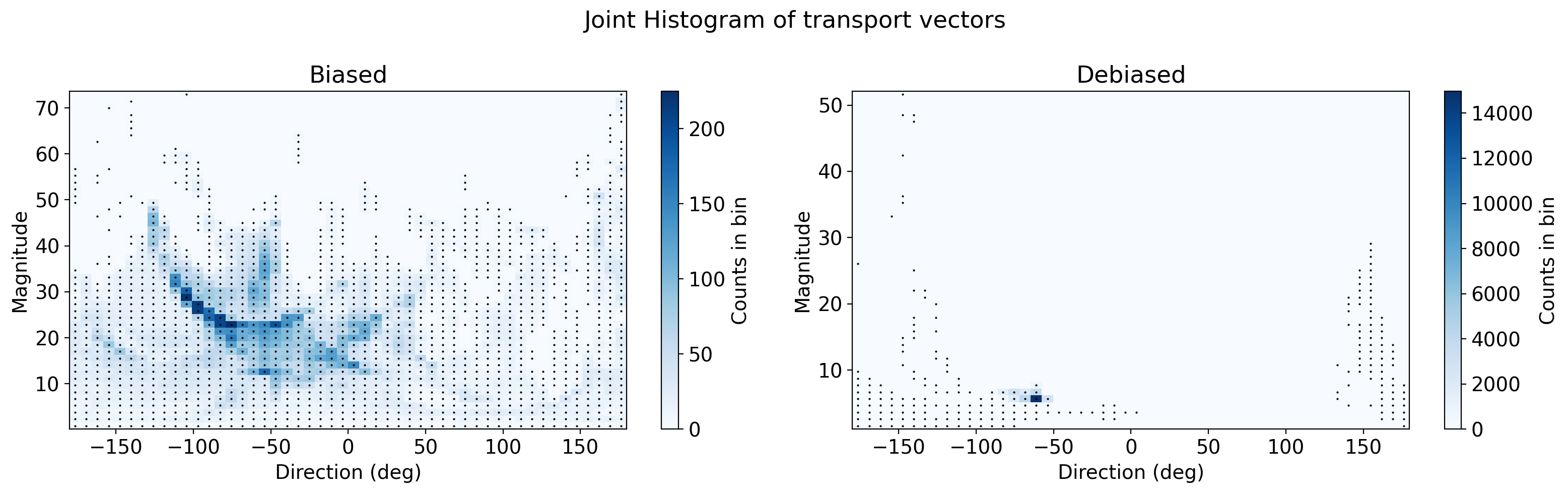}
    \caption{2D histogram of the magnitude and direction of the underlying transports vectors for the perturbed case fake000fake001, which has expected transport cost at \(-59^{\circ}\) and \(5.8\) points. The dotted bins, indicate those with non-zero mass. Left: biased transport vectors, which we know contract the support. Left: Debiased transport vector distribution which concentrates on the correct direction and magnitude. These results are based on TV with parameters \(\varepsilon=0.001L^2,\ \rho=L^2\).}
    \label{fig:2dhist}
\end{figure}

\begin{figure}
    \centering
    \begin{minipage}{0.5\linewidth}
        \centering
        \includegraphics[width=\linewidth]{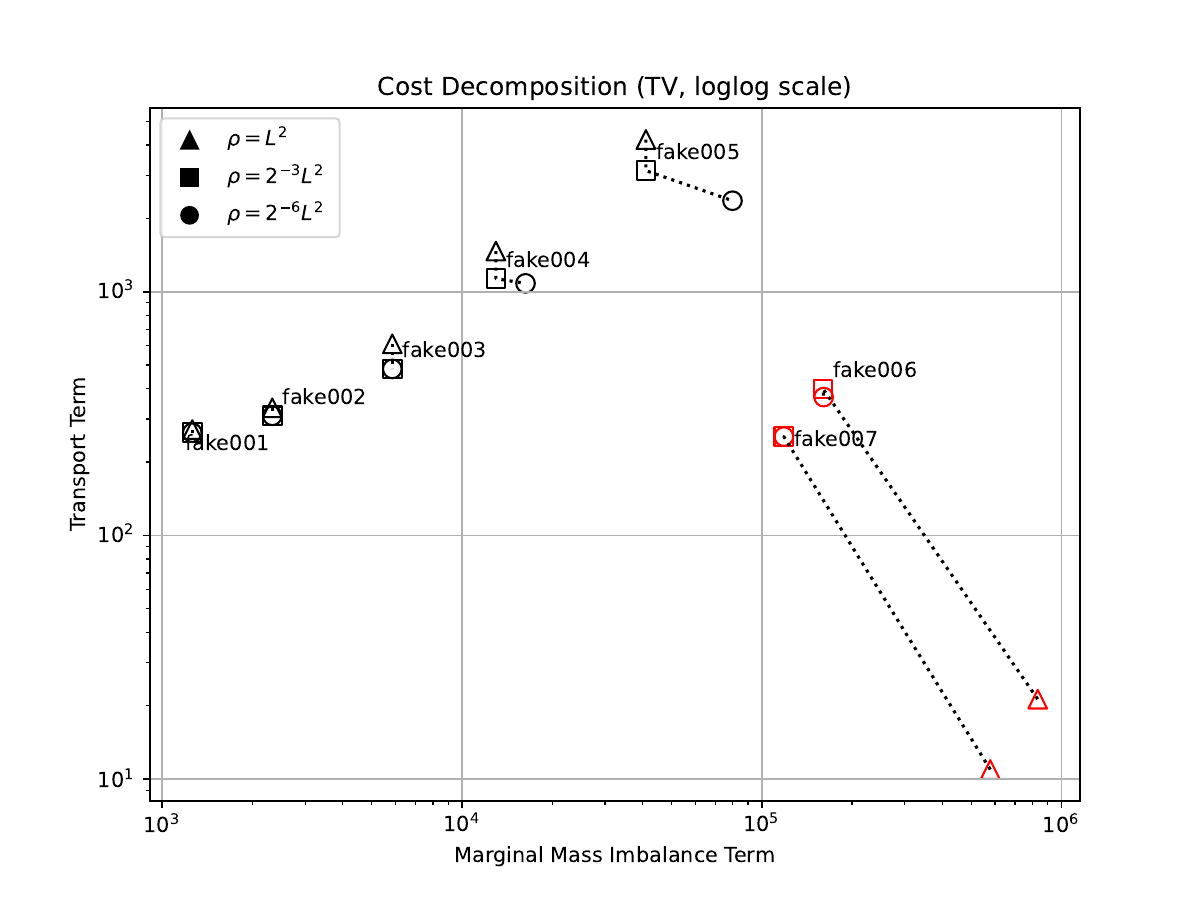}
    \end{minipage}%
    \begin{minipage}{0.5\linewidth}
        \centering
        \includegraphics[width=\linewidth]{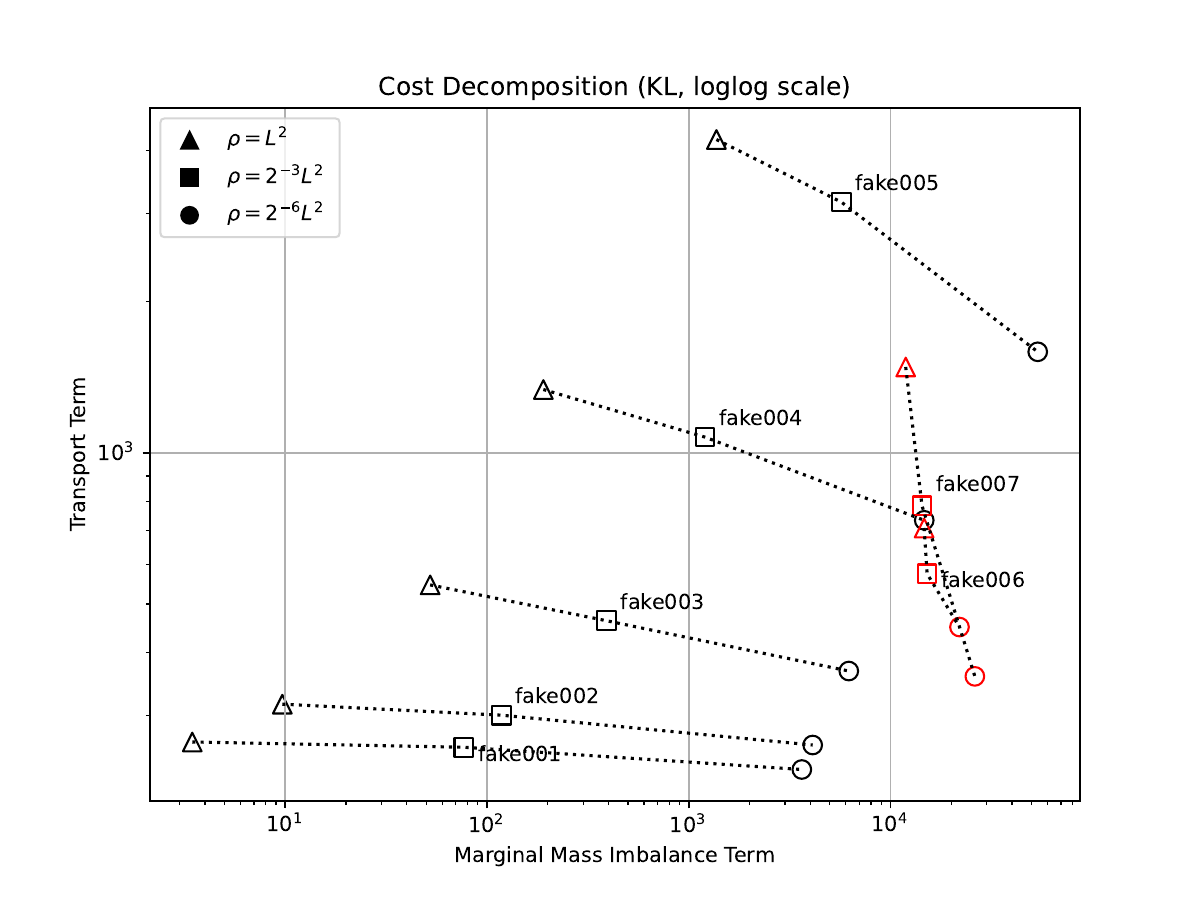}
    \end{minipage}
    \caption{
    Decomposition of cost, for perturbed case, into the transport and marginal mass imbalance term, or sum of the marginal penalty terms. These correspond to the first (without \(\rho\)  and  the regularisation term) and last two terms  of \ref{eq:uot_general_cost}. Left: TV penalty, Right: KL penalty. Whilst strictly all cases are unbalanced, fake006/7 have been coloured in red since they are the most extremely unbalanced.}
    \label{fig:perturbedspreadofcases}
\end{figure}

A further visualisation  following \citet{marzban_et_al_2009}, is a  2D histogram of the underlying transport vectors' direction and magnitude. This is a more powerful diagnostic tool than simple summary metrics like ATM and ATD which have been shown to oversimplify some cases.
Figure \ref{fig:2dhist}, illustrates that debiasing corrects the vectors  and concentrates them towards the true ATM, ATD; this is true with both penalties.
Two observations are key. First, the debiased case is able to successfully diagnose the correct angle and magnitude from the approximate mapping yet now with the inclusion of real intensities. Second, that debiasing, and hence the Sinkhorn divergence, is crucial to be able to draw conclusions.
Both have been highlighted throughout the above investigation.

The last of the perturbed cases are fake006 and fake007, which should be discernable from fake003.
The ATD maintains its trends being close to the expected angle, however the ATM is not well-preserved nor does it appear to scale with the change in mass (\(\sim +44\%, - 32\%\) respectively).
All the scores are penalised more than fake000fake003, and depending on the reach (and thus importance of transport) and flavour they score differently.
Whilst this may seem like complex behaviour, it allows selection of priorities by picking the appropriate formulation.
% That is, if imbalance is allowed and above a certain reach the forecast should be penalised more, by setting the parameter appropriately one can order these cases differently. 
Notably, fake006 always costs more than fake007 for TV penalisation which follows the size of mass imbalance.
However, KL penalisation changes with the reach, and for smaller reach (allowing cheaper mass destruction) fake006 can cost less than fake007. 

Examining these cases deeper, the decomposition in Figure \ref{fig:perturbedspreadofcases} does allow for clear trends in transport and mass balance where fake006/007 stand out since as \(\rho\) decreases, the transport cost goes up for KL and down for TV.  For KL, as we decrease the reach, more mass is allowed to be destroyed and hence there is less transport. 
Counter-intuitively for TV UOT demands mass balance in an unbalanced setting and the transport goes down. 
We hypothesis this to be to a bias effect, recalling the transport cost itself cannot be debiased.
In viewing the transport vectors there is an increase in the spread but not magnitude (Figure \ref{fig:fake007_2d_histograms}).
Unfortunately, from the scores alone the direction of imbalance is not diagnosable, i.e. the mass of fake006 is greater than fake000 which is greater than fake007.
However, by comparing terms \(D(\pi_0 | \mu_{O})\) and \(D(\pi_1 | \mu_{F})\) there is an indicated direction. Figure \ref{fig:dualdecomposition_0067} illustrates that when TV must destroy mass  it will do so to satisfy one constrain fully (giving one constraint approximately zero value, \(\sim \)1e-13). 
KL follows similar though less extreme trends.
This facilitates diagnosis of over and under forecasting through the ratio \(D(\pi_0 | \mu_{O})\)  relative to  \(D(\pi_1 | \mu_{F})\).  
For under forecasted fake007, \(D(\pi_0 | \mu_{O})\) >  \(D(\pi_1 | \mu_{F})\) since we have sufficient mass to fill the smaller  \(D(\pi_1 | \mu_{F})\), and vice versa. This ratio we call the marginal mass imbalance ratio, and if greater than one a forecaster is under-forecasting and below one over-forecasting. 
For the MesoVICT core case we explore the ratios use instead of their sum. 

Crucially, with real intensities and textured data there is still the possibility to  diagnose pure transport error and be robust against the double penalty problem.

\subsubsection{Spring 2005 cases}\label{section:spring}

\begin{figure}[h!]
    \centering
    \includegraphics[width=0.98\linewidth]{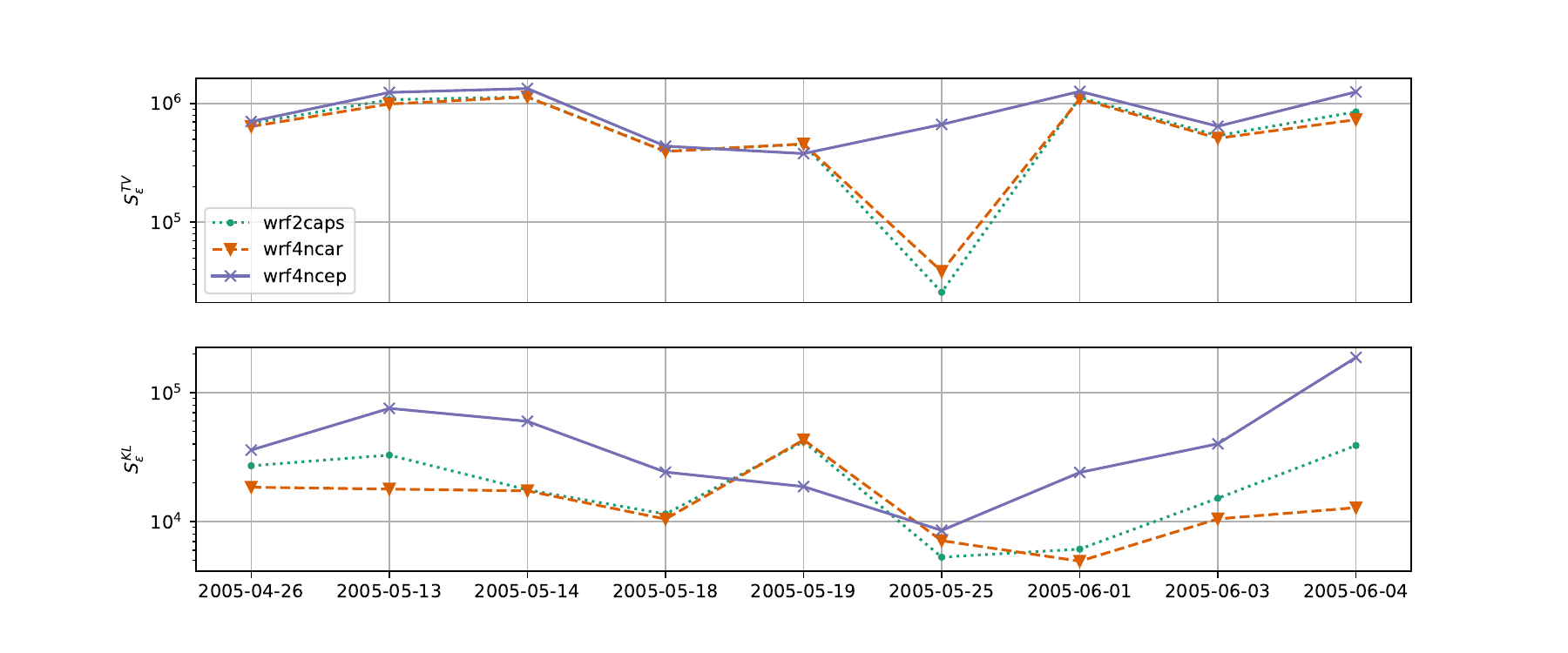}
    \caption{Spring 2005 cases' scores for 9 valid times. Note the difference in scale between \(\Sink^{TV}\) (top) and \(\Sink^{KL}\) (bottom), where  the models' rankings are mostly in agreement. However, there are clear anomalies at key dates: 2005-05-25 and 2005-05-19. Each model was at 24hr lead time, and interpolated on the coarser 4km grid. \(\varepsilon=0.001 L^2,\ \rho=L^2\). }
    \label{fig:spring_2005}
\end{figure}

The penultimate dataset considers 3 real model evaluations against observations.
These address \ref{q9} predominately, though \ref{q3} is also answered. 
Here the in-sample average gives \(M=200464\).
The experts' subjective evaluation ordered the models (wrf4ncar, wrf2caps, wrf4ncep), better to worse.
Across flavours, \(\Sink\) followed, on average (mean and median), in agreement with this ordering (see Figure \ref{fig:spring_2005}).
We notice that rankings for each day across various scores do not agree (Figure \ref{fig:sp2005_ranking}), i.e. we lack correlation with the subjective ranking.
This follows similar trends, such as for IWS \citep{gilleland_2010},  CRA \citep{ebert_gallus_2009}, and DAS \citep{keil_craig_2009}.
In fact, we have the advantage of agreeing on average with the expert assessment, however this is not a ground breaking statement due to the small number of experts sampled and small number of test cases considered.
Additionally, their scores did not themselves discriminate well a better or worse model configuration. 
Even if the Spearman rank correlation coefficient \citep{zwillinger1999crc} is calculated across the 9 days, Figure \ref{fig:heatmap-comparison}, there is still low and sometimes negative correlation. 
One conclusion that may be drawn is that the KL flavour does perform marginally better.
% However, on average agree with the overall ranking, with these heavy caveats. 

Notably, there are two large anomalies in the ranking with \(\Sink\): May 19, where wrf4ncar and wrf2caps perform worse than wrf4nceps, and May 25 where wrf4ncar performs worse than wrf2caps.
Otherwise, across the models the ranking was consistent, unlike the expert ranking.
Significantly, these two anomalies also correspond with two notable differences in the balance of mass. 
On May 19th there is a significant over forecast from all the models, shown in Figure 4 of \citet{ahijevychetal_2009}.
Then on May 25th wrf2caps and wrf4ncar under forecast, and both failed to pick up on a feature in the north-west, which wrf4ncep did forecast (see Figure \ref{fig:graphical_abstract}). 
This is reflected in the large ATM for wrf2caps and wrf4ncar (see Figure \ref{fig:spring_2005_all}) since the failure to predict this feature means this mass it is transported to the next closest feature, which is far away.
On average, \(\Sink^{KL/TV}\) follows subjective scoring and is able to pick up on anomalies which are diagnosable through further investigative techniques under the same methodology. 

\subsubsection{MesoVICT (VERA) Cases}\label{section:vera_mesovict}

To consolidate this methodology we investigate the MesoVICT Core Case (case 1).
Unfortunately, there is no expert evaluation hence this case is presented as a suggested evaluation strategy or operational use case. 
\(\varepsilon=0.005 L^2\)  given the grid (138 by 96 points), and M is (1958, 6405, 13361, 26685) depending on the model accumulation time. 

\begin{figure}
    \centering
    \includegraphics[width=0.995\linewidth]{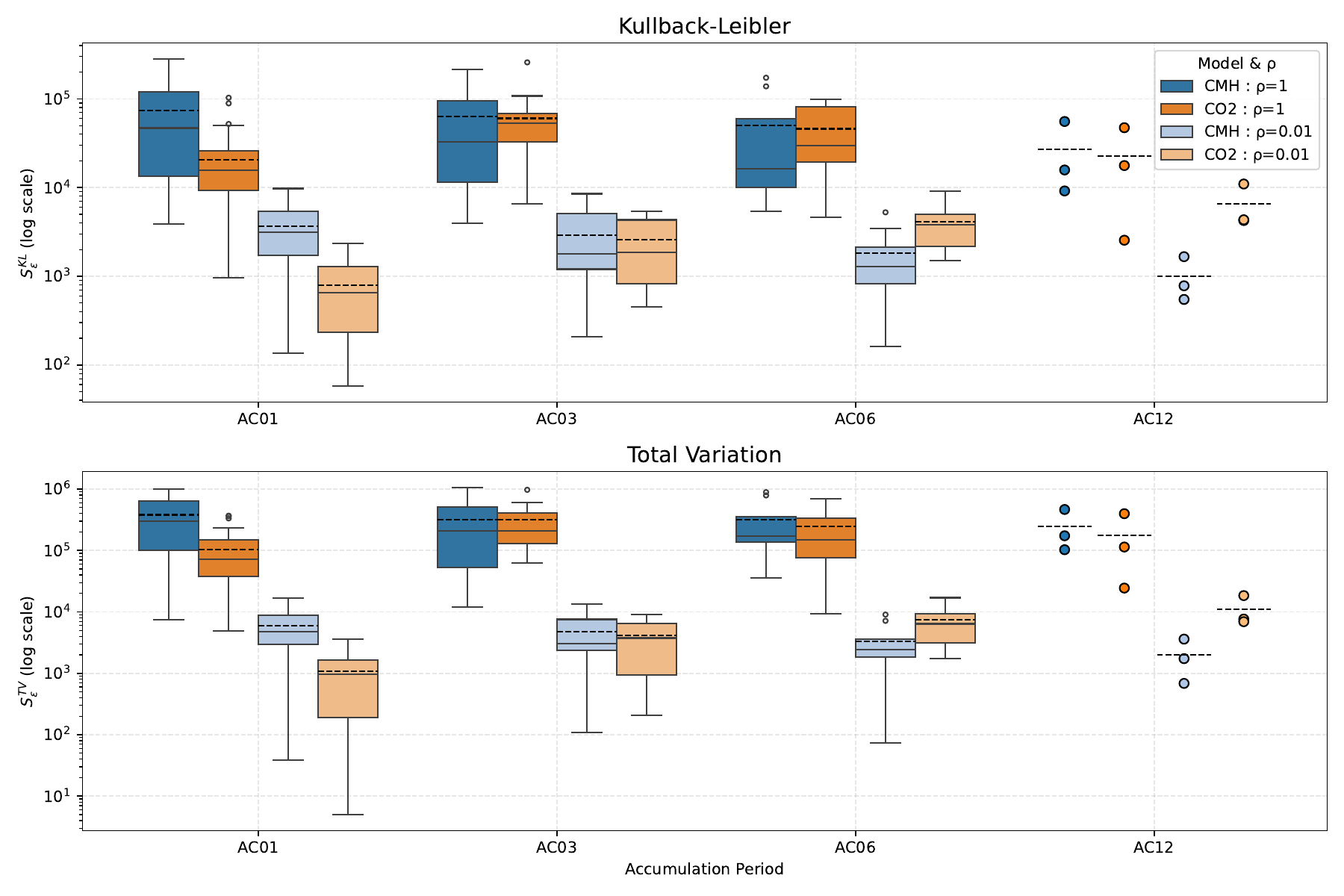}
    \caption{
    Box and whiskers plot showing the spread of the Sinkhorn divergence score for both flavours of penalty, and across accumulation time, (1, 3, 6 and 12 hours - left to right). The box and whiskers show the median (solid) line, mean (dashed) line, and interquartile range between the shaded regions. The whiskers show furthest non-anomalous point, being with-in 1.5 times the IQR. Otherwise they are considered anomalies and plotted separately. For AC12 since there were only 3 points to consider they are simply plotted and the mean given. A full box and whiskers would not be appropriate here. For each AC there are four boxes, corresponding to the two models and the two values of \(\rho\) tested (with a factor of \(L^2\)). \(\varepsilon=0.005 L^2\).}
    \label{fig:vera_ranking}
\end{figure}

Figure \ref{fig:vera_ranking} shows the spread of \(\Sink\) for the two models, across accumulation time and for \(\rho = 1 \) or \(0.01\). This information is plotted as a box and whiskers with the mean average added as the dotted line. 
For \(\Sink^{KL}\), CO2 outperforms CMH at AC01, whilst for larger accumulation periods, when timing errors are absorbed, the average performance starts to become comparable between the two models.
The median average performs better for CMH than CO2 for KL. 
For \(\Sink^{TV}\), these trends are similar though with less discernable difference in performance between models at AC03 and AC06. 
For lower \(\rho\), CMH starts to outperform CO2, suggesting it has accounted for its local transport and imbalance error through lag.

\begin{figure}
    \centering

        \includegraphics[width=\linewidth,trim= 70 70 70 70,clip]{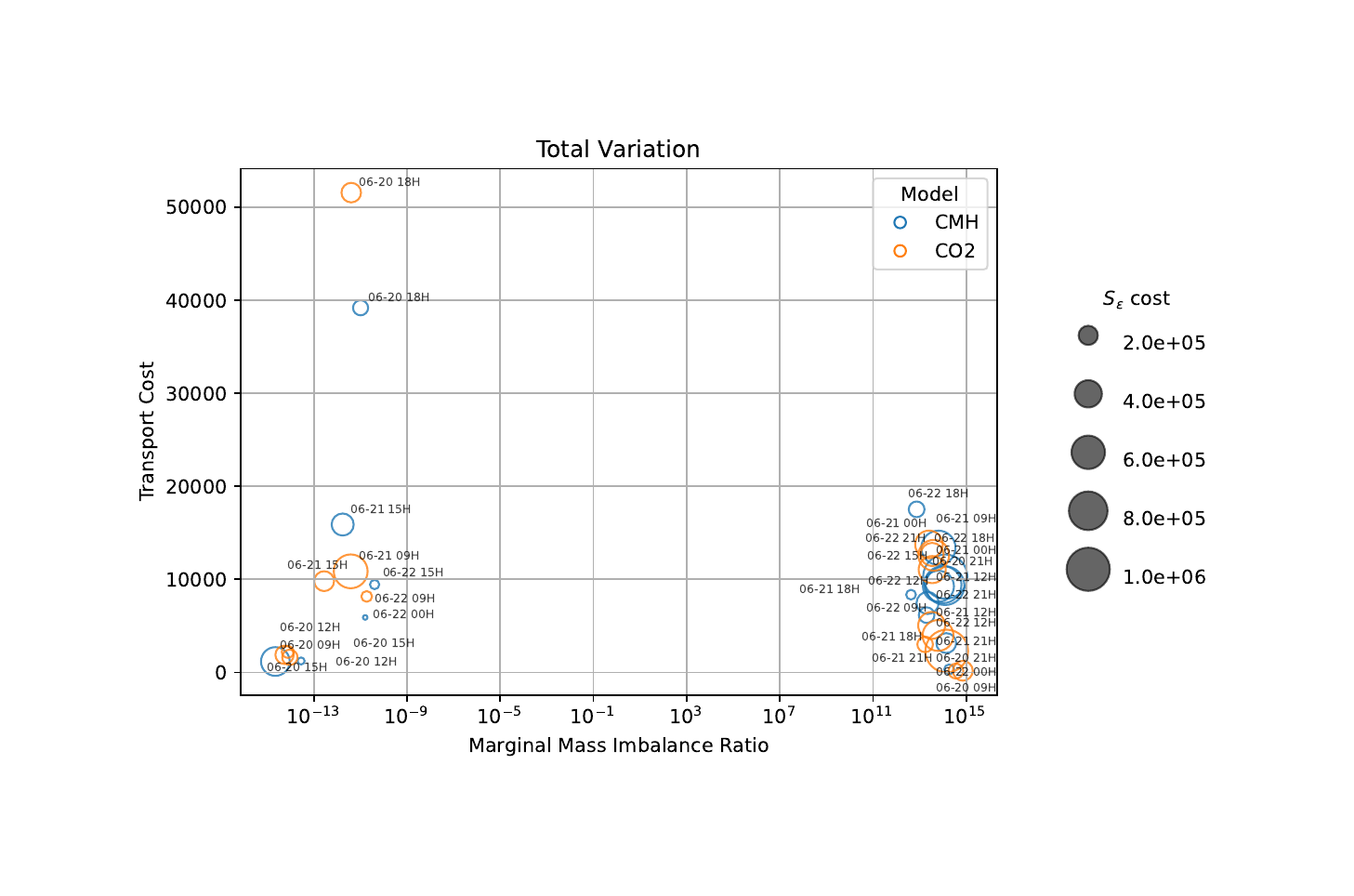}

        \includegraphics[width=\linewidth,trim= 70 70 70 70,clip]{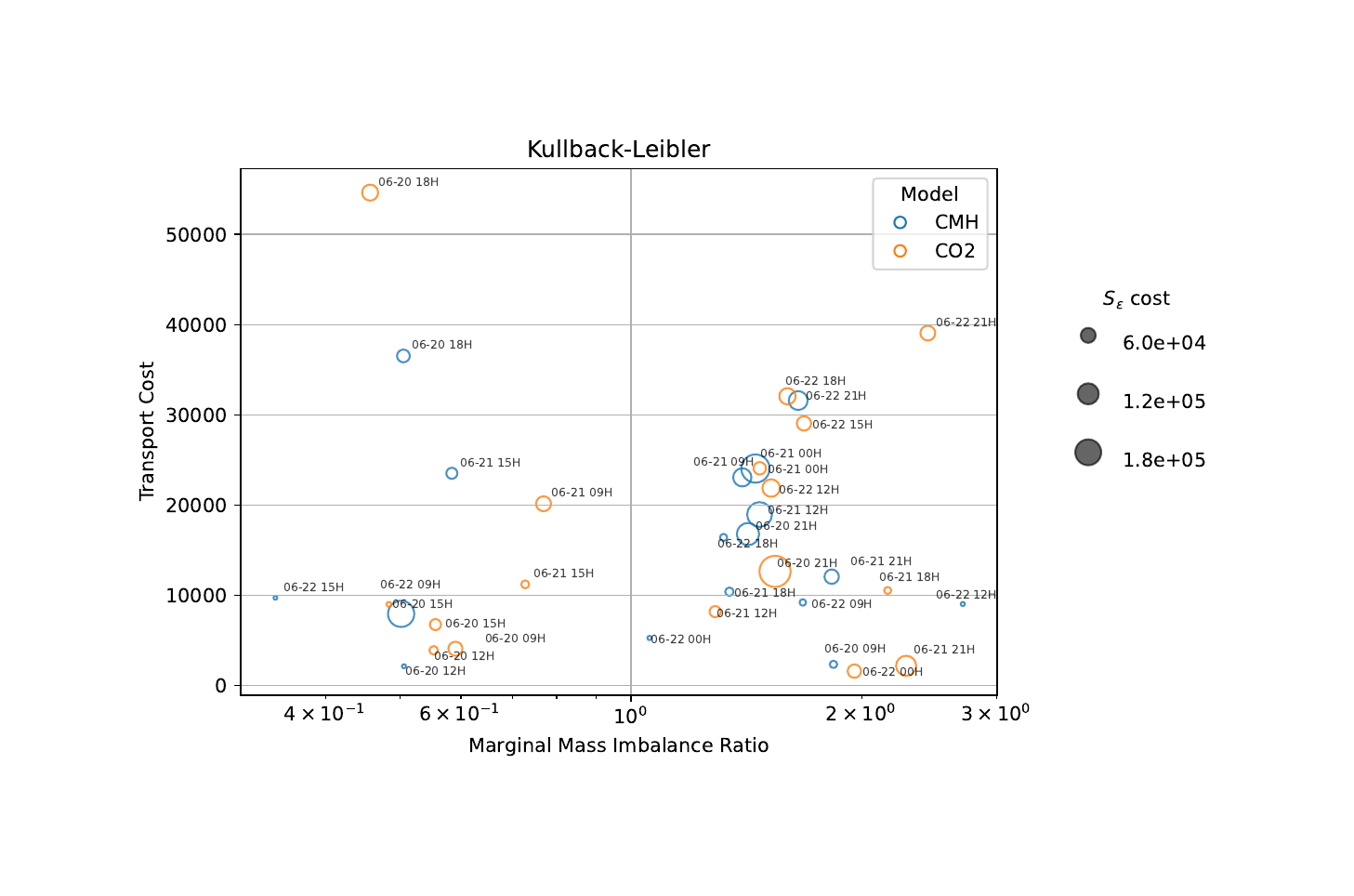}
    \caption{Decomposition of cost into the transport and marginal mass imbalance ratio, for the MesoVICT core case. This shows the decomposition for AC03, with \(\rho=L^2\). Notice CMH (blue) appears to have a higher transport cost across cases compared to CO2. However both over and under forecast for this short accumulation period. Case 06-20 18H AC03 stands out at having large transport distances in both penalities and models (From AC01, Figure \ref{fig:vera_decomposition_grid_rho1}, we see this was caused at the 17H valid time for CMH and accumulated for CO2).  Top: TV penalty, Bottom: KL penalty. }
    \label{fig:vera_decomposition_AC03}
\end{figure}

The time series for the spread of the \(\Sink\) and the ATM are shown in Figures \ref{fig:vera_ac0103_rho1} and \ref{fig:vera_ac0612_rho1} . 
Figure \ref{fig:vera_decomposition_AC03} presents a decomposition of the transport cost against the marginal mass imbalance ratio  (\(D(\pi_0 | \mu_{O})\) /  \(D(\pi_1 | \mu_{F})\)) for AC03. AC01 and AC06 are available in Figure \ref{fig:vera_decomposition_grid_rho1}, or Figure \ref{fig:vera_decomposition_grid_rho01} for \(\rho=0.01\). 

Beginning with the time series, consider AC06, with visualisations of the precipitation events available in the Supplementary Materials \ref{fig:vera_ac06_illustration_1},\ref{fig:vera_ac06_illustration_2}. 
This time resolution offers detail at a manageable quantity.
%Potentially AC12 is too course a time resolution, and AC01/AC03 offers too much information to start.
At 06-20 12H, CO2 over-forecasts and CMH outperforms it in its ATM and balance.
As the front travels from west to east, \(\Sink\) increases to a peak between 06-21 00H and 06-21 12H for CMH. 
At 06-21 00H both forecasts under-forecast significantly, although the ATM scores are reasonable since areal bias does not dominate.
In contrast, CO2 peaks then achieves its lowest \(\Sink\) at 06-21 12H. 
It improves its score by balancing total rain fall with the observation. 
In Figure \ref{fig:vera_ac06_illustration_1} a subtlety arises whereby, CO2 has achieved balance through a high intensity event entering from the northern boundary. 
Not matching spatial extend any better than CMH.
This stands out in the decomposition for being under forecasted with some transport error, but a small circle size. 
As the front moves easterly CMH improves its \(\Sink\) whilst CO2 worsens. 
It appears CO2 does not advect fast enough over the Alps, and this causes imbalance in mass which is penalised, as well as an increase in the ATM. 

% At 06-22 15H and 18H AC03 CMH over then under forecasts, relsulting in overall better balance at Ac06, than CMh, which under forcasts in both instanstes. this is coupled with an improvment int he ATM for CMH ATM and better \(\Sink\).  

Focusing on events with large transport error, 06-20-18H at AC03 (Figure \ref{fig:vera_ac03}) is highlighted in the decomposition figures. 
From the AC01 this was caused at the 17H valid time for CMH and accumulated for CO2.
Both models over-forecast and still miss much of the small scale structure seen in the observation.
Both also missed a high intensity feature in the south western region.
Overall, for AC01 and AC03, CMH has the higher transport error, yet  at longer accumulation times lag errors are absorbed.
This become more clear when local transport is prioritised through a smaller reach parameter (Figure \ref{fig:vera_decomposition_grid_rho01}).
It is clear that at AC01 CMH has more transport error. 
Then by AC06 CMH becomes cheaper and CO2 performs worse in both balance and transport. 
Validating how CO2 performed better for \(\Sink\) for AC01, AC03, but CMH for AC06. 

Significantly, the marginal mass imbalance ratio divides under- and over- forecasting, where both models have a tendency to under forecast. 
This may be due to the observation having smaller scale events than the forecasts. 
One strategy may be to remove these very low intensity events, yet as previously discussed we do not want to lose information.

Highlighted above is the sensitivity once again to total rain balance, which is penalised more in TV than KL. 
No strong trends were seen within the models, the 2.2KM lower resolution model CO2 did slightly out perform the 2.5 KM CMH, especially at short accumulation times.
However the resolution difference is small ,and for longer accumulations, and locally constrainted transport CMH started to perform better.
Although both tended to under forecast, with often high ATM.
This high ATM we suppose is due to small sale structure seen in the observation and not the forecast.

%% file: content/5_conclusion.tex
\section{Conclusions and Further Work}\label{section:conculsion}

Throughout the above discussion, UOT has been shown to possess many desirable properties as a spatial precipitation forecast verification method and associated toolbox. Our main contribution was the introduction of the Sinkhorn divergence to the forecast verification community. Overall UOT was shown to be an integrated methodology for tackling spatial and intensity errors simultaneously, whilst being parsimonious and diagnostic of known errors. For operational use, powerful diagrams and visualisation were demonstrated, with alignment with some subjective scores and intuition, though with the flexibility for users to select their priorities. 

Specifically, \(\Sink\) was shown to be robust to the double penalty problem and capable of diagnosing translational error. For the balanced cases, where the total area and intensity is equal in the two fields, both types of marginal penalties (TV and KL) were shown to be competent at diagnosing frequency error, or shape bias, and translational error. One major benefit in this setting is \(\Sink^{TV}\)'s geometric link with the reach parameter, or \(\rho\). Here, a user may define a strict distance above which a feature becomes disassociated and instead that mass is modified. This reflects a skill-scale similar to Cluster Analysis, which uses the number of clusters to define skill \citep{marzban_et_al_2009}.
In contrast, \(\Sink^{KL}\) did not possess such a strict relationship to the reach parameter, however it 
% did show effectiveness in the presence of spatial biases and noise. That is, \(\Sink^{KL}\) 
was shown to be more robust, both to small amounts of noise in the presence of translation error and to over-/under- forecasting in unbalanced regimes, i.e. in coping with spatial biases and noise. 

Unfortunately, UOT was demonstrated as not being able to diagnose rotational error. Rather, stretching and squeezing of a shape is undertaken, thus mimicking an aspect-ratio error correction. This is due to the convex nature of the optimal map. One potential solution is to separately optimise the rotation, like MODE \citep{davis_et_al_2006}, or even explicitly include rotation as part of the optimisation argument.
Further, UOT is not able to diagnose if a pair has one component being a subset of another. Such a scenario may be beneficial to differentiate, say, forecasts E3E11 verses E7E3, since depending on the user the value changes. While these pairs did receive different scores, UOT did not distinguish them based on subset relationships.

Some simpler geometric properties were exemplified, including insensitivity to zero points and, as a result, to the domain boundary, non-detrimental and predictable sensitivity to orientation on the regular grid, and it was vehemently demonstrated that UOT cannot be hedged.
UOT is performing transport and will move overlap to optimally assign regions of rain. Hence, if necessary, it will move overlap to reduce the overall cost. By preference, then, UOT favours pairs which have balanced total mass and minimal transport.

All these properties carried from the simple binary forecasts into the real textured data. Where through the real textured cases, it was also shown that by careful selection of the reach parameter, the importance associated to mass balance and translation can be changed.
Additionally, the Sinkhorn divergence was shown on average to follow subjective evaluation of real model performance, with a preference for KL penalisation.
Furthermore, through identification of anomalies in a time series, a more detailed exploration was enabled, calling attention to model runs which were both significantly imbalanced and missing features. A similar procedure, would be recommended for use in operation as demonstrated with the MesoVICT core case.

This investigation utilised diagrams of the decomposition of cost, simple times series, 2D histograms of the transport vectors, detailed depictions of the underlying transport vectors, and optimised marginals.
Further, in the research setting and given more computational expense a parameter sweep in the reach may be performed, with marginals returned spotlighting a scale of skill and at what reach certain features are no longer associated. 
% Or even the transport vectors found through the optimised marginal

Finally, while this current work is meant as an initial investigation into the potential behind precipitation forecast verification with the Sinkhorn divergence and its various tools, there is still much to consider, all contained within the same or similar framework.
Notably,  UOT is capable of having the observation-forecast pair on different grids. However, a thorough exploration of the effects and their influence on the costings would be needed for confidence to be given to the score.
There is also the possibility of looking into further decomposition of the costs, particularly in relation to the dual interpretation of the objective function. Between the primal and dual problem there are then 11 terms in total which may be investigated. 
Moreover, through numerical efficiency tricks, it is possible to capture information at successfully finer scales.
 
A further development would be turning the negatively oriented cost into a skill score. This could be against some random realisation over the domain, or climatology, or by considering bounds on the cost and using them to recentre the values.
Additionally, there are many marginal penalisation terms available beyond just KL and TV divergences, or as mentioned, the possibility of augmenting the objective function to consider rotation. Or even using UOT solely as a morphing technique before using a secondary evaluation - like OF methods. 
Related to this, is the possibility of examining quantities beyond precipitation, so long as they carry some density or mass, e.g. potential temperature and momentum.
One could also include topographical or location-based reach parameters, relating features which are geographically associated, but stopping (say over mountain ranges) features that should not be related.
 
Lastly, there is a strong capability for this methodology to extend into ensemble space. Underpinning all the mathematics is the extension to multi-marginal optimal transport \citep{benamou_2021, beier_etal_2022}, where each marginal would correspond to a member and one for the observation.
This would leverage conditional information on every member as well as granting access to a related problem of finding Wasserstein Barycentres, which are a transport-informed type of averaging, unlike the typical \(L_1\) average.

Recognising the importance of effective communication and user value, we welcome constructive feedback from readers on the visualisation choices presented here to help guide refinement in future studies and support the operational use of UOT.

\section{Supporting information}

The supplementary material, \ref{appendix:supplementary_material}, contains details of our implementation (Section \ref{appendix:implementation}), efficiency tricks used  (Section \ref{appendix:used_tricks}), and supporting figures for the idealised and real textured cases (Section \ref{appendix:extra_figures} and \ref{appendix:real_texture} respectively), as well as missing cases to match those examined in \citet{gilleland_et_al_2019} (Section \ref{appendix:all_cases}). 
Table \ref{table:binary_cases_here_2} describes the extra cases which are discussed in Section \ref{appendix:extra_figures}.

All figures are as follows:  
Figure \ref{fig:8case_2} studies reach sensitivity for C1 against displacement.  
Cases C6C12, C13C14, are in Figure \ref{fig:multiple_features}, E6E14, E2E10, in Figure \ref{fig:scaled_Cases}, transport vector illustration for C13C14, Figure \ref{fig:c13c14_transportvectors}.  
Cases C1C7, C1C8, in Figure \ref{fig:unbalanced_reach}, C1C9, E19E20, in Figure \ref{fig:subsets}.  
Cases E1E4, E2E4, C2C11, C1C6 in Figure \ref{fig:combined_multiple_features}.  
KL marginals for C1C6, C1C8 are in Figure \ref{fig:comparison_kl_marginal}.
Cases C6C7, C6C8, Figure \ref{fig:rho_balanced_reach}. 
TV marginal for C6C8 in Figure \ref{fig:c6c8_tv_marginal}, and C13C14 in Figure \ref{fig:c13c14_marginals}.  
 KL marginal for C6C8, C13C14, in Figure \ref{fig:c6c8_c13c14_kl_marginal}.  
C1C6, C1C7, C1C8 reach study, Figure \ref{fig:unbalanced_displcement_rho}.  
C1C9 2D histogram, Figure \ref{fig:c1c9_2d_hist}.  
E19E20 transport vectors, Figure \ref{fig:e19e20_transport_vectors}.
2D histograms, Figure \ref{fig:e19e20_kl_2d_hist}, \ref{fig:e20e19_kl_2d_hist}.  
E2E4 transport vectors, Figure \ref{fig:e2e4_transportvectors_tv}.  
E3E11, E7E3, E7E11 cost decomposition, Figure \ref{fig:ellipse_cost_distintergartion_across_rho}.  
Table of extreme and edge values, Table \ref{tab:extreme_edge_cases}.  
Cases P2P2, P2P5, P1P1, P1P1, Figure \ref{fig:new_new_paper_pcase_1}, S1S2, S1S3, H1H2, in Figure \ref{fig:scattered_hole}, C1C4, C1N3, C1N4, N1N2, in Figure \ref{fig:noise}.  
Decomposition for all ICP geometric cases Figure \ref{fig:cost_disintergration}.  

Perturbed real cases: Figure \ref{fig:perturbed_cases}, \ref{fig:perturbed_cases_smallrho} with \(\rho=1\) and \(2^{-6}\), and Figure \ref{fig:fake_cases_median_rho1} and \ref{fig:fake_cases_median_rhominus6} are demonstrations of Figure \ref{fig:perturbed_cases} and \ref{fig:perturbed_cases_smallrho} with a median average ATM and ATD.  
Figure \ref{fig:dualdecomposition_0067} shows how the marginal terms in the fake006/7 differ.
fake007 2D histogram is illustrated in Figure \ref{fig:fake007_2d_histograms}.  
For the Spring 2005, Figure \ref{fig:spring_2005_all} shows the ATM, complementing Figure \ref{fig:spring_2005}, and Figure \ref{fig:graphical_abstract} presents the graphical abstract image.  
The rank is shown in Figure \ref{fig:sp2005_ranking} per day, and across days in Figure \ref{fig:heatmap-comparison}.  

For the MesoVICT core cases;e time series across all accumulation times (Figure \ref{fig:vera_ac0103_rho1}, \ref{fig:vera_ac0612_rho1}, \ref{fig:vera_mass}); decompositions for \(\rho=1\) (Figure \ref{fig:vera_decomposition_grid_rho1}) and at  \(\rho=0.01\) (Figure \ref{fig:vera_decomposition_grid_rho01}); and finally, illustrations of the model and observation precipitation in Figure \ref{fig:vera_ac06_illustration_1}, \ref{fig:vera_ac06_illustration_2}, and \ref{fig:vera_ac03}.

Subsequent figures fill in missing cases: C3C5, C1C10, C3C4, C3C5, C3C4, Figure \ref{fig:new_new_paper_circles_0},  E4E8, E6E16, E2E17, E4E12, Figure \ref{fig:new_new_paper_ellipse_0}, E4E10, E4E14, E1E3, E1E13, Figure \ref{fig:new_new_paper_ellipse_1}, E5E7, E2E18, E2E6, E1E11, Figure \ref{fig:new_new_paper_ellipse_2}, E2E16, E1E14, Figure \ref{fig:new_new_paper_ellipse_3}, P1C1, P2P6, P1P5, P1P3, Figure \ref{fig:new_new_paper_pcase_0}, P2C1, P6P7, P1P4, Figure \ref{fig:new_new_paper_pcase_2}.

%% file: appendix/appendix.tex
\appendix
\renewcommand{\thesection}{S\arabic{section}}
\renewcommand{\thesubsection}{S\arabic{section}.\arabic{subsection}}
\renewcommand{\thesubsubsection}{S\arabic{section}.\arabic{subsection}.\arabic{subsubsection}}

\renewcommand{\thefigure}{S\arabic{figure}}
\renewcommand{\thetable}{S\arabic{table}}
\renewcommand{\theequation}{S\arabic{equation}}

\setcounter{section}{0}
\setcounter{figure}{0}
\setcounter{table}{0}
\setcounter{equation}{0}

\newpage

\section{Supplementary Material}\label{appendix:supplementary_material}
For full details on optimal transport (OT), its variations, and implementations, consult standard texts such as \citet{peyre_cuturi_2018, villani_2009, Santambrogio_2015}.

For the unbalanced Sinkhorn updates used in this work, with TV and KL marginal penalisation, see \citet{sejourne_et_al_2021} and Section 6 with-in, which defines the log-sum-exp approximation updates used in our implementation \citep{uot_own_implementation}. Their work builds up the unbalanced Sinkhorn divergence derived from the balanced Sinkhorn divergence introduced in \citet{feydy_sejourne_vialard_2018}.
They provide an overview of the dual-primal relationship that underpins most modern optimal transport theory and applications, however \citet{peyre_cuturi_2018, chizat_peyre_et_al_2016} similarly provide a discrete introduction of this for the balanced and unbalanced setting respectively. 
It is this dual-primal relationship that enables fast, parallelizable code and offers numerous diagnostic perspectives and tools, all within the same framework.

\subsection{Specifics of Implementation}\label{appendix:implementation}
A brief note on implementation. 
Code is available at \citet{uot_own_implementation}.
Following \citet{merigot_thibert_2020}, the mesh, number of Sinkhorn iterates and entropic parameter are linked for balanced OT. 
In particular, they guarantee error below, \(\eta > 0\), given iterations, \(k : k \gtrsim e^{\frac{C}{\varepsilon}}\log(\frac{1}{\eta})\) - for some constant \(C\).
This remains a heuristic for the unbalanced case, however it appears to work well.
In fact, often it is possible to achieve convergence before this, as the theoretical rates are not sharp.
Necessarily, \(\varepsilon\) itself is chosen from the grid resolution, hence it is not user defined; given the average distance realised in the cost function, \(\varepsilon\) is chosen to keep \(e^{-\frac{c(x, y)}{\varepsilon}}\) order 1.
In practise then, \(\varepsilon \sim 1/\sqrt{N}\), where \(N\) is the number of grid points (assuming length scale removed). This similarly defines the expected number of iterations to achieve \(\eta\) convergence, however since convergence is necessary to be able to assign belief in the score more iterations are added to our implementation if this error is not sufficiently low (similarly the algorithm may, and often does, terminate early).
  
Since reaching convergence is essential for this exploration, especially for the realisation of optimal transport vectors, a few different computational techniques and packages are discussed.
This should also aid replicability and utility of our findings.
Python libraries, such as POT \citep{pot_2021} and GeomLoss \citep{feydy_sejourne_vialard_2018, geomloss_webpage}, are available as off-the-shelf packages.
POT provides a wider range of methods; however, GeomLoss has been optimised for machine learning training.
GeomLoss does calculate a debiased cost too, though it does not prioritise convergence and instead relies on heuristics to provide an approximately optimal solution suitable for backpropagation and gradient descent.
POT, on the other hand, can take into account convergence and provides a range of methods in the unbalanced regime.
However, to our knowledge, it does not provide access to the debiased unbalanced Sinkhorn divergence.
  
Following GeomLoss, our own implementation relies on PyKEOps \citep{pykeops_2021, pykeops_webpage} for fast on-the-fly reductions (suitable for GPUs), while simultaneously taking into account convergence and providing the post-debiasing necessary for forming the Sinkhorn divergence.
Since efficiency was not the main goal of this work, rather demonstrating the capacity behind the method, not all possible numerical tricks were employed.
In fact, there is an abundance of well-studied tricks, particularly for the balanced OT cases, which may further be used heuristically in the unbalanced case. Such efficiency tricks include \(\varepsilon\)-scaling, multiscale damping, and kernel truncation \citep{schmitzer_2019}, capacity-constrained OT \citep{benamou_etal_2020}, and more recently, a method known as fast translation-invariant Sinkhorn which is specifically for UOT \citep{sejourne_vialard_pyre_2022}.
  
In this work, two straightforward efficiency methods were employed: \(\varepsilon\)-scaling, which solves a sequence of problems by progressively reducing \(\varepsilon\) in multiplicative steps (without fine-to-coarse gaining), and tensorisation for separable cost functions. Additionally, the Python package PyKeOps was utilised for on-the-fly computation, leveraging either CPU or GPU compute.

\subsubsection{Tensorisation and epsilon-annealing}\label{appendix:used_tricks}
First, tensoriation relies on the separation of the exponential of the cost function given regular Cartesian coordinates. Note this exponential is known as the kernel and is used in the Sinkhorn algorithm (See Remark 4.17 in \citet{peyre_cuturi_2018}). 
Given two  regular Cartesian grids, \((x_1, x_2), (y_1, y_2)\), and the square Euclidean distance one can write;
\begin{align}
   e^{\frac{(x_{1,i} - y_{1,j})^2 + (x_{2,i} - y_{2,j})^2}{2\varepsilon}} = e^{\frac{(x_{1,i} - y_{1,j})^2}{2\varepsilon}} \cdot e^{\frac{(x_{2,i} - y_{2,j})^2}{2\varepsilon}}.
\end{align}
Since the grid is regular \(x_{1,i} \in \{x_{1,k}\}_{k=1}^{N_1}\),  \(x_{2,i} \in \{x_{2,l}\}_{l=1}^{N_2}\), and similarly for \(y_{1,j} \in \{y_{1,s}\}_{s=1}^{M_1}\),  \(y_{2,j} \in \{y_{2,r}\}_{l=r}^{M_2}\), where \(N = N1\cdot N_2, M = M_1\cdot M_2\). 
This allows us to reduce the computation memory and only store these smaller matrices, where any summation can be performed via;
\begin{align}
    \sum_{i=1}^{N} e^{\frac{c_{i,j}}{\varepsilon}} = \sum_{k=1}^{N_1} e^{\frac{(x_{1,k} - y_{1,j})^2}{2\varepsilon}} \cdot \sum_{l=1}^{N_2} e^{\frac{(x_{2,l} - y_{2,j})^2}{2\varepsilon}}
\end{align}
which reduces the complexity from \(O(N) \) to \(O(max(N_1, N_2))\sim O(\sqrt{N}) \).
  
The second trick, relies on the convexity of the regularisation and that through improved convexity faster convergence in Sinkhorn is realised. 
By then stepping down in subsequent orders of \(\varepsilon\) but recycling the potentials from the higher values of \(\varepsilon\), a more stable and provably convergent algorithm  can be found \citep{schmitzer_2019, chizat_2024}. 
That is given a value you \(\varepsilon\) a sequence of problems is solved for \(\varepsilon_i \in \{1\cdot\omega, \ldots, \omega^{p}, \varepsilon\}\) where \(\omega < 1\) is some scaling parameter, and \(p\in\mathbb{Z}\) such that \(\omega^{p+1} \leq \varepsilon\). 
Previous optimal solutions are used as an initialisation to the subsequent one. 

\subsection{Extra Figures}\label{appendix:extra_figures}

The details below highlight particularly interesting points that aid interpretation of the methodology and may be of interest to some readers. Due to space limitations in the main text, these findings may have only been briefly referenced or omitted entirely.

\begin{table}
\centering
\footnotesize

\begin{minipage}{0.48\linewidth}
    \centering
    \begin{tabular}{c|p{5cm}}
    \hline
    \textbf{Case} & \textbf{Description} \\ \hline

    C12 & Two radius 20 circles centred at (120, 160) and (80, 40) \\
   E12 & E4 shifted (15, \(\rightarrow\)) and (20, \(\uparrow\)) \\ 
    E19 & Three vertical ovals centred at (100, 40), (100, 55),
and (125, 75), scaled by (40, 5),
(35, 5), and (25, 5), respectively \\ 
    E20 & E19 smoothed using a disk kernel with radius 12 \\ 
    P1 & Null field \\

    \end{tabular}
\end{minipage}
\hfill
\begin{minipage}{0.48\linewidth}
    \centering
    \begin{tabular}{c|p{5cm}}
    \hline
    \textbf{Case} & \textbf{Description} \\ \hline

    P3 & One point at (1, 1) \\ 
    P4 & One point at (200, 200) \\ 
    P5 & One point at (100, 100) \\ 
    P6 & Four points in the four corners\\ 
    P7 & Four points near the middle of each boundary
at (1, 100), (100, 1), (200, 100) and (100, 200) \\ 
    
    \end{tabular}
\end{minipage}
\caption{Table of idealised cases discussed in this extra figures section, this is still not exhaustive of all the cases. For full description of all possible cases, see \citet{gilleland_et_al_2019}. The domain in a 200 by 200 grid, indexing from (1,1) to (200,200), using coordinates (horizontal, vertical). For the ellipse cases, the ratio of the major vs minor axis is 5:1 with dimensions of the ellipses either 100 by 20 (large ellipse) or 25 by 5 (small ellipse). All ellipses are centred at (100, 100) unless translated. 
Notation; C circles, E ellipses, P points, S scattered, N noisy and H hole cases.}\label{table:binary_cases_here_2}
\end{table}

\subsubsection{Multiple Features}\label{appendix:multiple_features}

\begin{figure}[h!]
    \centering
     \begin{subfigure}{0.48\textwidth}
        \centering
        \includegraphics[width=\linewidth, trim= 70 70 30 70, clip]{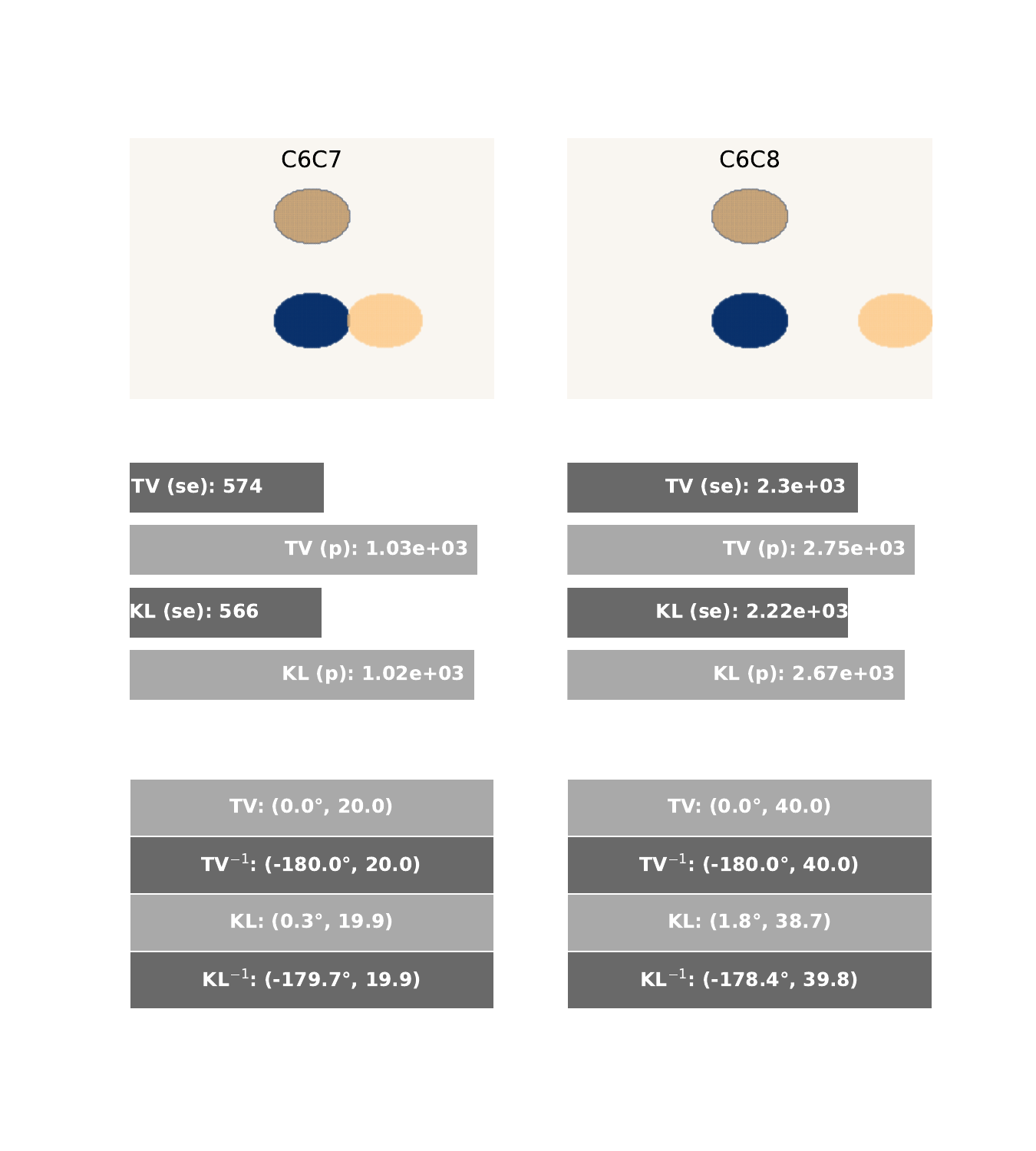}
        \caption{Demonstration that the logic of transportation holds when there is a hit and a miss present. The top four horizontal bars display \(\Sink^{TV}, \UOT^{TV}, \Sink^{KL}, \UOT^{KL}\), respectively. The lower table presents the mean (ATD, ATM) in both flavours, and with the forward and inverse vectors. The blue (darker) colour indicates observations, while the pale orange (lighter) represents forecasts. \(\rho = L^2\).}
        \label{fig:rho_balanced_reach}
    \end{subfigure}
    \hfill
    \begin{subfigure}{0.48\textwidth}
        \centering
       
        \includegraphics[width=0.65\linewidth, trim= 50  40 50 50, clip]{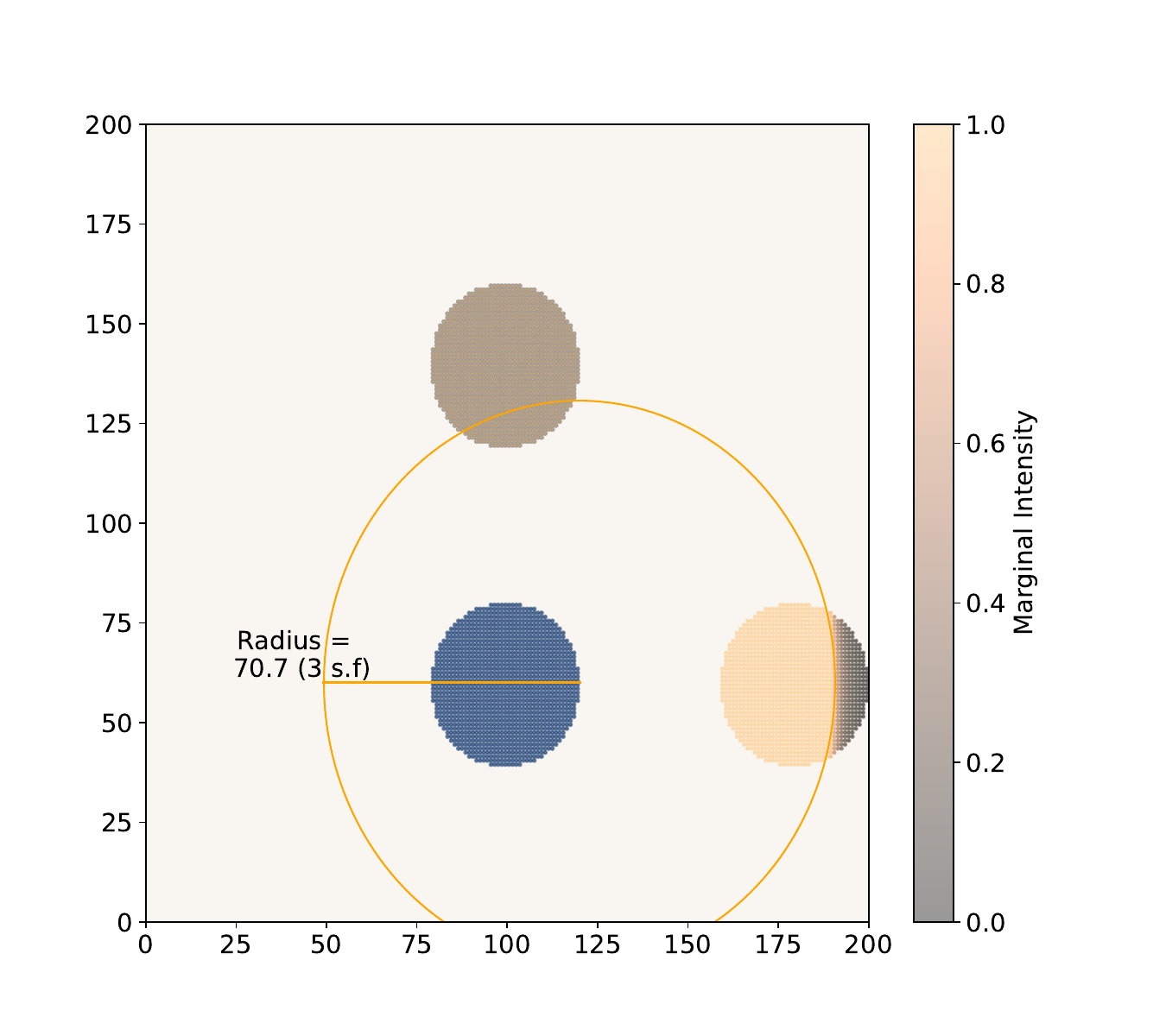}
        \subcaption{Case C6C8 with TV marginal penalty. The marginal \(\pi_1\) confirms there is a strict relationship with the reach parameter and it is maintained per feature. $\rho = 2^{-4} L^2$, the radius is equal to the reach.}
        \label{fig:c6c8_tv_marginal}
        
        \vspace{0.5em} 
 
        \includegraphics[width=0.65\linewidth, trim= 50  40 50 50, clip]{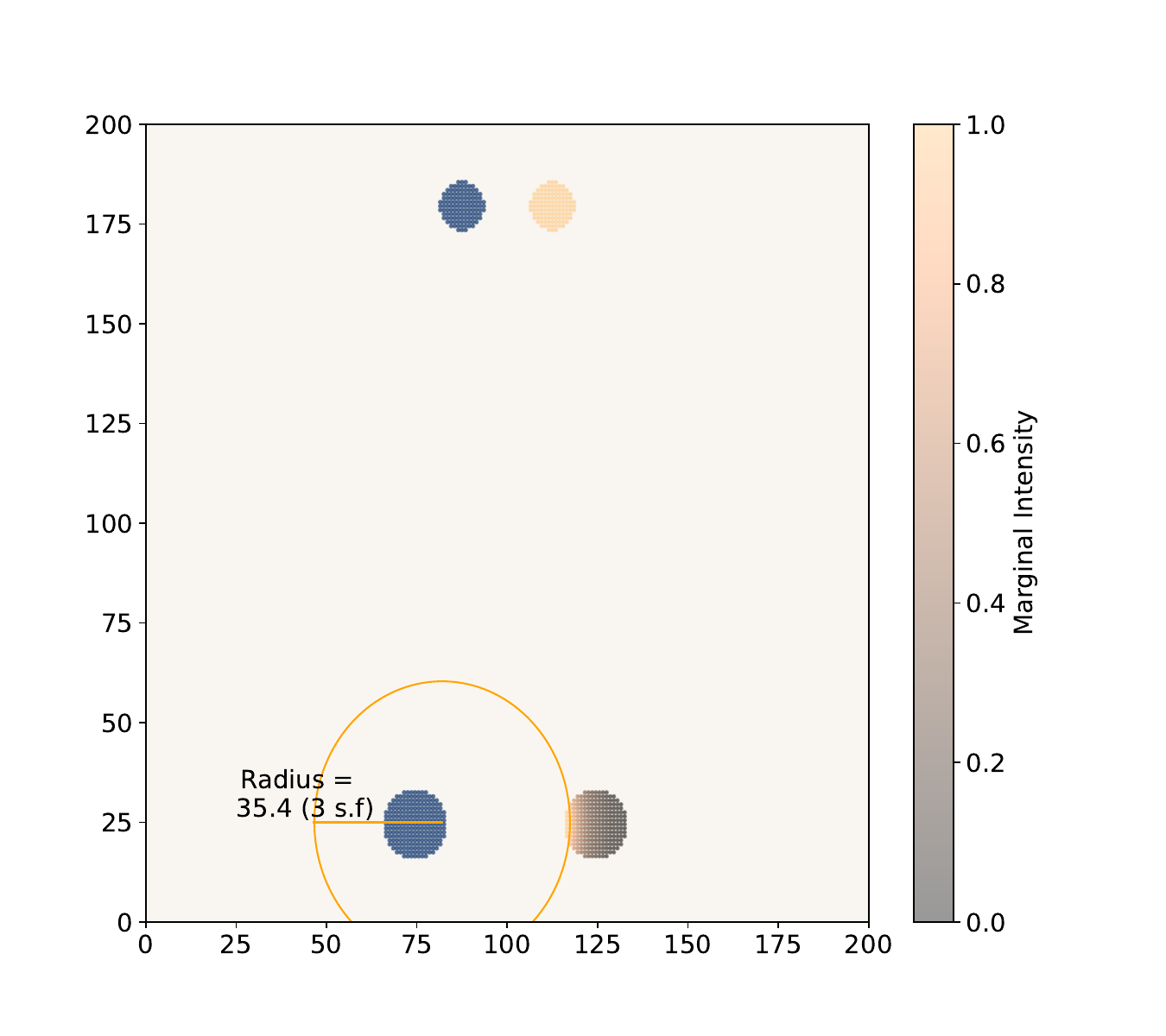}
        \subcaption{Case C13C14 with TV marginal penalty. The marginal \(\pi_1\) shows that a distance feature can be removed if deemed too far separated, and this is maintained per feature. Here, the radius is equal to the reach, $\rho = 2^{-6} L^2$.}
        \label{fig:c13c14_marginals}
    \end{subfigure}
  \caption{Figures exploring UOT's behaviour with multiple features present. (a) shows the scoring for cases with one hit and one miss which then moves away, and (b) and (c) show the returned marginals, \(\pi_1\), in blue (darker shades) with the observation in pale orange to grey shades. \(\varepsilon = 0.005L^2\).}
\end{figure}

Figure \ref{fig:rho_balanced_reach} considers cases C6C7 and C6C8 which each case containing one "hit" and one "miss", and are still in the balanced setting.
\(\Sink^{KL/TV}\) reports similar values to C1C2 and C1C3, since it is scoring the complete match as close to zero.
The \(\UOT^{KL/TV}\) cost fails to diagnosis this hit due to bias. 
Observe that the ATM is half the shift of the miss, due to the fact that only half the points require transporting.
The ATD remains close to 0 and 90 degrees.
Crucially, the full hit pays no transport or mass destruction fee (see Figure \ref{fig:c6c8_tv_marginal}), whilst the miss is penalised depending on \(\rho\) and analogously to C1C3, i.e. there is an equivalent \(\rho\) radius of influence dependence for the southern event.
It is important to note, that the hit's lack of cost is related to the zero transport required for the northern features, as opposed to the overlap. For example, C2C5 contains overlap, yet these points are shifted when the transport vectors are visualised. This is an important mathematical and behavioural distinction, from PAD, say, which a priori removes overlapped mass.
This behaviour arises due to the strict convexity of entropic regularisation of the plan, which makes longer transport distances more expensive.

\begin{figure}[h]
    \centering
     \begin{subfigure}[t]{0.45\textwidth}
        \centering
        \includegraphics[width=\linewidth, trim=70 70 30 70, clip]{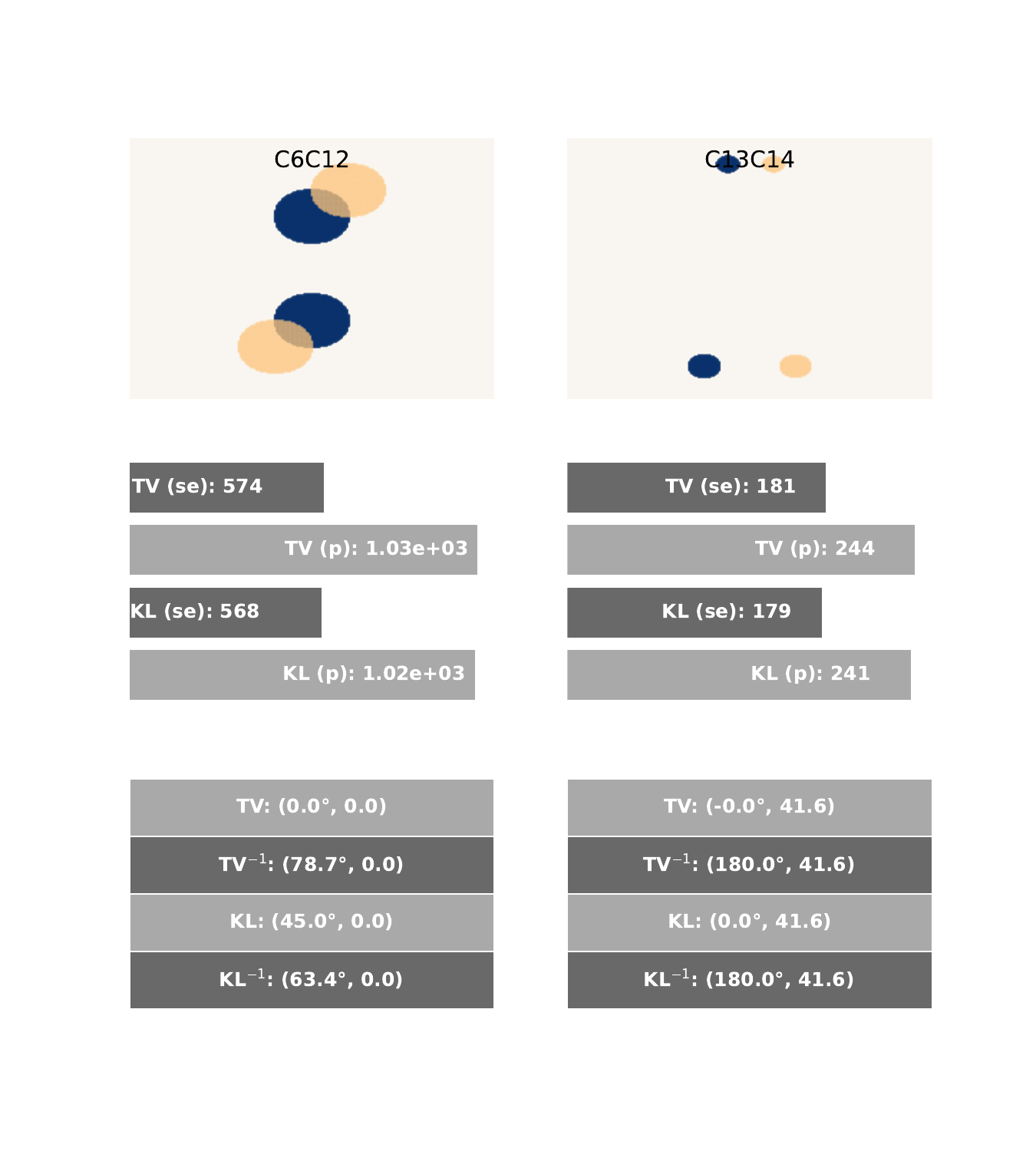}
        \caption{Balanced multiple features case, demonstrating that the logic of transportation holds when two transported features are present}
        \label{fig:multiple_features}
    \end{subfigure}\hfill
    \begin{subfigure}[t]{0.45\textwidth}
        \centering
        \includegraphics[width=\linewidth, trim=70 70 30 70, clip]{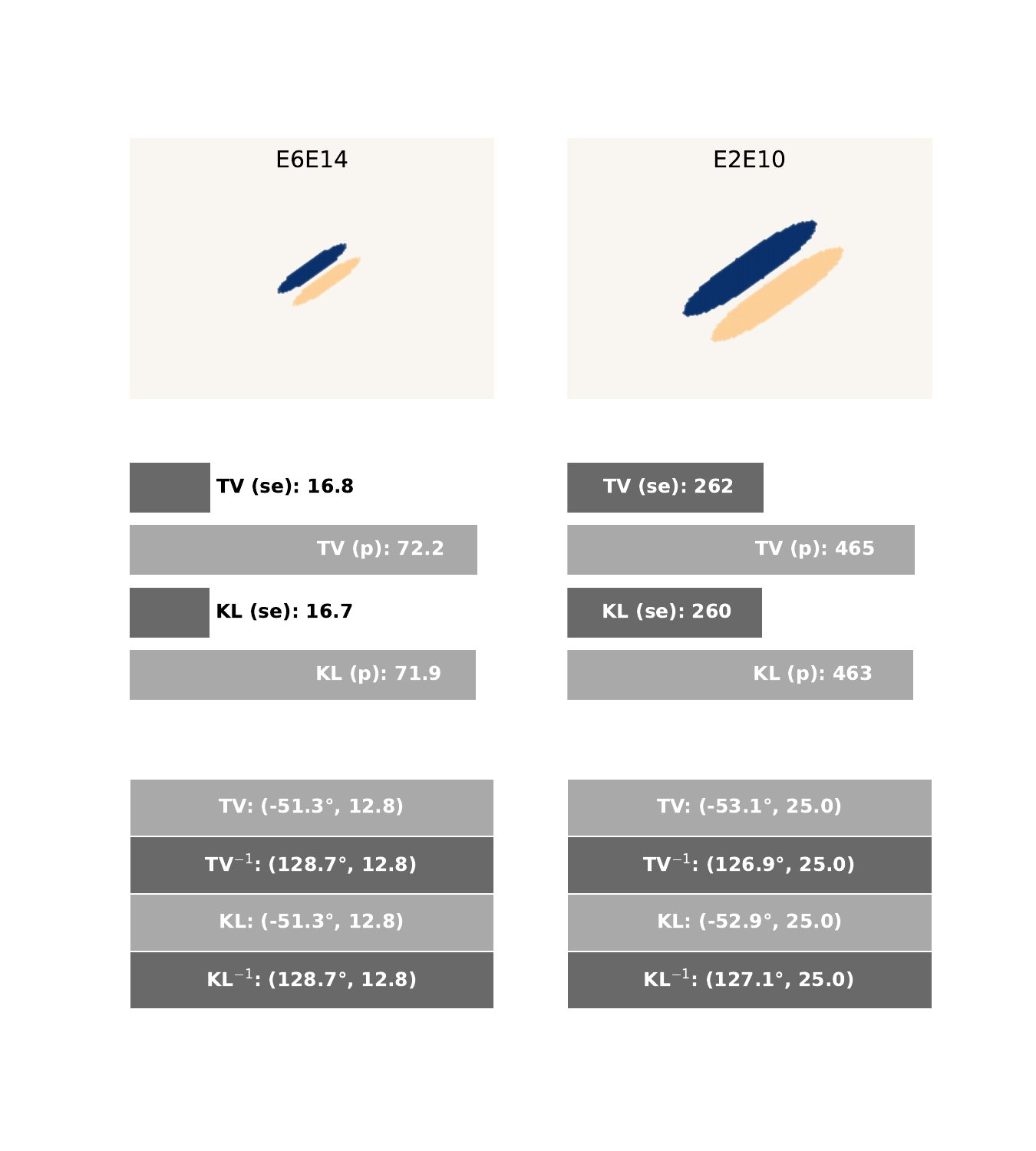}
        \caption{Balanced cases, yet with scaled size of events. The \(\Sink\) and ATM together allow powerful diagnostics of ranking, direction, and magnitude.}
        \label{fig:scaled_Cases}
    \end{subfigure}
    \caption{Figures exploring multiple features and scaling between different sizes of events. The top four horizfontal bars display; \(\Sink^{TV}, \UOT^{TV}, \Sink^{KL}, \UOT^{KL}\). The lower table presents the mean (ATD, ATM) in both flavours, and with the forward and inverse vectors. The blue (darker) colour indicates observations, while the pale orange (lighter) represents forecasts. \(\varepsilon = 0.005L^2,\ \rho = L^2\).}
    
\end{figure}

A further examination of the behaviour with multiple features is observed by the comparison between C6C12 and C1C2 (Figure \ref{fig:multiple_features}, and  Figure \ref{fig:translation}). 
Note these are comparable since the total north-south shift of C6C12 cancels and their west-east shift is equally 40 grid points.
Despite them having equivalent costs, they are discernable through the ATM and \(\UOT^{KL/TV}\).
C6C12 has each feature \(20\sqrt{2} \sim 28.3\) apart, though with equal and opposite magnitude.
This is reflected through zero average magnitude and arbitrary directions, which highlights a downfall of the ATM and ATD. 
Instead, considering the full distribution may offer more insight; this approach is taken in \citet{skok_lledo_2024}.
 
\begin{figure}[h]
    \centering
    \includegraphics[width=0.75\linewidth]{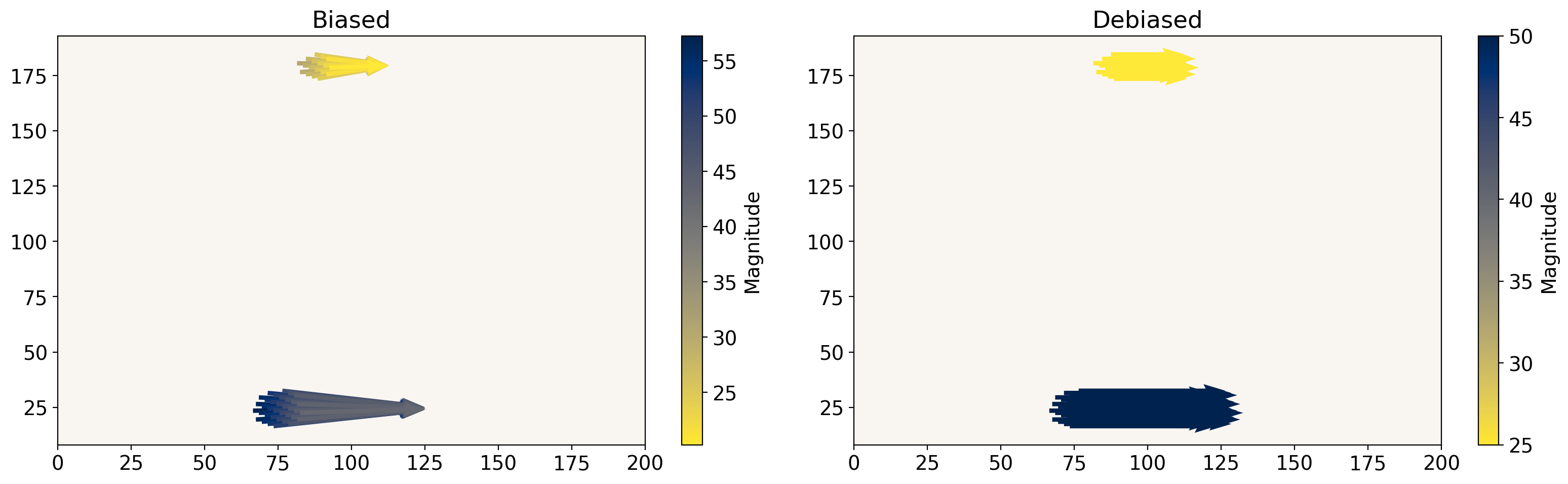}
    \caption{Biased vs debiased transport vector illustration with TV marginal penalty for the C13C14 case. Note, in the approximate map, the supports are still retained due to entropic regularisation, which assigns some mass to the supports and is defined everywhere by extension to the convex hull. Left: Biased UOT transport vectors, Right: Debiased UOT transport vectors. A regular sample of the vectors are shown to prevent overcrowding. \(\varepsilon = 0.005L^2,\ \rho = 2^{-4}L^2\).}
    \label{fig:c13c14_transportvectors}
\end{figure}

Additionally, through C13C14 (Figure \ref{fig:multiple_features}) it is possible to show that given sufficiently small reach parameter, one can dissociate features which are too far apart. 
Of course this is not directly reflected in the costs, however this is confirmed with visualisation of the approximate map (Figure \ref{fig:c13c14_transportvectors}, \ref{fig:c13c14_kl_marginal}). 
Moreover, by reducing the reach, a marginal with only the top feature retained is possible. The two southern events which are  too far separated are regarded as too expensive to transport and instead their associated mass is destroyed (Figure \ref{fig:c13c14_marginals}).
This demonstrates a dissociation of features which are too far apart, defined through the choice of \(\rho\).

Figure \ref{fig:multiple_features} illustrates two cases where multiple features are present in both fields, each exhibiting some form of symmetry. Case C6C12 shows the expected cost given the known total translation; however, the ATM cancels out because the translations are equal and opposite. Case C13C14 also shows the expected cost given the known translational error. Here, the ATM is an appropriately scaled version of the total transport, as the northern two features are mapped to each other, and similarly, the southern two features are mapped to each other. 
Figure \ref{fig:c13c14_transportvectors} shows this behaviour by viewing the transport vectors, clearly demonstrating that these distinct features remain separate in the transport.
However, if we suppose that is one wanted to separate influence one could calculate the barycentre using the optimised marginals.

\begin{figure}
    \centering
    \includegraphics[width=0.85\linewidth]{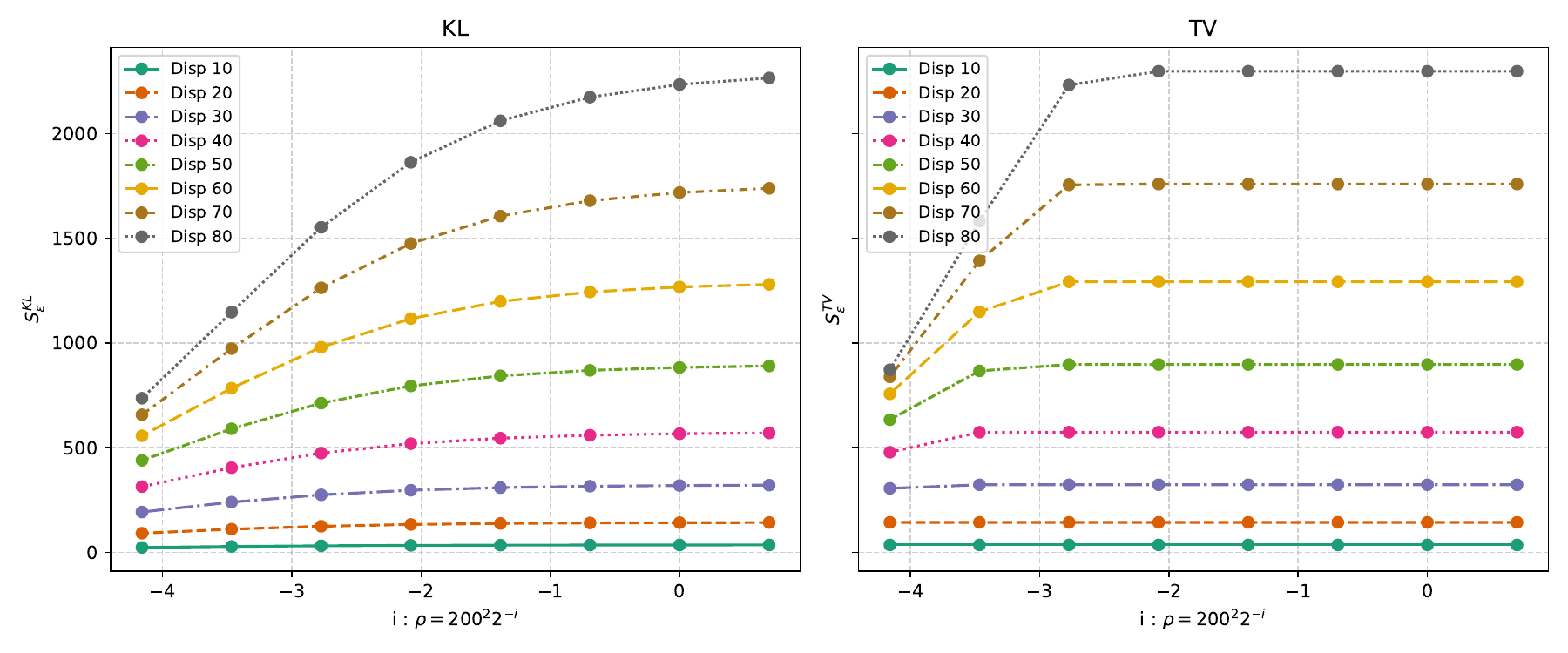}
    \caption{Reach sensitive study for case C1 at progressive 10 grid point increments from C1, thus in the balanced mass setting.
    Paired with Figure \ref{fig:rho_balanced_reach} both dependencies in displacement and reach are explored.
    For both flavours at low reach it becomes cheaper to destroy mass rather than transport, correcting  only local (small) displacement.
    Left: \(\Sink^{KL}\), Right: \(\Sink^{TV}\), Both: cost verse \(\rho\) values, with a curve per displacement case. \(\varepsilon = 0.005 L^2 \) .}
    \label{fig:8case_2}
\end{figure}

Continuing with the multi-feature cases, Figure \ref{fig:scaled_Cases} compares cases E6E14 and E2E10. Firstly, the debiased cost in both flavours correctly identifies the translational error. This is also reflected in the ATD and ATM. Secondly, they demonstrate that the cost is capable of scaling directly with the size of the event, i.e., the ellipses in E2 and E10 are 4 times larger, and there is approximately twice the transport. Then \(16.8 \cdot 4 \cdot 2^2 \sim 269\).

%%%%%%%%%%%%%%%%%%%
\begin{figure}[h]
    \centering
     \begin{subfigure}[t]{0.45\textwidth}
        \centering
        \includegraphics[width=\linewidth, trim= 70 70 30 70, clip]{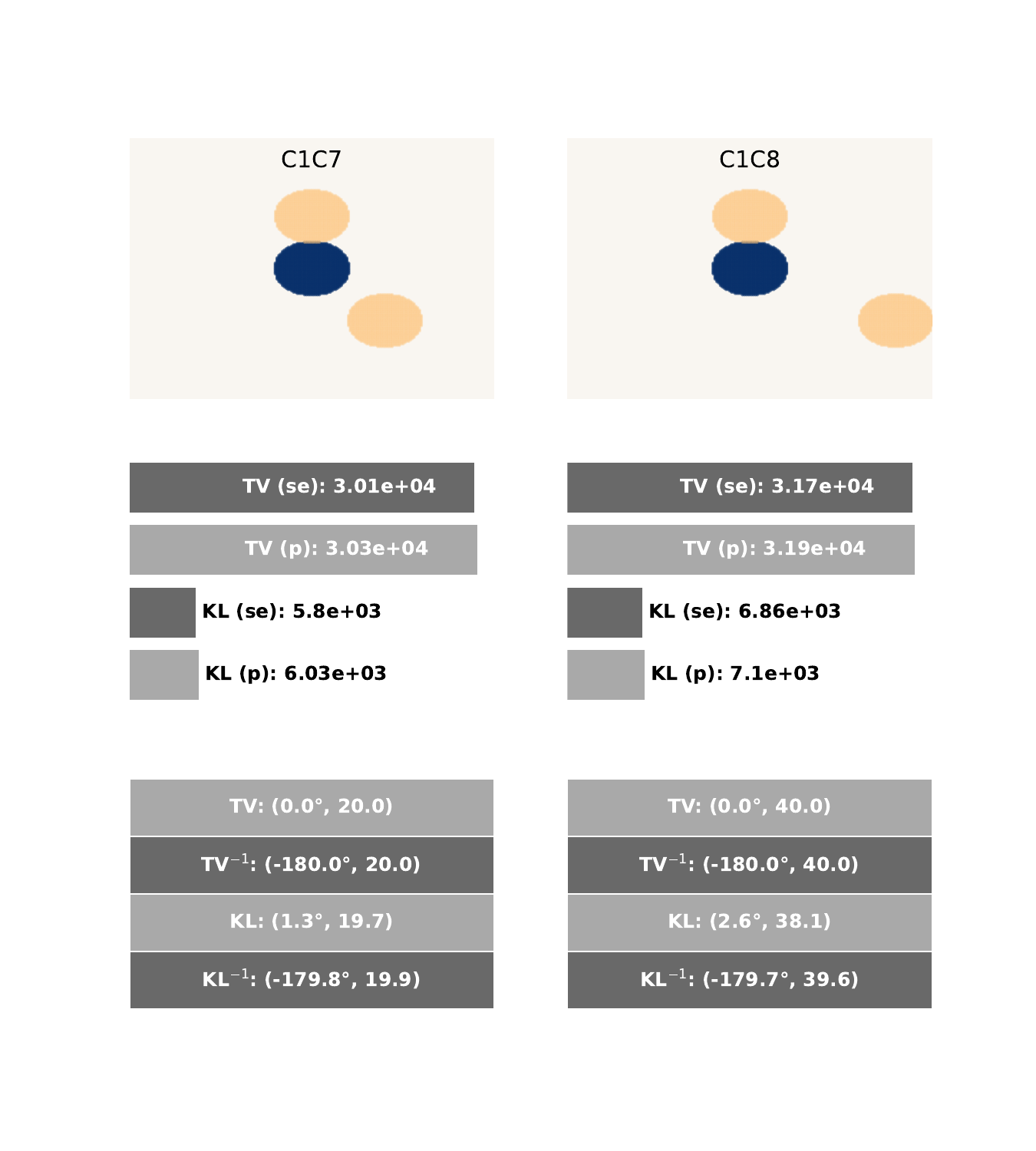}
        \caption{Illustration of cases C1C7, and C1C8, where there are two events, but the lower event moves away and thus becomes less associated with the closest event.}
        \label{fig:unbalanced_reach}
    \end{subfigure}\hfill
    \begin{subfigure}[t]{0.45\textwidth}
        \centering
        \includegraphics[width=\linewidth, trim= 70 70 30 70, clip]{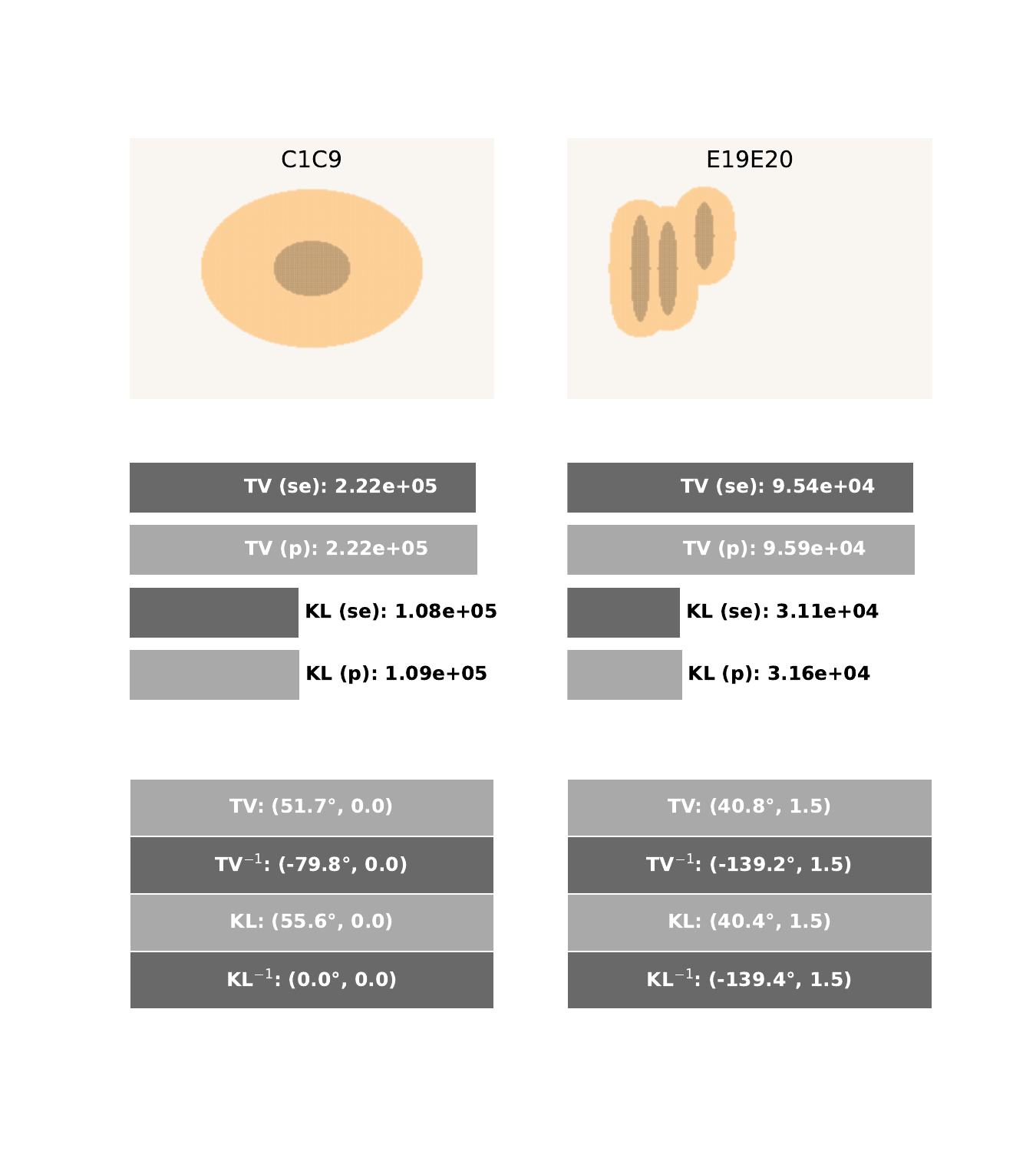}
        \caption{Testing UOT with subset events.
        However these score alone are not diagnostic of an event being a subset of the other.}
        \label{fig:subsets}
    \end{subfigure}
    \caption{Figures exploring multiple features and if UOT can diagnose subsetted events. The top four horizontal bars display; \(\Sink^{TV}, \UOT^{TV}, \Sink^{KL}, \UOT^{KL}\). The lower table presents the mean (ATD, ATM) in both flavours, and with the forward and inverse vectors. The blue (darker) colour indicates observations, while the pale orange (lighter) represents forecasts. \(\varepsilon = 0.005L^2,\ \rho = L^2\).}
    
\end{figure}

\begin{figure}[h!]
    \centering
     \begin{subfigure}[t]{0.45\textwidth}
        \centering
        \includegraphics[width=\linewidth, trim= 70 70 30 70, clip]{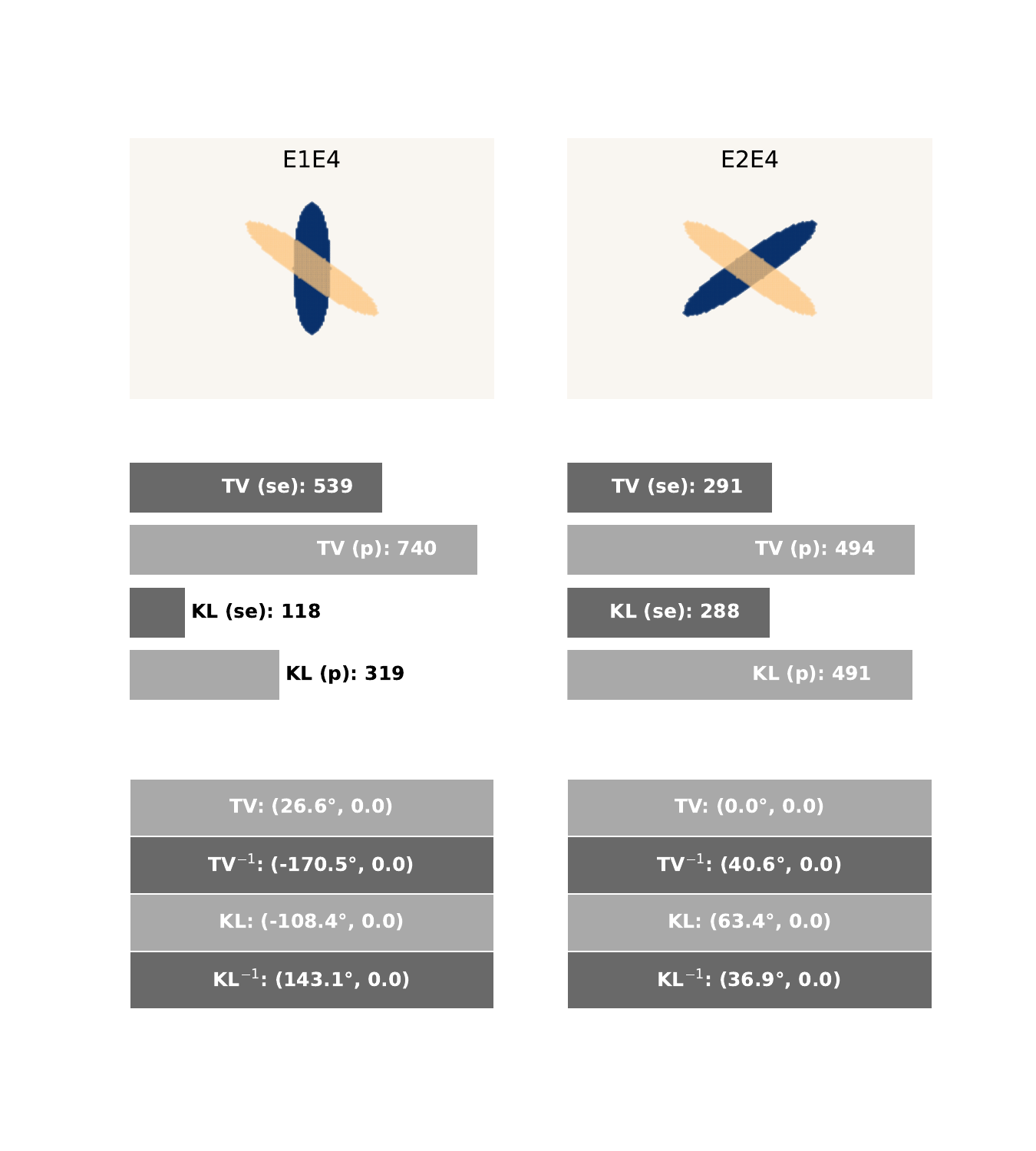}
        \caption{Rotated cases, E1E4 strictly is not balanced, though they are close in mass. UOT is not formulated to rotate shapes, instead stretches and squeezes them.}
        \label{fig:rotation_ellipse}
    \end{subfigure}\hfill
    \begin{subfigure}[t]{0.45\textwidth}
        \centering
        \includegraphics[width=\linewidth, trim= 70 70 30 70, clip]{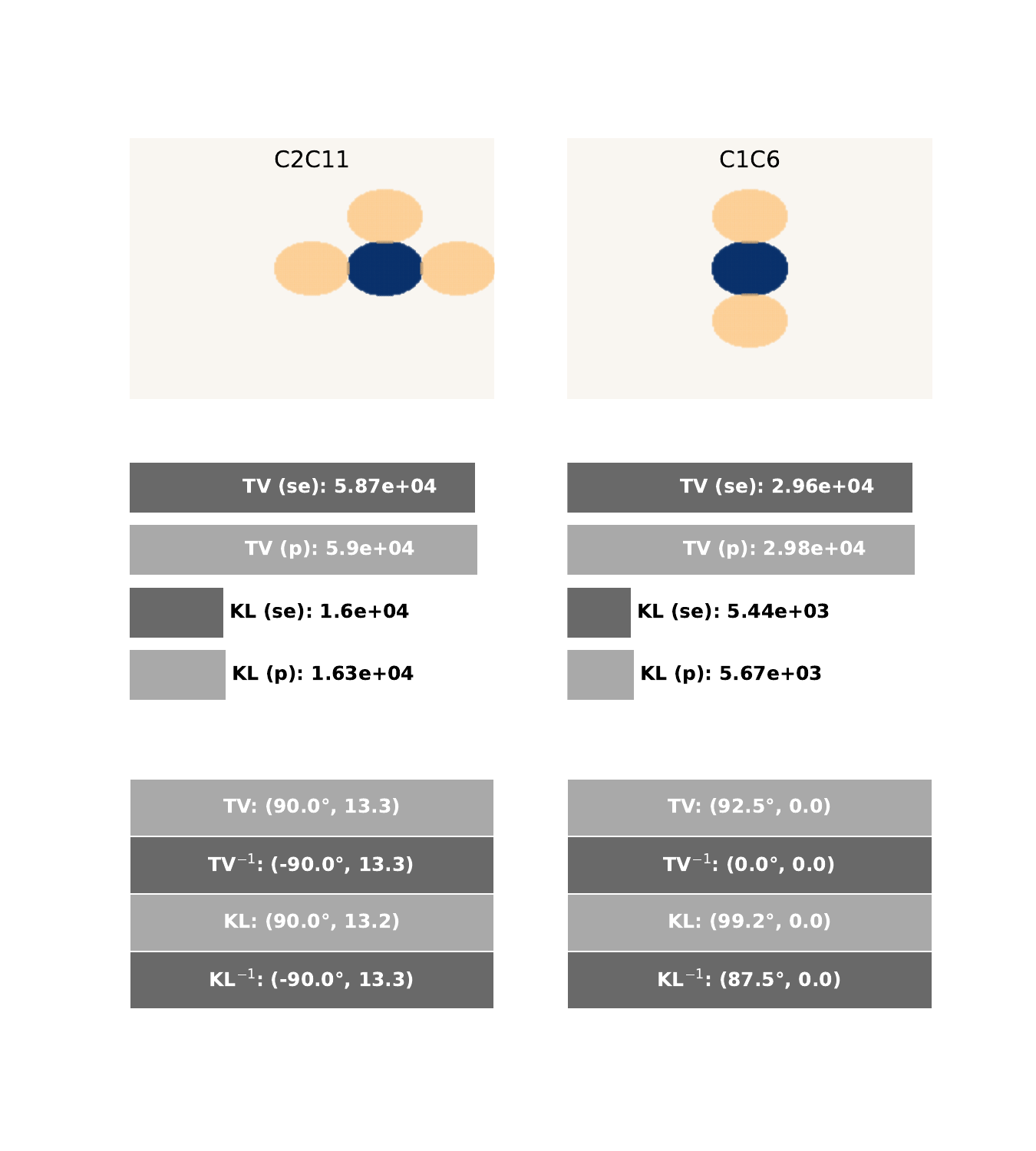}
        \caption{Over-/under-forecasting area extent cases, providing an unbalanced setup. In comparison to C1C2 both of these cases are penalised much more.}
        \label{fig:unbalanced_extent}
    \end{subfigure}
    \caption{ Figures exploring rotation and over-/under- forecasting of events providing an imbalance of total intensities. The top four horizontal bars display: \(\Sink^{TV}, \UOT^{TV}, \Sink^{KL}, \UOT^{KL}\). The lower table presents the mean (ATD, ATM) in both flavours, and with the forward and inverse vectors. The blue (darker) colour indicates observations, while the pale orange (lighter) represents forecasts. \(\varepsilon = 0.005L^2,\ \rho = L^2\).}
    \label{fig:combined_multiple_features}
\end{figure}

\begin{figure}[h]
    \centering
    \begin{subfigure}[b]{0.45\linewidth}
        \centering
        \includegraphics[width=\linewidth, trim= 50  40 50 50, clip]{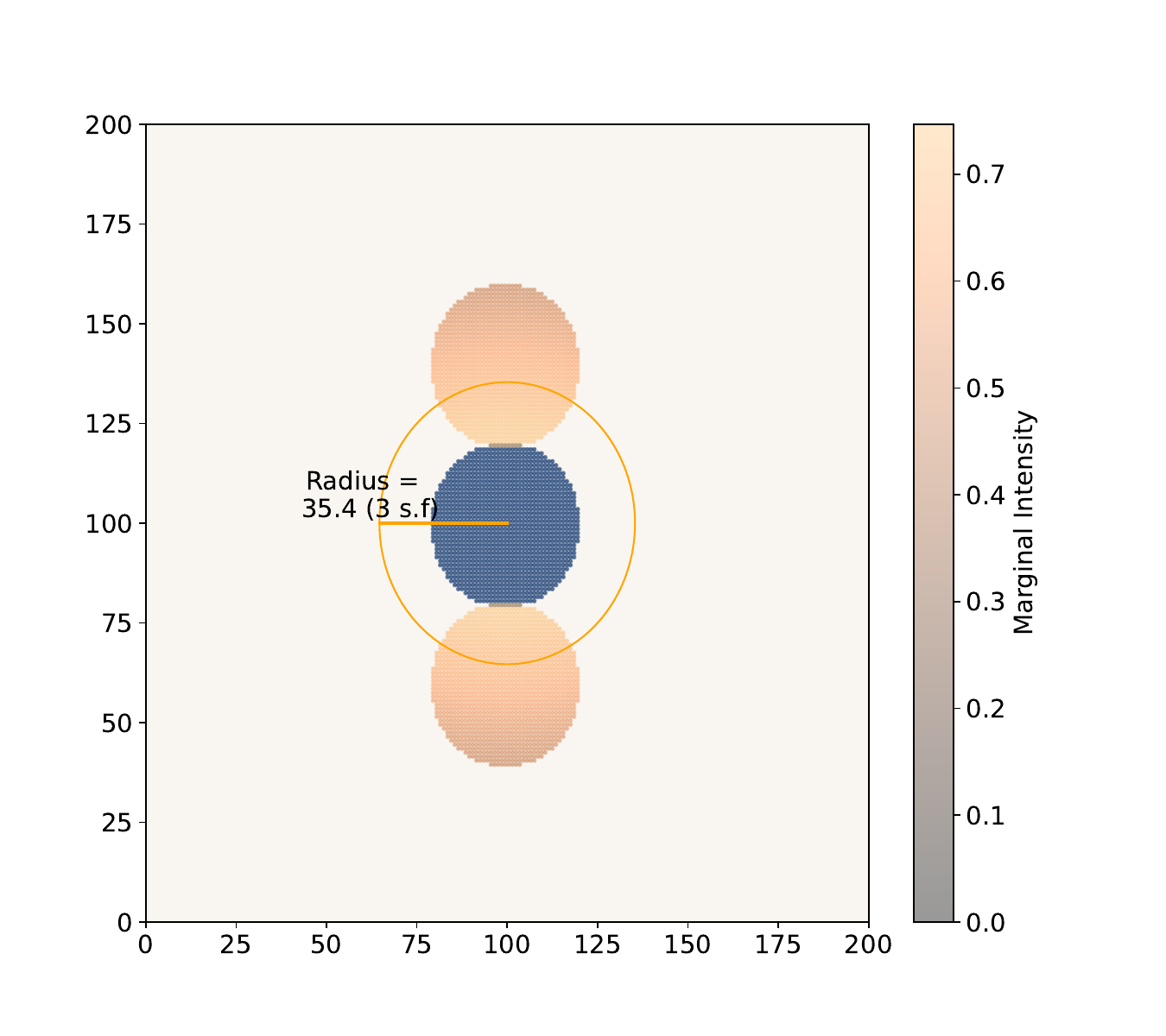}
        % \caption{}
        % \label{fig:c1c6_kl_marginal}
    \end{subfigure}
    \hfill
    \begin{subfigure}[b]{0.45\linewidth}
        \centering
        \includegraphics[width=\linewidth, trim= 50  40 50 50, clip]{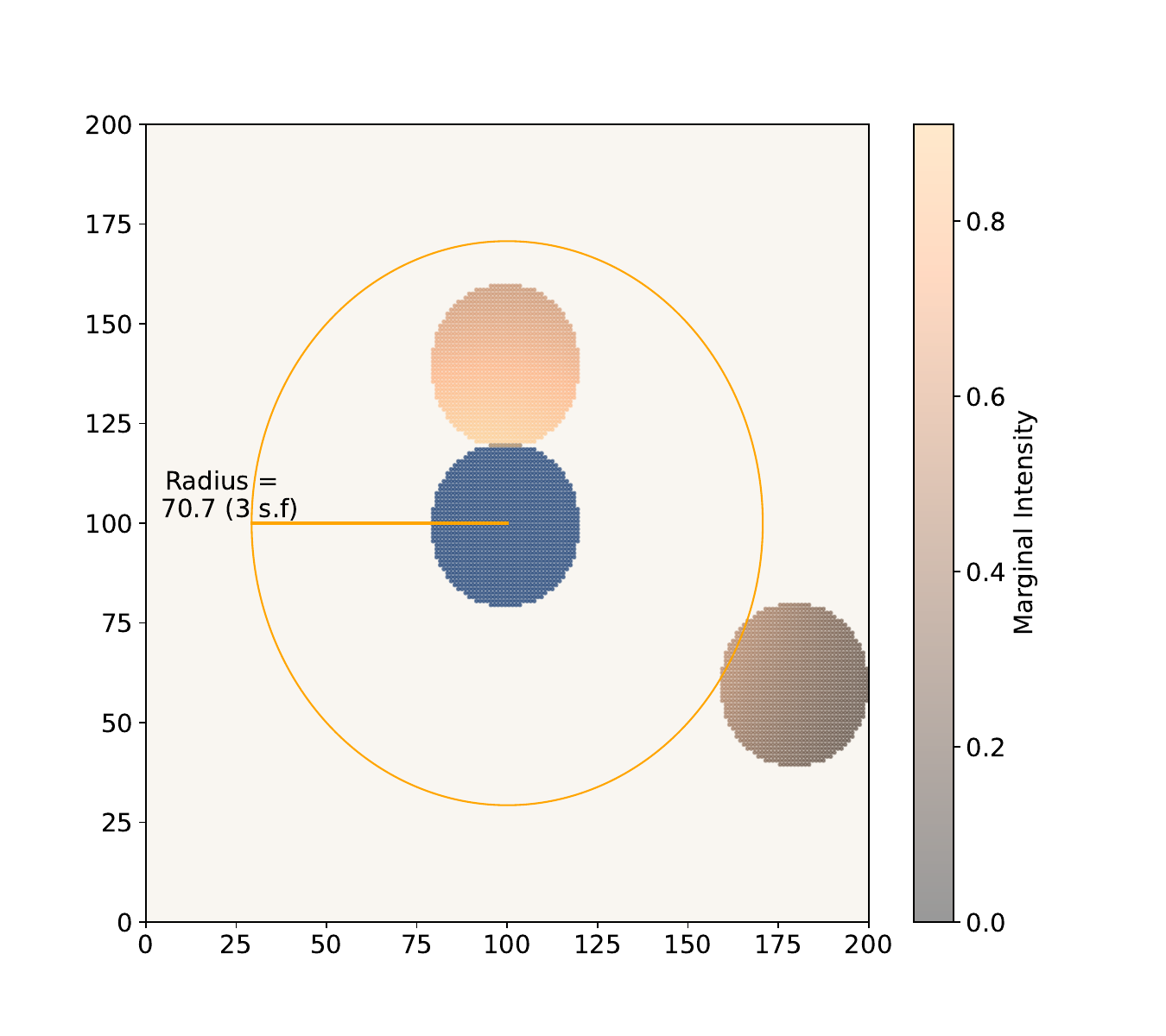}
        % \caption{}
        % \label{fig:c1c8_kl_marginal}
    \end{subfigure}
    \caption{Illustration of unbalanced displacement and \(\rho\) sensitivity with KL penalty. In both figures, the central observation C1 is shown in blue (darker), while the forecast is C6 (left) or C8 (right).
    Marginals \(\pi_1\) are plotted in pale orange to grey shades, not the forecast, highlighting regions where mass is destroyed. Left: there is again a relationship to the reach parameter although the mass is shared within half the reach and less strict than TV. The radius shown is equal to half the reach. Right: the closer event is assigned more intensity whilst the distance feature in this unbalanced setting is given very little. The radius is equal to the reach.  The figures were generated for \(\rho= 2^{-4} L^2, \varepsilon=0.005L^2\).}
    \label{fig:comparison_kl_marginal}
\end{figure}

 Figure \ref{fig:unbalanced_reach} then, explores the unbalanced setting \(\rho\) dependence. These results combine with C1C6, and Figure \ref{fig:comparison_kl_marginal}, which illustrate the returned marginals.
Notably, the strict reach interpretation is now ruined in this very unbalanced setting.
Instead, for the symmetric cases C1C6 the reach appears to have a halved influence, though when this symmetry is destroyed, mass is assigned to the closer feature. 
Moreover, the costing between KL and TV is now very different.
Here TV penalises the large mass imbalanced very harshly. 
The ATD and ATM then are suffering from a mean averaging downfall — they average out in the north-south shifts, and diagnoses only the west-east shift. This is a known issue and is overcome by viewing the whole spread of vectors direction and magnitude.

\begin{figure}[h]
    \centering
    \begin{subfigure}[b]{0.45\linewidth}
        \centering
        \includegraphics[width=\linewidth, trim= 50  40 50 50, clip]{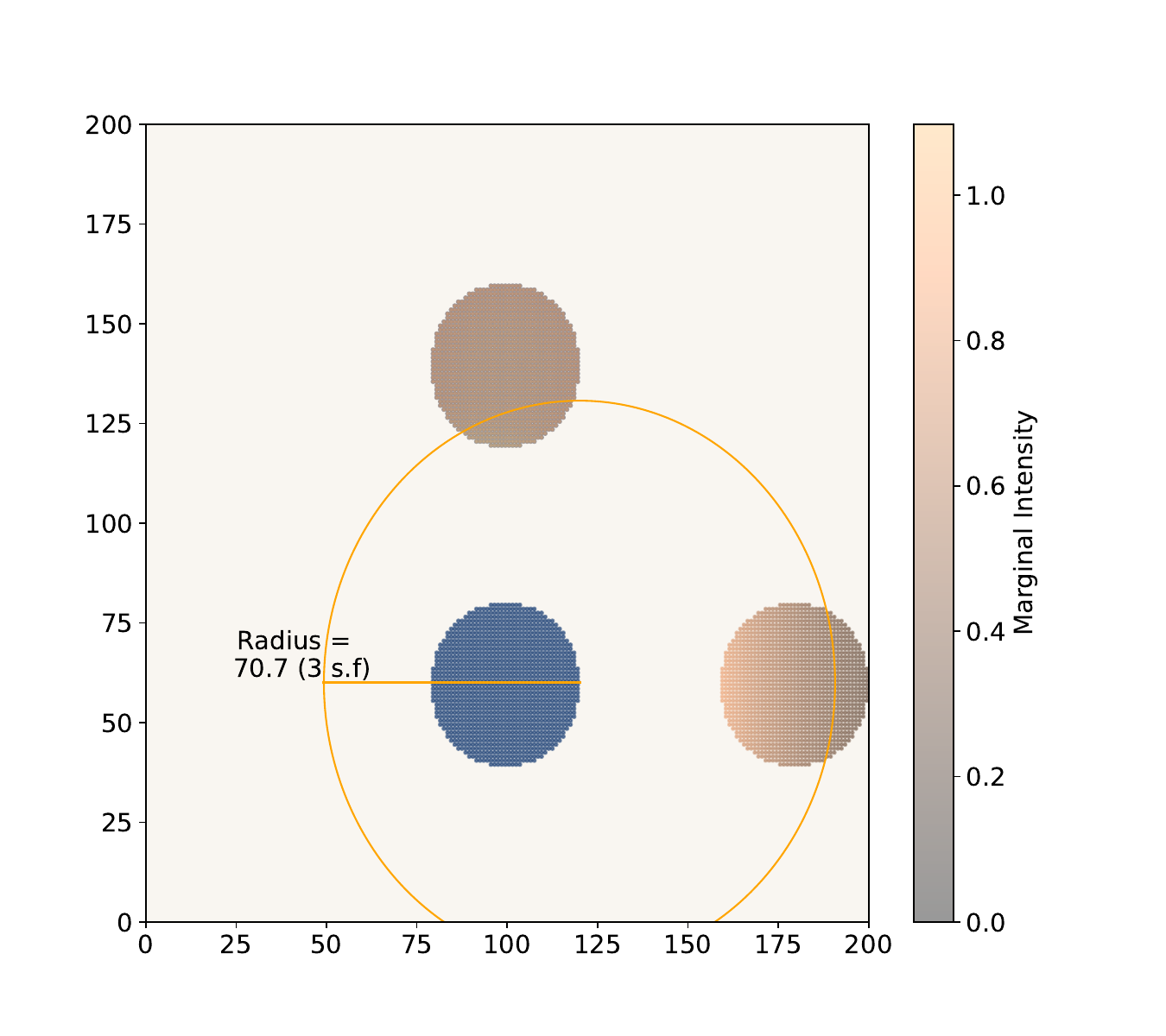}

    \end{subfigure}
    \hfill
    \begin{subfigure}[b]{0.45\linewidth}
        \centering
        \includegraphics[width=\linewidth, trim= 50  40 50 50, clip]{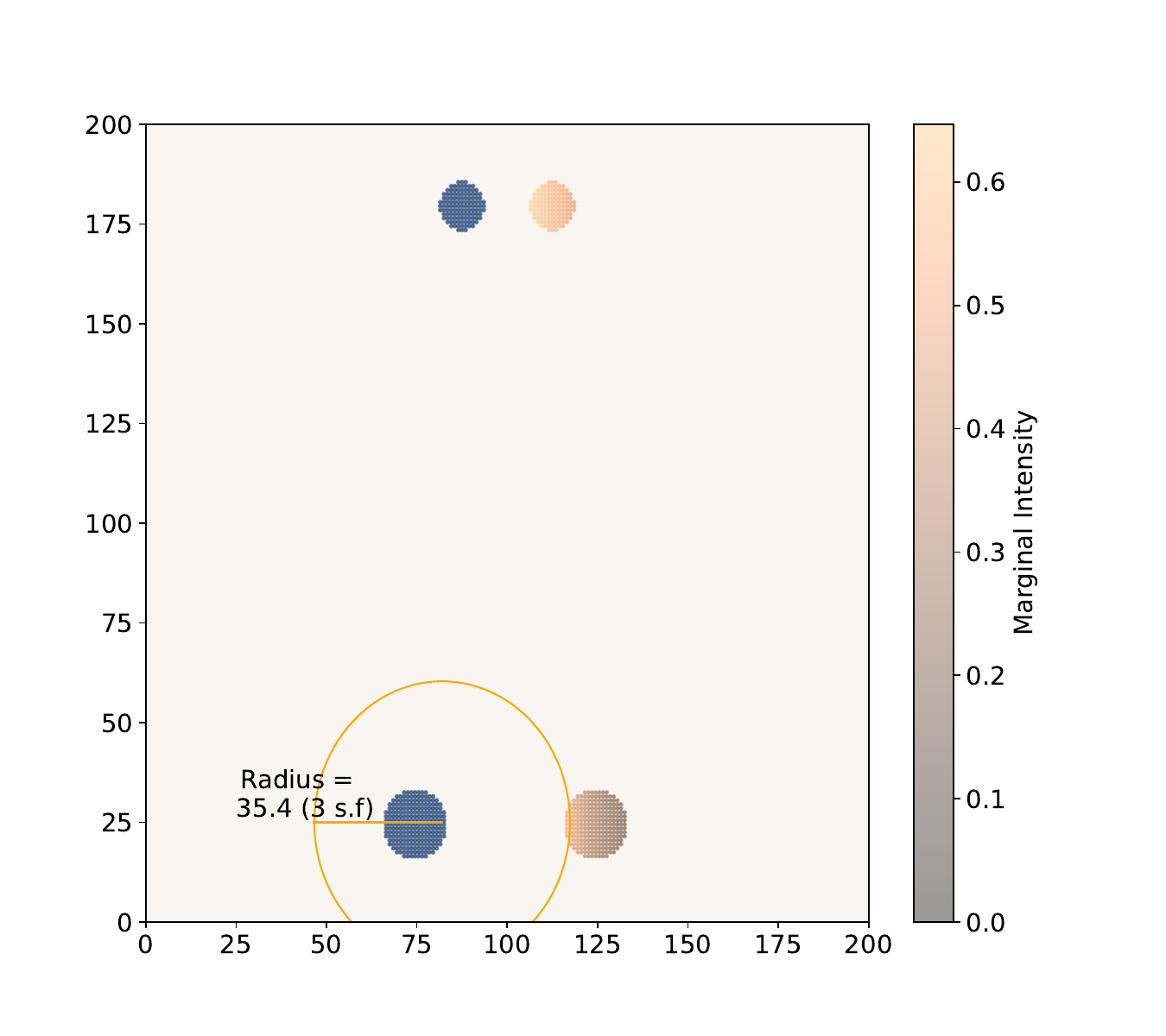}
       
    \end{subfigure}
    \caption{Illustration of unbalanced displacement and \(\rho\) sensitivity with KL penalty. 
    Marginals \(\pi_1\) are plotted in pale orange to grey shades, not the forecast, highlighting regions where mass is destroyed. Left: Case C6C8 with KL marginal penalty. The marginal \(\pi_1\) confirms there is a strict relationship with the reach parameter and it is maintained per feature. $\rho = 2^{-4} L^2$, the radius is equal to the reach. Right: Case C13C14 with KL marginal penalty. The marginal \(\pi_1\) shows that a distance feature can be removed if deemed too far separated, and this is maintained per feature. Here, the radius is equal to the reach, $\rho = 2^{-6} L^2$.}
    \label{fig:c6c8_c13c14_kl_marginal}
\end{figure}

\begin{figure}[h]
    \centering
    \includegraphics[width=0.85\linewidth]{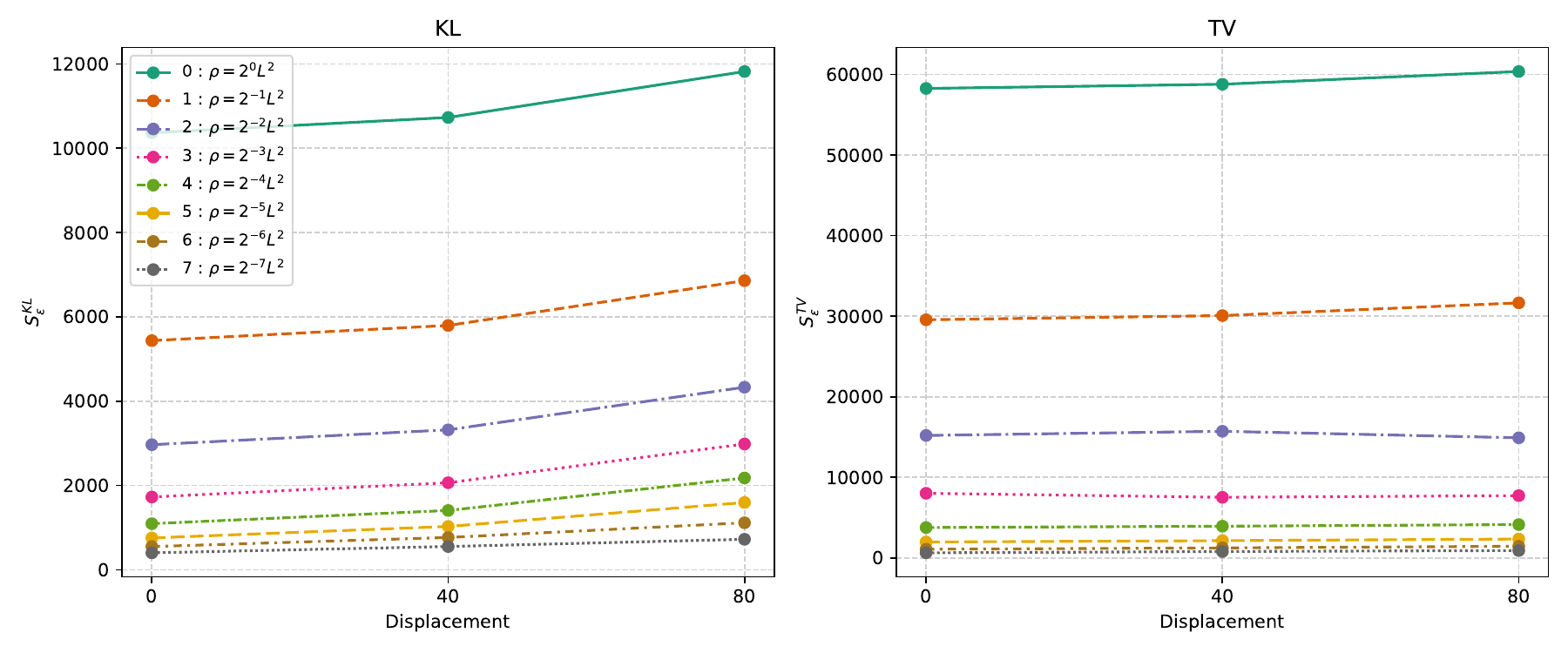}
    \includegraphics[width=0.85\linewidth]{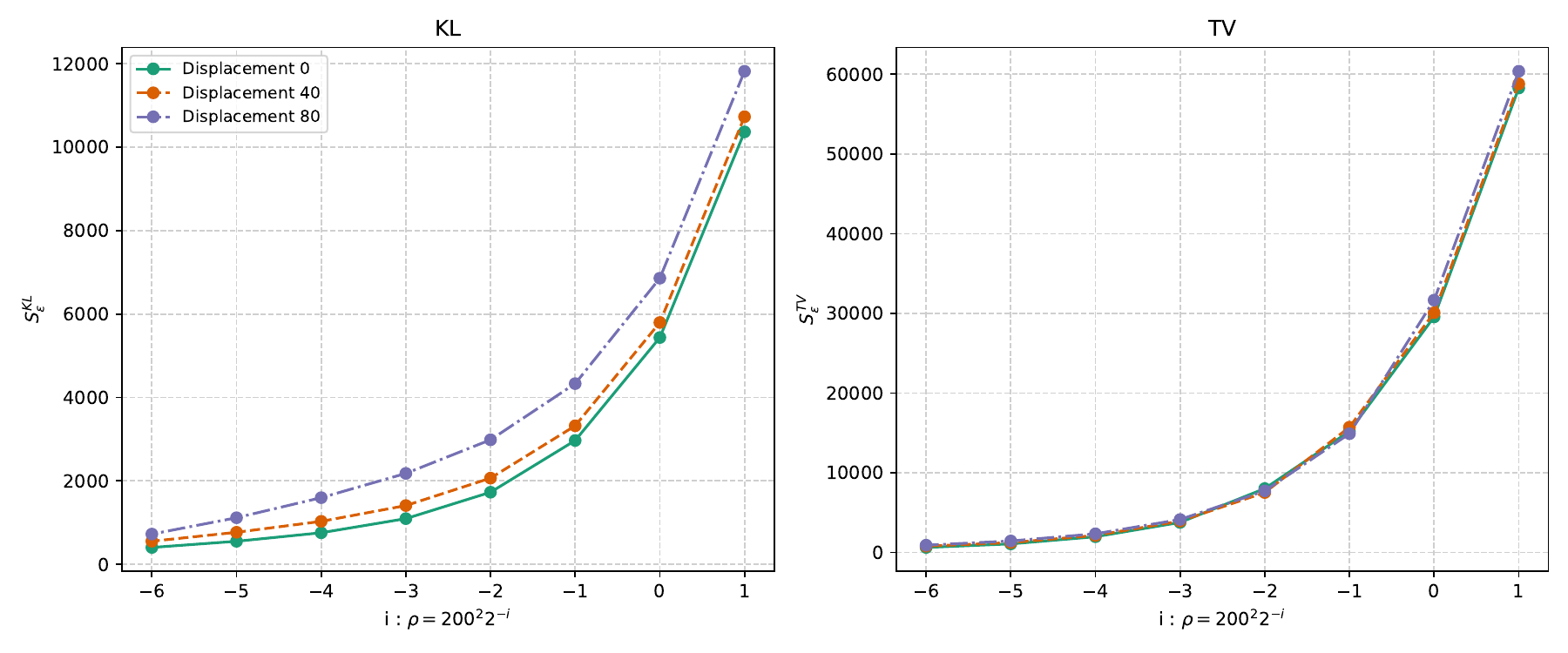}

    \caption{Reach sensitive study for case C1C6, C1C7, C1C8 which corresponds to a lower (southern) event being displaced away. This is now an unbalanced scenario and there is more extreme sensitivity to \(\rho\) than in the balanced setting. Notice that KL and TV have different scales. Compared to Figure \ref{fig:rho_balanced_reach} the relation is now swapped from dependence in displacement to \(\rho\), as imbalance costs dominate. Left: \(\Sink^{KL}\), Right: \(\Sink^{TV}\), Top: Cost verse displacement, with a curve per \(\rho\) values, Bottom: cost verse \(\rho\) values, with a curve per displacement case. \(\varepsilon = 0.005 L^2 \) . }
    \label{fig:unbalanced_displcement_rho}
\end{figure}

\subsubsection{Subsets}\label{appendix:subset_section}
\begin{figure}[h]
    \centering
    \includegraphics[width=0.85\linewidth]{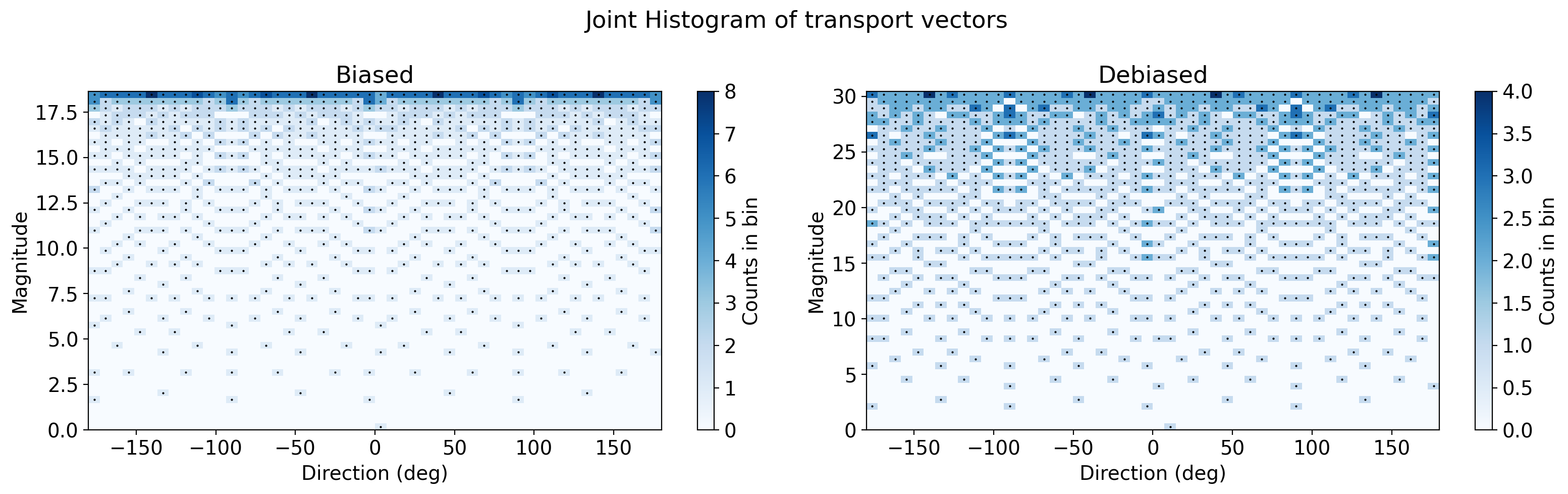}
    \caption{2D histogram of the magnitude and direction of the underlying transports vectors for the idealised case C1C9 which has a full \(360^{\circ}\) of necessary transport since one event is a subset of the other. The dotted bins, indicate those with non-zero mass. These results are based on TV marginal penalisation with parameters \(\varepsilon=0.005L^2, \rho=L^2\)} \label{fig:c1c9_2d_hist}
\end{figure}

\begin{figure}[h]
    \centering
    \includegraphics[width=0.75\linewidth]{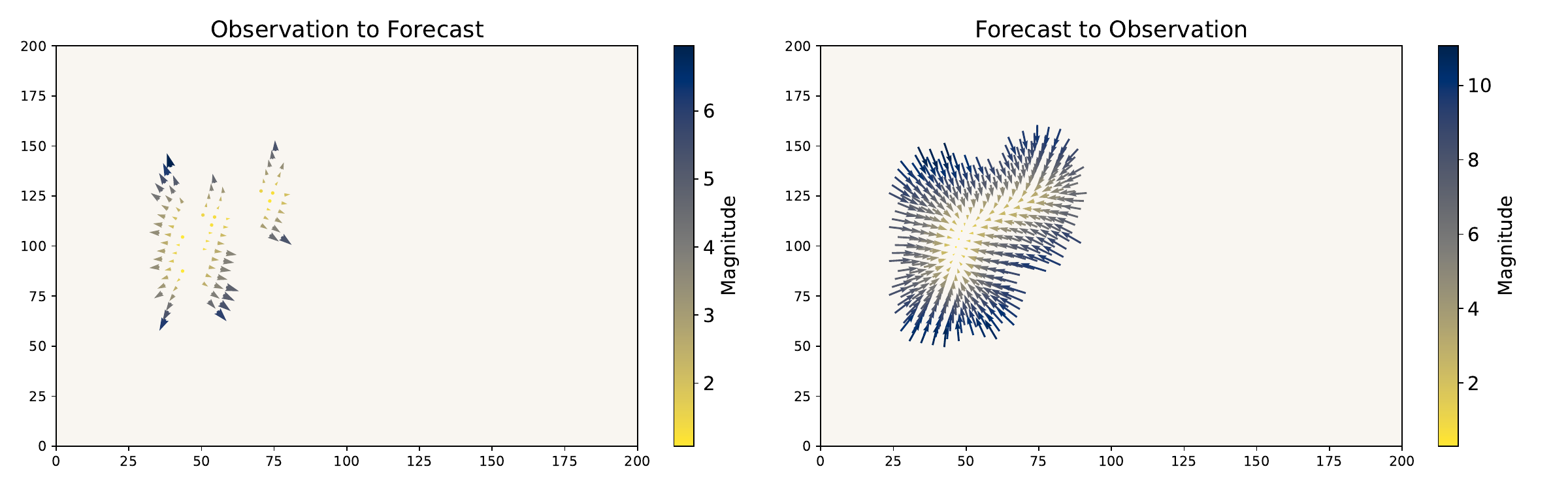}
    \includegraphics[width=0.75\linewidth]{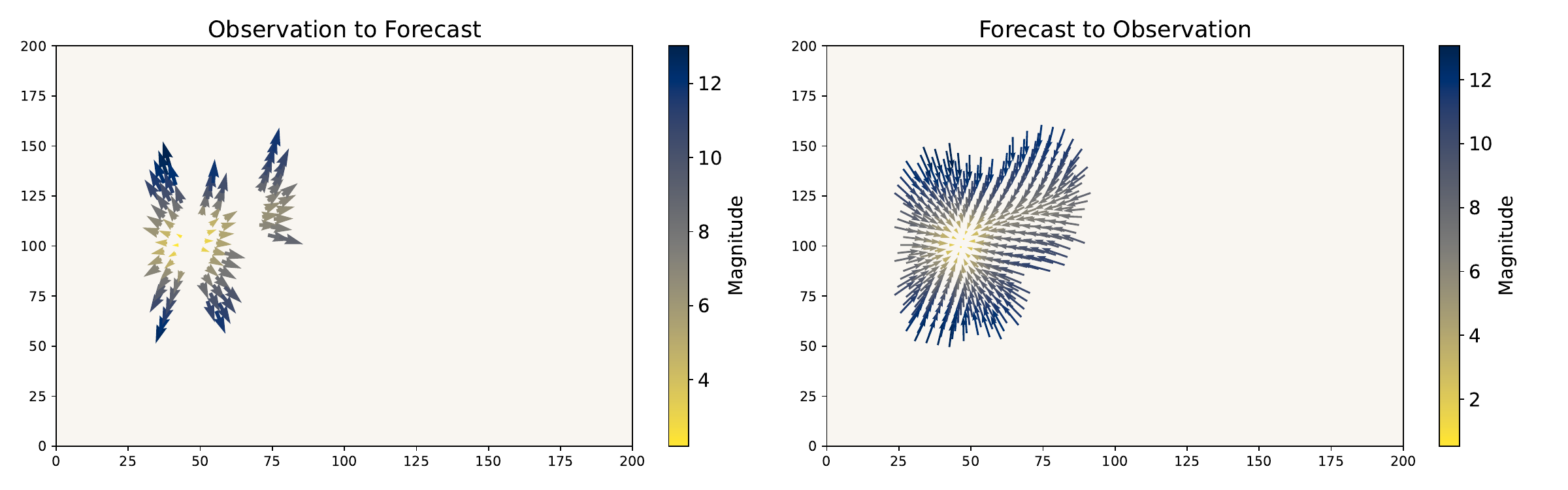}
    \caption{Illustration of the transports vectors for the idealised case E19E20 which has a full \(360^{\circ}\) of necessary transport since one event is a subset of the other. 
    Top: TV penalisation, Bottom: KL penalisation. Left: forward transport, right: inverse transport. \(\varepsilon=0.005L^2, \rho=L^2\)}
    \label{fig:e19e20_transport_vectors}
\end{figure}

\begin{figure}[h]
    \centering
    \includegraphics[width=0.75\linewidth]{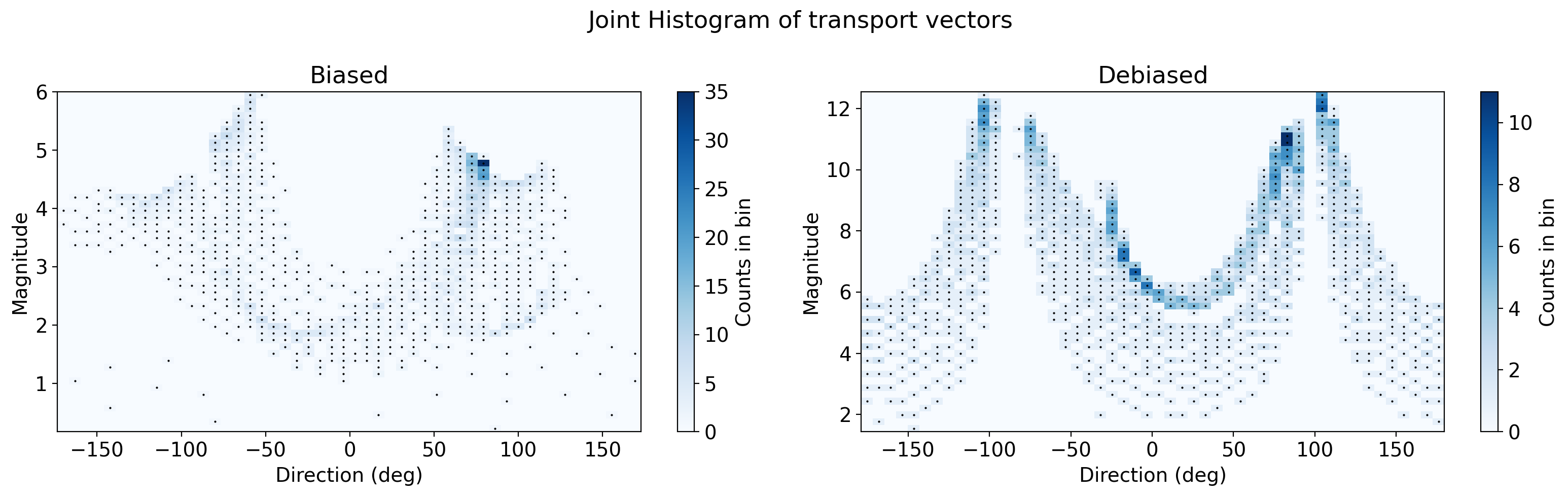}
    \caption{2D histogram of the magnitude and direction of the underlying transports vectors for the idealised case E19E20 which has a full \(360^{\circ}\) of necessary transport since one event is a subset of the other. The dotted bins, indicate those with non-zero mass. These results are based on KL marginal penalisation with parameters \(\varepsilon=0.005L^2, \rho=L^2\)}
    \label{fig:e19e20_kl_2d_hist}
\end{figure}

\begin{figure}[h]
    \centering
    \includegraphics[width=0.75\linewidth]{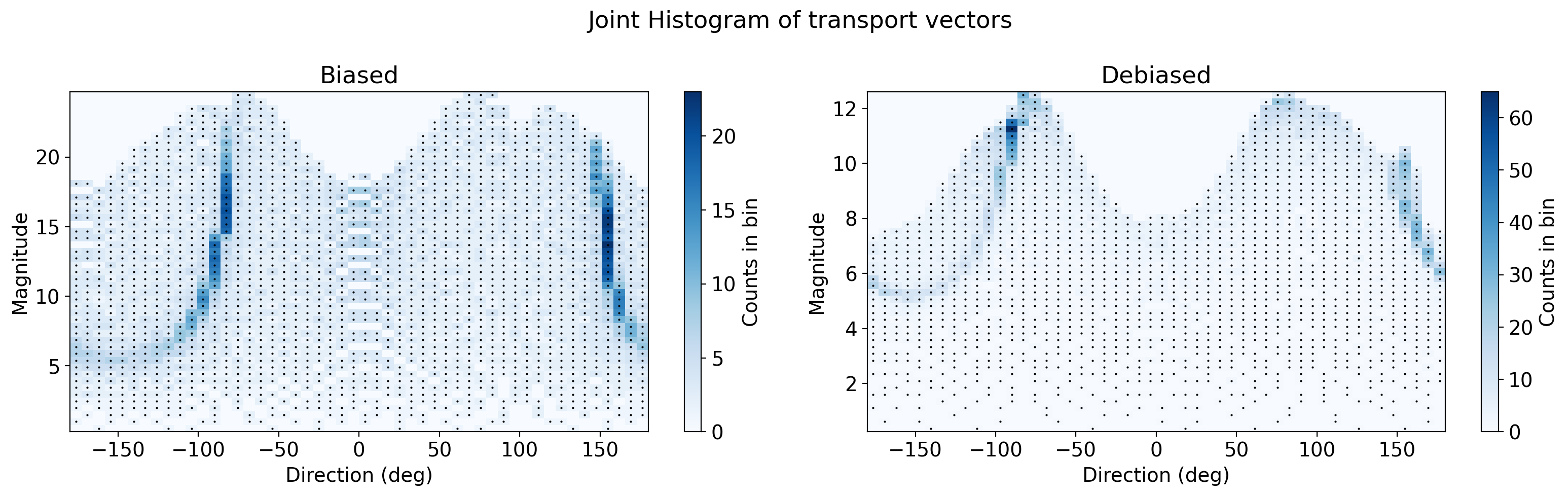}
    \caption{2D histogram of the magnitude and direction of the underlying transports vectors for the idealised case E20E19 (reverse of \ref{fig:e19e20_kl_2d_hist}) which has a full \(360^{\circ}\) of necessary transport since one event is a subset of the other. The dotted bins, indicate those with non-zero mass. These results are based on KL marginal penalisation with parameters \(\varepsilon=0.005L^2, \rho=L^2\)}
    \label{fig:e20e19_kl_2d_hist}
\end{figure}

\begin{figure}[h]
    \centering
    \includegraphics[width=0.75\linewidth, trim= 70 70 70 70, clip]{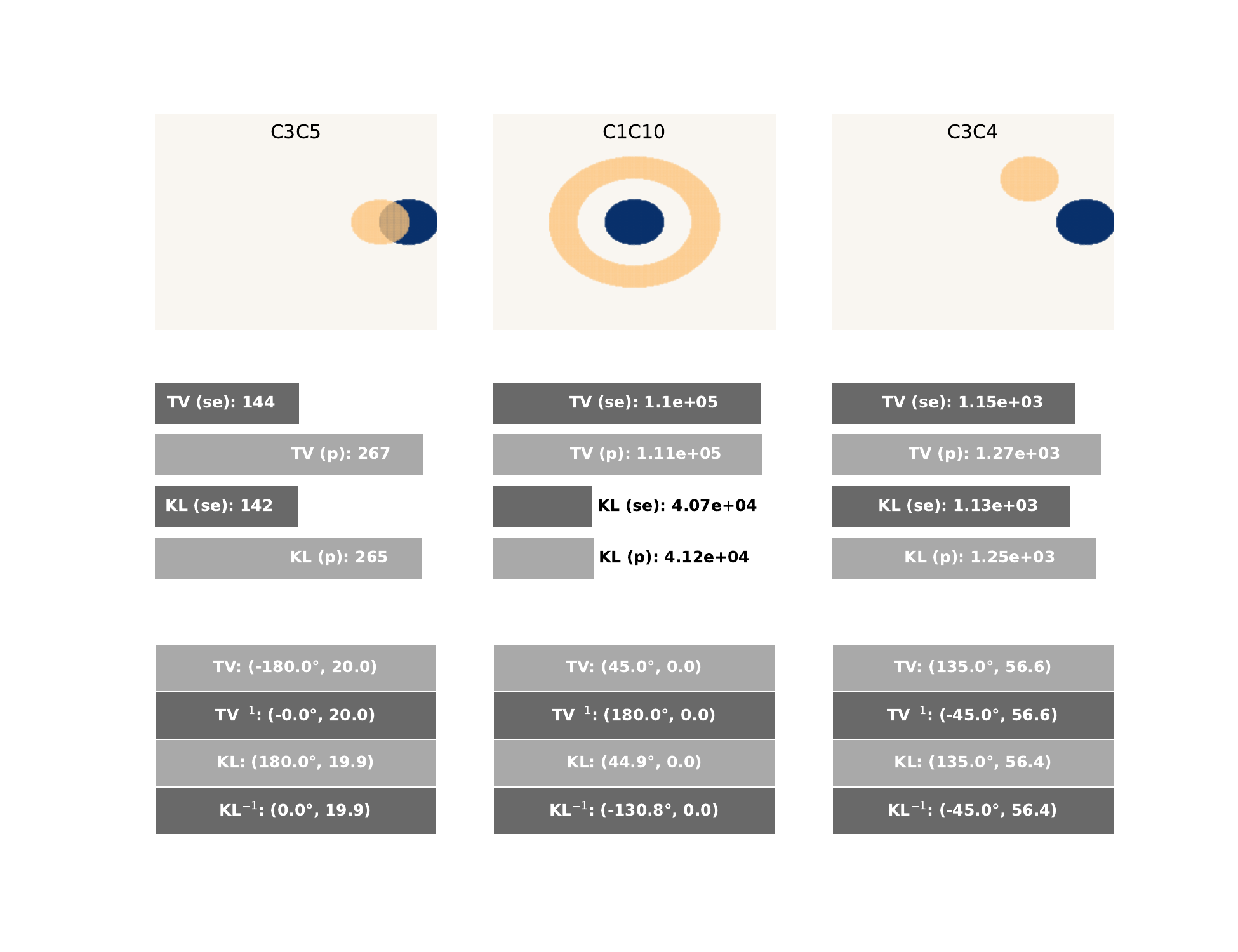}
    \caption{Cases C3C5, C1C10, C3C4. C3C5, C3C4 are pairs with C1C4, and C2C5 since they are reflected to shifted versions (relative to the boundary). C1C10 confirms the necessary but not sufficient condition of finding subsets of events, since the transport required here has a similar spread to C1C9 yet, C1 is not contained within C10.  The top four horizontal bars display; \(\Sink^{TV}, \UOT^{TV}, \Sink^{KL}, \UOT^{KL}\). The lower table presents the mean (ATD, ATM) in both flavours, and with the forward and inverse vectors. The blue (darker) colour indicates observations, while the pale orange (lighter) represents forecasts. \(\varepsilon = 0.005L^2,\ \rho = L^2\)}\label{fig:new_new_paper_circles_0}

\end{figure}

Figure \ref{fig:subsets} establishes UOT's behaviour when an event is a subset of another. 
Crucially, these are highly imbalanced cases and largely penalised for this large imbalanced.
Moreover, this alone does not allow for diagnosis of an event being a subset of another. 
The cost alone does not provide this information.
One potential route is through the 2D histogram of the transport vectors magnitude and directions; Figure \ref{fig:c1c9_2d_hist}, Figure \ref{fig:e19e20_kl_2d_hist} \& Figure \ref{fig:e20e19_kl_2d_hist}.
Observe for the simple cases C1C9 that the spread is uniform across a full \(360^{\circ}\), however for the non-symmetric case E19E20 this uniform spread is lost. 
However, there is still a spread across all angles.
This is a necessary but not sufficient condition of subset diagnosis. 
i.e. one could construct a case where the transport covers all angles but is not a subset of the other, e.g. C1C10 (Figure \ref{fig:new_new_paper_circles_0}) is case such a case. 
Otherwise, giving time for a deeper analysis, transport vectors can be visualised (See \ref{fig:e19e20_transport_vectors}). However in general without human input, UOT's score is not able to discern subsets.

\begin{figure}[h]
    \centering
    \includegraphics[width=0.75\linewidth]{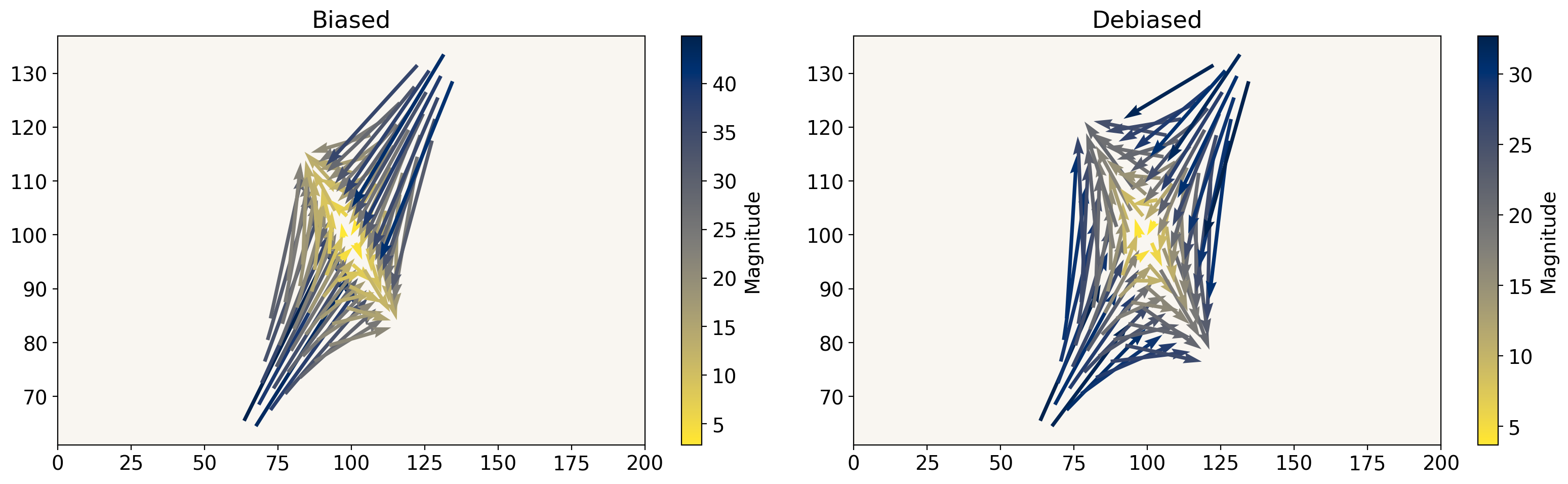}
    \caption{Biased vs debiased transport vector illustration with TV marginal penalty for the E2E4 case. Here the squeezing and stretching behaviour for rotated objects is demonstrated. Left: Biased UOT transport vectors, Right: Debiased UOT transport vectors. A regular sample of the vectors are shown to prevent overcrowding. \(\varepsilon = 0.005L^2,\ \rho = L^2\).}
    \label{fig:e2e4_transportvectors_tv}
\end{figure}

\subsubsection{Rotation and Scale for balanced events}\label{appendix:scale_rotation}

For a better understanding of how UOT treats rotation, consider Figure \ref{fig:e2e4_transportvectors_tv} \& Figure \ref{fig:e2e4_transportvectors_kl}, which illustrates the transport vectors for rotated objects E2E4 under TV and KL marginal penalisation respectively. 
Observer that the transport vectors squeeze and stretch the shape, not rotate. 
Instead, this follows an aspect-ratio correction, rather than rotation.
To examine this, E1E4 vs E2E4 is compared against the pure transport case of E1E9,  Figure \ref{fig:rotation_ellipse} and Figure \ref{fig:rotation_and_orientation} respectively.
E1E4 is scored significantly worse than E1E9 via \(\Sink^{TV}\) and better via \(\Sink^{KL}\).
For \(\Sink^{TV}\) most of the error appears in the marginal error terms, due to the imbalanced mass between E1 and E4.
In contrast, \(\Sink^{KL}\) keeps this penalty term lower.
For the \(90^{\circ}\) rotation, both flavours score worse than the pure transport, though now the mass is balanced so \(\Sink^{TV}\) scores this case better than E1E4, despite the rotation being larger. 
Notably, unlike CDST, \(\Sink^{TV/KL}\) is able to discern the two different rotations. CDST sees both rotations having the same centre of mass, causing a zero distance (this is also demonstrated in cases; C6C12, C1C6, C1C9 in Figures  \ref{fig:multiple_features}, \ref{fig:unbalanced_extent} and \ref{fig:subsets}).
However, the ATM cannot distinguish E1E4 vs E2E4, because the transport is equal and opposite.
On the other hand, \(\Sink^{KL}\) follows a similar trend to BDEL, where the scores indicate that E1E4 < E1E9 < E2E4.
This does not imply that UOT can take into account rotation: it cannot. 
Although, \(\Sink^{KL}\)  does appear to penalise the rotation in a way which may follow a more subjective evaluation. 
These cases further illustrate that \(\Sink^{TV}\)  is sensitive to mass imbalance. This shall be seen to be a major difference between the two flavours.
 
To complete the balanced cases, consider E6E14 and E2E10 (Figure \ref{fig:scaled_Cases}). 
Once more, the correct ATM and ATD are returned, with E6E14: \(\sqrt{8^2+10^2} = 12.8\), and E2E10 : \(\sqrt{15^2+20^2} = 25\).
\(\Sink^{TV/KL}\) ranks them in the expected order: E6E14 scores better (lower) than E2E10 since less effort is required to correct the transport error, and the areal extent is smaller. 
Moreover, the scores directly scale with spatial extent and transport error, since E2 and E10 are 4 times larger than E6 and E14, and there is approximately twice the transport error (in squared Euclidean cost), hence \(16.8 \cdot 4 \cdot 2^2 \sim 269\).
This establishes that in the balanced setting, the methodology follows a predictable scaling between pairs of events.

\subsubsection{E3E11, E7E3, E7E11 cost decomposition}\label{appendix:decomposition}

\begin{figure}[h]
    \centering
    \begin{minipage}{0.49\linewidth}
        \centering
        \includegraphics[width=\linewidth]{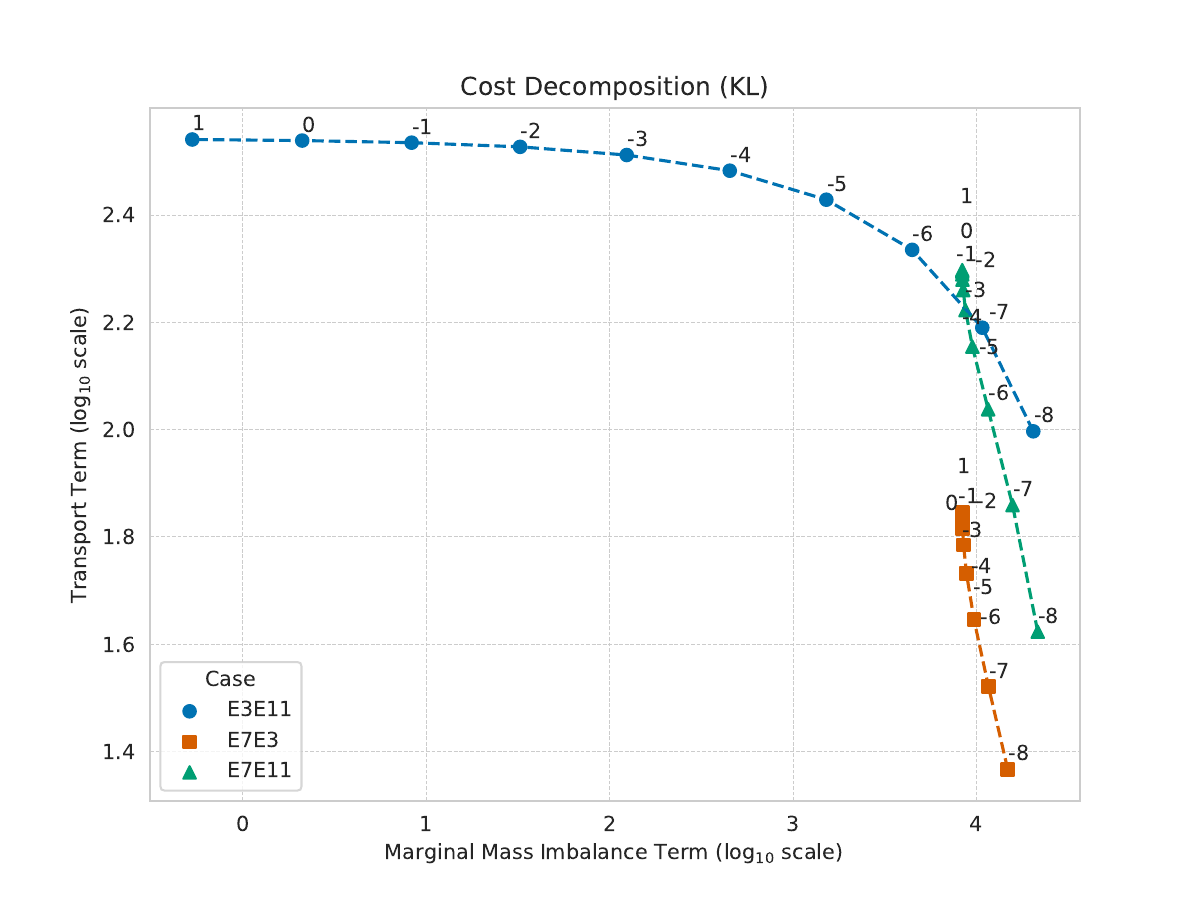}
    \end{minipage}
    \hfill
    \begin{minipage}{0.49\linewidth}
        \centering
        \includegraphics[width=\linewidth]{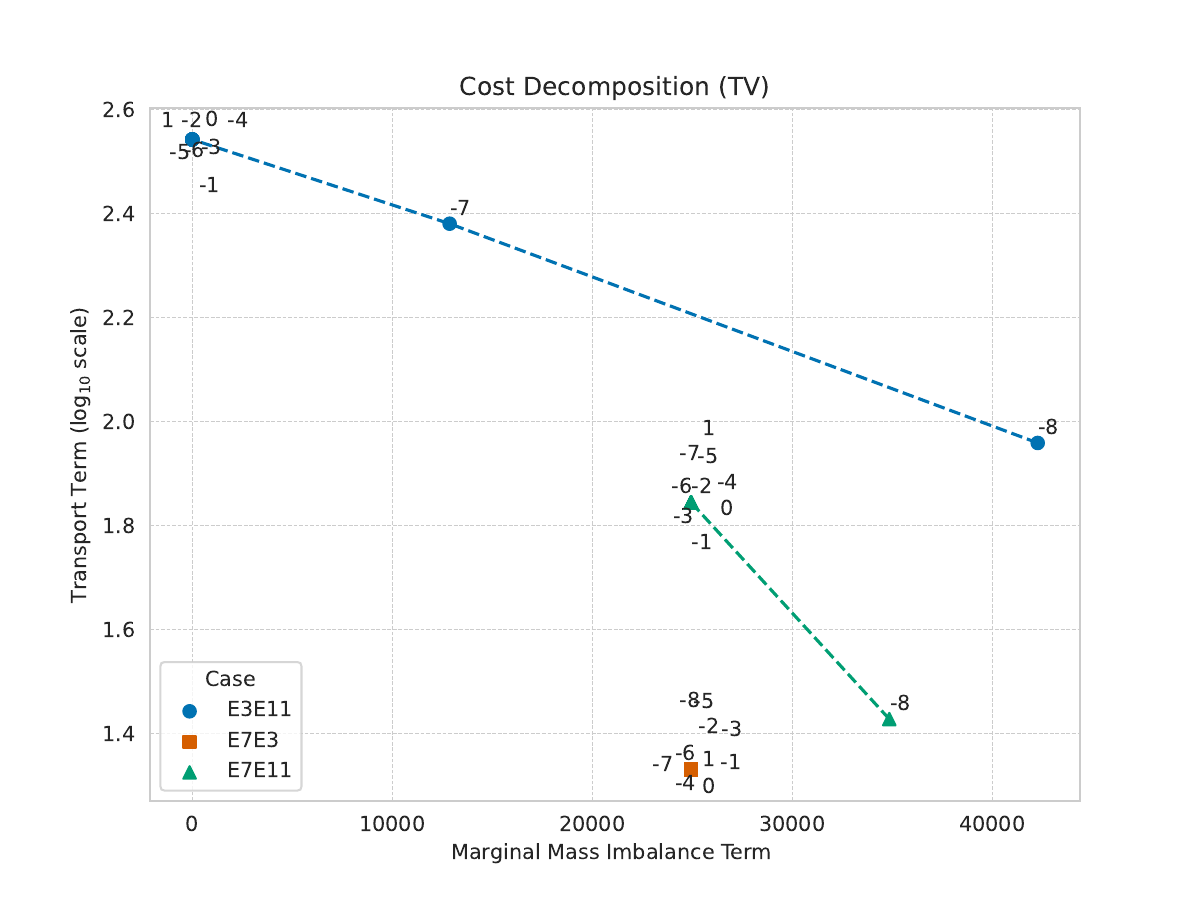}
    \end{minipage}
    \caption{Comparison of the Transport term (\(\sum_{i,j} c_{i,j} \pi_{i,j}\)) verse mass balanced, or penalty, terms \(D(\pi_0 | \mu_{O}) + D(\pi_1 | \mu_{F})\). Cases correspond to those examined in \ref{fig:hedging_forecast}, E3E11, E7E3, E7E11. Here the transport term is not debiased, since only \(\Sink\) is our debiased cost. The numbered nodes, correspond to different values of \(\rho : \rho = 200^2 2^{i}\) where \(i\) is the node number 1 to 8. Left: KL penalisation, right: TV penalisation. Notice, as the reach decreases the Transport term goes down, and balance up, since it is cheaper to destroy mass.   }
    \label{fig:ellipse_cost_distintergartion_across_rho}
\end{figure}

Additional insight is gained through Figure  \ref{fig:ellipse_cost_distintergartion_across_rho}. Here the cost is decomposed into its transport and marginal error, a form of mass imbalance, terms (The first, without \(\rho\), and two last terms in equation \ref{eq:uot_general_cost}).
Now E7E3 and E7E11 are discernable with equal marginal error terms, yet E7E11 has a larger transport term.
As the reach parameter changes the allowed transport is reduced and the marginal error increases, since its cheaper to destroy mass.
Again, behaviour in TV is sharper and more severe than in KL. Even for the smallest reach value tested, TV found it was not cheaper to destroy any more mass and instead opted to transport. However, there will be a point at which it becomes cheaper to start destroying mass as transport becomes even more localised (with a sufficiently small reach). KL offers similarly interpretable yet much more continuous, smooth results. 

\subsubsection{Hole, Edge and Extreme cases}\label{appendix:hole_edge}

The Hole case H1H2, compared against C1C2, determines how a method performs with increased spatial matching and compliments to an event.
\citet{gilleland_et_al_2019} made the observation that all the metrics  they studied (bar Zhu's) preferred H1H2.
UOT seems to offer similar behaviour, and \(\Sink^{KL/TV}\) scores lower for H1H2 but not for \(\UOT^{KL/TV}\). 
However, we suppose the reasoning to be different since most of the points do not require any transport. 
Whilst it looks as though the same transport is occurring, the compliment to C1C2 has the option to move the dips between the two sets. 
In this way, the average transport occurring is the same, but the total squared Euclidean costing is not the same\footnote{Another way to see this is to take any 1D slice in the x-direction that goes across the dip. Since the dip is present, it will move along to fill the new set and the points behind it move, recreating the dip. However, the dip is closer to start with and so requires less transport. For example, consider the 1D sets [O, X, O, O, O, O] transporting to [O, O, O, O, X, O]. This requires X to move 3 points which is 9 units in squared Euclidean costings. Whilst [X, O, X, X, X, X] transporting to [X, X, X, X, O, X] requires 3 single points to move but only 1 costing each.}. Considering the ATM, the returned score of 1.4 is correctly proportionally scaled against all those that do not move: \(1.4 \sim 40 \cdot \frac{1345}{200^2 - 1345}\) .

Then when comparing this full field with single point fields (P5, P3, P4 in Table \ref{tab:extreme_edge_cases}) the costs are large due to a lot of mass modification and small average transport.
These trends follow for P2P6. There are four points in P6, hence there is more mass available to distribute. The \(\Sink^{TV/KL}\) is lower in this case, but not significantly.
 
For P3P4 the transport cost is expected as \(199 \sqrt{2} \sim 281.4\) which is seen in the ATM.
In the balanced setting again \(\Sink^{TV}\) can attain the expected cost, whilst \(\Sink^{KL}\) is just below. 
When the reach is increased, both can achieve the expected score, whilst if it is reduced both eventually destroy  mass.
Notably, \(\UOT^{KL/TV}\) and \(\Sink^{KL/TV}\) are equal, since there is only one possible transport route between the two points.
Notice that P1P5 \(\Sink^{TV}\) measure the available mass (1 point) since 21 \(\sim 200^2 \cdot 1 / 1873.5\)

In summary, the edge cases demonstrate that much of the previous logic holds, where pure transport is diagnosed, and UOT remains sensitive for single points to the expected transport. However, this behaviour is  averaged with a larger number of points. Further, there are null cases which will need carefully defining.

\begin{table}[h]
    \centering
    \csvautobooktabularcenter{csvfiles/extreme.csv}
    \caption{Table of costings for the extreme and edge cases given \(\varepsilon = 0.005L^2\), \(\rho= L^2\). See Table \ref{table:binary_cases_here} for the description of the cases used here.}
    \label{tab:extreme_edge_cases}
\end{table}

\begin{figure}[h]
    \centering
    \includegraphics[width=1.0\linewidth, trim= 70 70 70 70, clip]{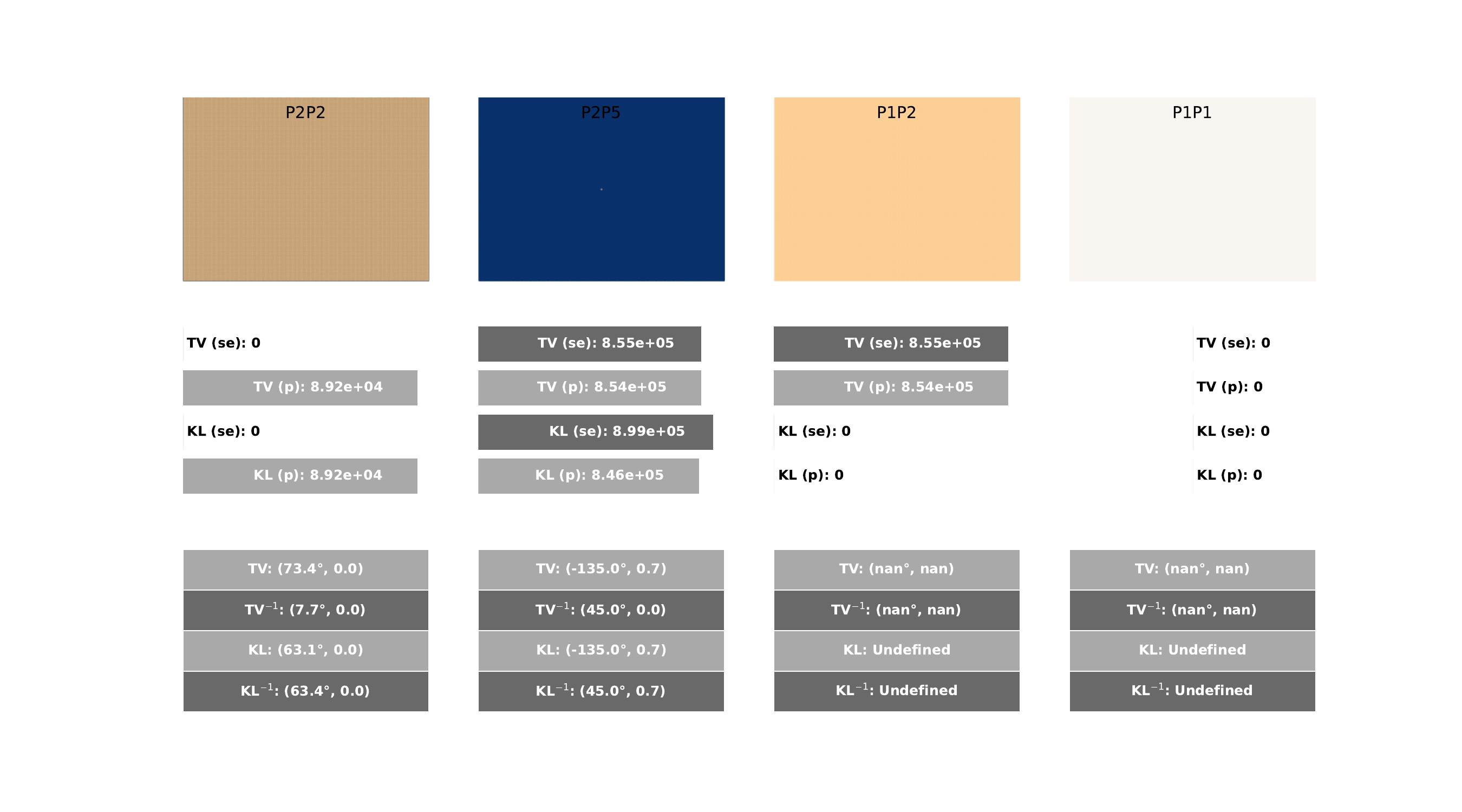}
    \caption{Cases P2P2, P2P5, P1P2, P1P1, recall that P1 is the null field, hence some quantities are undefined or zero. P2P2 returns the expected score for a perfect match of zero necessary transport. P1P2, P1P1 KL (se) and KL (p) are defined as zero. The top four horizontal bars display; \(\Sink^{TV}, \UOT^{TV}, \Sink^{KL}, \UOT^{KL}\). The lower table presents the mean (ATD, ATM) in both flavours, and with the forward and inverse vectors. The blue (darker) colour indicates observations, while the pale orange (lighter) represents forecasts. \(\varepsilon = 0.005L^2,\ \rho = L^2\)}\label{fig:new_new_paper_pcase_1}
\end{figure}

\begin{figure}[h]
    \centering
    \includegraphics[width=0.75\linewidth, trim= 70 70 70 70, clip]{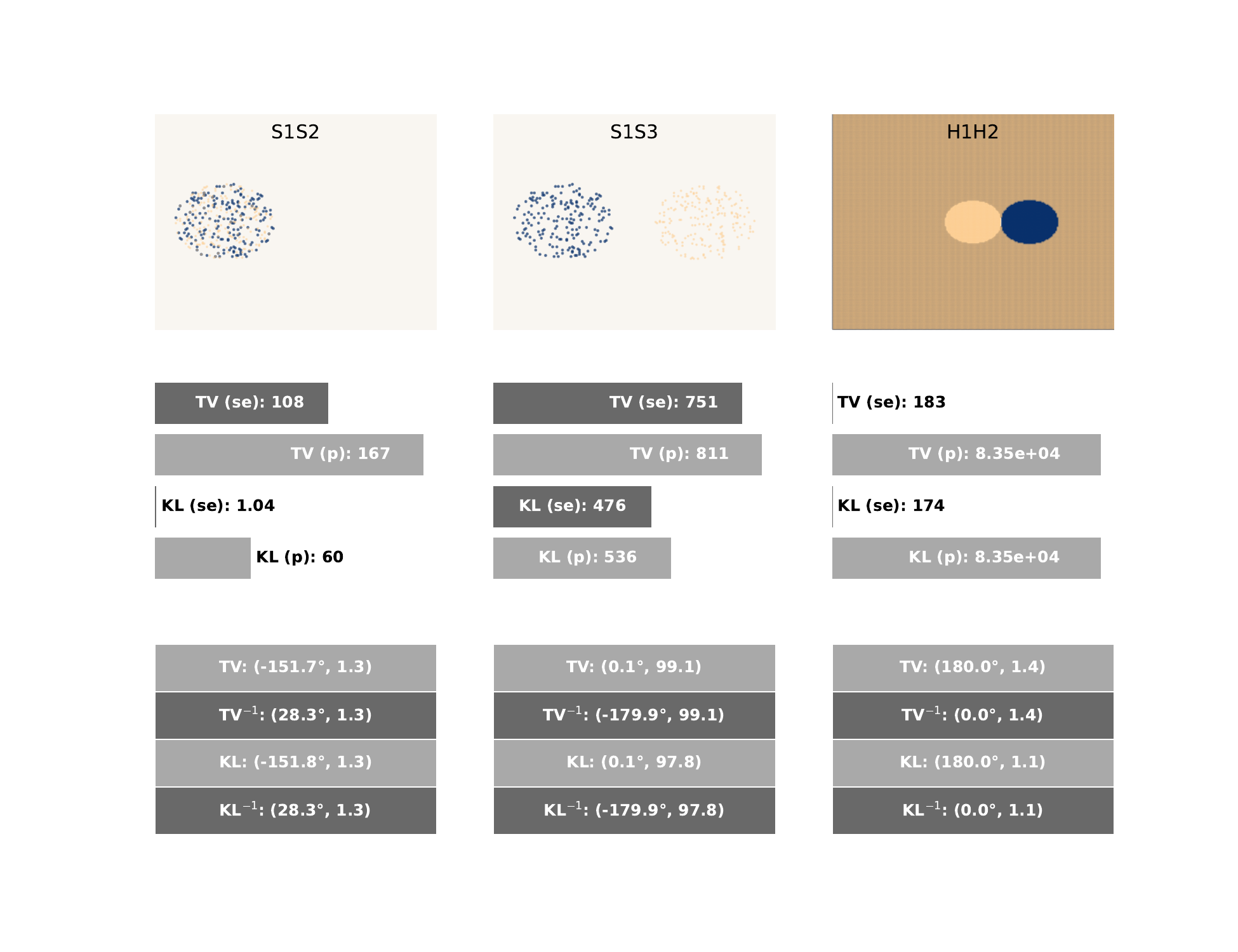}
    \caption{Scattered and compliment (hole) cases illustrating behaviour with small shower type events contained within an envelope and compliment logic or complete coverage with holes.
    The top four horizontal bars display; \(\Sink^{TV}, \UOT^{TV}, \Sink^{KL}, \UOT^{KL}\). The lower table presents the mean (ATD, ATM) in both flavours, and with the forward and inverse vectors. The blue (darker) colour indicates observations, while the pale orange (lighter) represents forecasts. \(\varepsilon = 0.005L^2,\ \rho = L^2\).}
    \label{fig:scattered_hole}
\end{figure}
\begin{figure}[h]
    \centering
    \includegraphics[width=1.0\linewidth, trim= 70 70 70 70, clip ]{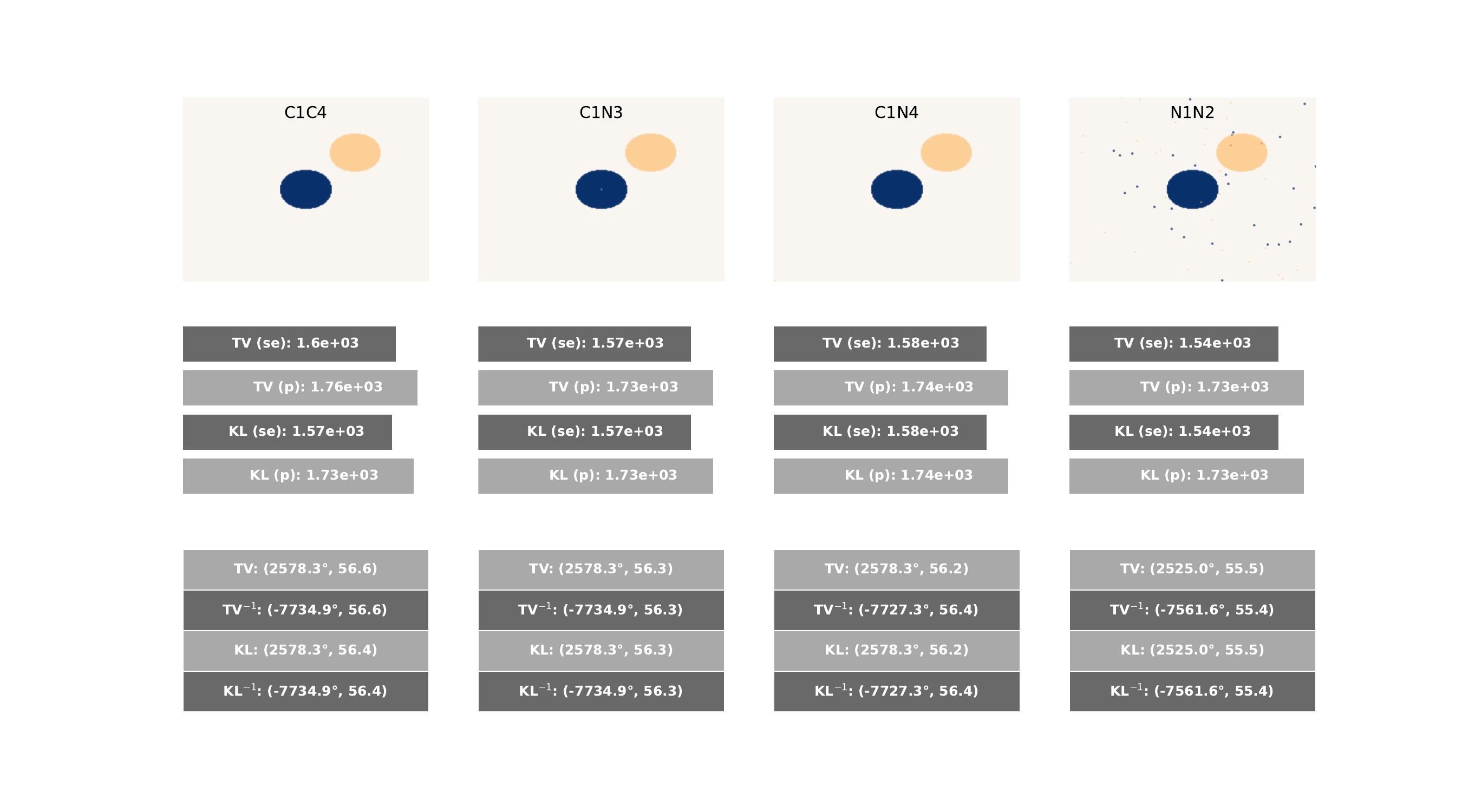}
    \caption{Noisy cases illustrating resistance to error in quality control or anomalous data depending on the type of marginal penalties. The top four horizontal bars display; \(\Sink^{TV}, \UOT^{TV}, \Sink^{KL}, \UOT^{KL}\). The lower table presents the mean (ATD, ATM) in both flavours, and with the forward and inverse vectors. The blue (darker) colour indicates observations, while the pale orange (lighter) represents forecasts. \(\varepsilon = 0.005L^2,\ \rho = L^2\).}
    \label{fig:noise}
\end{figure}

\begin{figure}[h]
    \centering
    \begin{minipage}{0.49\linewidth}
        \centering
        \includegraphics[width=\linewidth]{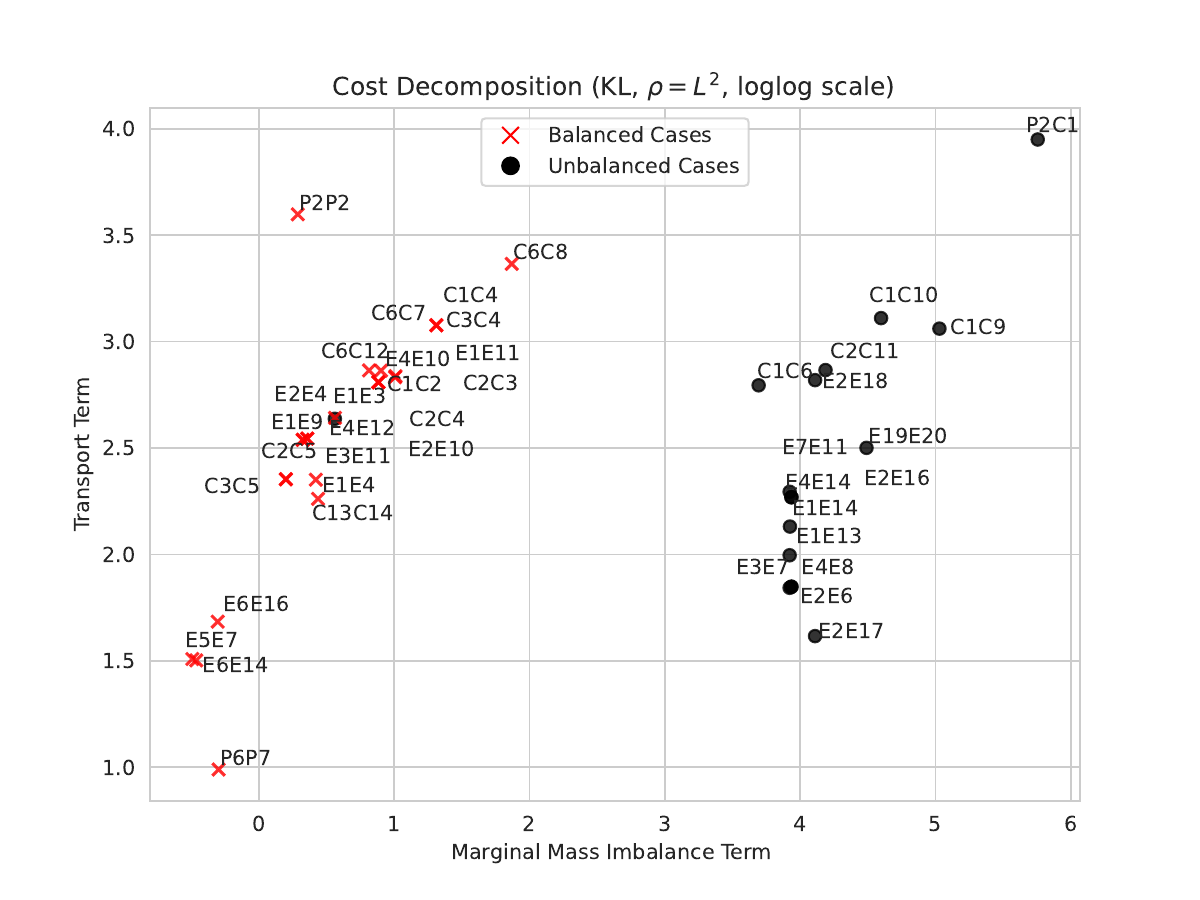}
    \end{minipage}
    \hfill
    \begin{minipage}{0.49\linewidth}
        \centering
        \includegraphics[width=\linewidth]{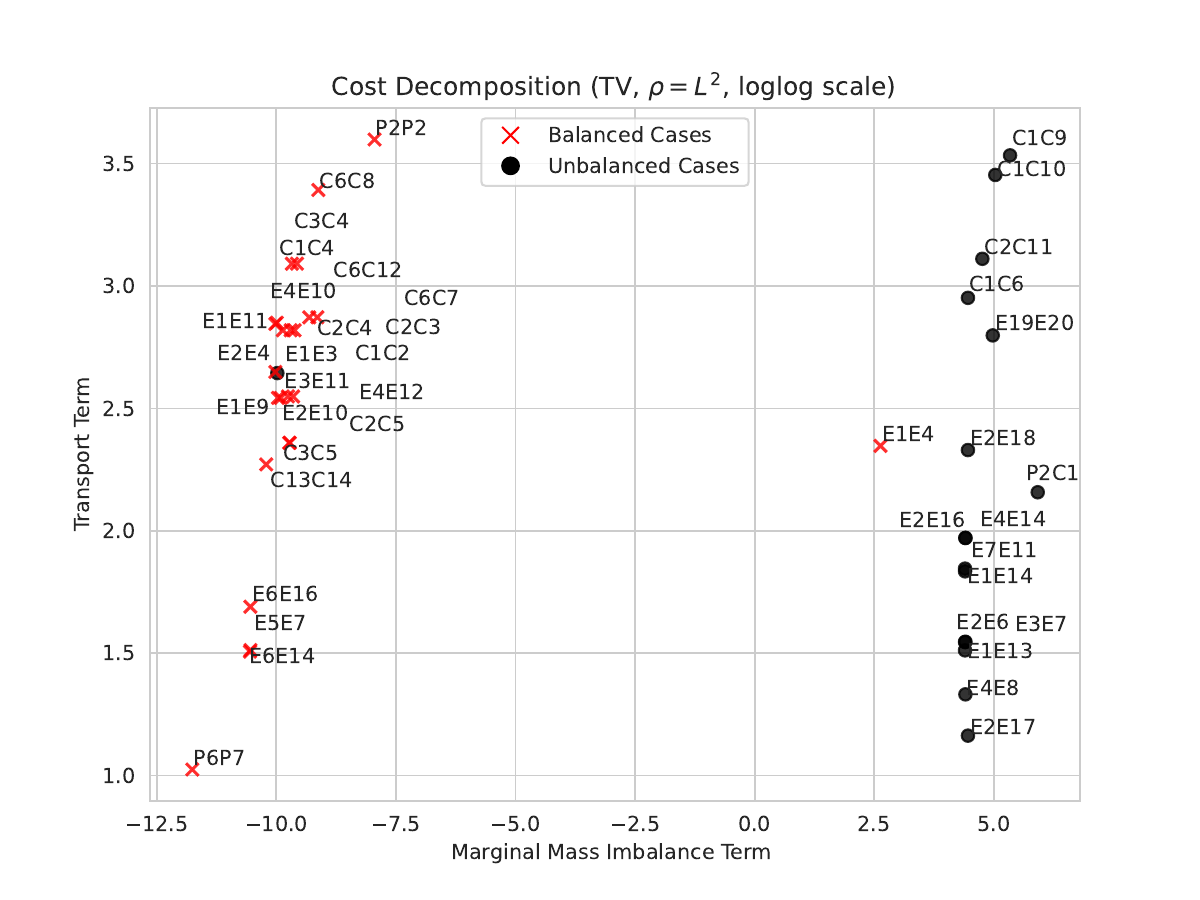}
    \end{minipage}
    \caption{Comparison of the Transport term (\(\sum_{i,j} c_{i,j} \pi_{i,j}\)) verse mass balanced, or penalty, terms \(D(\pi_0 | \mu_{O}) + D(\pi_1 | \mu_{F})\). Cases correspond to those in the original ICP study \cite{gilleland_et_al_2019.} Notice that these terms are not debiased, which can be seen as case P2P2 has some non-zero transport cost, yet achieves zero for \(\Sink^{KL/TV}\).
    Balanced and unbalanced cases are connected by colour. E2E4 is almost balanced and so is E1E4 hence these appear on the opposite cluster depending on the penalty.
    Left: KL penalisation, Right: TV penalisation. Notice E1E4 falls in the balanced cluster within KL but not TV, due to TV being more sensitive to mass imbalanced (by eye the shapes are similar and differ only by a few grid points). We can see the KL penalisation is more smooth in the balance term, whilst TV returns close to zero score for all balanced cases. }
    \label{fig:cost_disintergration}
\end{figure}

\subsection{Real textured cases}\label{appendix:real_texture}

\begin{figure}[h]
    \centering
    \includegraphics[width=0.8\linewidth]{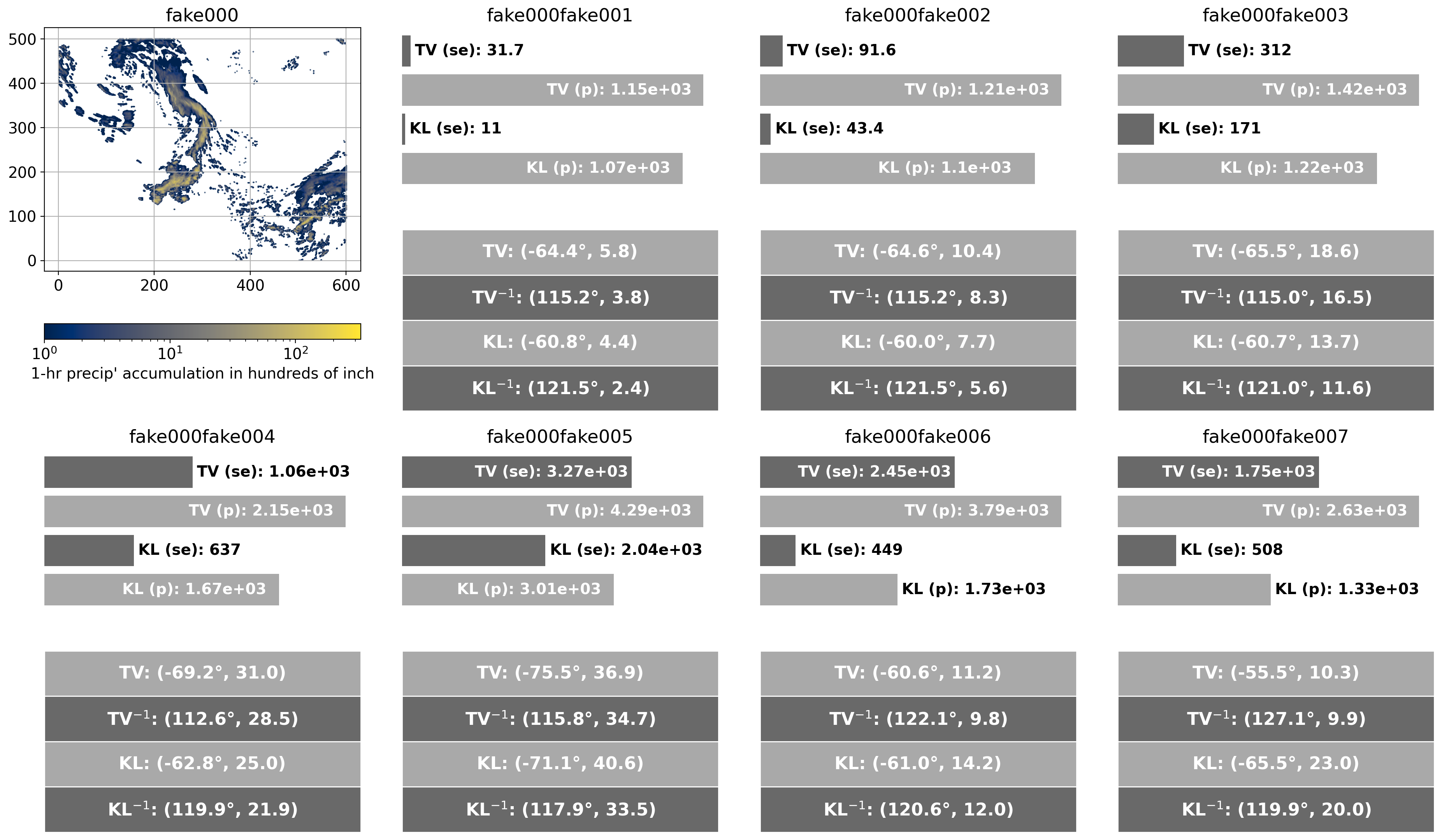}
    \caption{Perturbed cases with real textured intensities and smaller \(\rho\). All cases are a changed version of fake000, shown in the top left. The top four horizontal bars display; \(\Sink^{TV}, \UOT^{TV}, \Sink^{KL}, \UOT^{KL}\). The lower table presents the mean (ATD, ATM) in both flavours, and with the forward and inverse vectors. The blue (darker) colour indicates observations, while the pale orange (lighter) represents forecasts. \(\varepsilon = 0.001, \rho = 2^{-6}L^2\), reach \(\sim 106\).}
    \label{fig:perturbed_cases_smallrho}
\end{figure}

\begin{figure}[h]
    \centering
    \includegraphics[width=0.8\linewidth]{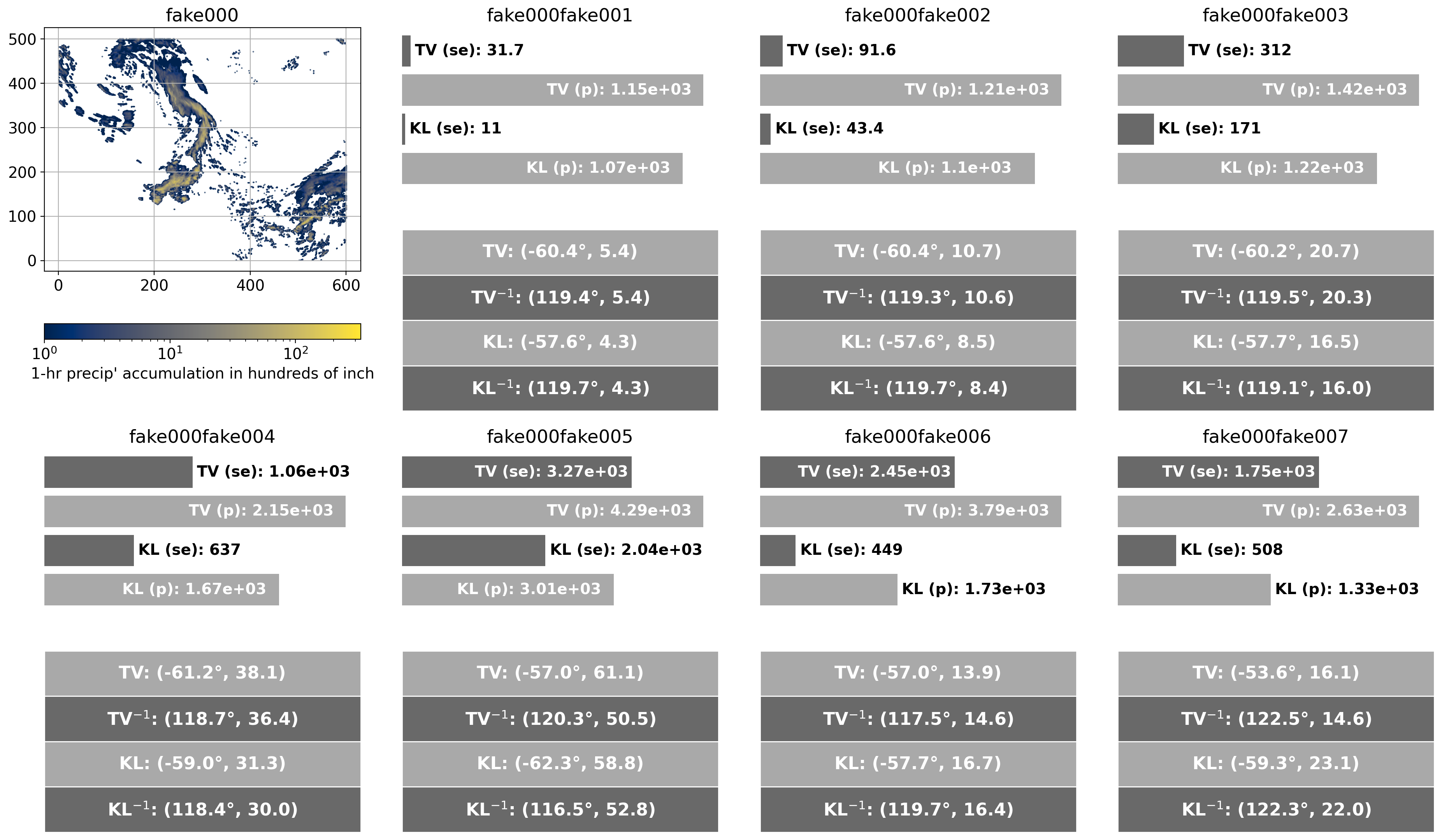}
    \caption{Median ATM version of \ref{fig:perturbed_cases_smallrho}.
    The top four horizontal bars display; \(\Sink^{TV}, \UOT^{TV}, \Sink^{KL}, \UOT^{KL}\). The lower table presents the median ATM and ATD in both flavours, and with the forward and inverse vectors. The blue (darker) colour indicates observations, while the pale orange (lighter) represents forecasts. \(\varepsilon = 0.001, \rho = 2^{-6}L^2\), reach \(\sim 106\).}
    \label{fig:fake_cases_median_rhominus6}
\end{figure}

Within the perturbed cases, Figure \ref{fig:perturbed_cases_smallrho} demonstrates the returned cost for the perturbed fake cases, now with a small reach parameter.
This means that the transport is kept more local.
This has the effect that TV is able to diagnose more reasonably to the expected score as we k ow transport is not lag,e and smaller reach keep transport local.
KL still performed better in this imbalanced setting. 
Additionally, both improve their ATM and ATD since only local transport is allowed. 
However, in the more perturbed cases, which now go above the reach, the ATD and ATM struggle to align with the expected displacement: \(5.8\), approximately doubling sequentially \((11.7, 23.3, 46.6, 93.3)\), and in the forward direction at a maintained \(-59^{\circ}\).
However, in Figure \ref{fig:fake_cases_median_rho1} \& Figure \ref{fig:fake_cases_median_rhominus6}, the ATM and ATD are restored, closer to the expected values, simply by changing the averaging from a mean average to the median.

\begin{figure}[h]
    \centering
    \includegraphics[width=0.8\linewidth]{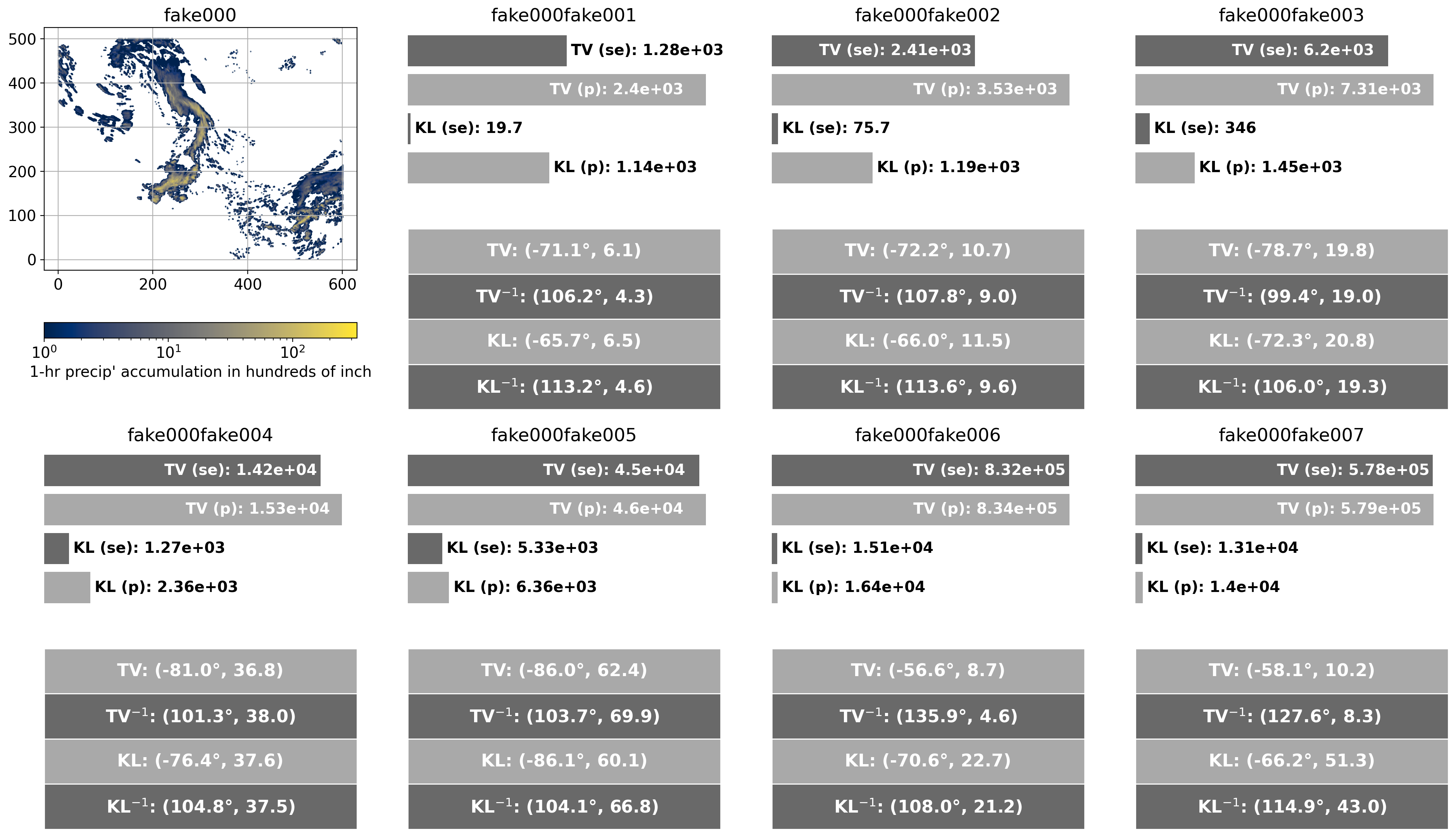}
    \caption{Perturbed cases with real textured intensities. All cases are a changed version of fake000, shown in the top left. The top four horizontal bars display; \(\Sink^{TV}, \UOT^{TV}, \Sink^{KL}, \UOT^{KL}\). The lower table presents the mean (ATD, ATM) in both flavours, and with the forward and inverse vectors. The blue (darker) colour indicates observations, while the pale orange (lighter) represents forecasts. \(\varepsilon = 0.001, \rho = L^2\), reach \(\sim 849\).}
    \label{fig:perturbed_cases}
\end{figure}

\begin{figure}[h]
    \centering
    \includegraphics[width=0.8\linewidth]{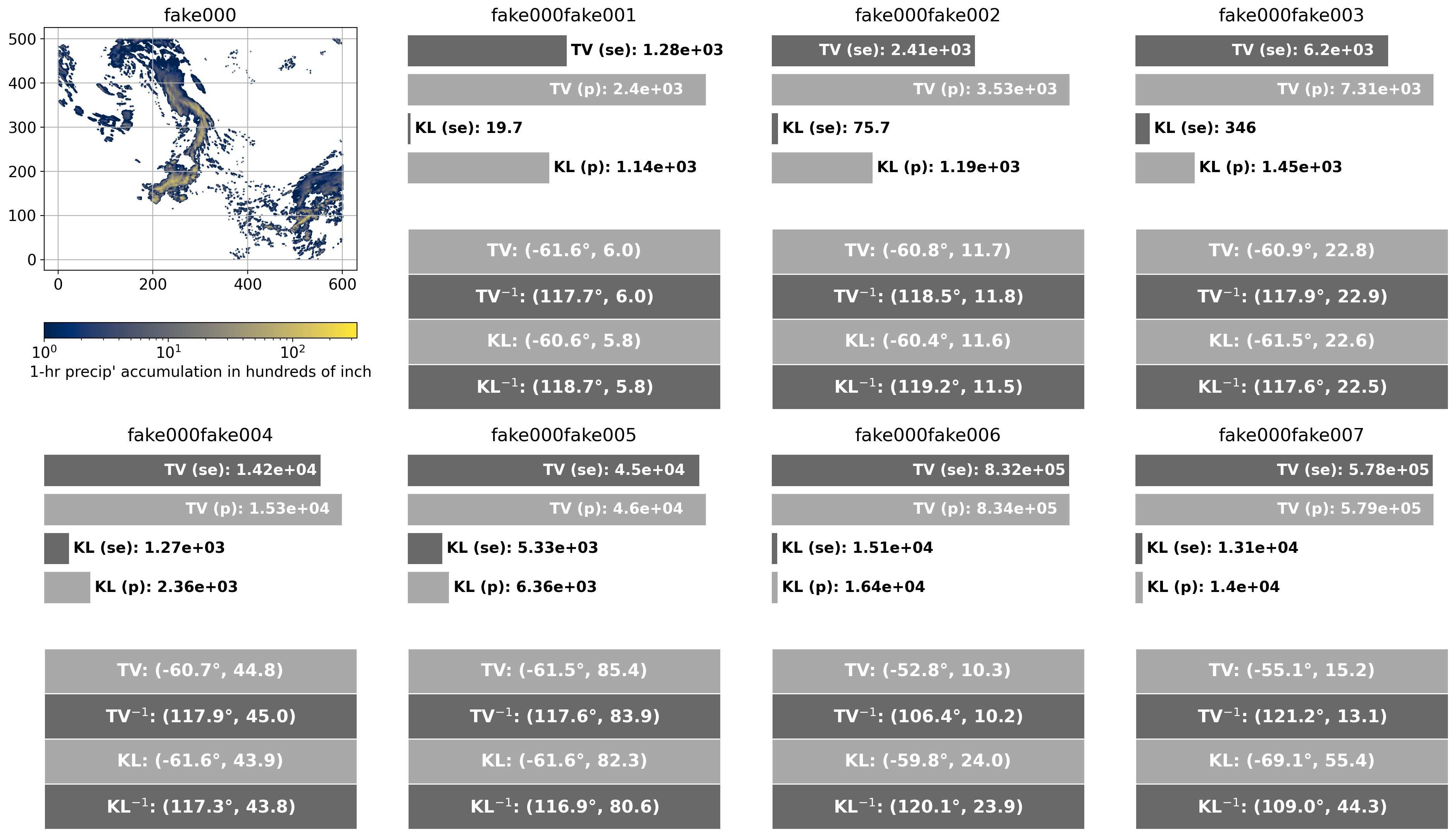}
    \caption{Median ATM version of \ref{fig:perturbed_cases}.
    The top four horizontal bars display; \(\Sink^{TV}, \UOT^{TV}, \Sink^{KL}, \UOT^{KL}\). The lower table presents the median ATM and ATD in both flavours, and with the forward and inverse vectors. The blue (darker) colour indicates observations, while the pale orange (lighter) represents forecasts. \(\varepsilon = 0.001, \rho = L^2\), reach \(\sim 849\).}
    \label{fig:fake_cases_median_rho1}
\end{figure}

\begin{figure}
    \centering
    \includegraphics[width=0.5\linewidth]{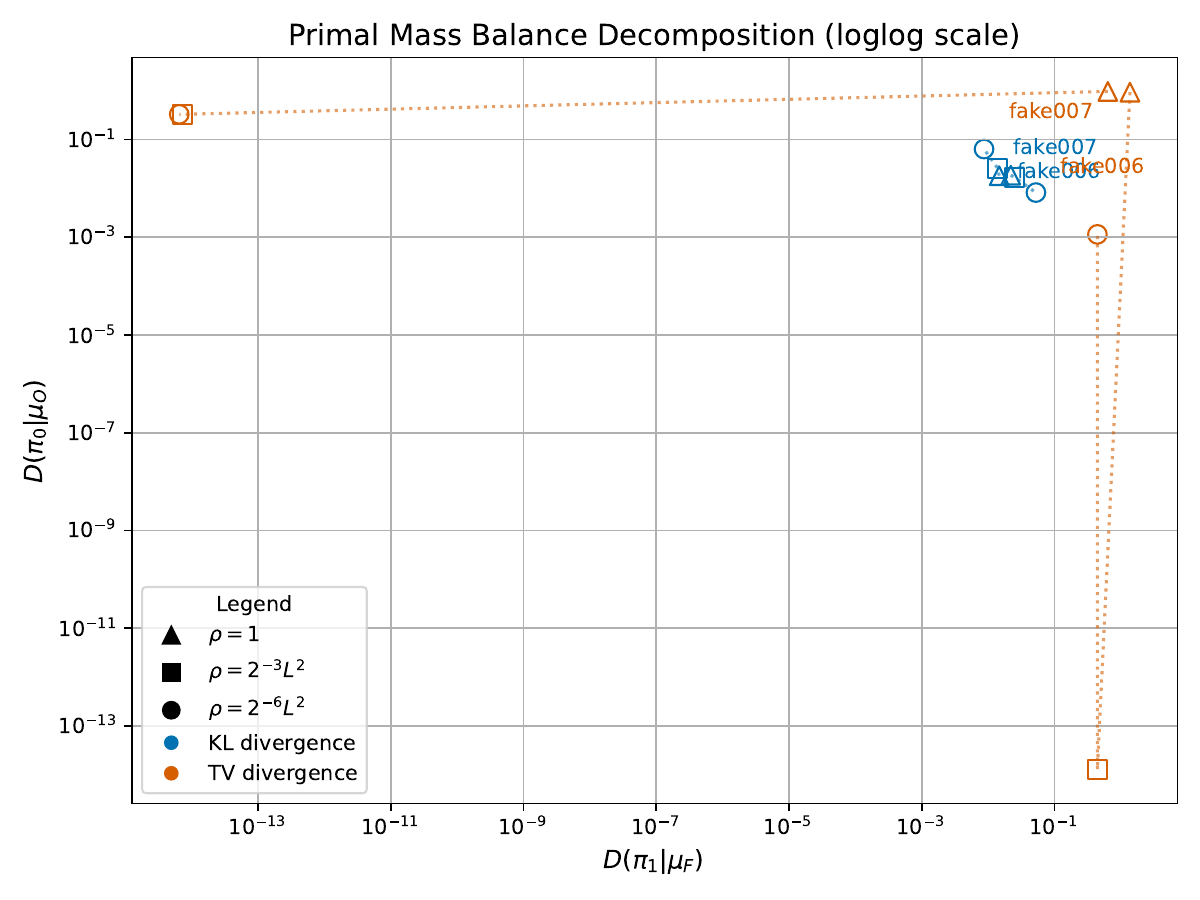}
    \caption{Primal mass balance decomposition (of marginal mass imbalance error,  \ref{eq:uot_general_cost}) of the cost for unbalanced cases fake000fake006 and fake000fake007. Using the primal disintegration of transport vs mass balance alone the direction of imbalance was not possible, i.e. are we over or under forecasting. Yet through the ratio of the marginal penalties paid a direction is highlighted, recalling fake006 is over forecasting and fake007 under.}
    \label{fig:dualdecomposition_0067}
\end{figure}

\begin{figure}
    \centering
    \includegraphics[width=0.75\linewidth]{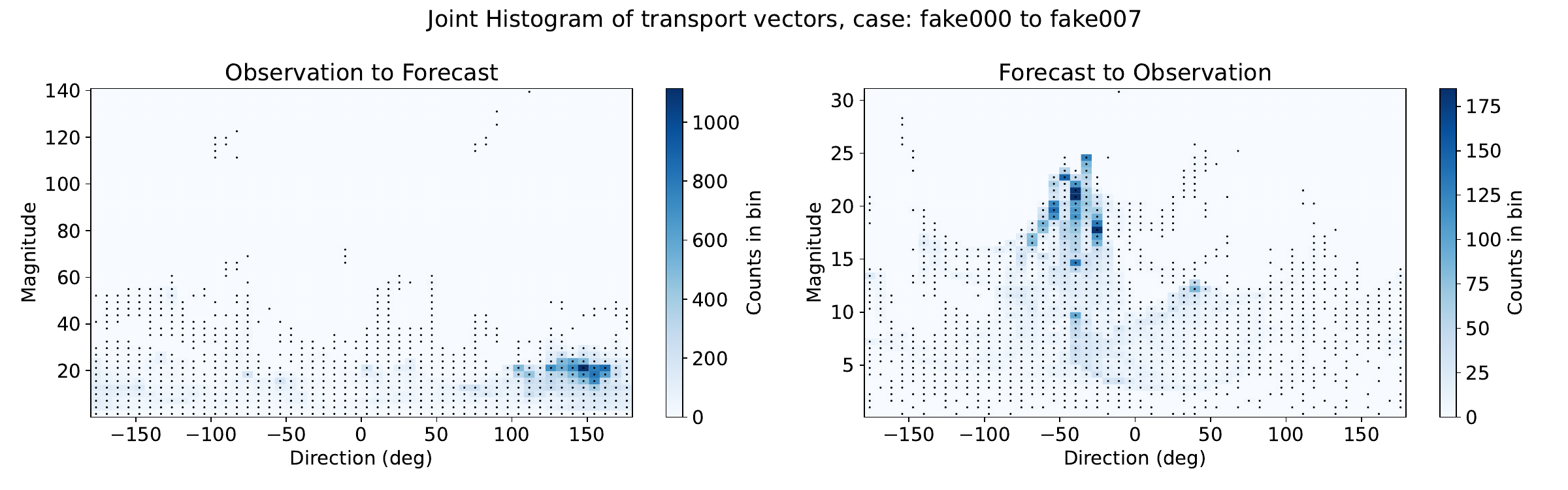}
    \includegraphics[width=0.75\linewidth]{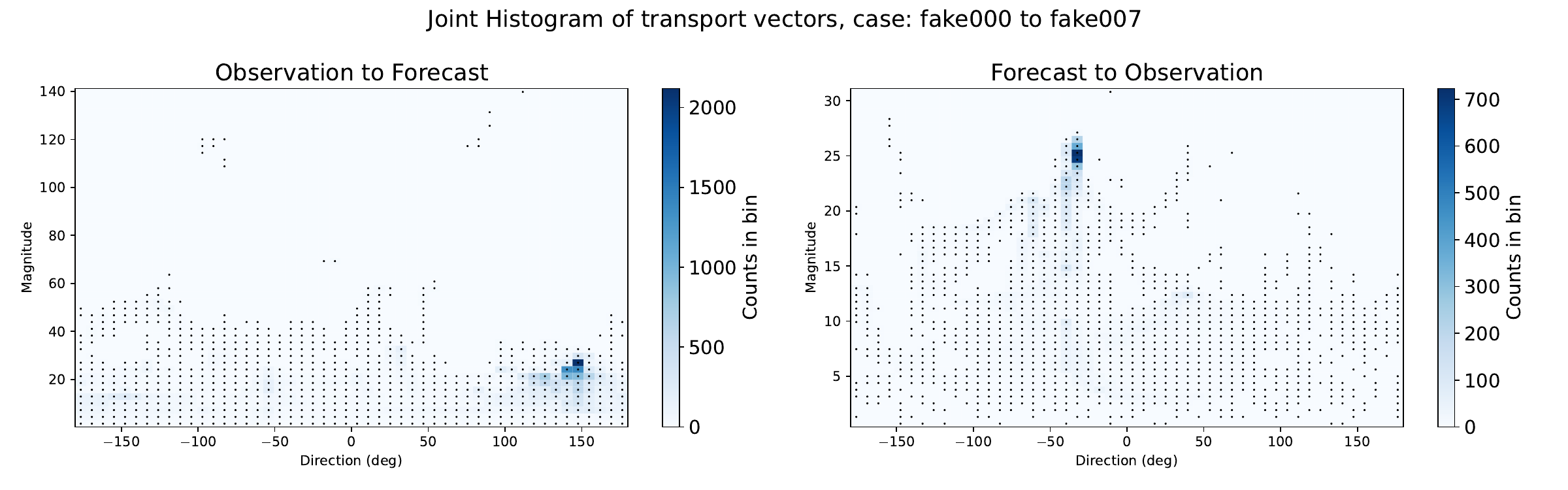}
    \caption{2D histogram of the magnitude and direction of the underlying transports vectors for the perturbed case fake000fake007, which is a more unbalanced version of fake000fake003.
    To understand the counter-intuitive behaviour in Figure \ref{fig:perturbedspreadofcases}, we examine the transport with \(\rho=1\) and \(2^{-6}\), top and bottom  respectively.
     The dotted bins, indicate those with non-zero mass. These results are based on TV marginal penalisation with parameters \(\varepsilon=0.001L^2, \rho=L^2\)}
    \label{fig:fake007_2d_histograms}
\end{figure}

\begin{figure}[h]
    \centering
    \includegraphics[width=0.98\linewidth]{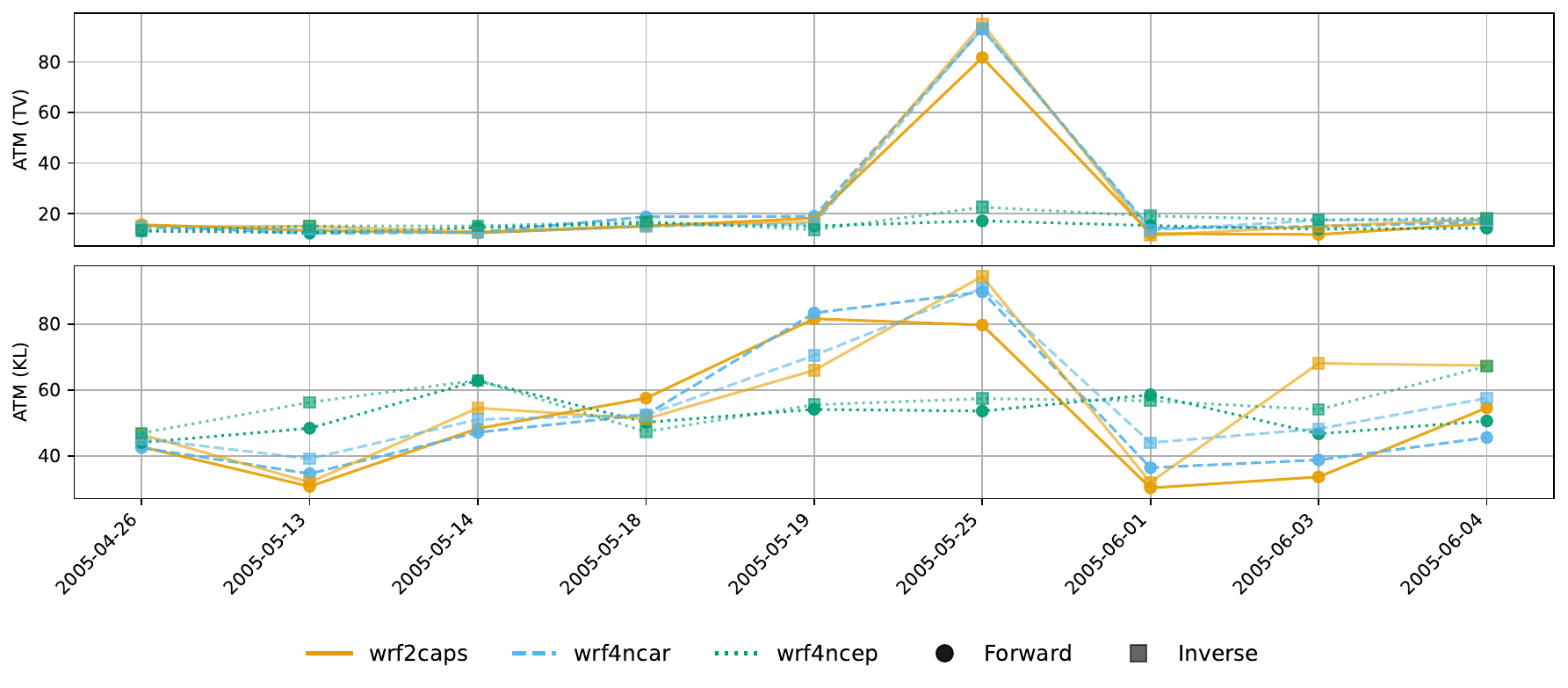}
    \caption{ Spring 2005 cases over 9 times stamps, demonstrating the ATM across flavours and time. There is again a clear anomaly on 2005-05-25 in both flavours, and an anomaly on the 2005-05-19 in the KL penalty. 
    Each model was at 24hr lead time, and interpolated on the coarser 4km grid. Top: TV penalty, Bottom: KL penalty. The key at the bottom of  the figure allows for comparison of models and forward and backwards ATM. \(\varepsilon=0.001L^2, \rho=L^2\). }
    \label{fig:spring_2005_all}
\end{figure}

\begin{figure}[h]
    \centering
    \includegraphics[width=0.75\linewidth]{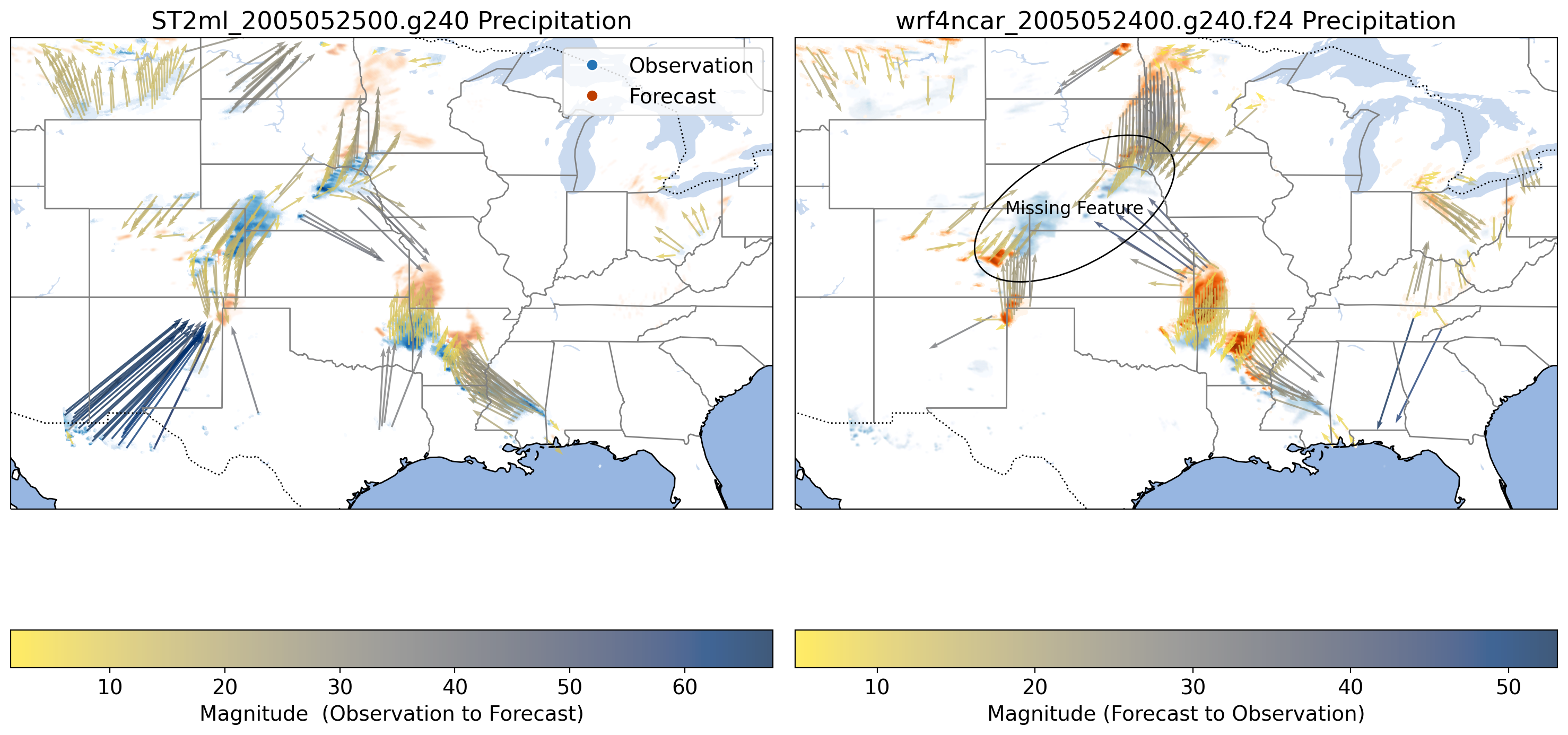}
    \includegraphics[width=0.75\linewidth]{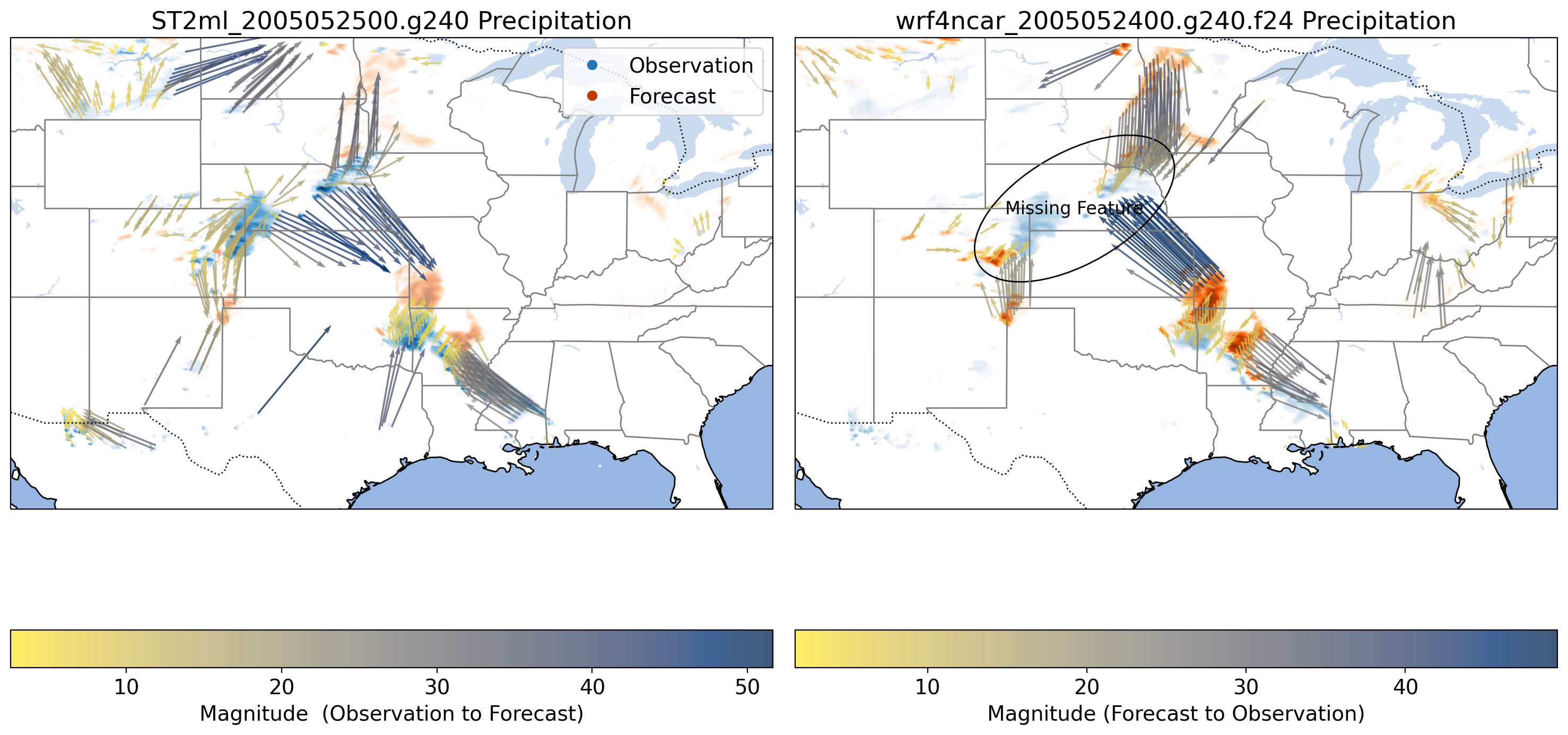}
    \caption{Illustration of the transport vectors with KL (Top rows) and TV (bottom rows) marginal penalisation, found on May25th between the observations ST2ml\_2005052500.g240 (blues shades) and forecast wrf4ncar\_2005052400.g240.f24 (orange shades). Notice the events that are well-matched; those in the north-west corner and those around Arkansas. There is a large missing feature above Nebraska in the forecast, which is contained in the ellipse (RHS). Left: Observation to forecast transport vectors with forecast intensities superimposed in lighter shades, Right: Forecast to observation transport vectors with observation intensities superimposed in lighter shades. \(\varepsilon = 0.001L^2, \rho =1L^2\)}\label{fig:graphical_abstract}
\end{figure}

Turning to the Spring 2005 dataset, Figure \ref{fig:spring_2005_all} shows ATM data to pair with Figure \ref{fig:spring_2005}. Observe the same anomalies for found on May 25th and May 19th.  
Then Figure \ref{fig:graphical_abstract} illustrates the transport vectors with KL marginal penalisation, found on May25th between the observations ST2ml\_2005052500.g240 (blues shades) and forecast wrf4ncar\_2005052400.g240.f24 (orange shades). 
notice the events that are well-matched; those in the north-west corner and those around Arkansas.
However, within the forecast there is a large missing feature above Nebraska and on the Mexico boarder.
This in turns has lead to large transport required to correct the missing features, and it is found in the forward and inverse map. 

\begin{figure}
    \centering
    \begin{minipage}{0.5\linewidth}
        \centering
        \includegraphics[width=\linewidth]{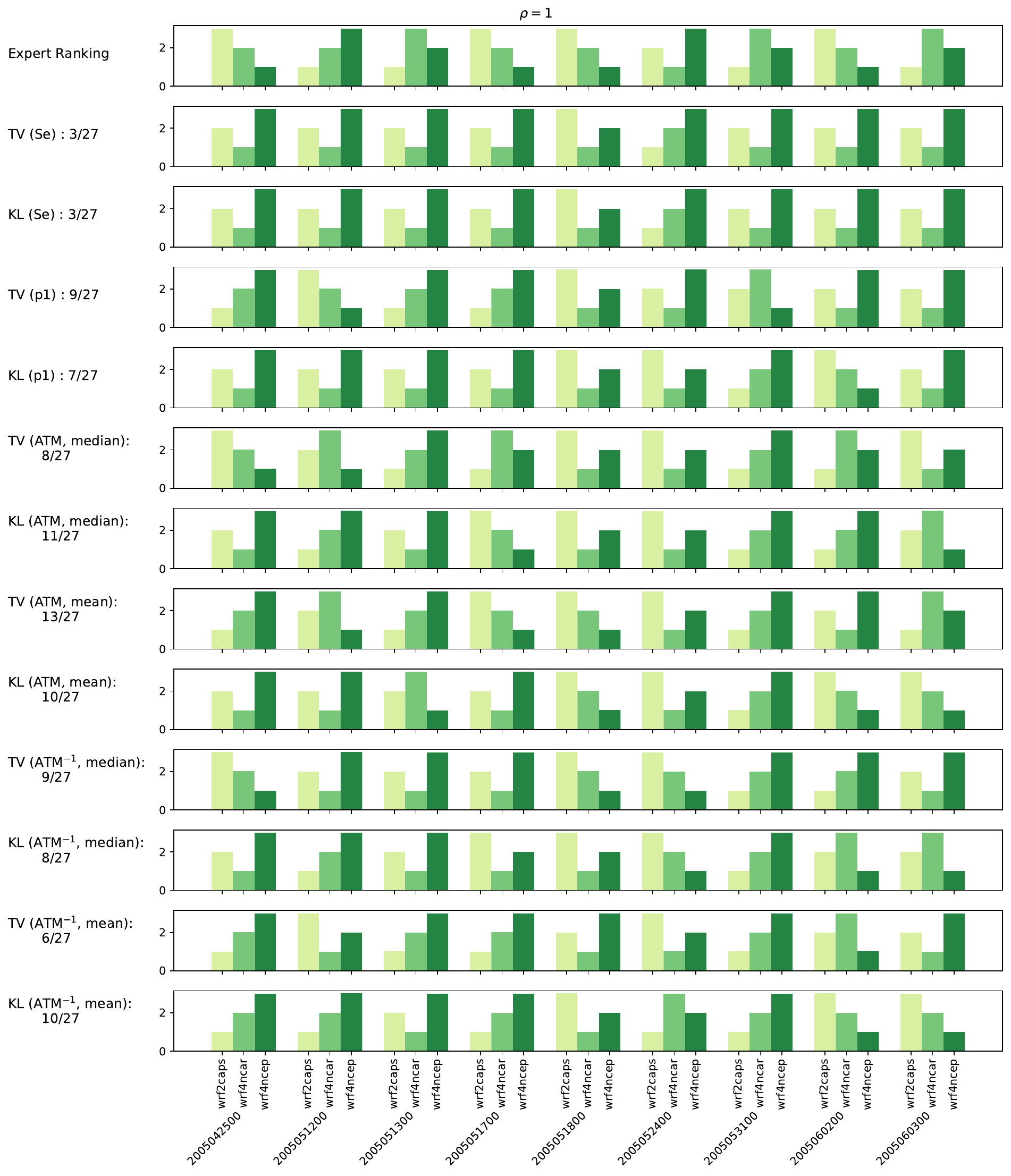}
    \end{minipage}%
    \begin{minipage}{0.5\linewidth}
        \centering
        \includegraphics[width=\linewidth]{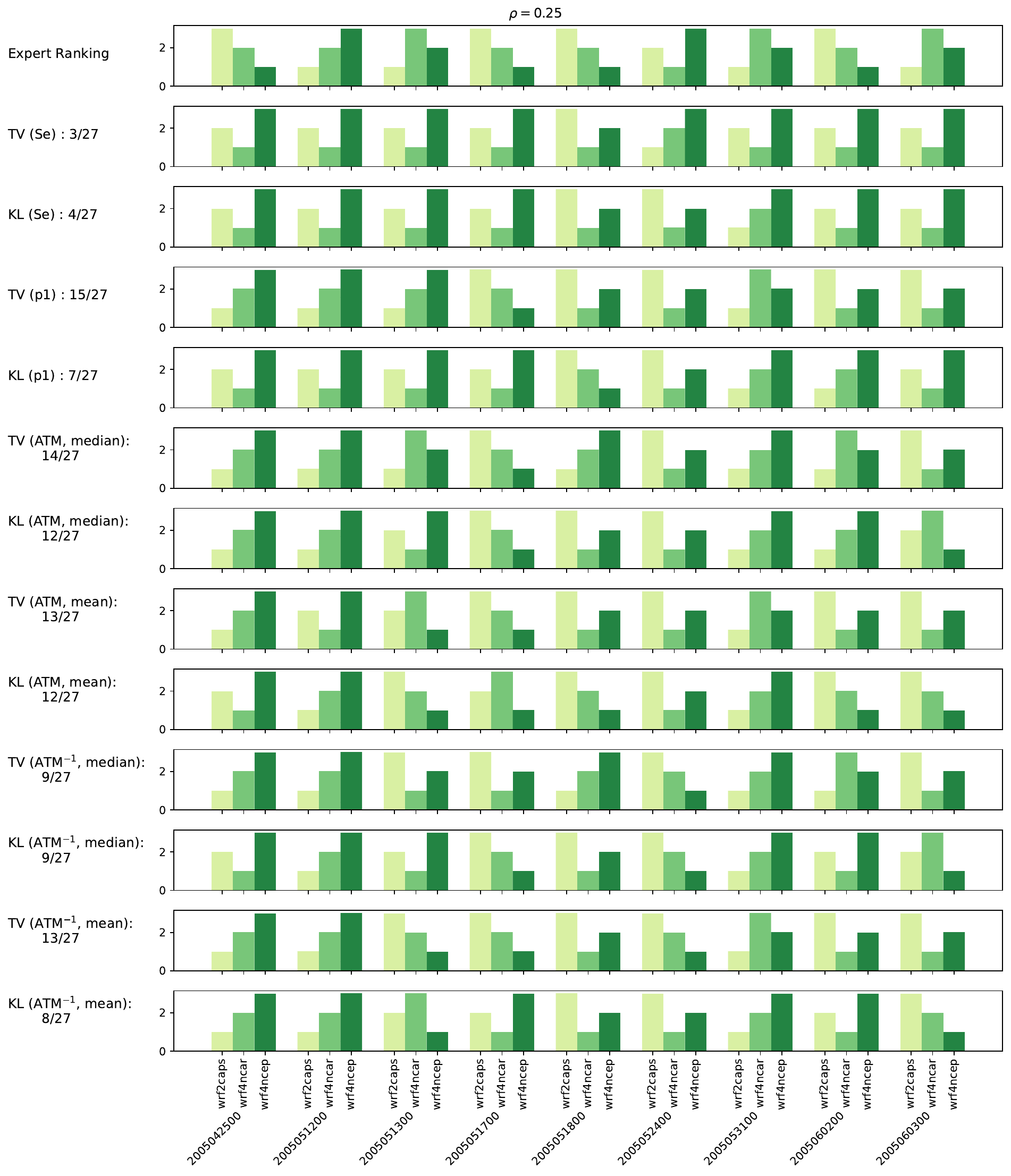}
    \end{minipage}
    \caption{Exploration of ranking for each model each day, against the expert subjective assessment. The hashed-in bars, correspond to matches to the experts, and the faction (/27) is the total number of agreements. The left correspondence to the larger reach value and the right the lower. \(\rho = 1, 0.25\) respectively. In taking a more local transport approach, and thus smaller reach, we do gain some more agreement with the expert score. However mainly for the ATM and  not significantly, nor in all cases.}
    \label{fig:sp2005_ranking}
\end{figure}

\begin{figure}[h]
    \centering
    \begin{subfigure}[b]{0.48\linewidth}
        \centering
        \includegraphics[width=\linewidth]{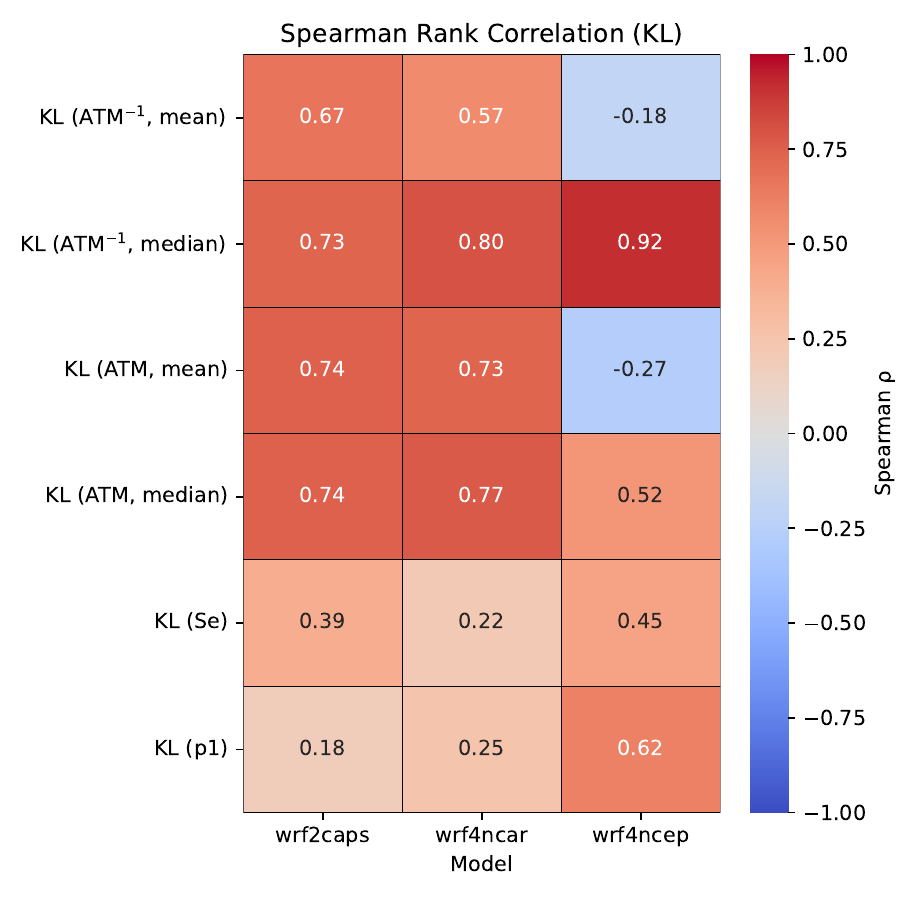}
    \end{subfigure}
    \hfill
    \begin{subfigure}[b]{0.48\linewidth}
        \centering
        \includegraphics[width=\linewidth]{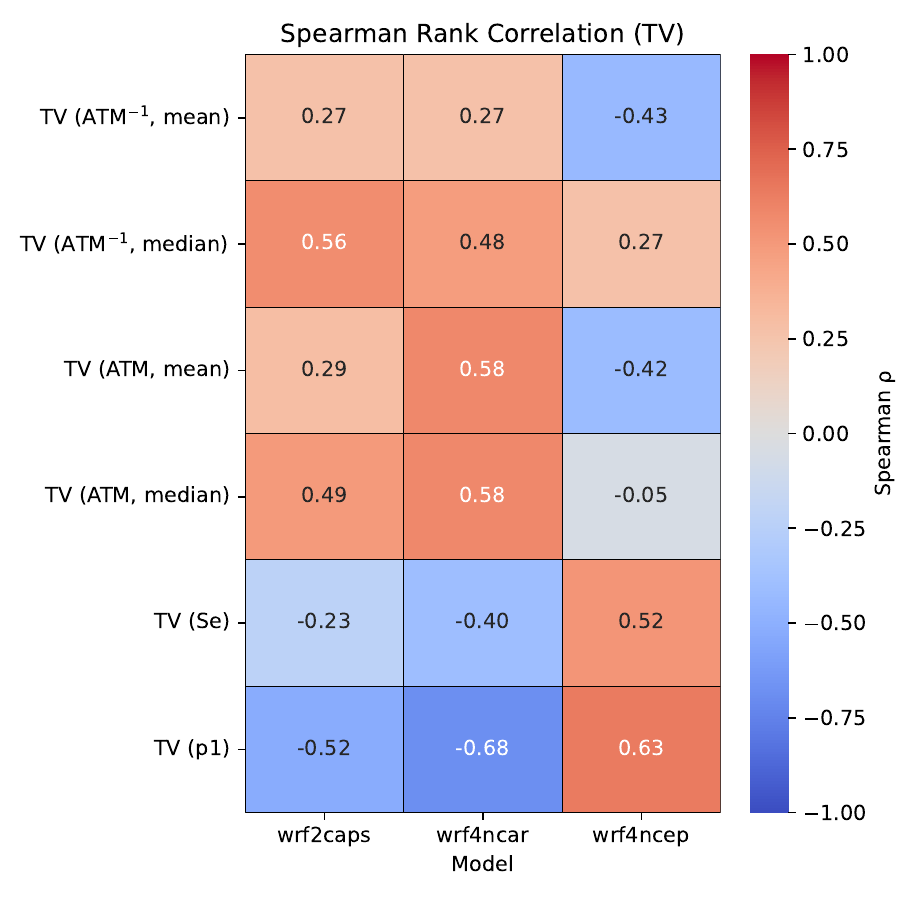}
        
    \end{subfigure}
    \caption{Comparison of Spearman rank correlations for KL  (left )and TV (right ) scores across models. The columns correspond to each model,  with agreement in ranking across the 9 days being tested. Down the rows, we examine: inverse mean ATM, inverse median ATM, forward mean ATM, forward median ATM, \(\Sink\), cost term  \((\sum_{i,j}c_{i,j}\pi_{i,j})\). }
    \label{fig:heatmap-comparison}
\end{figure}

Figures \ref{fig:sp2005_ranking} and \ref{fig:heatmap-comparison} show the ranking against expert evaluation and are discussed in the main report. Spearman’s rank correlation coefficient has not been calculated for the model day by day ranking; instead, a visual demonstration is presented, as only three ranks are being tested.


For the MesoVICT data preparation, each model provides output on similar, but not identical, gridded domains. Therefore, the VERA data is clipped to match each model's respective limited area before any comparisons are made. Additionally, to place the data on a regular grid (allowing use of the tensorisation trick; see Supplementary Materials \ref{appendix:used_tricks}), some edge points are filled with zero weights to achieve a rectangular, regular grid.
Since only paired valid times, in both models, are compared and as the CMH is ran for 18 hours there is a 6 hour gap each day seen on the plots. 
This also means the initialisation time between the models was different. 

\begin{figure}
    \centering
    \includegraphics[width=0.95\linewidth]{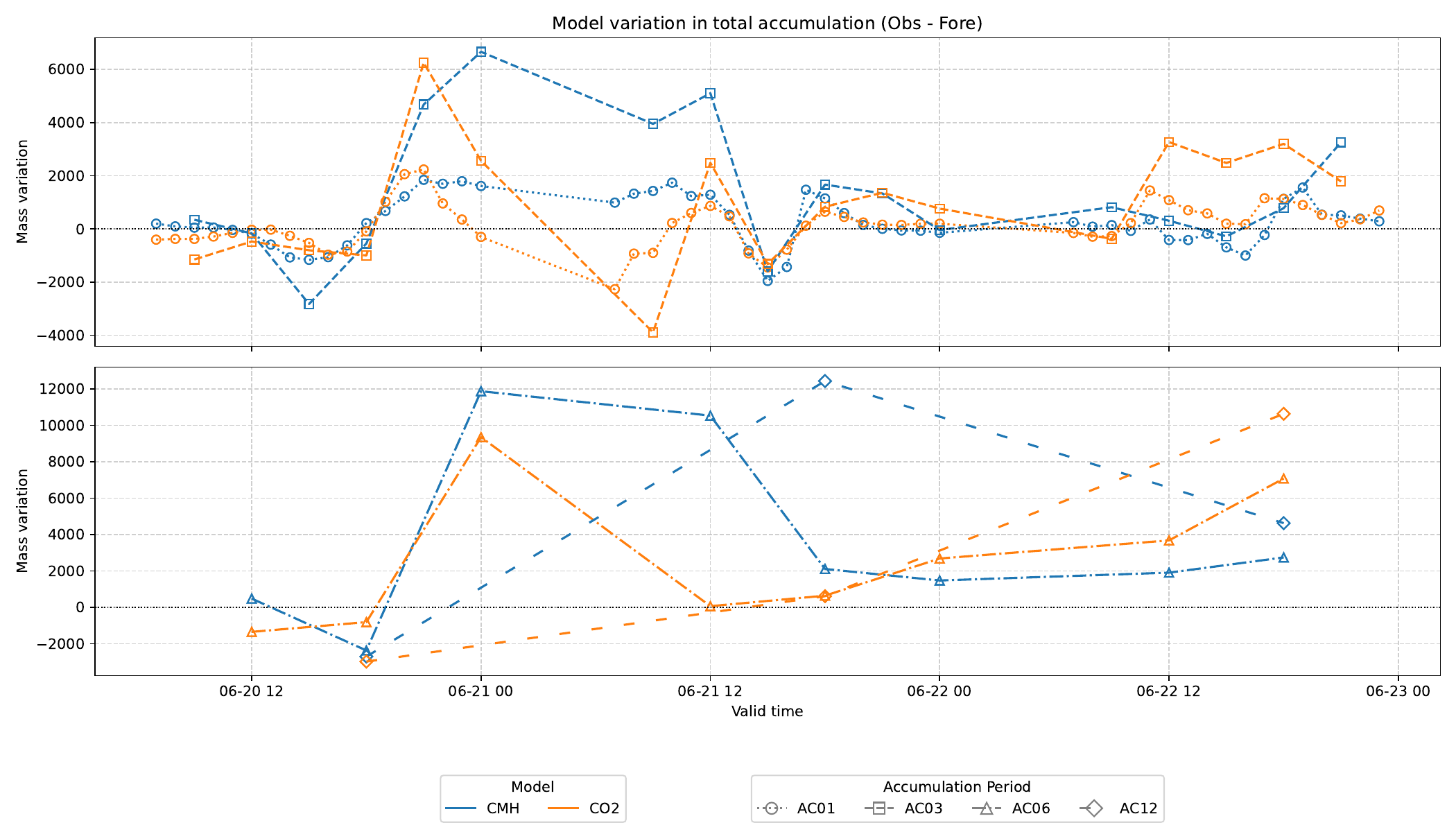}
    \caption{Illustration of mass balance over time for both models, and across different accumulation periods. Here, we take the difference between observed and forecast total accumulated rainfall. Positive values indicate under-forecasting, while negative values indicate over-forecasting. As seen in the decomposition diagrams, these models predominantly under-forecast.}
    \label{fig:vera_mass}
\end{figure}

\begin{figure}[h]
    \centering
    \begin{minipage}{0.49\linewidth}
        \includegraphics[width=\linewidth,trim= 70 70 70 70,clip]{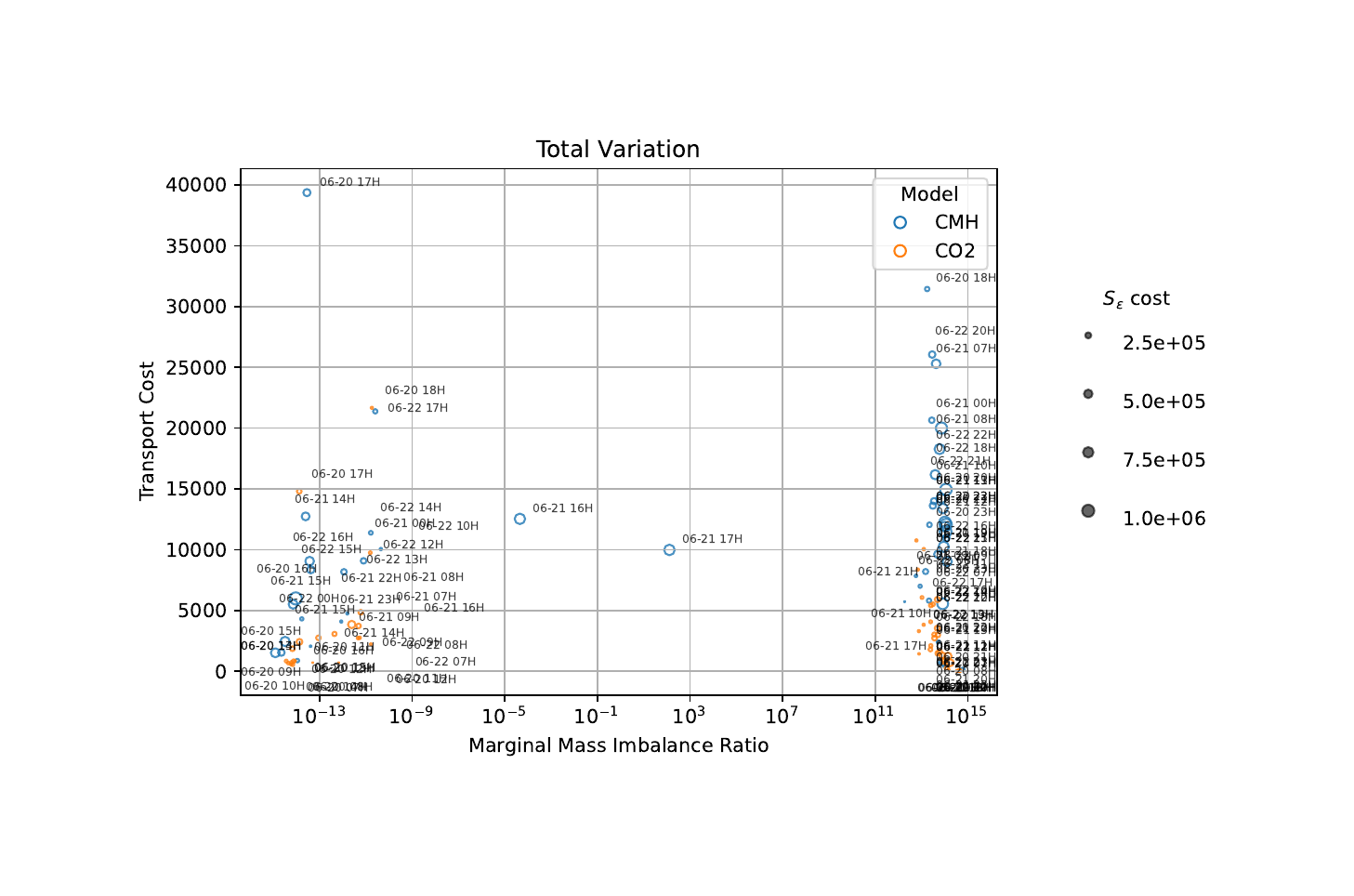}
    \end{minipage}
    \begin{minipage}{0.49\linewidth}
        % \subcaption*{(b)}
        \includegraphics[width=\linewidth,trim= 70 70 70 70,clip]{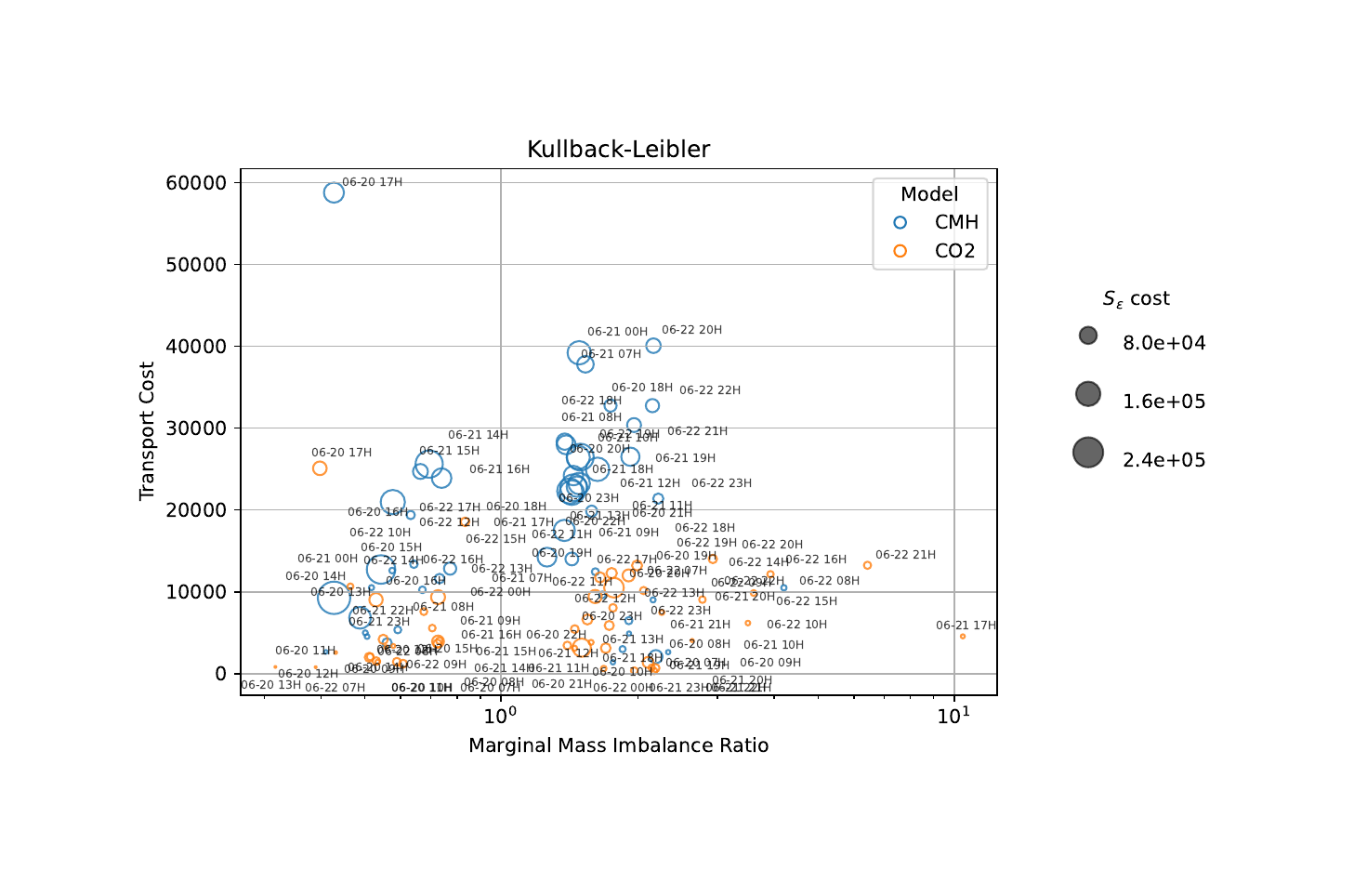}
    \end{minipage} \\
    \begin{minipage}{0.49\linewidth}
        % \subcaption*{(c)}
        \includegraphics[width=\linewidth,trim= 70 70 70 70,clip]{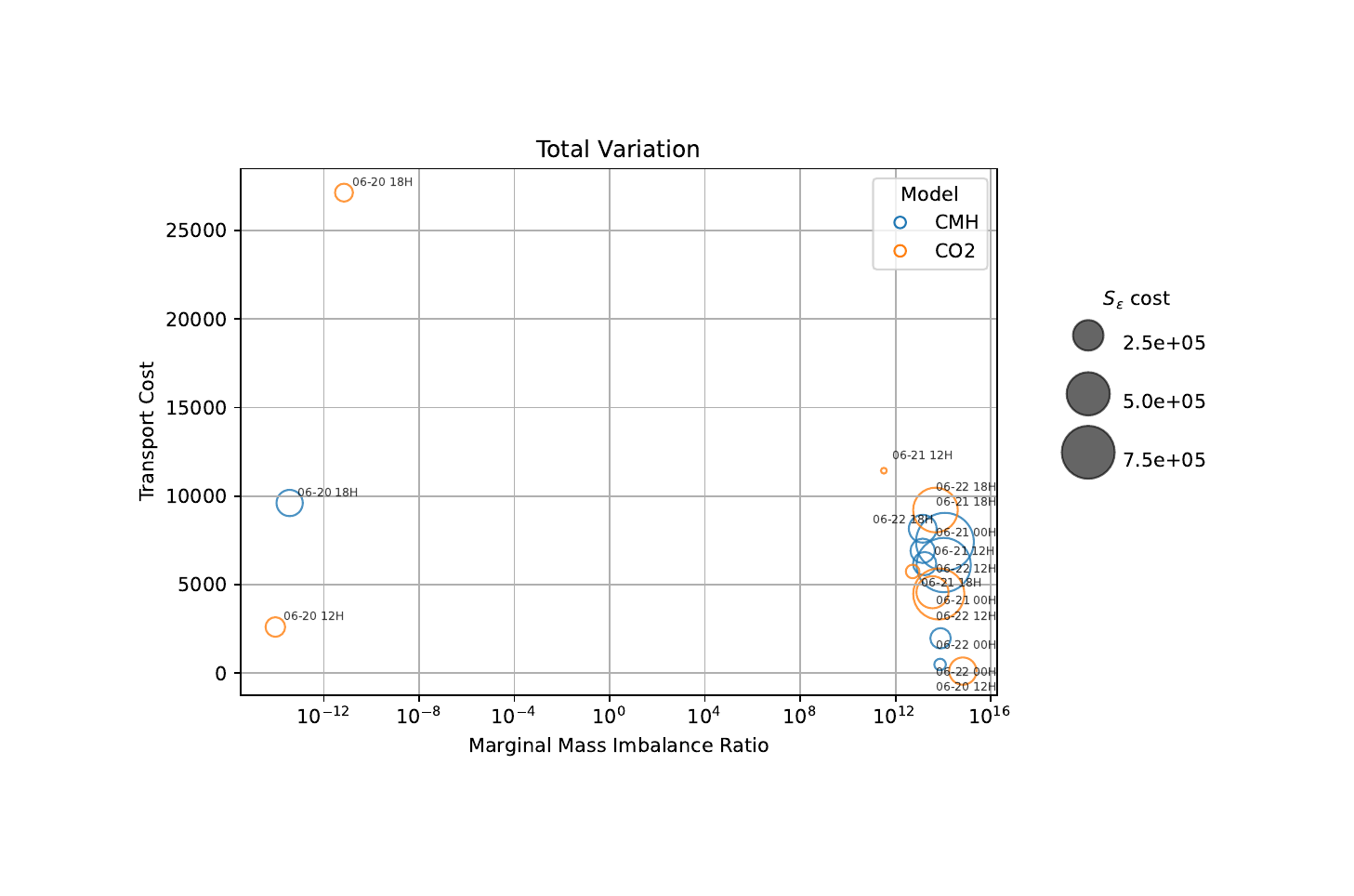}
    \end{minipage}
    \begin{minipage}{0.49\linewidth}
        % \subcaption*{(d)}
        \includegraphics[width=\linewidth,trim= 70 70 70 70,clip]{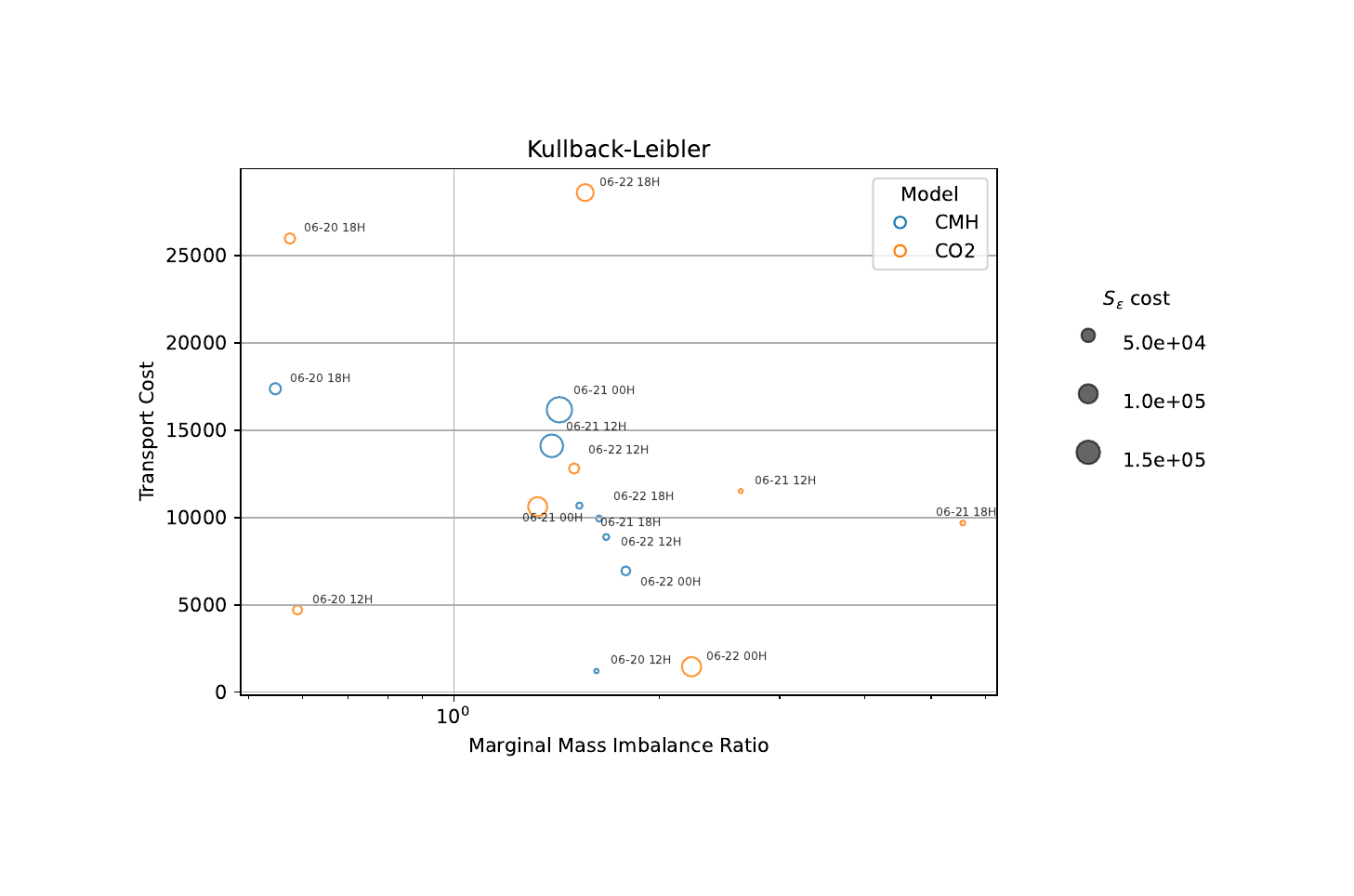}
    \end{minipage}
        \caption{Decomposition of cost into the transport and marginal mass imbalance ratio, for the MesoVICT core case. This shows the decomposition for AC01 (top) and AC06 (bottom), with \(\rho=L^2\). AC03 is available in Figure \ref{fig:vera_decomposition_AC03}. For AC01 notice that CMH (blue) appears to have a higher transport cost across cases compared to CO2. However both over and under forecast for this short accumulation period. By AC06 this trends seem to have gone, assuming lag errors have been absorbed, and it is seem we tend to under-forecast given that the marginal mass imbalance ratio is above 1.  Top: TV penalty, Bottom: Kl penalty. Here \(\varepsilon=0.005L^2, \rho=L^2\). }
    \label{fig:vera_decomposition_grid_rho1}
\end{figure}

\begin{figure}[h]
    \centering
    \begin{minipage}{0.49\linewidth}
        % \subcaption*{(e)}
        \includegraphics[width=\linewidth, trim= 70 70 70 70,clip]{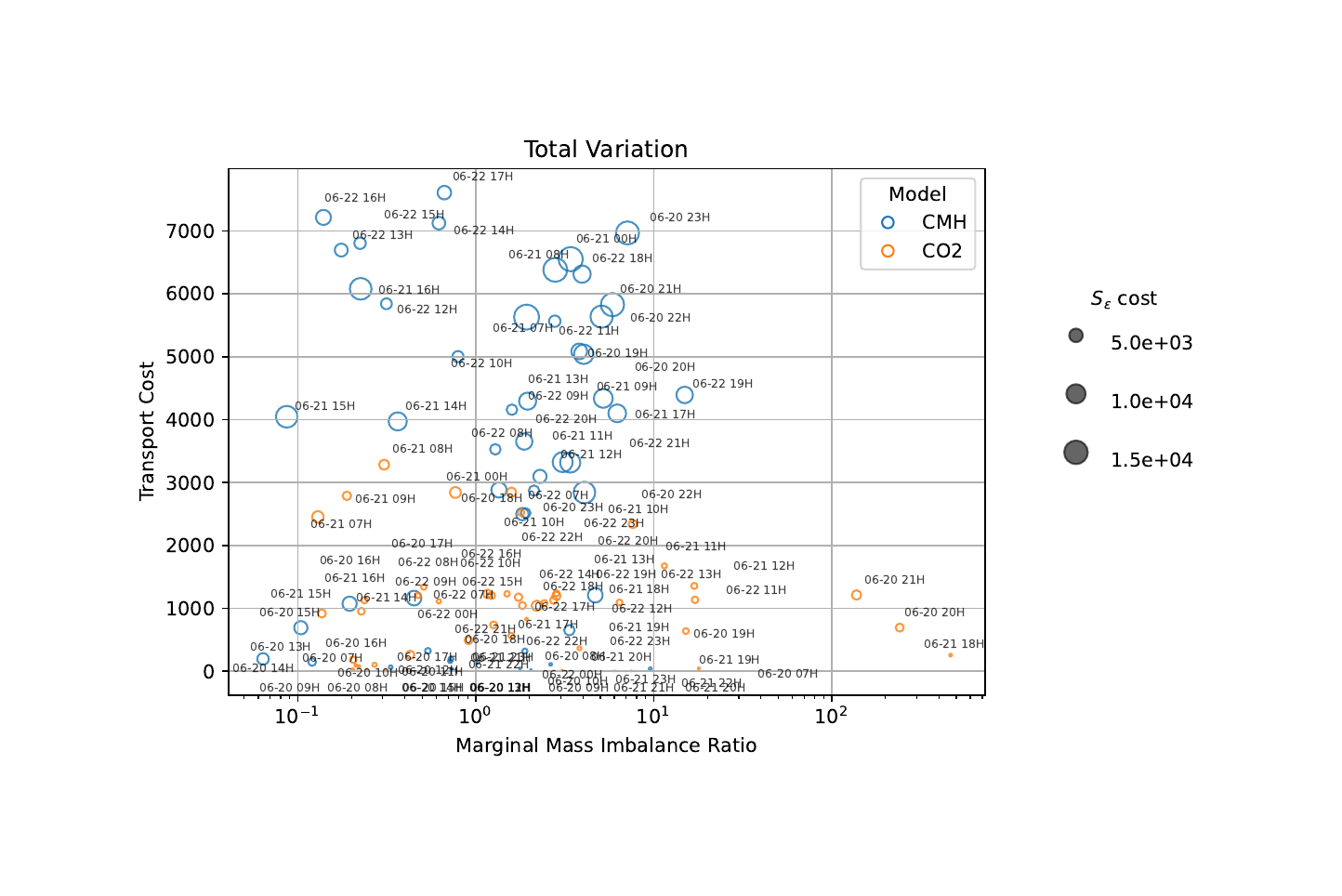}
    \end{minipage}
    \begin{minipage}{0.49\linewidth}
        % \subcaption*{(f)}
        \includegraphics[width=\linewidth,trim= 70 70 70 70,clip]{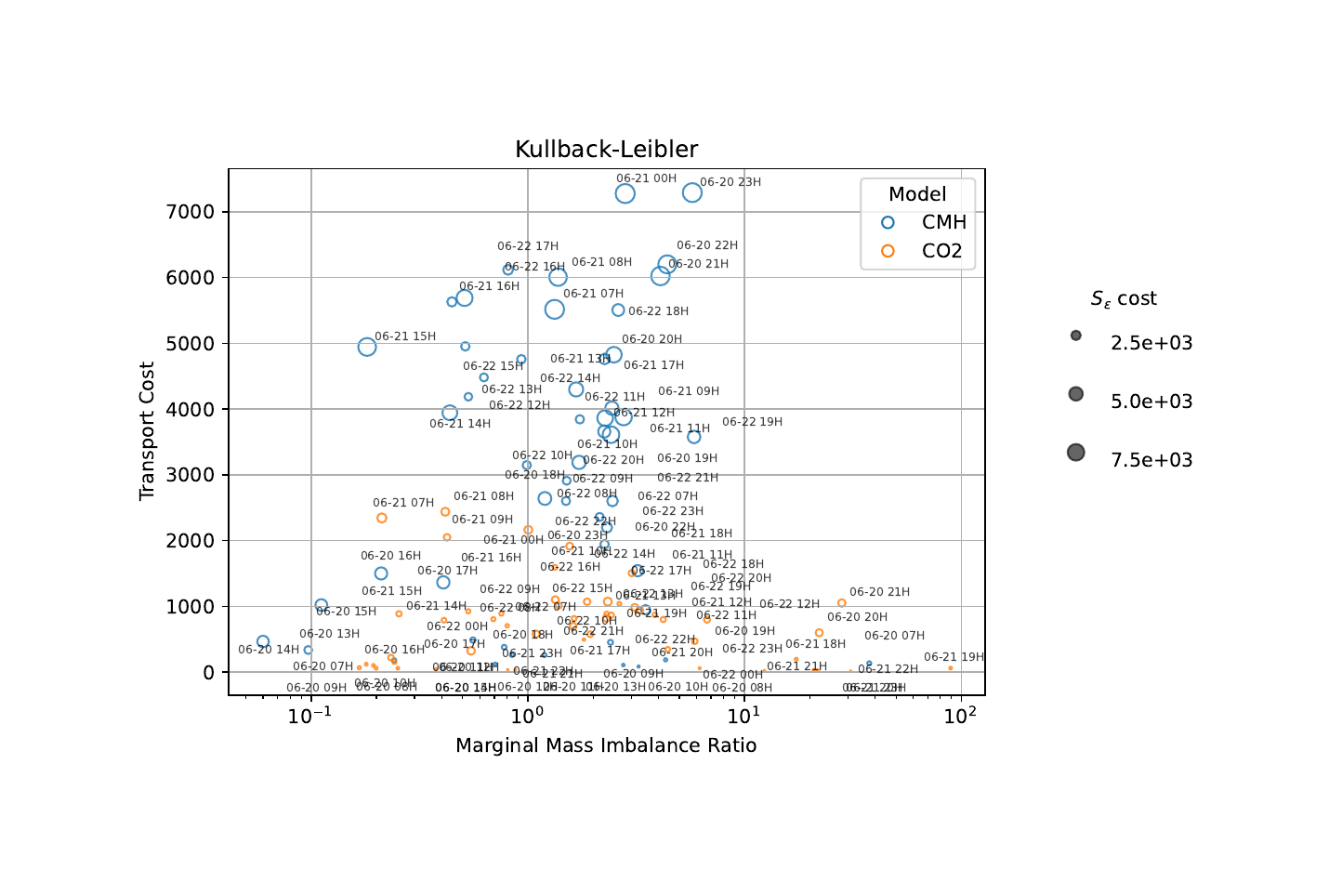}
    \end{minipage} \\
    \begin{minipage}{0.49\linewidth}
        % \subcaption*{(g)}
        \includegraphics[width=\linewidth,trim= 70 70 70 70,clip]{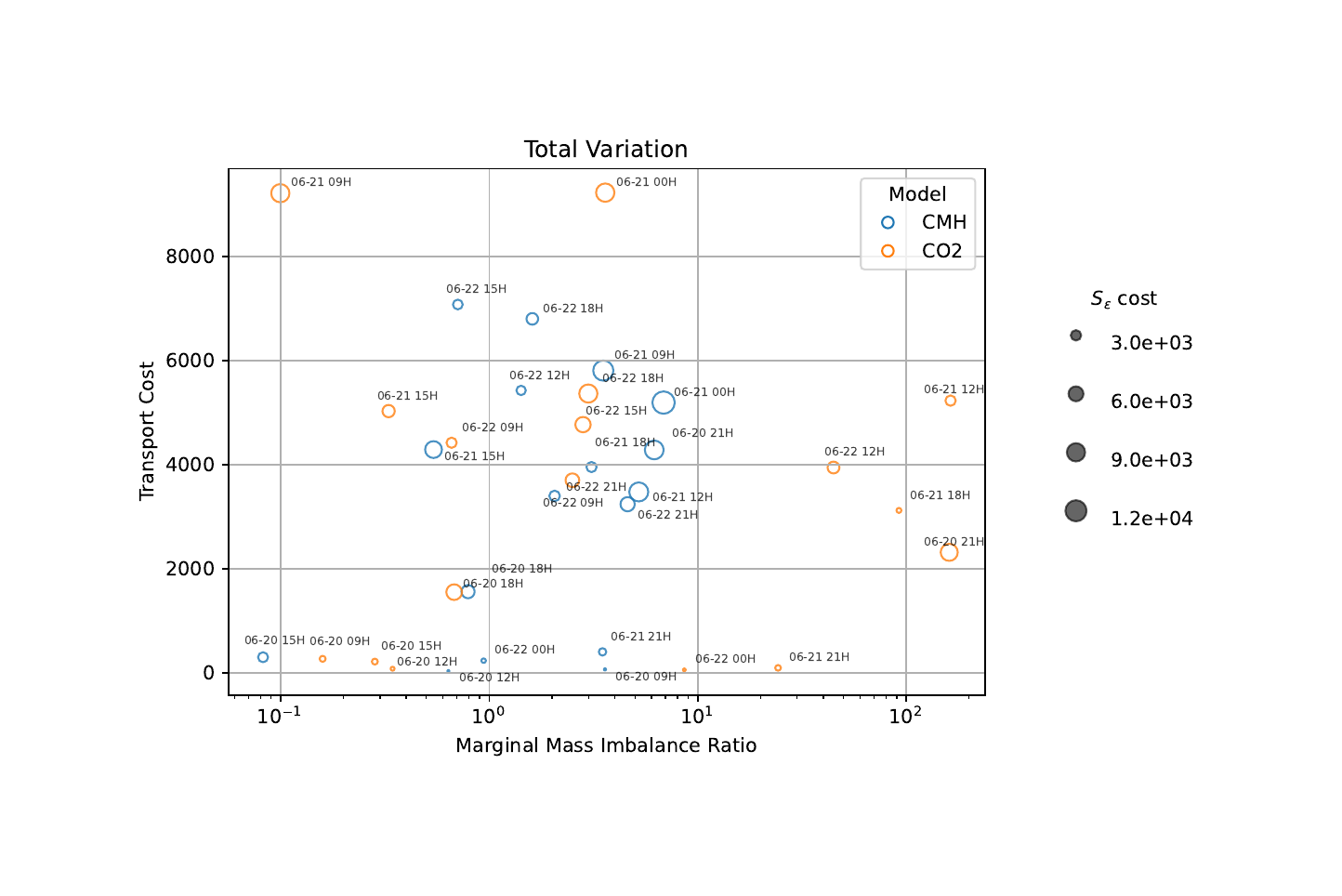}
    \end{minipage}
    \begin{minipage}{0.49\linewidth}
        % \subcaption*{(h)}
        \includegraphics[width=\linewidth,trim= 70 70 70 70,clip]{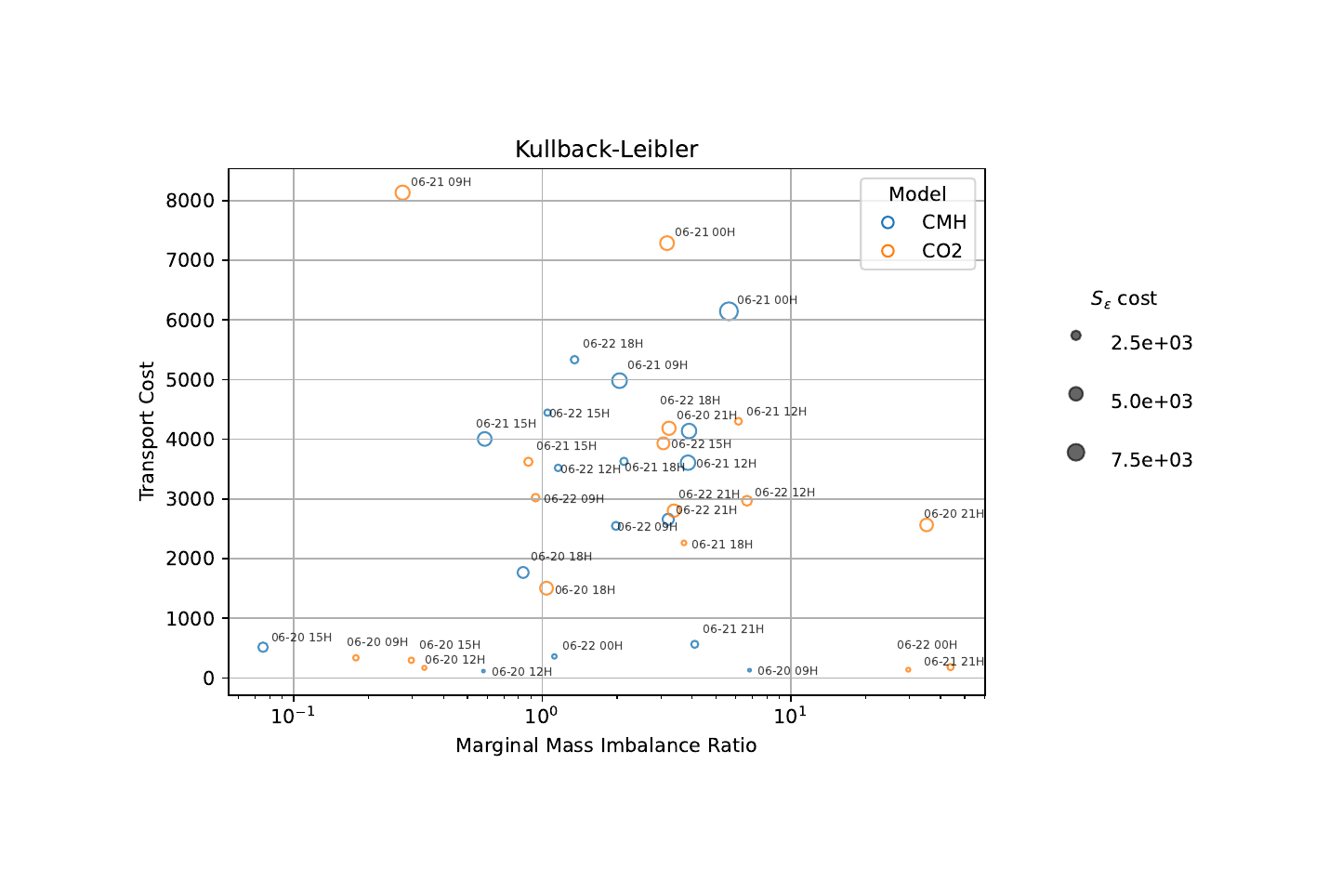}
    \end{minipage} \\
    \begin{minipage}{0.49\linewidth}
        % \subcaption*{(i)}
        \includegraphics[width=\linewidth,trim= 70 70 70 70,clip]{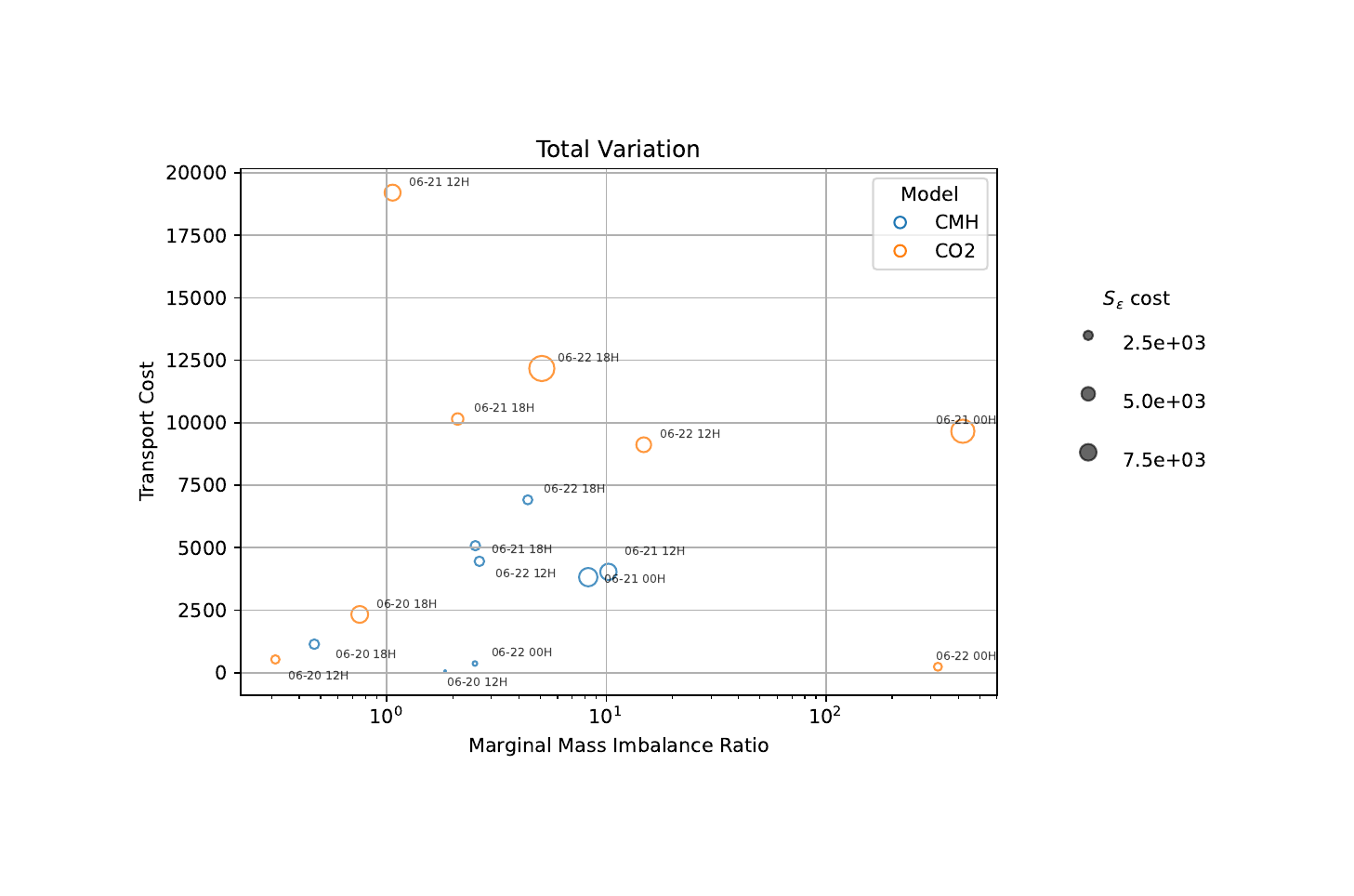}
    \end{minipage}
    \begin{minipage}{0.49\linewidth}
        % \subcaption*{(j)}
        \includegraphics[width=\linewidth,trim= 70 70 70 70,clip]{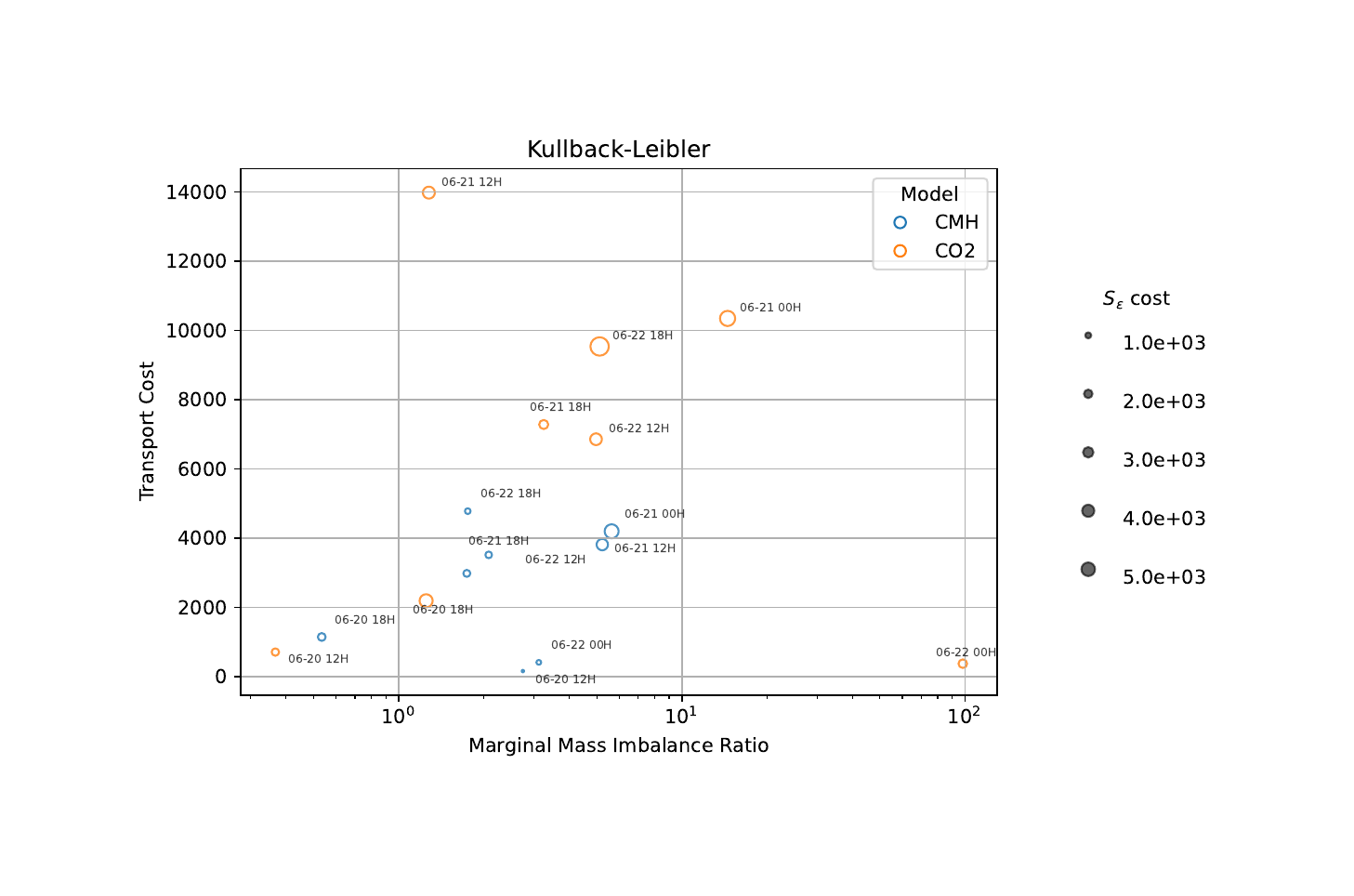}
    \end{minipage} 
    \caption{Decomposition of cost into the transport and marginal mass imbalance ratio, for the MesoVICT core case, with smaller reach parament thus propitiating local transport. Many of the same trends remain, however now that local transport is prioritised, it becomes even more clear that at AC01 CMH has more transport error, though by the AC06 is appears CMH becomes more consistent and CO2 performs worse in balance and transport.  Top: TV penalty, Bottom: Kl penalty. Here \(\varepsilon=0.005L^2, \rho=10^{-2}L^2\). }
    \label{fig:vera_decomposition_grid_rho01}
\end{figure}

\begin{figure}[h]
    \centering
    \includegraphics[width=\textwidth, trim=0 70 0 70, clip]{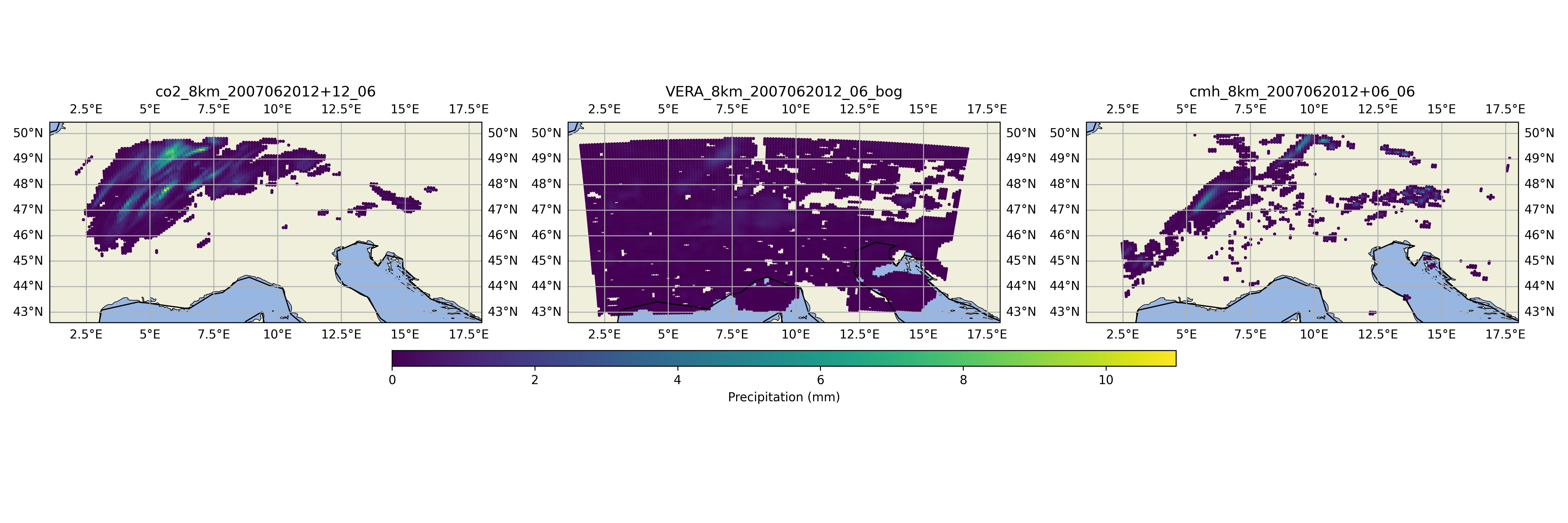}
    \includegraphics[width=\textwidth, trim=0 70 0 70, clip]{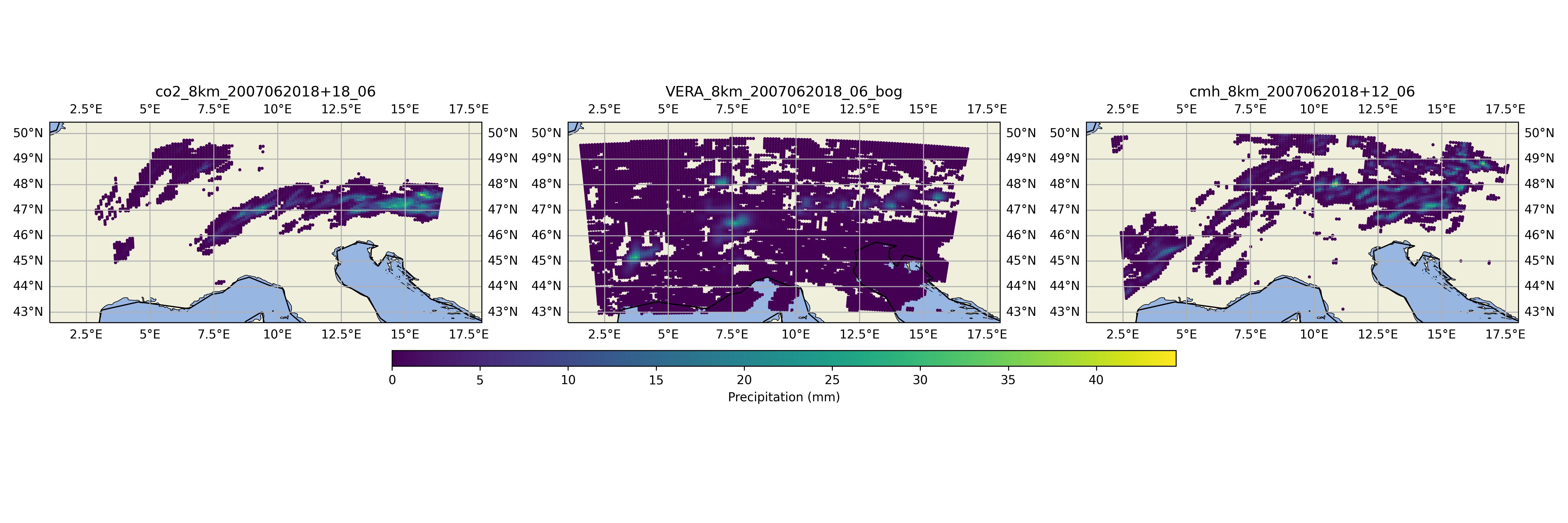}
    \includegraphics[width=\textwidth, trim=0 70 0 70, clip]{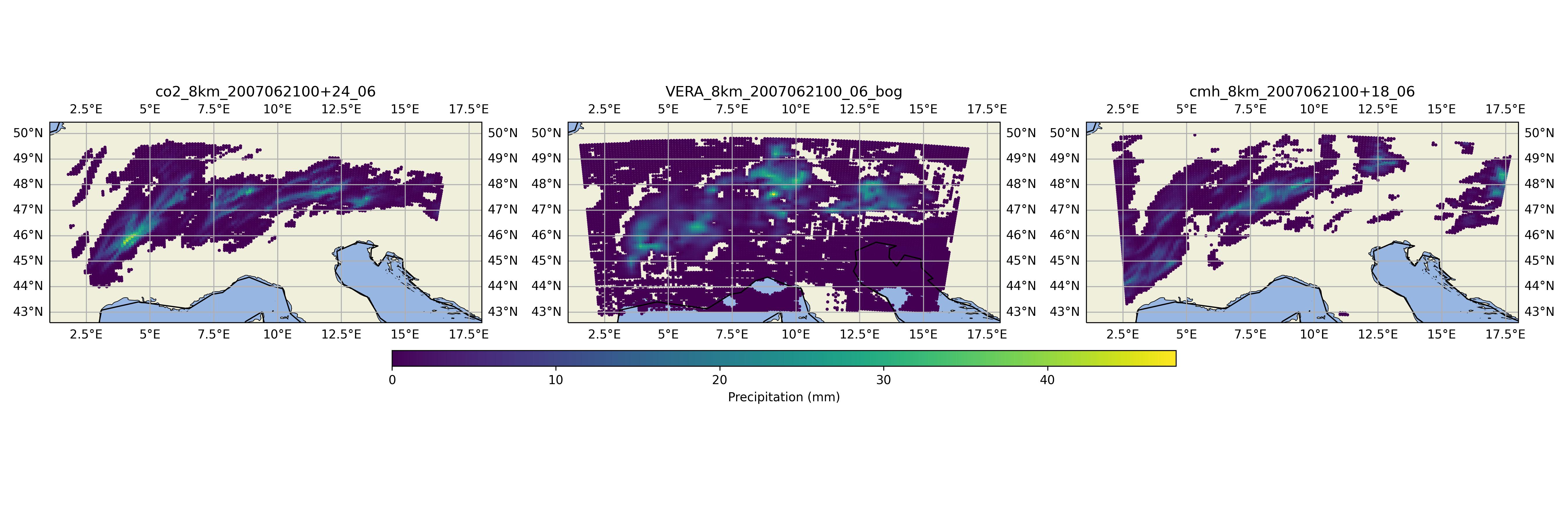}
    \includegraphics[width=\textwidth, trim=0 70 0 70, clip]{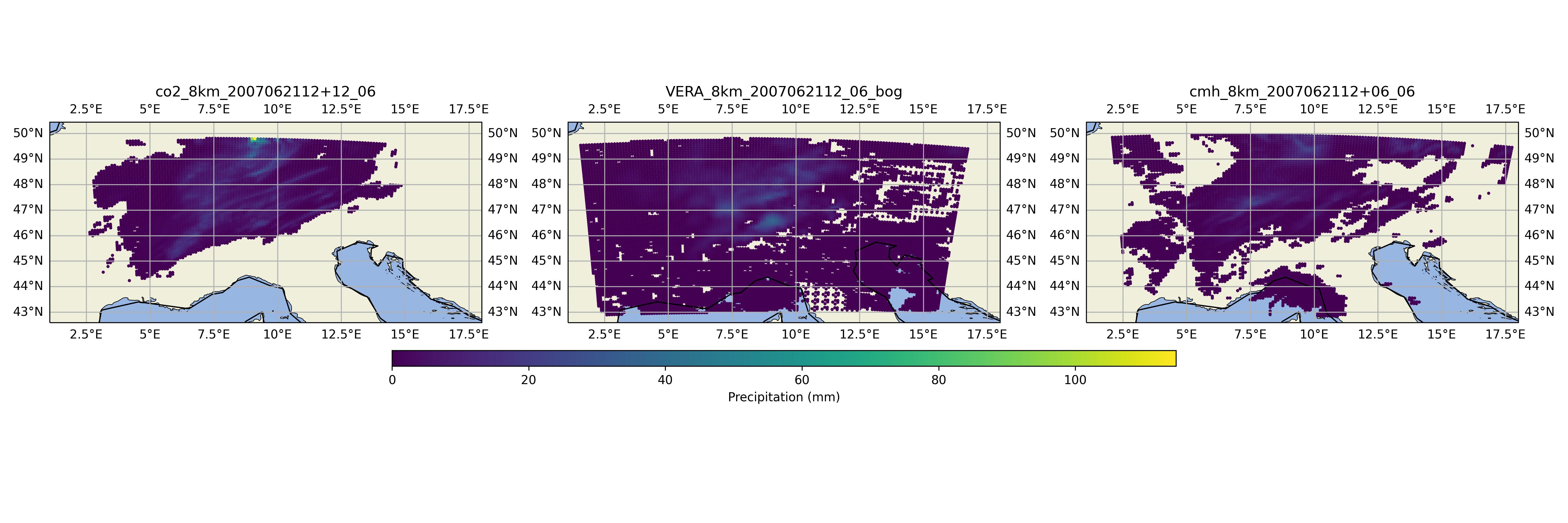}
        \caption{AC06 comparison in chronological order; 06-20 12H, 06-20 18H, 06-21 00H, 06-21 12H. Left: CO2, Middle: VERA, Right: CMH, zero precipitation points are not plotted. }
    \label{fig:vera_ac06_illustration_1}
\end{figure}

\begin{figure}[h]
    \includegraphics[width=\textwidth, trim=0 70 0 70, clip]{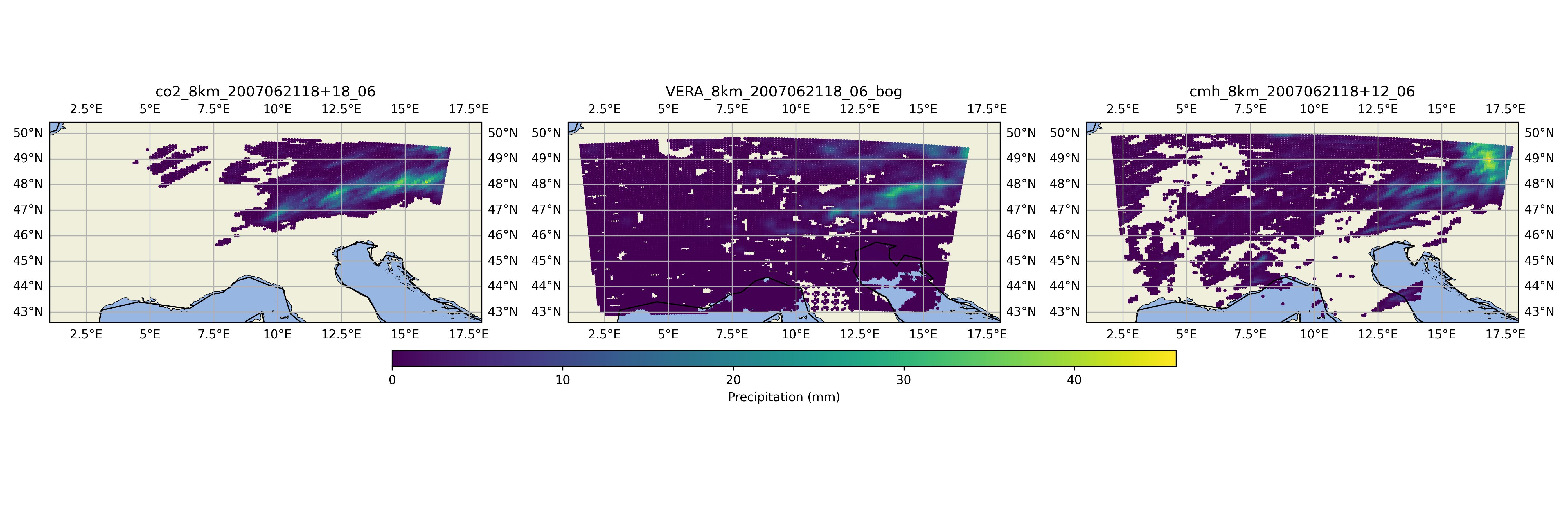}
    \includegraphics[width=\textwidth, trim=0 70 0 70, clip]{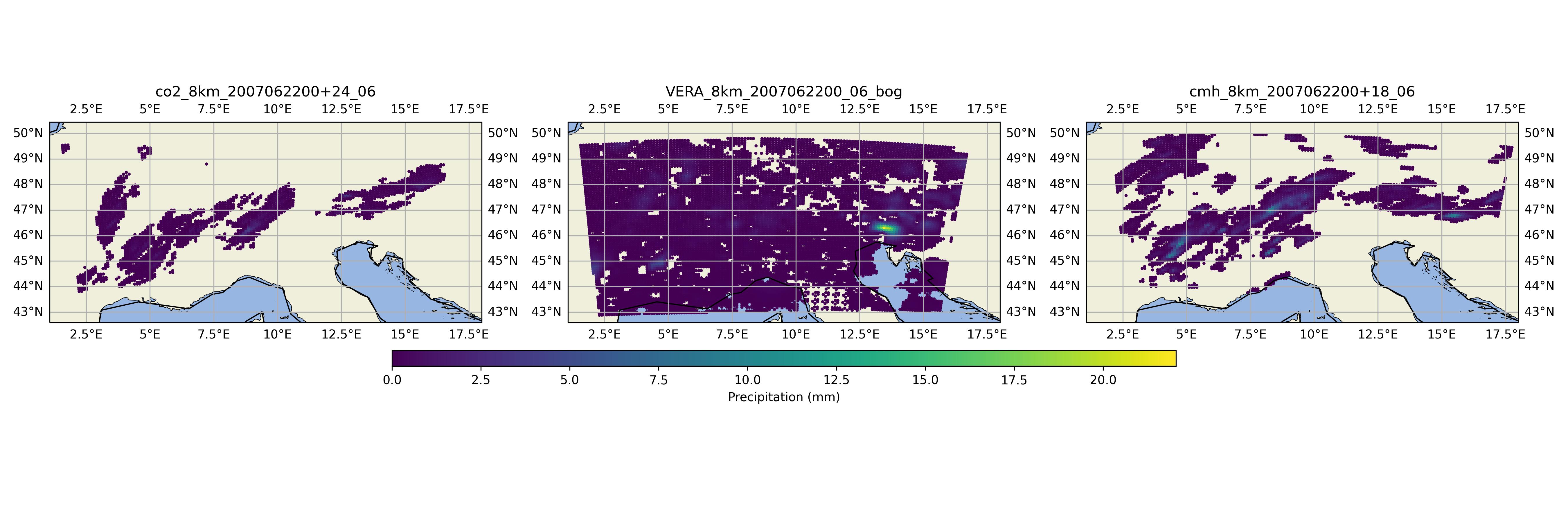}
    \includegraphics[width=\textwidth, trim=0 70 0 70, clip]{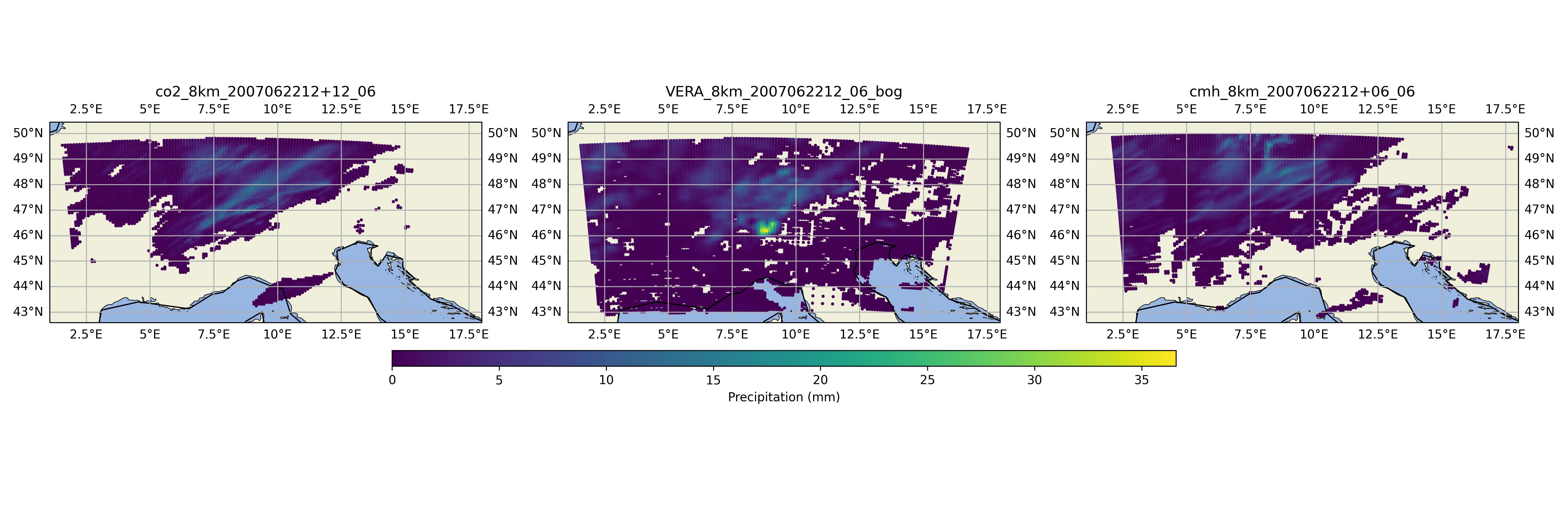}
    \includegraphics[width=\textwidth, trim=0 70 0 70, clip]{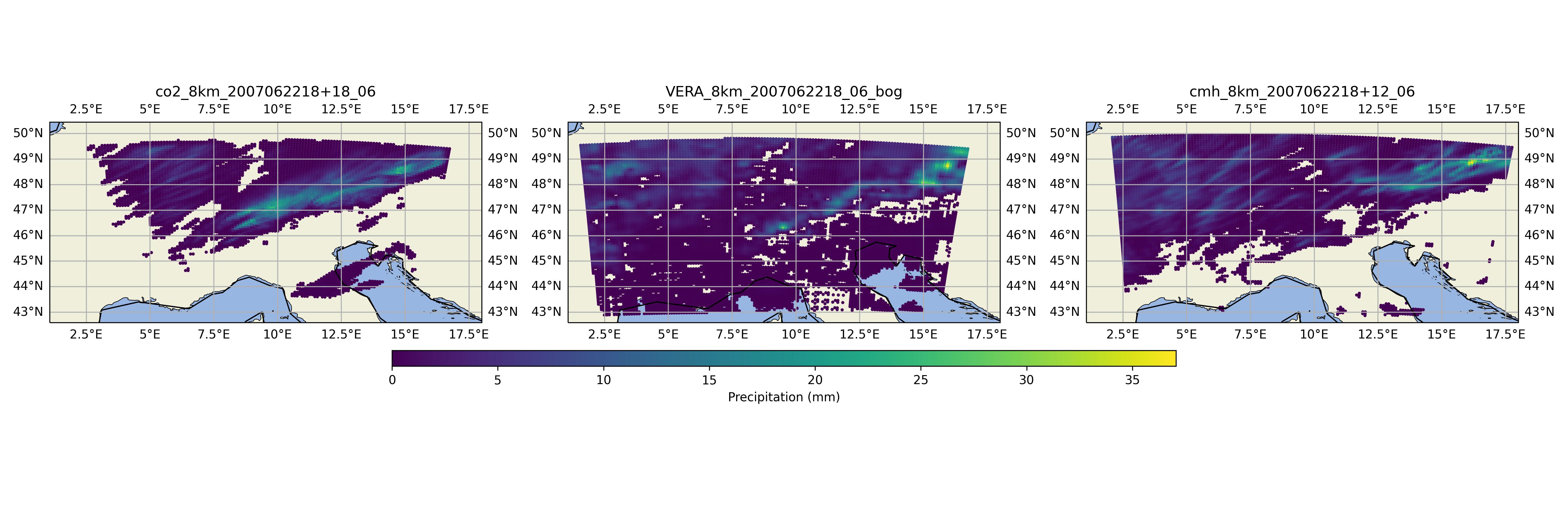}
    \caption{AC06 comparison in chronological order, 06-21 18H, 06-22 00H, 06-22 12H, 06-22 18H. Left: CO2, Middle: VERA, Right: CMH, zero precipitation points are not plotted. }
    \label{fig:vera_ac06_illustration_2}
\end{figure}

\begin{figure}[h]
    \includegraphics[width=\textwidth, trim=0 70 0 70, clip]{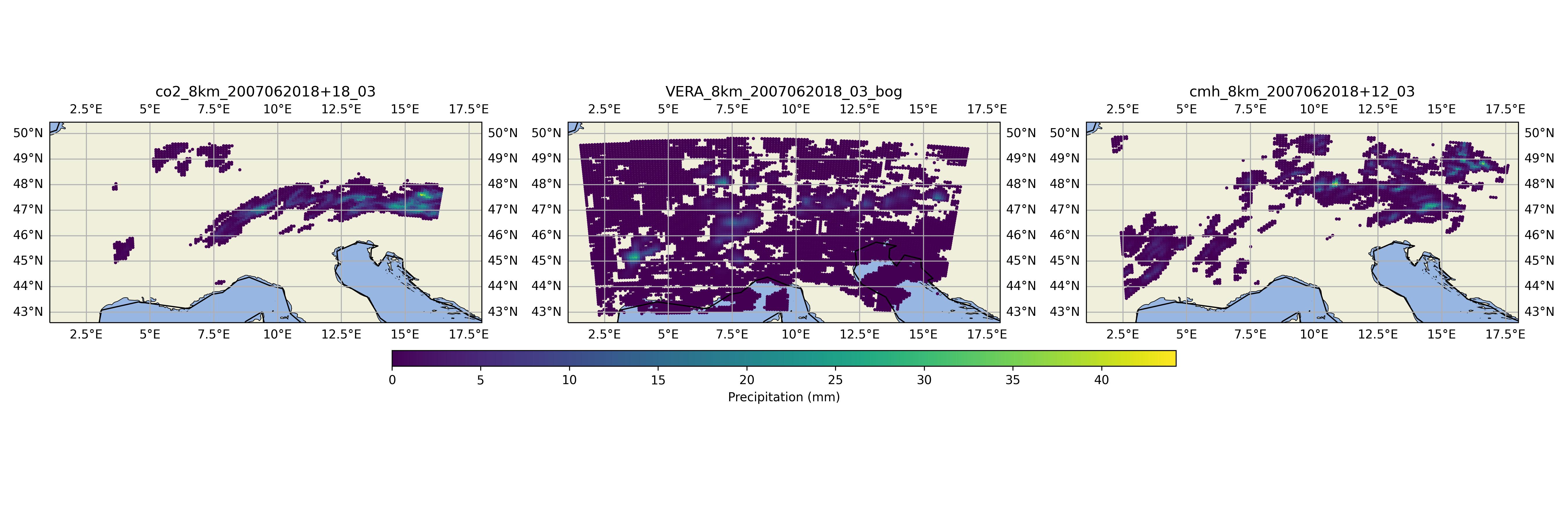}
    \caption{AC03 anomalous outlier 06-20 18H for large transport cost. Left: CO2, Middle: VERA, Right: CMH, zero precipitation points are not plotted. }
    \label{fig:vera_ac03}
\end{figure}

Lastly, there are many figures for the analysis of the MesoVICT core case 1, with the main analysis the report. Available are time series across all accumulation times (Figure \ref{fg:vera_ac0103_rho1}, \ref{fg:vera_ac0612_rho1}, and \ref{fig:vera_mass}); decompositions for \(\rho=1\) (Figure \ref{fig:vera_decomposition_grid_rho1}) and at the smaller, more localised value of \(\rho=0.01\) (Figure \ref{fig:vera_decomposition_grid_rho01}); and finally, illustrations of the model and observation precipitation in Figure \ref{fig:vera_ac06_illustration_1}, \ref{fig:vera_ac06_illustration_2}, and \ref{fig:vera_ac03}. 

\begin{figure}
    \centering
    \includegraphics[width=0.8\linewidth]{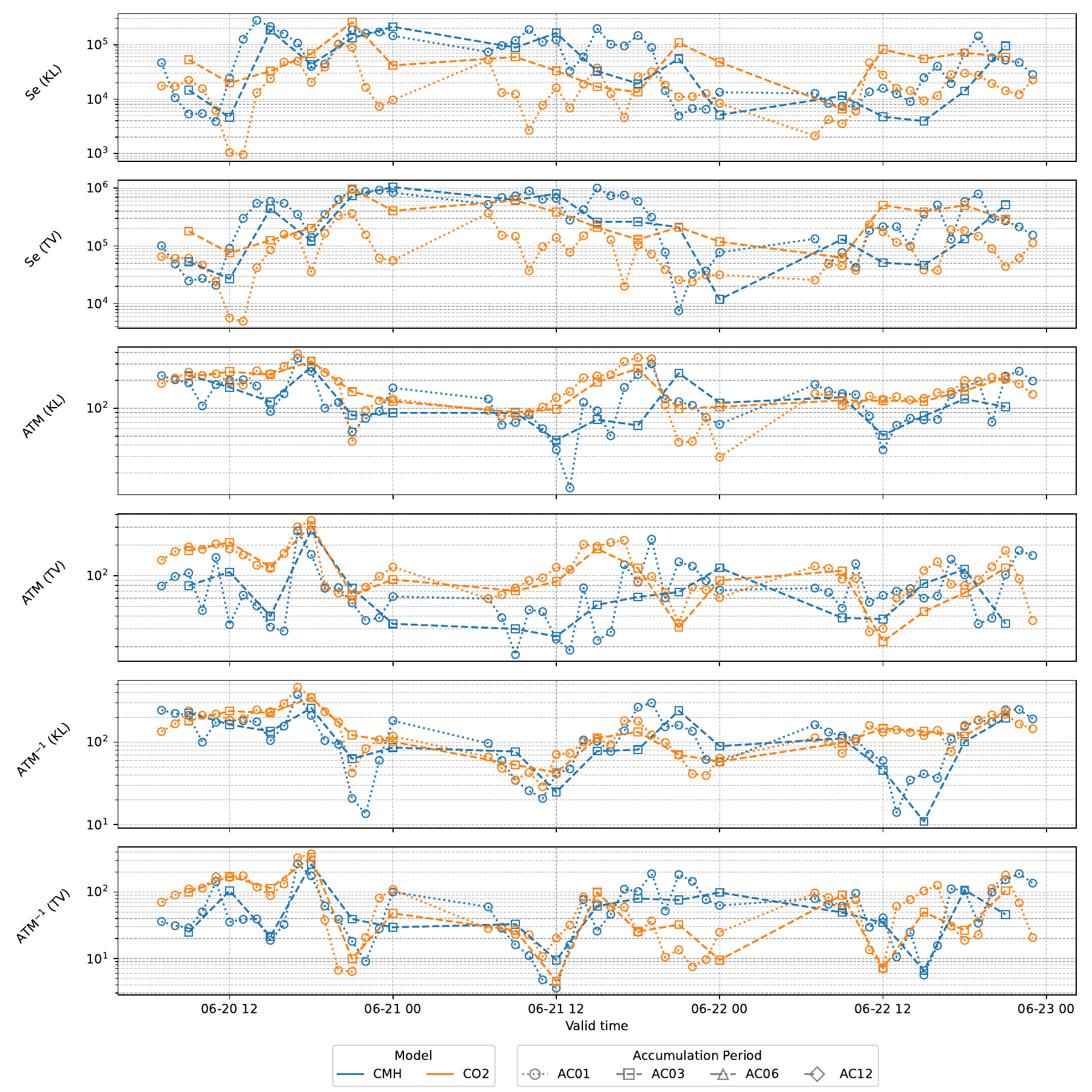}
    \caption{Time series plots for various scores for the AC01 and AC03 MesoVICT data in both models. Down in the rows we examine: \(\Sink^{KL}, \Sink^{TV},\) mean forward ATM (KL), mean forward ATM (TV), mean inverse ATM (KL), mean inverse ATM (TV). Given the short accumulation time, there are large variations in the scores; however, CO2 does appear less volatile whilst CMH, especially at AC01, appears more so. Though it improves with longer accumulation as timing errors are absorbed.   Here \(\varepsilon=0.005L^2, \rho=L^2\).
    }
    \label{fig:vera_ac0103_rho1}
\end{figure}

\begin{figure}
    \centering
    \includegraphics[width=0.8\linewidth]{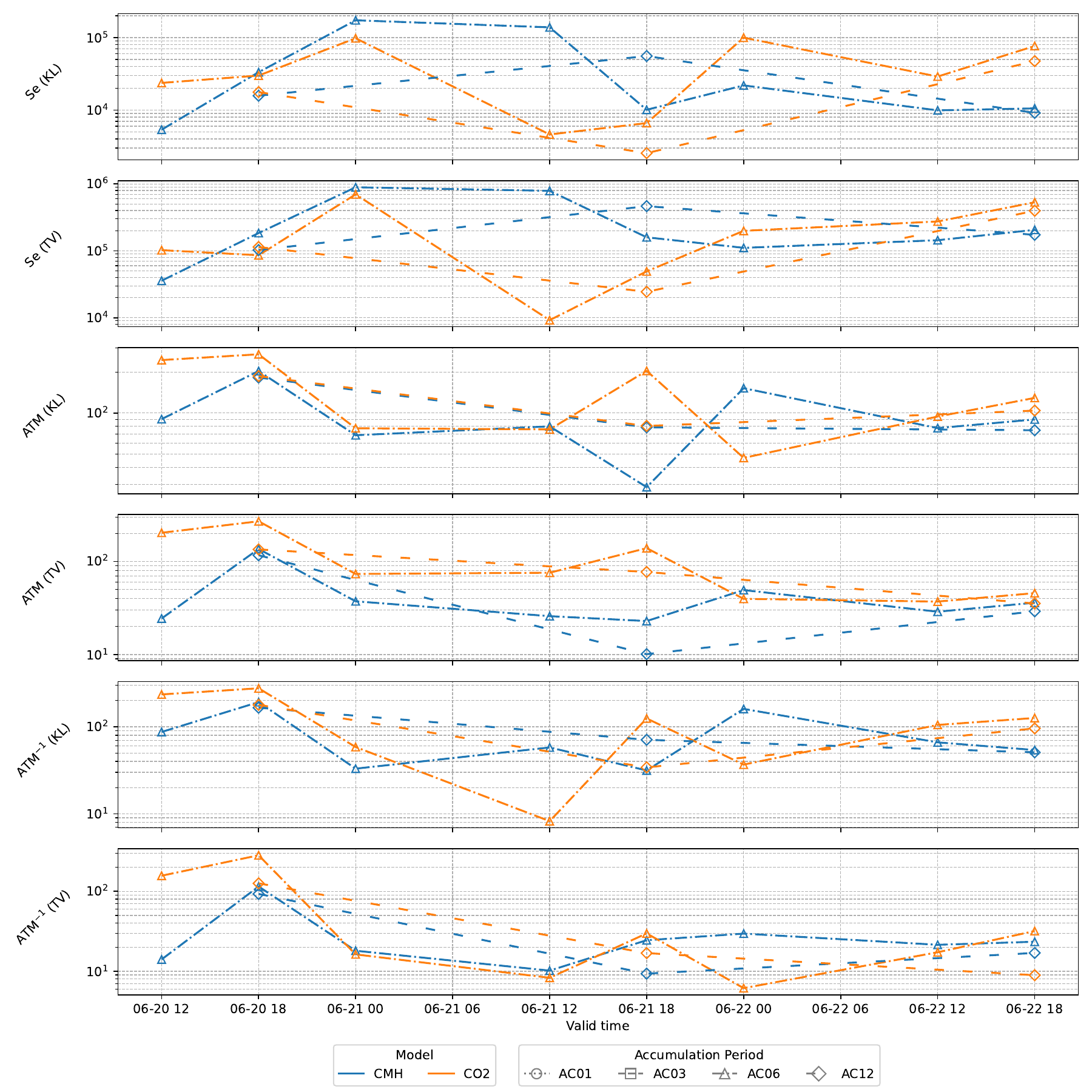}
    \caption{Time series plots for various scores for the AC06 and AC12 MesoVICT data in both models. Down in the rows we examine: \(\Sink^{KL}, \Sink^{TV},\) mean forward ATM (KL), mean forward ATM (TV), mean inverse ATM (KL), mean inverse ATM (TV). Given the long accumulation time, variations are smaller, thouhg there are still clear trends. Tiwth CO2 performing better thatn CMH as the front moves across the domain (round 06-21 12H). At later times, CMH does appear to start outperforming CO2. Here \(\varepsilon=0.005L^2, \rho=L^2\).}
    \label{fig:vera_ac0612_rho1}
\end{figure}


\subsection{Missing cases}\label{appendix:all_cases}
To avoid cherry-picking results the full array of testing cases are presented to align with \citet{gilleland_et_al_2019}. Together with the results above all cases included in the initial study are repeated here for UOT in two flavours. The top four horizontal bars display; \(\Sink^{TV}, \UOT^{TV}, \Sink^{KL}, \UOT^{KL}\). The lower table presents the mean (ATD, ATM) in both flavours, and with the forward and inverse vectors. The blue (darker) colour indicates observations, while the pale orange (lighter) represents forecasts. \(\varepsilon = 0.005L^2,\ \rho = L^2\).

\begin{figure}[h]
    \centering
    \includegraphics[width=1.0\linewidth, trim= 70 70 70 70, clip]{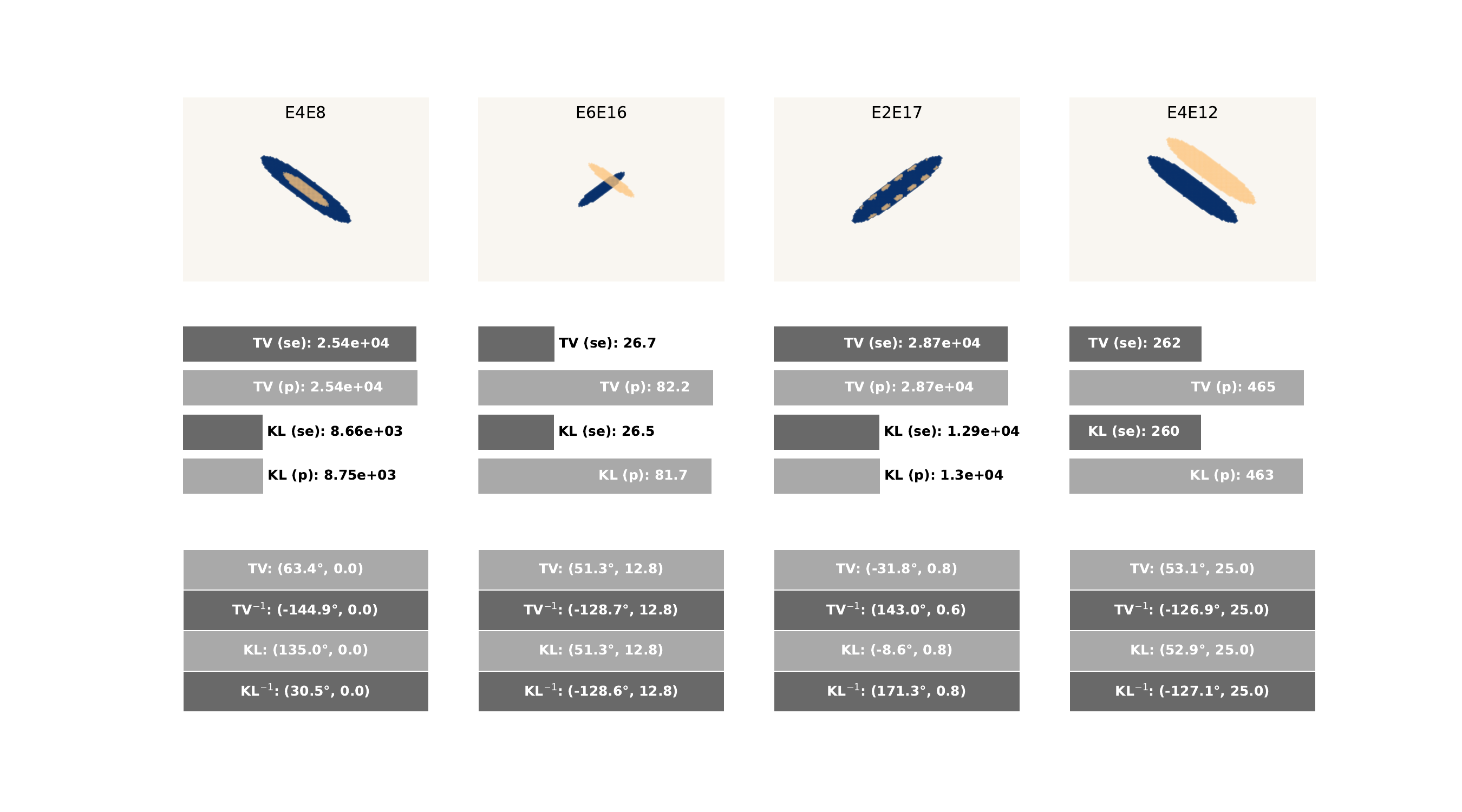}
    \caption{Cases E4E8, E6E16, E2E17, E4E12, which demonstrate different behaviours on scales, scatter showers in an envelope and simple transportation. All are considered more realistic shapes for complex terrain than the circular cases. The top four horizontal bars display; \(\Sink^{TV}, \UOT^{TV}, \Sink^{KL}, \UOT^{KL}\). The lower table presents the mean (ATD, ATM) in both flavours, and with the forward and inverse vectors. The blue (darker) colour indicates observations, while the pale orange (lighter) represents forecasts. \(\varepsilon = 0.005L^2,\ \rho = L^2\)}\label{fig:new_new_paper_ellipse_0}

\end{figure}

\begin{figure}[h]
    \centering
    \includegraphics[width=1.0\linewidth, trim= 70 70 70 70, clip]{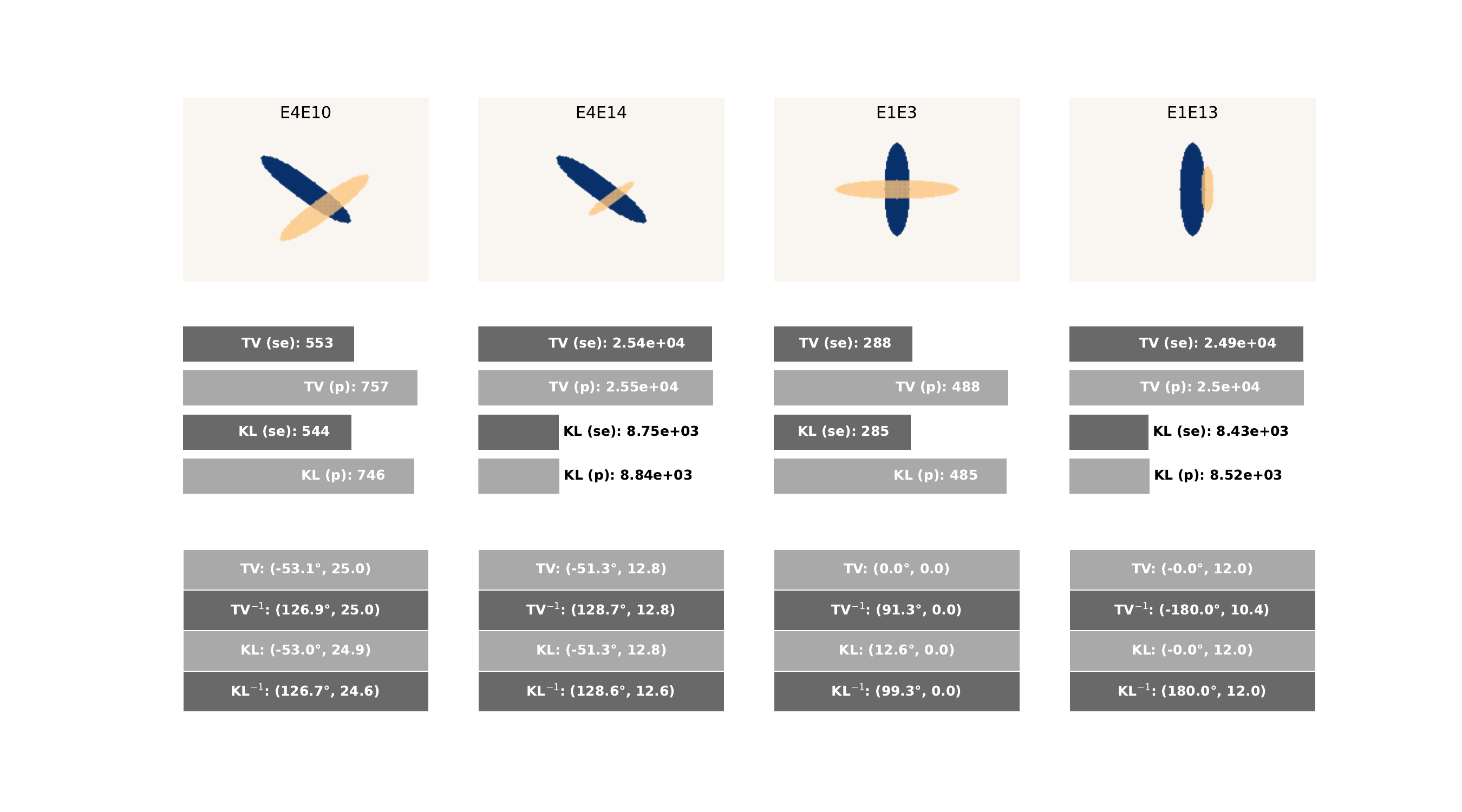}
    \caption{Cases E4E10, E4E14, E1E3, E1E13, which demonstrate different behaviours on scales, rotation, or aspect ratio errors and combinations of both. All are considered more realistic shapes for complex terrain than the circular cases. The top four horizontal bars display; \(\Sink^{TV}, \UOT^{TV}, \Sink^{KL}, \UOT^{KL}\). The lower table presents the mean (ATD, ATM) in both flavours, and with the forward and inverse vectors. The blue (darker) colour indicates observations, while the pale orange (lighter) represents forecasts. \(\varepsilon = 0.005L^2,\ \rho = L^2\)}\label{fig:new_new_paper_ellipse_1}

\end{figure}

\begin{figure}[h]
    \centering
    \includegraphics[width=1.0\linewidth, trim= 70 70 70 70, clip]{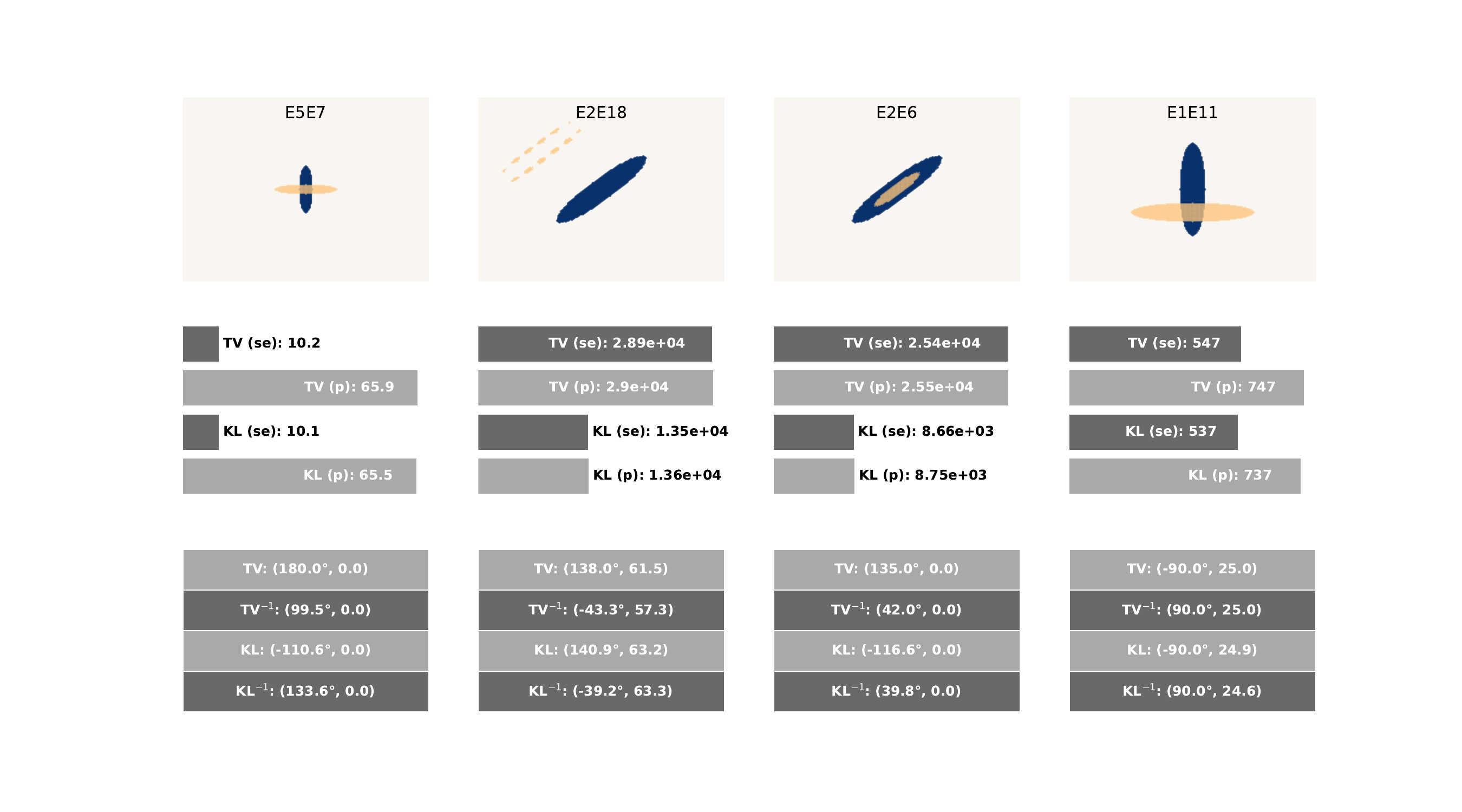}
    \caption{Cases E5E7, E2E18, E2E6, E1E11, which demonstrate different behaviours on scales, rotation with translation as well as missed scatter showers. All are considered more realistic shapes for complex terrain than the circular cases. The top four horizontal bars display; \(\Sink^{TV}, \UOT^{TV}, \Sink^{KL}, \UOT^{KL}\). The lower table presents the mean (ATD, ATM) in both flavours, and with the forward and inverse vectors. The blue (darker) colour indicates observations, while the pale orange (lighter) represents forecasts. \(\varepsilon = 0.005L^2,\ \rho = L^2\)}\label{fig:new_new_paper_ellipse_2}

\end{figure}

\begin{figure}[h]
    \centering
    \includegraphics[width=0.5\linewidth, trim= 70 70 70 70, clip]{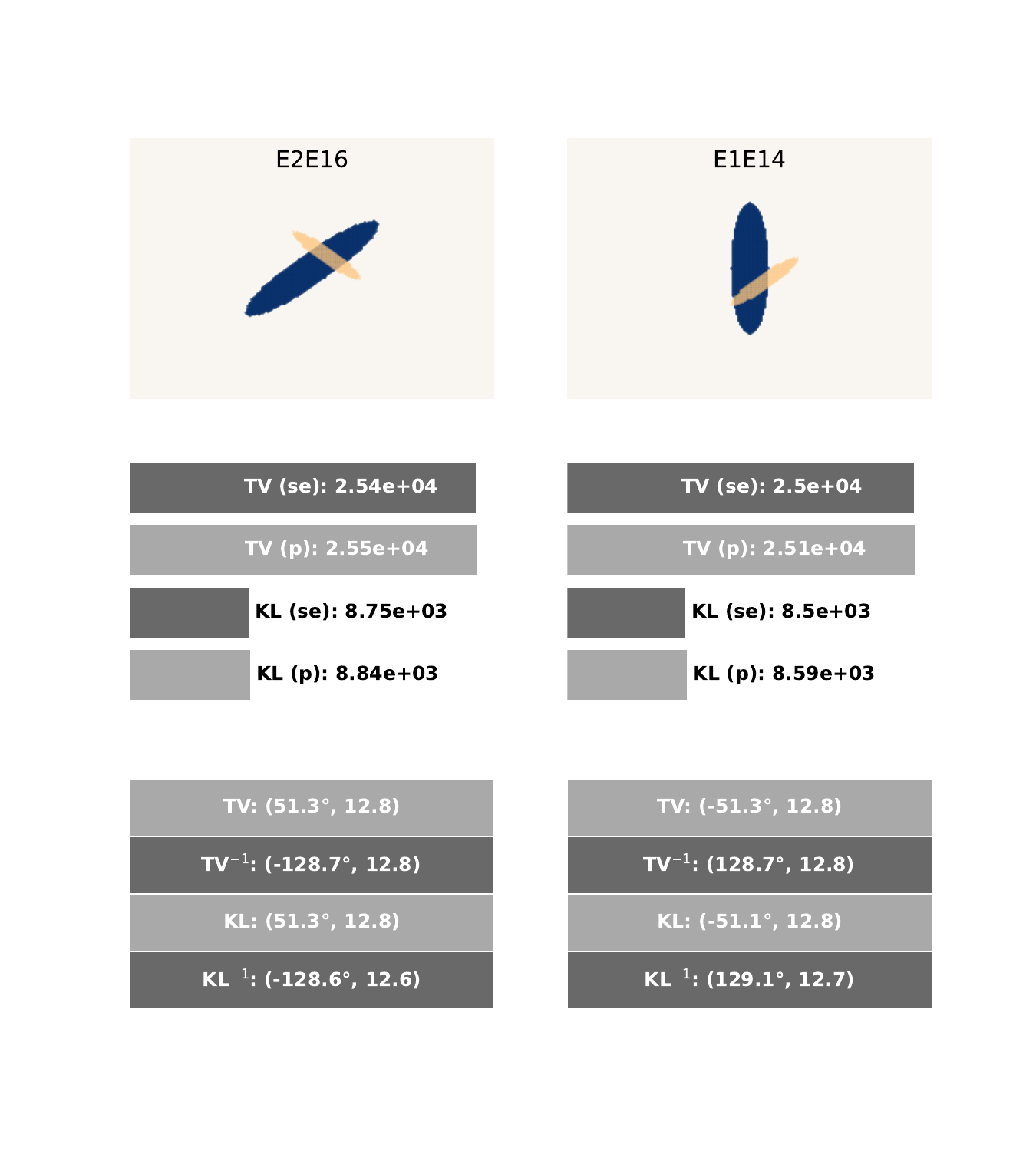}
    \caption{Cases E2E16, E1E14, which demonstrate different behaviours with scales rotation and transitional error combined. All are considered more realistic shapes for complex terrain than the circular cases. The top four horizontal bars display; \(\Sink^{TV}, \UOT^{TV}, \Sink^{KL}, \UOT^{KL}\). The lower table presents the mean (ATD, ATM) in both flavours, and with the forward and inverse vectors. The blue (darker) colour indicates observations, while the pale orange (lighter) represents forecasts. \(\varepsilon = 0.005L^2,\ \rho = L^2\)}\label{fig:new_new_paper_ellipse_3}
\end{figure}

\begin{figure}[h]
    \centering
    \includegraphics[width=1.0\linewidth, trim= 70 70 70 70, clip]{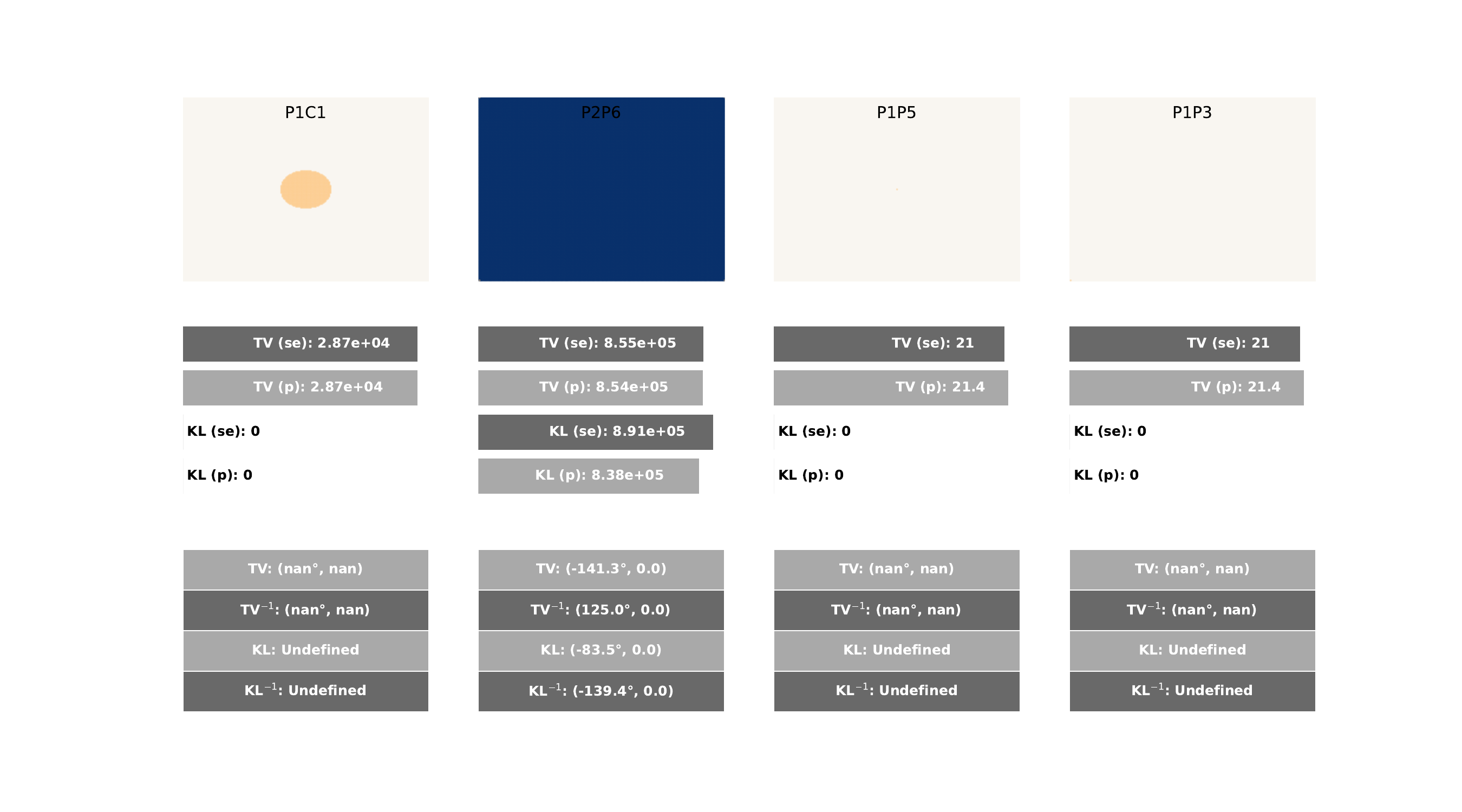}
    \caption{Cases P1C1, P2P6, P1P5, P1P3, recall that P1 is the null field, hence some quantities are undefined or zero. P1C1, P1P5, P1P3 KL (se) and KL (p) are defined as zero. The top four horizontal bars display; \(\Sink^{TV}, \UOT^{TV}, \Sink^{KL}, \UOT^{KL}\). The lower table presents the mean (ATD, ATM) in both flavours, and with the forward and inverse vectors. The blue (darker) colour indicates observations, while the pale orange (lighter) represents forecasts. \(\varepsilon = 0.005L^2,\ \rho = L^2\)}\label{fig:new_new_paper_pcase_0}
\end{figure}

\begin{figure}[h]
    \centering
    \includegraphics[width=0.75\linewidth, trim= 70 70 70 70, clip]{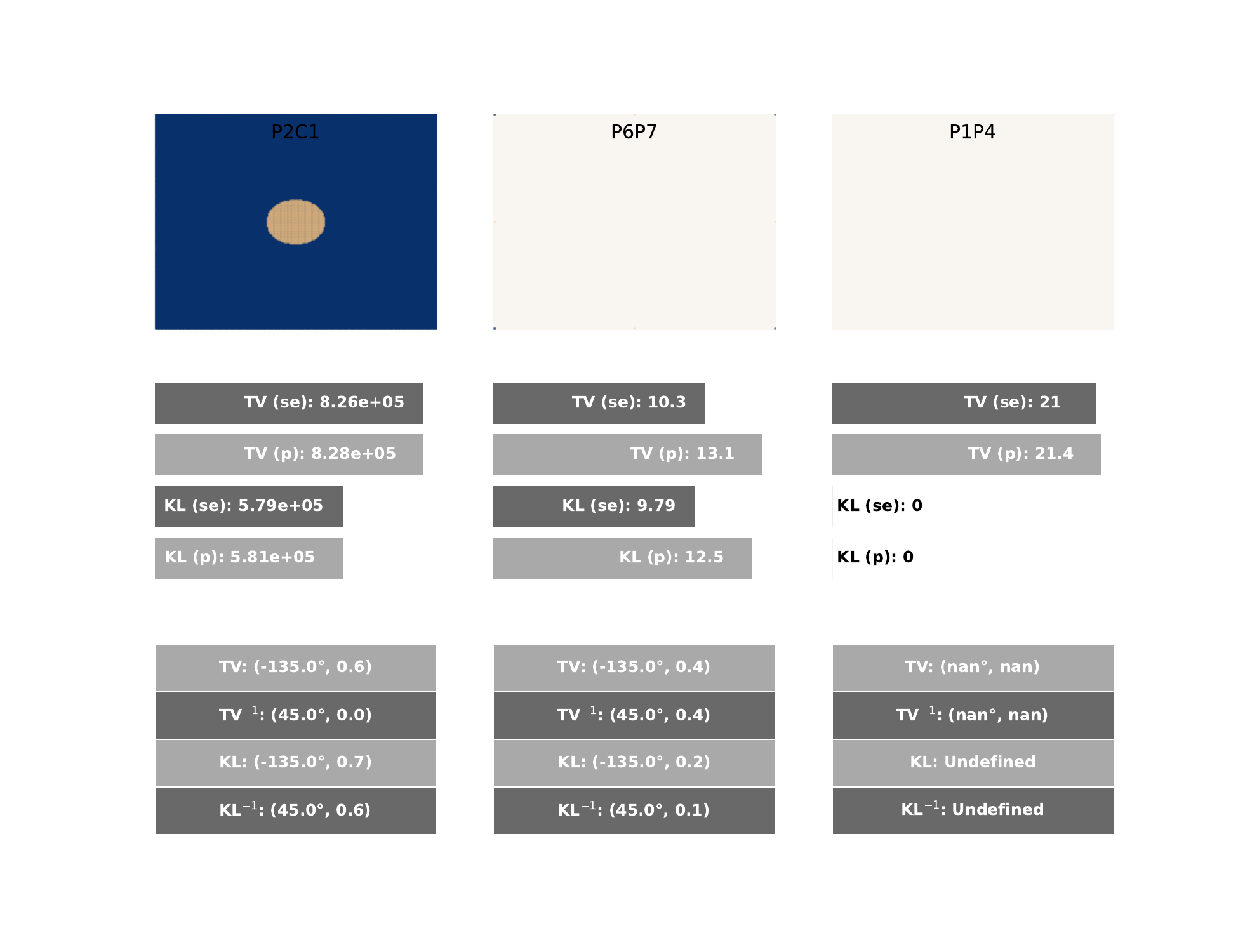}
    \caption{Cases P2C1, P6P7, P1P4, recall that P1 is the null field, hence some quantities are undefined or zero. P1P4, KL (se) and KL (p) are defined as zero. The top four horizontal bars display; \(\Sink^{TV}, \UOT^{TV}, \Sink^{KL}, \UOT^{KL}\). The lower table presents the mean (ATD, ATM) in both flavours, and with the forward and inverse vectors. The blue (darker) colour indicates observations, while the pale orange (lighter) represents forecasts. \(\varepsilon = 0.005L^2,\ \rho = L^2\)}\label{fig:new_new_paper_pcase_2}
\end{figure}

%% file: main.bbl
\begin{thebibliography}{75}
\expandafter\ifx\csname natexlab\endcsname\relax\def\natexlab#1{#1}\fi
\providecommand{\url}[1]{\texttt{#1}}
\providecommand{\href}[2]{#2}
\providecommand{\path}[1]{#1}
\providecommand{\DOIprefix}{doi:}
\providecommand{\ArXivprefix}{arXiv:}
\providecommand{\URLprefix}{URL: }
\providecommand{\Pubmedprefix}{pmid:}
\providecommand{\doi}[1]{\href{http://dx.doi.org/#1}{\path{#1}}}
\providecommand{\Pubmed}[1]{\href{pmid:#1}{\path{#1}}}
\providecommand{\bibinfo}[2]{#2}
\ifx\xfnm\relax \def\xfnm[#1]{\unskip,\space#1}\fi
%Type = Article
\bibitem[{Ahijevych et~al.(2009)Ahijevych, Gilleland, Brown and Ebert}]{ahijevychetal_2009}
\bibinfo{author}{Ahijevych, D.}, \bibinfo{author}{Gilleland, E.}, \bibinfo{author}{Brown, B.G.}, \bibinfo{author}{Ebert, E.E.}, \bibinfo{year}{2009}.
\newblock \bibinfo{title}{{Application of Spatial Verification Methods to Idealized and NWP-Gridded Precipitation Forecasts}}.
\newblock \bibinfo{journal}{Weather and Forecasting} \bibinfo{volume}{24}, \bibinfo{pages}{1485--1497}.
\newblock \DOIprefix\doi{ 10.1175/2009WAF2222298.1}.
%Type = Article
\bibitem[{Beier et~al.(2022)Beier, von Lindheim, Neumayer and Steidl}]{beier_etal_2022}
\bibinfo{author}{Beier, F.}, \bibinfo{author}{von Lindheim, J.}, \bibinfo{author}{Neumayer, S.}, \bibinfo{author}{Steidl, G.}, \bibinfo{year}{2022}.
\newblock \bibinfo{title}{{Unbalanced Multi-marginal Optimal Transport}}.
\newblock \bibinfo{journal}{Journal of Mathematical Imaging and Vision} , \bibinfo{pages}{1--20}\DOIprefix\doi{ 10.1007/S10851-022-01126-7/FIGURES/5}.
%Type = Article
\bibitem[{Benamou(2003)}]{benamou_2003}
\bibinfo{author}{Benamou, J.D.}, \bibinfo{year}{2003}.
\newblock \bibinfo{title}{Numerical resolution of an “unbalanced” mass transport problem}.
\newblock \bibinfo{journal}{ESAIM: Mathematical Modelling and Numerical Analysis} \bibinfo{volume}{37}, \bibinfo{pages}{851–868}.
\newblock \DOIprefix\doi{ 10.1051/m2an:2003058}.
%Type = Misc
\bibitem[{Benamou(2021)}]{benamou_2021}
\bibinfo{author}{Benamou, J.D.}, \bibinfo{year}{2021}.
\newblock \bibinfo{title}{{Optimal transportation, modelling and numerical simulation}}.
\newblock \DOIprefix\doi{ 10.1017/S0962492921000040}.
%Type = Article
\bibitem[{Benamou and Brenier(1998)}]{benamou_brenier_1998}
\bibinfo{author}{Benamou, J.D.}, \bibinfo{author}{Brenier, Y.}, \bibinfo{year}{1998}.
\newblock \bibinfo{title}{{Weak Existence for the Semigeostrophic Equations Formulated as a Coupled Monge–Amp{\`{e}}re/Transport Problem}}.
\newblock \bibinfo{journal}{SIAM Journal on Applied Mathematics} \bibinfo{volume}{58}, \bibinfo{pages}{1450--1461}.
\newblock \DOIprefix\doi{ 10.1137/S0036139995294111}.
%Type = Article
\bibitem[{Benamou and Brenier(2000)}]{benamou_brenier_2000}
\bibinfo{author}{Benamou, J.D.}, \bibinfo{author}{Brenier, Y.}, \bibinfo{year}{2000}.
\newblock \bibinfo{title}{{A computational fluid mechanics solution to the Monge-Kantorovich mass transfer problem}}.
\newblock \bibinfo{journal}{Numerische Mathematik 2000 84:3} \bibinfo{volume}{84}, \bibinfo{pages}{375--393}.
\newblock \DOIprefix\doi{ 10.1007/S002110050002}.
%Type = Article
\bibitem[{Benamou and Martinet(2020)}]{benamou_etal_2020}
\bibinfo{author}{Benamou, J.D.}, \bibinfo{author}{Martinet, M.}, \bibinfo{year}{2020}.
\newblock \bibinfo{title}{{Capacity Constrained Entropic Optimal Transport, Sinkhorn Saturated Domain Out-Summation and Vanishing Temperature}} .
%Type = Article
\bibitem[{Berman(2020)}]{berman_2020}
\bibinfo{author}{Berman, R.J.}, \bibinfo{year}{2020}.
\newblock \bibinfo{title}{{The Sinkhorn algorithm, parabolic optimal transport and geometric Monge–Amp{\`{e}}re equations}}.
\newblock \bibinfo{journal}{Numerische Mathematik} \bibinfo{volume}{145}, \bibinfo{pages}{771--836}.
\newblock \DOIprefix\doi{ 10.1007/s00211-020-01127-x}.
%Type = Article
\bibitem[{Bocquet et~al.(2024)Bocquet, Vanderbecken, Farchi, Dumont Le~Brazidec and Roustan}]{bocquet_et_all_2024}
\bibinfo{author}{Bocquet, M.}, \bibinfo{author}{Vanderbecken, P.J.}, \bibinfo{author}{Farchi, A.}, \bibinfo{author}{Dumont Le~Brazidec, J.}, \bibinfo{author}{Roustan, Y.}, \bibinfo{year}{2024}.
\newblock \bibinfo{title}{{Bridging classical data assimilation and optimal transport: the 3D-Var case}}.
\newblock \bibinfo{journal}{Nonlinear Processes in Geophysics} \bibinfo{volume}{31}, \bibinfo{pages}{335--357}.
\newblock \DOIprefix\doi{ 10.5194/npg-31-335-2024}.
%Type = Article
\bibitem[{Brown et~al.(2011)Brown, Gilleland and Ebert}]{brown_et_al_2012}
\bibinfo{author}{Brown, B.G.}, \bibinfo{author}{Gilleland, E.}, \bibinfo{author}{Ebert, E.E.}, \bibinfo{year}{2011}.
\newblock \bibinfo{title}{{Forecasts of Spatial Fields}}.
\newblock \bibinfo{journal}{Forecast Verification} , \bibinfo{pages}{95--117}\DOIprefix\doi{ 10.1002/9781119960003.CH6}.
%Type = Article
\bibitem[{Casati et~al.(2004)Casati, Ross and Stephenson}]{casati_ross_stephenson_2004}
\bibinfo{author}{Casati, B.}, \bibinfo{author}{Ross, G.}, \bibinfo{author}{Stephenson, D.B.}, \bibinfo{year}{2004}.
\newblock \bibinfo{title}{{A new intensity-scale approach for the verification of spatial precipitation forecasts}}.
\newblock \bibinfo{journal}{Meteorol. Appl} \bibinfo{volume}{11}, \bibinfo{pages}{141--154}.
\newblock \DOIprefix\doi{ 10.1017/S1350482704001239}.
%Type = Misc
\bibitem[{Charlier et~al.()Charlier, Feydy and Glaun{\`{e}}s}]{pykeops_webpage}
\bibinfo{author}{Charlier, B.}, \bibinfo{author}{Feydy, J.}, \bibinfo{author}{Glaun{\`{e}}s, J.A.}, .
\newblock \bibinfo{title}{{Kernel Operations on the GPU, with autodiff, without memory overflows — KeOps}}.
\newblock \URLprefix \url{ https://www.kernel-operations.io/keops/index.html}.
%Type = Article
\bibitem[{Charlier et~al.(2021)Charlier, Feydy, Glaun{\`{e}}s, Collin and Durif}]{pykeops_2021}
\bibinfo{author}{Charlier, B.}, \bibinfo{author}{Feydy, J.}, \bibinfo{author}{Glaun{\`{e}}s, J.A.}, \bibinfo{author}{Collin, F.D.}, \bibinfo{author}{Durif, G.}, \bibinfo{year}{2021}.
\newblock \bibinfo{title}{{Kernel Operations on the GPU, with Autodiff, without Memory Overflows}}.
\newblock \bibinfo{journal}{Journal of Machine Learning Research} \bibinfo{volume}{22}, \bibinfo{pages}{1--6}.
\newblock \URLprefix \url{ http://jmlr.org/papers/v22/20-275.html}.
%Type = Article
\bibitem[{Chizat(2024)}]{chizat_2024}
\bibinfo{author}{Chizat, L.}, \bibinfo{year}{2024}.
\newblock \bibinfo{title}{{Annealed Sinkhorn for Optimal Transport: convergence, regularization path and debiasing}} \URLprefix \url{ http://arxiv.org/abs/2408.11620}.
%Type = Article
\bibitem[{Chizat et~al.(2016)Chizat, Peyr{\'{e}}, Schmitzer and Vialard}]{chizat_peyre_et_al_2016}
\bibinfo{author}{Chizat, L.}, \bibinfo{author}{Peyr{\'{e}}, G.}, \bibinfo{author}{Schmitzer, B.}, \bibinfo{author}{Vialard, F.X.}, \bibinfo{year}{2016}.
\newblock \bibinfo{title}{{Scaling Algorithms for Unbalanced Transport Problems}} \DOIprefix\doi{ 10.48550/arxiv.1607.05816}.
%Type = Inproceedings
\bibitem[{Cuturi(2013)}]{cuturi_2013}
\bibinfo{author}{Cuturi, M.}, \bibinfo{year}{2013}.
\newblock \bibinfo{title}{{Sinkhorn Distances: Lightspeed Computation of Optimal Transport}}, in: \bibinfo{editor}{Burges, C.J.}, \bibinfo{editor}{Bottou, L.}, \bibinfo{editor}{Welling, M.}, \bibinfo{editor}{Ghahramani, Z.}, \bibinfo{editor}{Weinberger, K.Q.} (Eds.), \bibinfo{booktitle}{Advances in Neural Information Processing Systems}, \bibinfo{publisher}{Curran Associates, Inc.}
\newblock \URLprefix \url{ https://proceedings.neurips.cc/paper_files/paper/2013/file/af21d0c97db2e27e13572cbf59eb343d-Paper.pdf}.
%Type = Article
\bibitem[{Davis et~al.(2006)Davis, Brown and Bullock}]{davis_et_al_2006}
\bibinfo{author}{Davis, C.A.}, \bibinfo{author}{Brown, B.}, \bibinfo{author}{Bullock, R.}, \bibinfo{year}{2006}.
\newblock \bibinfo{title}{{Object-Based Verification of Precipitation Forecasts. Part I: Methodology and Application to Mesoscale Rain Areas}}.
\newblock \bibinfo{journal}{Monthly Weather Review} \bibinfo{volume}{134}, \bibinfo{pages}{1772--1784}.
\newblock \DOIprefix\doi{ 10.1175/MWR3145.1}.
%Type = Article
\bibitem[{Davis et~al.(2009)Davis, Brown, Bullock and Halley-Gotway}]{davis_et_al_2009}
\bibinfo{author}{Davis, C.A.}, \bibinfo{author}{Brown, B.G.}, \bibinfo{author}{Bullock, R.}, \bibinfo{author}{Halley-Gotway, J.}, \bibinfo{year}{2009}.
\newblock \bibinfo{title}{{The Method for Object-Based Diagnostic Evaluation (MODE) Applied to Numerical Forecasts from the 2005 NSSL/SPC Spring Program}}.
\newblock \bibinfo{journal}{Weather and Forecasting} \bibinfo{volume}{24}, \bibinfo{pages}{1252--1267}.
\newblock \DOIprefix\doi{ 10.1175/2009WAF2222241.1}.
%Type = Article
\bibitem[{Dorninger et~al.(2018)Dorninger, Gilleland, Casati, Mittermaier, Ebert, Brown and Wilson}]{dorninger_et_al_2018}
\bibinfo{author}{Dorninger, M.}, \bibinfo{author}{Gilleland, E.}, \bibinfo{author}{Casati, B.}, \bibinfo{author}{Mittermaier, M.P.}, \bibinfo{author}{Ebert, E.E.}, \bibinfo{author}{Brown, B.G.}, \bibinfo{author}{Wilson, L.J.}, \bibinfo{year}{2018}.
\newblock \bibinfo{title}{{The Setup of the MesoVICT Project}}.
\newblock \bibinfo{journal}{Bulletin of the American Meteorological Society} \bibinfo{volume}{99}, \bibinfo{pages}{1887--1906}.
\newblock \DOIprefix\doi{ 10.1175/BAMS-D-17-0164.1}.
%Type = Techreport
\bibitem[{Dorninger et~al.(2013)Dorninger, Mittermaier, Gilleland, Ebert, Brown and Wilson}]{dorninger_et_al_2013}
\bibinfo{author}{Dorninger, M.}, \bibinfo{author}{Mittermaier, M.P.}, \bibinfo{author}{Gilleland, E.}, \bibinfo{author}{Ebert, E.E.}, \bibinfo{author}{Brown, B.G.}, \bibinfo{author}{Wilson, L.J.}, \bibinfo{year}{2013}.
\newblock \bibinfo{title}{MesoVICT: Mesoscale Verification Inter-Comparison over Complex Terrain}.
\newblock \bibinfo{type}{Technical Report} \bibinfo{number}{NCAR/TN-505+STR}. NCAR.
\newblock \DOIprefix\doi{ 10.5065/D6416V21}.
%Type = Article
\bibitem[{Ebert and McBride(2000)}]{ebert_mcbridge_2000}
\bibinfo{author}{Ebert, E.}, \bibinfo{author}{McBride, J.}, \bibinfo{year}{2000}.
\newblock \bibinfo{title}{Verification of precipitation in weather systems: determination of systematic errors}.
\newblock \bibinfo{journal}{Journal of Hydrology} \bibinfo{volume}{239}, \bibinfo{pages}{179--202}.
\newblock \DOIprefix\doi{ https://doi.org/10.1016/S0022-1694(00)00343-7}.
%Type = Article
\bibitem[{Ebert(2008)}]{ebert_2008}
\bibinfo{author}{Ebert, E.E.}, \bibinfo{year}{2008}.
\newblock \bibinfo{title}{{Fuzzy verification of high-resolution gridded forecasts: a review and proposed framework}}.
\newblock \bibinfo{journal}{Meteorological Applications} \bibinfo{volume}{15}, \bibinfo{pages}{51--64}.
\newblock \DOIprefix\doi{ 10.1002/met.25}.
%Type = Article
\bibitem[{Ebert and Gallus(2009)}]{ebert_gallus_2009}
\bibinfo{author}{Ebert, E.E.}, \bibinfo{author}{Gallus, W.A.}, \bibinfo{year}{2009}.
\newblock \bibinfo{title}{Toward better understanding of the contiguous rain area (cra) method for spatial forecast verification}.
\newblock \bibinfo{journal}{Weather and Forecasting} \bibinfo{volume}{24}, \bibinfo{pages}{1401 -- 1415}.
\newblock \DOIprefix\doi{ 10.1175/2009WAF2222252.1}.
%Type = Misc
\bibitem[{{Eric Gilleland}(n.d.)}]{eric_reference_list}
\bibinfo{author}{{Eric Gilleland}}, \bibinfo{year}{n.d.}
\newblock \bibinfo{title}{{Spatial Forecast Verification Reference List}}.
\newblock \URLprefix \url{ https://projects.ral.ucar.edu/icp/references.html}.
%Type = Article
\bibitem[{Farchi et~al.(2016)Farchi, Bocquet, Roustan, Mathieu and Qu{\'{e}}rel}]{farchi_et_al_2016}
\bibinfo{author}{Farchi, A.}, \bibinfo{author}{Bocquet, M.}, \bibinfo{author}{Roustan, Y.}, \bibinfo{author}{Mathieu, A.}, \bibinfo{author}{Qu{\'{e}}rel, A.}, \bibinfo{year}{2016}.
\newblock \bibinfo{title}{{Using the Wasserstein distance to compare fields of pollutants: application to the radionuclide atmospheric dispersion of the Fukushima-Daiichi accident}}.
\newblock \bibinfo{journal}{Tellus B: Chemical and Physical Meteorology} \DOIprefix\doi{ 10.3402/tellusb.v68.31682}.
%Type = Misc
\bibitem[{Feydy and Roussillon(n.d.)}]{geomloss_webpage}
\bibinfo{author}{Feydy, J.}, \bibinfo{author}{Roussillon, P.}, \bibinfo{year}{n.d.}
\newblock \bibinfo{title}{{Geometric Loss functions between sampled measures, images and volumes — GeomLoss}}.
\newblock \URLprefix \url{ https://www.kernel-operations.io/geomloss/index.html}.
%Type = Inproceedings
\bibitem[{Feydy et~al.(2019)Feydy, S{\'{e}}journ{\'{e}}, Vialard, Amari, Trouve and Peyr{\'{e}}}]{feydy_sejourne_vialard_2018}
\bibinfo{author}{Feydy, J.}, \bibinfo{author}{S{\'{e}}journ{\'{e}}, T.}, \bibinfo{author}{Vialard, F.X.}, \bibinfo{author}{Amari, S.i.}, \bibinfo{author}{Trouve, A.}, \bibinfo{author}{Peyr{\'{e}}, G.}, \bibinfo{year}{2019}.
\newblock \bibinfo{title}{{Interpolating between Optimal Transport and MMD using Sinkhorn Divergences}}, in: \bibinfo{editor}{Chaudhuri, K.}, \bibinfo{editor}{Sugiyama, M.} (Eds.), \bibinfo{booktitle}{Proceedings of the Twenty-Second International Conference on Artificial Intelligence and Statistics}, \bibinfo{publisher}{PMLR}. pp. \bibinfo{pages}{2681--2690}.
\newblock \URLprefix \url{ https://proceedings.mlr.press/v89/feydy19a.html}.
%Type = Article
\bibitem[{Feyeux et~al.(2018)Feyeux, Vidard and Nodet}]{feyuex_vidard_nodet_2018}
\bibinfo{author}{Feyeux, N.}, \bibinfo{author}{Vidard, A.}, \bibinfo{author}{Nodet, M.}, \bibinfo{year}{2018}.
\newblock \bibinfo{title}{{Optimal transport for variational data assimilation}}.
\newblock \bibinfo{journal}{Nonlinear Processes in Geophysics} \bibinfo{volume}{25}, \bibinfo{pages}{55--66}.
\newblock \DOIprefix\doi{ 10.5194/npg-25-55-2018}.
%Type = Article
\bibitem[{Flamary et~al.(2021)Flamary, Courty, Gramfort, Alaya, Boisbunon, Chambon, Chapel, Corenflos, Fatras, Fournier, Gautheron, Gayraud, Janati, Rakotomamonjy, Redko, Rolet, Schutz, Seguy, Sutherland, Tavenard, Tong and Vayer}]{pot_2021}
\bibinfo{author}{Flamary, R.}, \bibinfo{author}{Courty, N.}, \bibinfo{author}{Gramfort, A.}, \bibinfo{author}{Alaya, M.Z.}, \bibinfo{author}{Boisbunon, A.}, \bibinfo{author}{Chambon, S.}, \bibinfo{author}{Chapel, L.}, \bibinfo{author}{Corenflos, A.}, \bibinfo{author}{Fatras, K.}, \bibinfo{author}{Fournier, N.}, \bibinfo{author}{Gautheron, L.}, \bibinfo{author}{Gayraud, N.T.H.}, \bibinfo{author}{Janati, H.}, \bibinfo{author}{Rakotomamonjy, A.}, \bibinfo{author}{Redko, I.}, \bibinfo{author}{Rolet, A.}, \bibinfo{author}{Schutz, A.}, \bibinfo{author}{Seguy, V.}, \bibinfo{author}{Sutherland, D.J.}, \bibinfo{author}{Tavenard, R.}, \bibinfo{author}{Tong, A.}, \bibinfo{author}{Vayer, T.}, \bibinfo{year}{2021}.
\newblock \bibinfo{title}{{POT: Python Optimal Transport}}.
\newblock \bibinfo{journal}{Journal of Machine Learning Research} \bibinfo{volume}{22}, \bibinfo{pages}{1--8}.
\newblock \URLprefix \url{ http://jmlr.org/papers/v22/20-451.html}.
%Type = Misc
\bibitem[{Francis(n.d.)}]{uot_own_implementation}
\bibinfo{author}{Francis, J.}, \bibinfo{year}{n.d.}
\newblock \bibinfo{title}{{UOTforSpatialVerification: Unbalanced Optimal transport implementation with total variation or Kullback Leiber marginal penalisation for Precipitation forecast verification.}}
\newblock \URLprefix \url{ https://github.com/Jacob-Francis/UOTforSpatialVerification}.
%Type = Article
\bibitem[{Gilleland(2013)}]{gilleland_2013}
\bibinfo{author}{Gilleland, E.}, \bibinfo{year}{2013}.
\newblock \bibinfo{title}{Testing competing precipitation forecasts accurately and efficiently: The spatial prediction comparison test}.
\newblock \bibinfo{journal}{Monthly Weather Review} \bibinfo{volume}{141}, \bibinfo{pages}{340 -- 355}.
\newblock \DOIprefix\doi{ 10.1175/MWR-D-12-00155.1}.
%Type = Article
\bibitem[{Gilleland(2021)}]{gilleland_2021}
\bibinfo{author}{Gilleland, E.}, \bibinfo{year}{2021}.
\newblock \bibinfo{title}{Novel measures for summarizing high-resolution forecast performance}.
\newblock \bibinfo{journal}{Advances in Statistical Climatology, Meteorology and Oceanography} \bibinfo{volume}{7}, \bibinfo{pages}{13--34}.
\newblock \DOIprefix\doi{ 10.5194/ascmo-7-13-2021}.
%Type = Article
\bibitem[{Gilleland(2022)}]{cran_spatialvx}
\bibinfo{author}{Gilleland, E.}, \bibinfo{year}{2022}.
\newblock \bibinfo{title}{Comparing spatial fields with {SpatialVx}: Spatial forecast verification in {R}}.
\newblock \bibinfo{journal}{Unpublished} \bibinfo{volume}{0}, \bibinfo{pages}{1--69}.
\newblock \DOIprefix\doi{ 10.5065/4px3-5a05}.
%Type = Article
\bibitem[{Gilleland et~al.(2009)Gilleland, Ahijevych, Brown, Casati and Ebert}]{gillelandetal_2009a}
\bibinfo{author}{Gilleland, E.}, \bibinfo{author}{Ahijevych, D.}, \bibinfo{author}{Brown, B.G.}, \bibinfo{author}{Casati, B.}, \bibinfo{author}{Ebert, E.E.}, \bibinfo{year}{2009}.
\newblock \bibinfo{title}{{Intercomparison of Spatial Forecast Verification Methods}}.
\newblock \bibinfo{journal}{Weather and Forecasting} \bibinfo{volume}{24}, \bibinfo{pages}{1416--1430}.
\newblock \DOIprefix\doi{ 10.1175/2009WAF2222269.1}.
%Type = Article
\bibitem[{Gilleland et~al.(2010a)Gilleland, Ahijevych, Brown and Ebert}]{gilleland_et_al_2010}
\bibinfo{author}{Gilleland, E.}, \bibinfo{author}{Ahijevych, D.A.}, \bibinfo{author}{Brown, B.G.}, \bibinfo{author}{Ebert, E.E.}, \bibinfo{year}{2010}a.
\newblock \bibinfo{title}{{Verifying forecasts spatially}}.
\newblock \bibinfo{journal}{Bulletin of the American Meteorological Society} \bibinfo{volume}{91}, \bibinfo{pages}{1365--1373}.
\newblock \DOIprefix\doi{ 10.1175/2010BAMS2819.1}.
%Type = Article
\bibitem[{Gilleland et~al.(2010b)Gilleland, Lindström and Lindgren}]{gilleland_2010}
\bibinfo{author}{Gilleland, E.}, \bibinfo{author}{Lindström, J.}, \bibinfo{author}{Lindgren, F.}, \bibinfo{year}{2010}b.
\newblock \bibinfo{title}{Analyzing the image warp forecast verification method on precipitation fields from the icp}.
\newblock \bibinfo{journal}{Weather and Forecasting} \bibinfo{volume}{25}, \bibinfo{pages}{1249 -- 1262}.
\newblock \DOIprefix\doi{ 10.1175/2010WAF2222365.1}.
%Type = Article
\bibitem[{Gilleland et~al.(2020)Gilleland, Skok, Brown, Casati, Dorninger, Mittermaier, Roberts and Wilson}]{gilleland_et_al_2019}
\bibinfo{author}{Gilleland, E.}, \bibinfo{author}{Skok, G.}, \bibinfo{author}{Brown, B.G.}, \bibinfo{author}{Casati, B.}, \bibinfo{author}{Dorninger, M.}, \bibinfo{author}{Mittermaier, M.P.}, \bibinfo{author}{Roberts, N.}, \bibinfo{author}{Wilson, L.J.}, \bibinfo{year}{2020}.
\newblock \bibinfo{title}{{A Novel Set of Geometric Verification Test Fields with Application to Distance Measures}}.
\newblock \bibinfo{journal}{Monthly Weather Review} \bibinfo{volume}{148}, \bibinfo{pages}{1653--1673}.
\newblock \DOIprefix\doi{ 10.1175/MWR-D-19-0256.1}.
%Type = Inproceedings
\bibitem[{{Gramfort A} et~al.(2015){Gramfort A}, Peyr{\'{e}} and {Cuturi M}}]{gramfort_peyre_cuturi_2015}
\bibinfo{author}{{Gramfort A}}, \bibinfo{author}{Peyr{\'{e}}, G.}, \bibinfo{author}{{Cuturi M}}, \bibinfo{year}{2015}.
\newblock \bibinfo{title}{{Fast Optimal Transport Averaging of Neuroimaging Data}}, in: \bibinfo{editor}{Ourselin, S.}, \bibinfo{editor}{Alexander, D.}, \bibinfo{editor}{Westin, C.}, \bibinfo{editor}{Cardoso, M.} (Eds.), \bibinfo{booktitle}{Information Processing in Medical Imaging}, \bibinfo{publisher}{Springer International Publishing}, \bibinfo{address}{Cham}. pp. \bibinfo{pages}{261--272}.
%Type = Inproceedings
\bibitem[{Grauman and Darrell(2005)}]{grauman_darrell_2005}
\bibinfo{author}{Grauman, K.}, \bibinfo{author}{Darrell, T.}, \bibinfo{year}{2005}.
\newblock \bibinfo{title}{{The pyramid match kernel: discriminative classification with sets of image features}}, in: \bibinfo{booktitle}{Tenth IEEE International Conference on Computer Vision (ICCV'05) Volume 1}, pp. \bibinfo{pages}{1458--1465}.
\newblock \DOIprefix\doi{ 10.1109/ICCV.2005.239}.
%Type = Article
\bibitem[{Hou et~al.(2014)Hou, Kakar, Neeck, Azarbarzin, Kummerow, Kojima, Oki, Nakamura and Iguchi}]{hou_et_al_2014}
\bibinfo{author}{Hou, A.Y.}, \bibinfo{author}{Kakar, R.K.}, \bibinfo{author}{Neeck, S.}, \bibinfo{author}{Azarbarzin, A.A.}, \bibinfo{author}{Kummerow, C.D.}, \bibinfo{author}{Kojima, M.}, \bibinfo{author}{Oki, R.}, \bibinfo{author}{Nakamura, K.}, \bibinfo{author}{Iguchi, T.}, \bibinfo{year}{2014}.
\newblock \bibinfo{title}{{The Global Precipitation Measurement Mission}}.
\newblock \bibinfo{journal}{Bulletin of the American Meteorological Society} \bibinfo{volume}{95}, \bibinfo{pages}{701--722}.
\newblock \DOIprefix\doi{ 10.1175/BAMS-D-13-00164.1}.
%Type = Article
\bibitem[{Hyun et~al.(2022)Hyun, Mishra, Follett, Jonsson, Kulk, Forget, Racault, Jackson, Dutkiewicz, M{\"{u}}ller and Bien}]{hyan_et_al_2022}
\bibinfo{author}{Hyun, S.}, \bibinfo{author}{Mishra, A.}, \bibinfo{author}{Follett, C.L.}, \bibinfo{author}{Jonsson, B.}, \bibinfo{author}{Kulk, G.}, \bibinfo{author}{Forget, G.}, \bibinfo{author}{Racault, M.F.}, \bibinfo{author}{Jackson, T.}, \bibinfo{author}{Dutkiewicz, S.}, \bibinfo{author}{M{\"{u}}ller, C.L.}, \bibinfo{author}{Bien, J.}, \bibinfo{year}{2022}.
\newblock \bibinfo{title}{{Ocean movers distance: using optimal transport for analysing oceanographic data}}.
\newblock \bibinfo{journal}{Proceedings of the Royal Society A} \bibinfo{volume}{478}.
\newblock \DOIprefix\doi{ 10.1098/RSPA.2021.0875}.
%Type = Article
\bibitem[{Kain et~al.(2008)Kain, Weiss, Bright, Baldwin, Levit, Carbin, Schwartz, Weisman, Droegemeier, Weber and Thomas}]{Kain_et_al_2008}
\bibinfo{author}{Kain, J.S.}, \bibinfo{author}{Weiss, S.J.}, \bibinfo{author}{Bright, D.R.}, \bibinfo{author}{Baldwin, M.E.}, \bibinfo{author}{Levit, J.J.}, \bibinfo{author}{Carbin, G.W.}, \bibinfo{author}{Schwartz, C.S.}, \bibinfo{author}{Weisman, M.L.}, \bibinfo{author}{Droegemeier, K.K.}, \bibinfo{author}{Weber, D.B.}, \bibinfo{author}{Thomas, K.W.}, \bibinfo{year}{2008}.
\newblock \bibinfo{title}{{Some Practical Considerations Regarding Horizontal Resolution in the First Generation of Operational Convection-Allowing NWP}}.
\newblock \bibinfo{journal}{Weather and Forecasting} \bibinfo{volume}{23}, \bibinfo{pages}{931--952}.
\newblock \DOIprefix\doi{ 10.1175/WAF2007106.1}.
%Type = Article
\bibitem[{Kantorovitch(1958)}]{kantoro_1958}
\bibinfo{author}{Kantorovitch, L.}, \bibinfo{year}{1958}.
\newblock \bibinfo{title}{{On the Translocation of Masses}}.
\newblock \bibinfo{journal}{Management Science} \bibinfo{volume}{5}, \bibinfo{pages}{1--4}.
\newblock \URLprefix \url{ http://www.jstor.org/stable/2626967}.
%Type = Article
\bibitem[{Keil and Craig(2007)}]{keil_craig_2007}
\bibinfo{author}{Keil, C.}, \bibinfo{author}{Craig, G.C.}, \bibinfo{year}{2007}.
\newblock \bibinfo{title}{{A Displacement-Based Error Measure Applied in a Regional Ensemble Forecasting System}}.
\newblock \bibinfo{journal}{Monthly Weather Review} \bibinfo{volume}{135}, \bibinfo{pages}{3248--3259}.
\newblock \DOIprefix\doi{ 10.1175/MWR3457.1}.
%Type = Article
\bibitem[{Keil and Craig(2009)}]{keil_craig_2009}
\bibinfo{author}{Keil, C.}, \bibinfo{author}{Craig, G.C.}, \bibinfo{year}{2009}.
\newblock \bibinfo{title}{{A Displacement and Amplitude Score Employing an Optical Flow Technique}}.
\newblock \bibinfo{journal}{Weather and Forecasting} \bibinfo{volume}{24}, \bibinfo{pages}{1297 -- 1308}.
\newblock \DOIprefix\doi{ 10.1175/2009WAF2222247.1}.
%Type = Article
\bibitem[{Kolouri et~al.(2017)Kolouri, Park, Thorpe, Slepcev and Rohde}]{kolouri_et_al_2017}
\bibinfo{author}{Kolouri, S.}, \bibinfo{author}{Park, S.R.}, \bibinfo{author}{Thorpe, M.}, \bibinfo{author}{Slepcev, D.}, \bibinfo{author}{Rohde, G.K.}, \bibinfo{year}{2017}.
\newblock \bibinfo{title}{{Optimal Mass Transport: Signal processing and machine-learning applications}}.
\newblock \bibinfo{journal}{IEEE Signal Processing Magazine} \bibinfo{volume}{34}, \bibinfo{pages}{43--59}.
\newblock \DOIprefix\doi{ 10.1109/MSP.2017.2695801}.
%Type = Article
\bibitem[{Liero et~al.(2016)Liero, Mielke and Savar{\'{e}}}]{liero_mielke_savare_2016}
\bibinfo{author}{Liero, M.}, \bibinfo{author}{Mielke, A.}, \bibinfo{author}{Savar{\'{e}}, G.}, \bibinfo{year}{2016}.
\newblock \bibinfo{title}{{Optimal transport in competition with reaction: The Hellinger-Kantorovich distance and geodesic curves}}.
\newblock \bibinfo{journal}{SIAM Journal on Mathematical Analysis} \bibinfo{volume}{48}, \bibinfo{pages}{2869--2911}.
\newblock \DOIprefix\doi{ 10.1137/15M1041420}.
%Type = Inproceedings
\bibitem[{Lin and Mitchell(2005)}]{lin_mitchell_2005}
\bibinfo{author}{Lin, Y.}, \bibinfo{author}{Mitchell, K.E.}, \bibinfo{year}{2005}.
\newblock \bibinfo{title}{{The NCEP Stage II/IV Hourly Precipitation Analyses: Development and Applications}}, in: \bibinfo{booktitle}{Proceedings of the 19th Conference on Hydrology}, \bibinfo{publisher}{American Meteorological Society}, \bibinfo{address}{San Diego, CA, USA}.
%Type = Article
\bibitem[{Liu and Frank(2022)}]{liu_frank_2022}
\bibinfo{author}{Liu, X.}, \bibinfo{author}{Frank, J.}, \bibinfo{year}{2022}.
\newblock \bibinfo{title}{{Ensemble data assimilation using optimal control in the Wasserstein metric}}.
\newblock \bibinfo{journal}{Journal of Computational Science} \bibinfo{volume}{65}, \bibinfo{pages}{101895}.
\newblock \DOIprefix\doi{ 10.1016/J.JOCS.2022.101895}.
%Type = Inproceedings
\bibitem[{Lucas and Kanade(1981)}]{lucas_kanade_1981}
\bibinfo{author}{Lucas, B.D.}, \bibinfo{author}{Kanade, T.}, \bibinfo{year}{1981}.
\newblock \bibinfo{title}{{An Iterative Image Registration Technique with an Application to Stereo Vision}}, in: \bibinfo{booktitle}{IJCAI'81: 7th international joint conference on Artificial intelligence}, \bibinfo{address}{Vancouver, Canada}. pp. \bibinfo{pages}{674--679}.
\newblock \URLprefix \url{ https://hal.science/hal-03697340}.
%Type = Article
\bibitem[{Magyar and Sambridge(2023)}]{magyar_sambridge_2023}
\bibinfo{author}{Magyar, J.C.}, \bibinfo{author}{Sambridge, M.}, \bibinfo{year}{2023}.
\newblock \bibinfo{title}{{Hydrological objective functions and ensemble averaging with the Wasserstein distance}}.
\newblock \bibinfo{journal}{Hydrology and Earth System Sciences} \bibinfo{volume}{27}, \bibinfo{pages}{991--1010}.
\newblock \DOIprefix\doi{ 10.5194/hess-27-991-2023}.
%Type = Article
\bibitem[{Marzban and Sandgathe(2010)}]{marzban_sandgathe_2010}
\bibinfo{author}{Marzban, C.}, \bibinfo{author}{Sandgathe, S.}, \bibinfo{year}{2010}.
\newblock \bibinfo{title}{{Optical Flow for Verification}}.
\newblock \bibinfo{journal}{Weather and Forecasting} \bibinfo{volume}{25}, \bibinfo{pages}{1479--1494}.
\newblock \DOIprefix\doi{ 10.1175/2010WAF2222351.1}.
%Type = Article
\bibitem[{Marzban et~al.(2009)Marzban, Sandgathe, Lyons and Lederer}]{marzban_et_al_2009}
\bibinfo{author}{Marzban, C.}, \bibinfo{author}{Sandgathe, S.}, \bibinfo{author}{Lyons, H.}, \bibinfo{author}{Lederer, N.}, \bibinfo{year}{2009}.
\newblock \bibinfo{title}{{Three spatial verification techniques: Cluster analysis, variogram, and optical flow}}.
\newblock \bibinfo{journal}{Weather and Forecasting} \bibinfo{volume}{24}, \bibinfo{pages}{1457--1471}.
\newblock \DOIprefix\doi{ 10.1175/2009WAF2222261.1}.
%Type = Article
\bibitem[{Mass et~al.(2002)Mass, Ovens, Westrick and Colle}]{mass_et_al_2002}
\bibinfo{author}{Mass, C.F.}, \bibinfo{author}{Ovens, D.}, \bibinfo{author}{Westrick, K.}, \bibinfo{author}{Colle, B.A.}, \bibinfo{year}{2002}.
\newblock \bibinfo{title}{{DOES INCREASING HORIZONTAL RESOLUTION PRODUCE MORE SKILLFUL FORECASTS?: The Results of Two Years of Real-Time Numerical Weather Prediction over the Pacific Northwest}}.
\newblock \bibinfo{journal}{Bulletin of the American Meteorological Society} \bibinfo{volume}{83}, \bibinfo{pages}{407--430}.
\newblock \DOIprefix\doi{ 10.1175/1520-0477(2002)083<0407:DIHRPM>2.3.CO;2}.
%Type = Article
\bibitem[{M{\'{e}}rigot and Thibert(2020)}]{merigot_thibert_2020}
\bibinfo{author}{M{\'{e}}rigot, Q.}, \bibinfo{author}{Thibert, B.}, \bibinfo{year}{2020}.
\newblock \bibinfo{title}{{Optimal transport: discretization and algorithms}}.
\newblock \bibinfo{journal}{Handbook of Numerical Analysis} \bibinfo{volume}{22}, \bibinfo{pages}{133--212}.
\newblock \DOIprefix\doi{ 10.1016/bs.hna.2020.10.001}.
%Type = Article
\bibitem[{Mittermaier et~al.(2013)Mittermaier, Roberts and Thompson}]{mittermaier_2013}
\bibinfo{author}{Mittermaier, M.}, \bibinfo{author}{Roberts, N.}, \bibinfo{author}{Thompson, S.A.}, \bibinfo{year}{2013}.
\newblock \bibinfo{title}{{A long-term assessment of precipitation forecast skill using the Fractions Skill Score}}.
\newblock \bibinfo{journal}{Meteorological Applications} \bibinfo{volume}{20}, \bibinfo{pages}{176--186}.
\newblock \DOIprefix\doi{ 10.1002/MET.296}.
%Type = Misc
\bibitem[{Monge(1781)}]{monge1781}
\bibinfo{author}{Monge, G.}, \bibinfo{year}{1781}.
\newblock \bibinfo{title}{{M\'{e}moire sur la th\'{e}orie des d\'{e}blais et des remblais}}.
\newblock \URLprefix \url{ https://search.worldcat.org/title/Memoire-sur-la-theorie-des-deblais-et-des-remblais/oclc/51928110}.
%Type = Article
\bibitem[{Murphy(1993)}]{murphy_1992}
\bibinfo{author}{Murphy, A.H.}, \bibinfo{year}{1993}.
\newblock \bibinfo{title}{{What Is a Good Forecast? An Essay on the Nature of Goodness in Weather Forecasting}}.
\newblock \bibinfo{journal}{Weather and Forecasting} \bibinfo{volume}{8}, \bibinfo{pages}{281--293}.
\newblock \DOIprefix\doi{ https://doi.org/10.1175/1520-0434(1993)008<0281:WIAGFA>2.0.CO;2}.
%Type = Misc
\bibitem[{{NCAR}(n.d.)}]{ral_icp_webpage}
\bibinfo{author}{{NCAR}}, \bibinfo{year}{n.d.}
\newblock \bibinfo{title}{{RAL | JNT | Forecast Evaluation and Applied Statistics | Spatial Forecast Methods Inter-Comparison Project (ICP/MesoVICT)}}.
\newblock \URLprefix \url{ https://projects.ral.ucar.edu/icp/}.
%Type = Article
\bibitem[{Nishizawa(2024)}]{nishizawa_2024}
\bibinfo{author}{Nishizawa, S.}, \bibinfo{year}{2024}.
\newblock \bibinfo{title}{{Extracting Latent Variables From Forecast Ensembles and Advancements in Similarity Metric Utilizing Optimal Transport}}.
\newblock \bibinfo{journal}{Journal of Geophysical Research: Machine Learning and Computation} \bibinfo{volume}{1}.
\newblock \DOIprefix\doi{ 10.1029/2023JH000112}.
%Type = Article
\bibitem[{Papayiannis et~al.(2018)Papayiannis, Galanis and Yannacopoulos}]{papayiannis_et_al_2018}
\bibinfo{author}{Papayiannis, G.I.}, \bibinfo{author}{Galanis, G.N.}, \bibinfo{author}{Yannacopoulos, A.N.}, \bibinfo{year}{2018}.
\newblock \bibinfo{title}{{Model aggregation using optimal transport and applications in wind speed forecasting}}.
\newblock \bibinfo{journal}{Environmetrics} \bibinfo{volume}{29}.
\newblock \DOIprefix\doi{ https://doi.org/10.1002/env.2531}.
%Type = Article
\bibitem[{Peyr{\'{e}} and Cuturi(2019)}]{peyre_cuturi_2018}
\bibinfo{author}{Peyr{\'{e}}, G.}, \bibinfo{author}{Cuturi, M.}, \bibinfo{year}{2019}.
\newblock \bibinfo{title}{{Computational Optimal Transport: With Applications to Data Science}}.
\newblock \bibinfo{journal}{Foundations and Trends{\textregistered} in Machine Learning} \bibinfo{volume}{11}, \bibinfo{pages}{355--607}.
\newblock \DOIprefix\doi{ 10.1561/2200000073}.
%Type = Article
\bibitem[{Rubner et~al.(2000)Rubner, Tomasi and Guibas}]{rubner_et_al_2000}
\bibinfo{author}{Rubner, Y.}, \bibinfo{author}{Tomasi, C.}, \bibinfo{author}{Guibas, L.J.}, \bibinfo{year}{2000}.
\newblock \bibinfo{title}{The earth mover's distance as a metric for image retrieval}.
\newblock \bibinfo{journal}{International journal of computer vision} \bibinfo{volume}{40}, \bibinfo{pages}{99--121}.
\newblock \DOIprefix\doi{ 10.1023/A:1026543900054}.
%Type = Article
\bibitem[{Santambrogio(2015)}]{Santambrogio_2015}
\bibinfo{author}{Santambrogio, F.}, \bibinfo{year}{2015}.
\newblock \bibinfo{title}{{Optimal Transport for Applied Mathematicians}} \bibinfo{volume}{87}.
\newblock \DOIprefix\doi{ 10.1007/978-3-319-20828-2}.
%Type = Article
\bibitem[{Schmitzer(2019)}]{schmitzer_2019}
\bibinfo{author}{Schmitzer, B.}, \bibinfo{year}{2019}.
\newblock \bibinfo{title}{{Stabilized sparse scaling algorithms for entropy regularized transport problems}}.
\newblock \bibinfo{journal}{SIAM Journal on Scientific Computing} \bibinfo{volume}{41}, \bibinfo{pages}{A1443--A1481}.
\newblock \DOIprefix\doi{ 10.1137/16M1106018}.
%Type = Article
\bibitem[{S{\'{e}}journ{\'{e}} et~al.(2019)S{\'{e}}journ{\'{e}}, Feydy, Vialard, Trouv{\'{e}} and Peyr{\'{e}}}]{sejourne_et_al_2021}
\bibinfo{author}{S{\'{e}}journ{\'{e}}, T.}, \bibinfo{author}{Feydy, J.}, \bibinfo{author}{Vialard, F.X.}, \bibinfo{author}{Trouv{\'{e}}, A.}, \bibinfo{author}{Peyr{\'{e}}, G.}, \bibinfo{year}{2019}.
\newblock \bibinfo{title}{{Sinkhorn Divergences for Unbalanced Optimal Transport}} \URLprefix \url{ http://arxiv.org/abs/1910.12958}.
%Type = Article
\bibitem[{S{\'{e}}journ{\'{e}} et~al.(2022)S{\'{e}}journ{\'{e}}, Vialard and Peyr{\'{e}}}]{sejourne_vialard_pyre_2022}
\bibinfo{author}{S{\'{e}}journ{\'{e}}, T.}, \bibinfo{author}{Vialard, F.X.}, \bibinfo{author}{Peyr{\'{e}}, G.}, \bibinfo{year}{2022}.
\newblock \bibinfo{title}{{Faster Unbalanced Optimal Transport: Translation invariant Sinkhorn and 1-D Frank-Wolfe}}.
\newblock \bibinfo{journal}{Proceedings of Machine Learning Research} \bibinfo{volume}{151}, \bibinfo{pages}{4995--5021}.
\newblock \URLprefix \url{ https://arxiv.org/abs/2201.00730v1}.
%Type = Article
\bibitem[{Skok(2022)}]{skok_2022}
\bibinfo{author}{Skok, G.}, \bibinfo{year}{2022}.
\newblock \bibinfo{title}{{A New Spatial Distance Metric for Verification of Precipitation}}.
\newblock \bibinfo{journal}{Applied Sciences} \bibinfo{volume}{12}.
\newblock \DOIprefix\doi{ 10.3390/app12084048}.
%Type = Article
\bibitem[{Skok(2023)}]{skok_2023}
\bibinfo{author}{Skok, G.}, \bibinfo{year}{2023}.
\newblock \bibinfo{title}{{Precipitation attribution distance}}.
\newblock \bibinfo{journal}{Atmospheric Research} \bibinfo{volume}{295}, \bibinfo{pages}{106998}.
\newblock \DOIprefix\doi{ 10.1016/J.ATMOSRES.2023.106998}.
%Type = Article
\bibitem[{Skok and Lled{\'{o}}(2024)}]{skok_lledo_2024}
\bibinfo{author}{Skok, G.}, \bibinfo{author}{Lled{\'{o}}, L.}, \bibinfo{year}{2024}.
\newblock \bibinfo{title}{{Spatial verification of global precipitation forecasts}} \URLprefix \url{ https://arxiv.org/abs/2407.20624v1}.
%Type = Article
\bibitem[{Steinacker et~al.(2000)Steinacker, H{\"a}berli and P{\"o}ttschacher}]{steinacker_et_al_2000}
\bibinfo{author}{Steinacker, R.}, \bibinfo{author}{H{\"a}berli, C.}, \bibinfo{author}{P{\"o}ttschacher, W.}, \bibinfo{year}{2000}.
\newblock \bibinfo{title}{A transparent method for the analysis and quality evaluation of irregularly distributed and noisy observational data}.
\newblock \bibinfo{journal}{Monthly Weather Review} \bibinfo{volume}{128}, \bibinfo{pages}{2303--2316}.
\newblock \DOIprefix\doi{ 10.1175/1520-0493(2000)128<2303:ATMFTA>2.0.CO;2}.
%Type = Article
\bibitem[{Vanderbecken et~al.(2023)Vanderbecken, Dumont Le~Brazidec, Farchi, Bocquet, Roustan, Potier and Broquet}]{vanderbecken_et_al_2023}
\bibinfo{author}{Vanderbecken, P.J.}, \bibinfo{author}{Dumont Le~Brazidec, J.}, \bibinfo{author}{Farchi, A.}, \bibinfo{author}{Bocquet, M.}, \bibinfo{author}{Roustan, Y.}, \bibinfo{author}{Potier, E.}, \bibinfo{author}{Broquet, G.}, \bibinfo{year}{2023}.
\newblock \bibinfo{title}{{Accounting for meteorological biases in simulated plumes using smarter metrics}}.
\newblock \bibinfo{journal}{Atmospheric Measurement Techniques} \bibinfo{volume}{16}, \bibinfo{pages}{1745--1766}.
\newblock \DOIprefix\doi{ 10.5194/amt-16-1745-2023}.
%Type = Article
\bibitem[{Villani(2009)}]{villani_2009}
\bibinfo{author}{Villani, C.}, \bibinfo{year}{2009}.
\newblock \bibinfo{title}{{Optimal Transport: Old and New}} \bibinfo{volume}{338}.
\newblock \DOIprefix\doi{ 10.1007/978-3-540-71050-9}.
%Type = Article
\bibitem[{Wernli et~al.(2008)Wernli, Paulat, Hagen and Frei}]{wernli_2008}
\bibinfo{author}{Wernli, H.}, \bibinfo{author}{Paulat, M.}, \bibinfo{author}{Hagen, M.}, \bibinfo{author}{Frei, C.}, \bibinfo{year}{2008}.
\newblock \bibinfo{title}{{SAL—A Novel Quality Measure for the Verification of Quantitative Precipitation Forecasts}}.
\newblock \bibinfo{journal}{Monthly Weather Review} \bibinfo{volume}{136}, \bibinfo{pages}{4470 -- 4487}.
\newblock \DOIprefix\doi{ 10.1175/2008MWR2415.1}.
%Type = Book
\bibitem[{Zwillinger and Kokoska(1999)}]{zwillinger1999crc}
\bibinfo{author}{Zwillinger, D.}, \bibinfo{author}{Kokoska, S.}, \bibinfo{year}{1999}.
\newblock \bibinfo{title}{CRC standard probability and statistics tables and formulae}.
\newblock \bibinfo{publisher}{Crc Press}.

\end{thebibliography}
